\documentclass{article}
\usepackage[english]{babel}
\usepackage{amsfonts}
\usepackage{amscd}

\title
{\bf  Skands and coskands\\ (The non-founded set theory with individuals and its model in the Field of all Conway numbers)}

\author{Ju. T. Lisica\thanks {{\it Dedicated to my son Andrew Yu. Lisitsa.} 
 }\\
St. Tikhon's Orthodox  Humanities University \\
Likhov Lane 6, p. 1, Moscow, 127051,\\
e-mail: jutlisica@yandex.ru}

\date{}
\begin{document}
\maketitle

\begin{abstract}
The basic one in this work is the axiomatic set theory $NBG$ (von Neumann-Bernays-G{\"o}del), which is a first-order theory with its own axioms, including in particular the axiom of choice ${\bf AC}$ and the axiom of regularity ${\bf RA}$. The universal class ${\bf V}$ of all sets in this theory exactly coincides with the class of all founded sets, i.e., such $X\in{\bf V}$ that {\it does not exist} an infinitely descending $\in$-sequence $X\ni X_1\ni X_2\ni...\ni X_n\ni...$ of sets $X_n$, $n=1,2,3,...\,\,$. 

In the first part of the paper, a new concept of {\it skand} is introduced -- a random aggregate, or \grqq decreasing\grqq\, tuple composed of founded sets, e.g., $X=\{1,\{2,\{3,\{...\,\,\,...\}\}\}\}$, and the theory of $NBG^-=NBG-{\bf RA}$, i.e., the theory of $NBG$ without the axiom of regularity ${\bf RA}$, to which is added the new axiom ${\bf SEA}$ of the existence of infinite-length skands and the pseudo-founding axiom  ${\bf PFA}$. These new axioms are a negation of the axiom of regularity and are thus less restrictive than the axiom of regularity ${\bf RA}$ in the sense that they admit the existence of non-founded sets, and the axiom of regularity excludes the existence of such sets. At the same time of course the axiom of extensionality ${\bf EA}$ is replaced by a more accurate axiom of extensionality ${\bf EEA}$, since it takes into account the equality of new objects. 

In the second part of the paper, a new concept of {\it coskand} is introduced, which is dual to a notion of skand and is a random aggregate, or \grqq increasing\grqq\, tuple composed of founded sets and the theory of $NBG$ and actually is a theory  $NBG[\cal U]$ with individuals as limiting coskands, e.g., $X=...\{3,\{2,\{1,\{0\}\}\}\}...\,\,$.

New mathematical objects as skands and coskands provide an opportunity  to build more complex formal theories of non-founded sets with usual individuals and in addition with individuals like coskands as well as with coskands based over non-founded sets, not only over founded sets.

As a consequence of this theory is the following preferred formulas for contradiction and existence or nonexistence:
$$(\exists x)(\forall y[y\in x\Leftrightarrow y\notin y])\Rightarrow(\neg\exists x)(\forall y[y\in x\Leftrightarrow y\notin y])\,\, {\bf vs}\,\, x\in x\Leftrightarrow x\notin x$$

\begin{center}
and 

$\neg\exists X\forall x[x\in X]\,\,{\bf vs}\,\,\exists X\neg\exists\{X\}$.
\end{center}

This new approach makes it possible to expand the understanding of the theory of non-founded sets with individuals but also to enrich the theory with new mathematical objects that were previously considered  {\it undefinable notions} like  \grqq trans-definite sets\grqq\, (see \cite{l315}, p. 311) and in particular to bypass paradoxical Russell-type sets, which as a result leads to new consistent non-mathematical objects, previously conventionally called {\it sets without boundaries}, which go beyond not only the universal of the proper class {\bf V}, but in mathematics itself, but having their own interesting and diverse, say, \grqq geography\grqq.

Models of all these formal set theories are constructed in the Field of all Conway numbers ${\bf No}$. A digest of the main constructions of Conway number theory and the latest new results of this theory is given.

The difference between infinities (say, $+\infty_\omega$) and infinite numbers (say, $\omega$)  in the structure of the Conway numbers fields is explained. 
 
The paper also considers previously unnoticed patterns of the infinite sequences $(\alpha)_{0<\alpha<\lambda}$  of ordinal numbers $0<\alpha<\lambda$ ($\lambda$ is a limit ordinal) and their reciprocals $(\frac{1}{\alpha})_{0<\alpha<\lambda}$ as \grqq pulsars\grqq\, or pulsating \grqq jumps\grqq\, and their corresponding periodic \grqq compressions\grqq,  as well as various algebraic structures of infinitesimal Conway numbers.

Examples of exponential functions $y=\omega^x$ and $y=e^x$, $x\in{\bf No}$, are given, for the first of which the argument coincides with the value of the function at some points, i.e., $\omega^{x_0}=x_0$, and for the second, the argument can be greater than the value of the function, i.e., $e^{x_1}<x_1$, a kind of anti-Zeno paradox: the tortoise (input) catches up with Achilles (output) in the first case and the tortoise overtakes him in the second case.
\end{abstract}

\bigskip

\begin{center}
{\bf 1. Introduction}
\end{center}

Axiomatic set theories, more precisely, formal logistic systems that are first-order theories with logical axioms and proper axioms of set theories (which, together with the definition of equality of sets, become first-order theories with equality) are well known and carefully developed. Such, for example, are set theory $ZF$ with an infinite scheme (system) of Zermelo-Frenkel axioms and set theory $NBG$ with a finite system of von Neumann-Bernays-G\"{o}del axioms. Both of these systems can include the axiom of choice ${\bf CA}$ and the axiom of regularity ${\bf RA}$. Here we will assume that they are included in the system of axioms of these theories, and through $ZF^-$ and $NBG^-$ denote them without the axiom of regularity ${\bf RA}$, respectively. The latter, with the existing axiom of choice, is equivalent to the axiom of foundation ${\bf FA}$ (see: Proposition 4.39 in \cite{l7}).

Therefore, we will talk here about the axiom  ${\bf FA}$, which says that for any set 
$X\in{\bf V}$ from the universe ${\bf V}$ of all sets  an infinitely descending $\in$-sequence $X\ni X_1\ni X_2\ni...$ of sets $X_n\in{\bf V}$, $n=1,2,3,...\,\,$ {\it does not exist}.

Note also that working with the formal theory of $NBG^-$, some authors (see: \cite{l31}, \cite{l94}) in those places where one has to use the axiom of choice {\bf CA} (for any set of $X$ there is such a choice function $F:2^X\setminus\emptyset\rightarrow X$, that $F(Y)\in Y$) for any element $Y\in 2^X\setminus\emptyset$), use the Neumann axiom {\bf N}, which states that the universe ${\bf V}$ of all sets of the theory is bijective to the class ${\bf On}$ of all ordinals, that in  the absence of the axiom of foundation, ${\bf FA}$ is a stronger axiom than even the axiom of global choice {\bf E} (there is a unitary operator $\sigma$ that maps any nonempty set $X$ to its element $\sigma(X)\in X$, see: \cite{l31}, \cite{l95}). In this work, we will not need such a strengthening of the axiom of choice, with the exception of one mention of the strong axiom of plenitude of individuals ${\bf SAP}$, but it is not excluded in the further development of the theory of non-founded sets.

The axiom of the foundation ${\bf FA}$ is very important, for example, in the axiom of the pair {\bf PA}: for any two different sets $X$ and $Y$, there is a single set $\{X,Y\}$, the elements of which are exactly $X$ and $Y$ and only them. At the same time, without the foundation axiom, it is impossible to prove that the set $\{X,Y\}$ is different from $X$ or from $Y$, since this, generally speaking, is incorrect: in non-founded set theories $ZF^-$ or $NBG^-$, i.e. without the foundation axiom, there are such different sets of $X$ and $Y$ that $\{X,Y\}=X$ or $\{X,Y\}=Y$; for example, if $X=\{0\}=\{\{\}\}\stackrel{def}{=}1$ and $Y=\{1,Y\}$, then $\{X,Y\}=\{1,Y\}=Y$.

(We will see below that such a set exists.) On the contrary, in $ZF$ and $NBG$ it is always
 $\{X,Y\}\not=X$ and $\{X,Y\}\not=Y$, since the disordered pair $\{X,Y\}$ of the founded sets $X$ and $Y$ is always founded and therefore is not a proper element, this means that it does not coincide with either $X$ or $Y$, since any family of founded sets always does not coincide with any of its elements (see: \cite{l14}, (3.7), p. 94). 
 
This directly follows from the fact that for two founded sets for two founded sets $X$ and $Y$ the ordered pair $Z=<X,Y>\stackrel{def}{=}\{\{X\},\{X,Y\}\}$ is also a founded set and $Z\not=X$, $Z\not=Y$, $Z\notin X$ and $Z\notin Y$ (see: \cite{l98}, p. 25). For non-founded sets, this is generally not the case, for example, if $X=\{X\}=\{\{...\}\}$, and $Y=\{X,Y\}=\{X,\{X,Y\}\}$, then $Z=<X,Y>=\{\{X\},\{X,Y\}\}=\{X,\{X,Y\}\}=Y$ and $Z\in Y$.

Despite this advantage of accuracy and distinctness in the axiom of a pair and the definition of an ordered pair, simplifying the definitions of ordinal and cardinal numbers, (\cite{l144}, pp. 110, 129 and p. 123) and, for example, that the generalized continuum hypothesis in this form: $2^{\aleph_\alpha}=\aleph_{\alpha+1}$ for any ordinal $\alpha$ -- entails the axiom of choice (without the axiom of regularity, this is not the case \cite{l14}, p. 103) and many other advantages, the axiom of regularity (and with the axiom of choice ${\bf CA}$ the axiom of foundation ${\bf FA}$) in a sense greatly narrows not only the \grqq class\grqq\, of all sets \grqq of naive\grqq\, Cantor's set theory (Cantor's \grqq inconsistent multiplicity\grqq), which has no denotation, but also the class ${\bf V}^-$ of all sets of the theory of $NBG^-$.

In this case, the so-called reflexive sets with respect to the membership relation $\in$ are excluded, i.e., when the set $X$ is a proper element: $X\in X$ (note that the only predicative letter $A^2_2$ of similar formal set theories, which is usually abbreviated to the symbol of the membership relation $\in$, when writing $X\in Y$ for $A^2_2(X,Y)$ and $X\notin Y$ for $\neg A^2_2(X,Y)$, does not exclude the reflexive relation $A^2_2(X,X)$, i.e. $X\in X$), and the so-called extraordinary Mirimanoff sets \cite{l2} are also excluded, i.e., sets $X$ that have infinitely descending $\in$-sequences $X\ni X_1\ni X_2\ni...\ni X_n\ni...\,\,$. Note that the axiom of foundation is not restrictive, since the multiplicity of all founded sets in the naive Cantor theory also remains inconsistent, and in the formal theory $NBG$ is its proper class, so it is more correct to say that this axiom is not restrictive, but narrows the field of reasoning of any set theory, although, as in this work It is the basis for further extensions of the universal class of all sets of formal theories excluding the axiom of foundation.

The so-called Russell \grqq principle of the vicious circle\grqq\, was also excluded. However, von Neumann \cite{l901}, introducing the axiom of regularity ${\bf RA}$ in 1925, pursued exactly this goal: to prohibit such sets, since it was believed that they lead to antinomies and paradoxes.

 But it immediately became clear that mathematics was deprived of the field of modeling such natural processes as cyclicity, the study of orthodox statements such as \grqq statements of a liar\grqq\, and many problems of linguistics, later computer science, graph theory, game theory, flow theory, etc., in which cyclic models lie outside the universe of founded sets.
 
 The uncertainty and ambiguity of the relation $X\in X$ was ignored, since the set $X$ is not defined by its elements in $X\setminus\{X\}$, i.e., by all elements other than $X$: with the same elements, there is its proper class of different sets containing themselves as therefore, the equality $X=Y$ of such sets coincides with the tautology $X=Y\Leftrightarrow X=Y$, therefore, $X$ is not defined as a set, since according to Cantor: \grqq A set is a collection of $M$ {\it defined} and {\it distinguishable objects} $m$\grqq\, \cite{l1}; and $X$ in a reflexive set, its element equal to the set itself is undefinable and indistinguishable, therefore, the set itself is undefinable.

 A similar uncertainty and ambiguity arises for the Mirimanoff sets $X\ni X_1\ni X_2\ni...$ in particular for $X\ni X\ni X\ni...$  the irregular set $X=\{\{\{...\}\}\}$, the only element of which is the set $X$ itself: there are a lot of such different sets - they make up their proper class, if one does not introduce an additional axiom of uniqueness that would identify all the elements of this proper class, see: \cite{l93}.
 
 Already in the first half of the 20th century, many authors offered their own \smallbreak
 \parindent=0 cm
 \grqq anti-founded axioms\grqq\, ${\bf AFA}$ and simply added their ${\bf AFA}$-universes of sets designed to solve these problems, for example, \cite{l91}, \cite{l31}, \cite{l90}, \cite{l92}, \cite{l94} , etc. All these theories are consistent if $ZF^-$ or $NBG^-$ are consistent. Nevertheless, they themselves differ in the definition of equality in the axiom of extensionality.

\parindent=0,5 cm

In this work, we will construct a new natural theory of unfounded sets that is different from the ones mentioned above. In some sense, it will be the most minimal and very similar to the theory of founded sets, in particular, the class of all its sets contains the class of all founded sets as a subclass, and the sets that are not contained in this class are non-founded but can be said to be \grqq pseudo-founded\grqq.

Here, for certainty, we will choose the formal theory of $NBG$ as the basis for constructing our theory of non-founded sets (it would be easy to base the theory of $ZF$ in the same way) and abandon the axiom of foundation ${\bf FA}$, i.e. consider the theory of $NBG^-=NBG-{\bf FA}$. Next, we define a new mathematical object, which we will call a {\it skand}: a random aggregate (or tuple) of founded sets or their absence. Two new axioms about the existence of skands ${\bf ESA}$ and pseudo-founded ${\bf PFA}$ are restrictive axioms, but less rigid than the axiom of foundation. Thus, the resulting theory is such an extension of $NBG^-$ that it naturally expands the class ${\bf V}\stackrel{def}{=}{\bf V}^{(0)}$ of all the founded sets of the theory of $NBG$ underlying the new theory to the class 
${\bf V}^{(1)}\subset{\bf V}^-$, with ${\bf V}^{(0)}\subset{\bf V}^{1}$.

This is the first step of inductive theory building $NBG^{(\nu)}$ when $\nu=1$.

A similar sequential extension is constructed in the same natural way for any ordinal number $\nu>1$.
 
 The first application of the theory of skands is the natural definition of coskands, which, in the case of their infinite length, are either sets of a new theory, which is an extension of the theory of $NBG^-$, or new individuals of this set theory.

Another application of this non-founded set theory is new proofs of the results of set theory used by \grqq paradoxical sets\grqq, such as the diagonal Cantor set $K\stackrel{def}{=}\{x\notin\varphi(x)\}$, where 
$\varphi:X\rightarrow 2^X$ the assumed bijection; as the Russell set ${\bf R}\stackrel{def}{=}\{X\in {\bf R}\Leftrightarrow X\notin X\}$; as the Mirimanoff set ${\bf W}$ of all founded sets; as a Zermelo set $M$ containing as its elements all its subsets; as a Hilbert set $H$ constructed naturally from natural numbers. All the results in the theory of $NBG^-$ obtained with the help of these paradoxical sets -- that they have no correspondence under the assumed bijection, or that they are proper classes, or that they do not exist -- are bypassed by the corresponding \grqq paradoxes\grqq. Moreover, the same method proves that the set ${\bf S}\stackrel{def}{=}\{X\in{\bf S}\Leftrightarrow X\in X\}$ of all reflexive sets is not a set, but is its proper class, although the formal substitution instead of The $X$ of the assumed set ${\bf S}$, which allows the first-order theory, does not lead to a contradictory formula, as in the Russell paradox. Thus, the supposed dual \grqq sets\grqq\, ${\bf R}$ and ${\bf S}$ in set theory become equal and are their proper classes.

Models for the theory of skands and the theory of coskands have been constructed in the universal class ${\bf No}$ of all Conway numbers. This universal class has the algebraic structure of a linearly ordered field and can therefore be useful for models and even non-standard models of other formal theories. The paper also considers previously unnoticed patterns of the infinite sequences $(\alpha)_{0<\alpha<\lambda}$  of ordinal numbers $0<\alpha<\lambda$ ($\lambda$ is a limit ordinal) and their reciprocals $(\frac{1}{\alpha})_{0<\alpha<\lambda}$ as \grqq pulsars\grqq\, or pulsating \grqq jumps\grqq\, and their corresponding periodic \grqq compressions\grqq,  as well as various algebraic structures of infinitesimal Conway numbers.

The Conway function $y=\omega^x$, $x\in{\bf No}$, is used to define the exponential and logarithmic functions $y=e^x$, $x\in{\bf No}$, and $y=\ln x$, $x\in{\bf D}^+_\Omega\subset{\bf No}$,  which, unlike the Conway function $y=\omega^x$, are continuous in the linearly ordered topology in ${\bf No}$ and have unexpected properties: for the Conway function $y=\omega^x$ the argument coincides with the value of the function at some points, i.e., $\omega^{\varepsilon_0}=\varepsilon_0$, and for the function $y=e^x$ the argument can be greater than the value of the function, i.e., $e^{\omega\cdot\varepsilon_0}=\varepsilon_0<\omega\cdot\varepsilon_0$, a kind of anti-Zeno paradox: the tortoise (argument) catches up with Achilles (value of function) in the first case and the tortoise overtakes him in the second case.

\begin{center}
{\bf 2. Necessary information from the theory of sets with individuals, prerequisites and notation}
\end{center}

 As already mentioned, for the sake of certainty, we will consider here the set theory $NBG$, which contains the axiom of choice ${\bf CA}$ and the axiom of foundation ${\bf FA}$, which is equivalent to the axiom of regularity {\bf RA}. This theory will be the basis of all our constructions. Since further constructions will use the theory of sets with individuals, then for the uniformity of the description we consider the theory of sets with some fixed set or proper class of ${\cal U}$ individuals (they also say atoms or proto-elements), which, generally speaking, can be empty (the theory of \grqq pure sets\grqq, considered above in Annotations and in the Introduction). Let's denote this theory by $NBG[{\cal U}]$. We borrow the information we need here from the monograph \cite{l7}, pp. 297-304. Similarly, we can consider the theory of Zermelo-Frenkel sets with individuals $ZF[{\cal U}]$.)
$$
$$
 
 {\bf 1.} Set theory $NBG^-[{\cal U}]$ is a formal first-order logical theory with equality, which is constructed as follows: the variables of this theory are lowercase bold Latin letters ${\bf x}$, ${\bf y}$, ${\bf z}$,...; there is one binary predicative letter $A^2_2$ and one unary predicative letter $A^1_1$, which are abbreviated: $A^2_2({\bf x,y})$ to ${\bf x}\in{\bf y}$, $\neg A^2_2({\bf x,y})$ to ${\bf x}\notin{\bf y}$ and $A^1_1({\bf x})$ to ${\bf Cls}({\bf x})$, and reads like \grqq${\bf x}$ there is a class\grqq. The usual entries of formulas in \grqq pure\grqq\, set theory $NBG^-$, as $(\forall X){\cal B}(X)$ and $(\exists X){\cal B}(X)$ in $NBG^-[{\cal U}]$, written here as $(\forall {\bf x})({\bf Cls}({\bf x})\Rightarrow{\cal B}({\bf x}))$ and $(\exists {\bf x})({\bf Cls}({\bf x})\wedge{\cal B}({\bf x}))$ accordingly (${\cal B}(X)$ means that $X$ has the property ${\cal B}$); notation ${\bf Pr}({\bf x})$ is given by the formula ${\bf Cls}({\bf x})\wedge\neg {\bf M}({\bf x})$, and reads \grqq${\bf x}$ has its proper class\grqq; designation ${\bf Ur}({\bf x})$ is given by the formula $\neg {\bf Cls}({\bf x})$, and reads \grqq ${\bf x}$ there is an individual\grqq. 
 
 Thus, in the theory of $NBG^-[{\cal U}]$, two concepts are distinguished: {\it classes} and {\it individuals}, as objects of the theory, and classes have two types: {\it sets}  and {\it proper classes}. Denote by ${\bf El}({\bf x})$ formula ${\bf M}({\bf x})\vee {\bf Ur}({\bf x})$ and reads \grqq ${\bf x}$ there is an element\grqq. In this interpretation, sets and individuals are objects that are elements, or members of classes; proper classes are also objects, but are not elements of any other objects -- sets or proper classes.

The definition of the equality of two theoretical objects $NBG^-[{\cal U}]$, i.e. ${\bf x}={\bf y}$ is an abbreviation of the following formula:
\begin{equation}
 \label{f1947}
 [{\bf Cls}({\bf x})\wedge {\bf Cls}({\bf y})\wedge(\forall{\bf z})({\bf z}\in{\bf x}\Leftrightarrow{\bf z}\in{\bf y})]\vee[{\bf Ur}({\bf x})\wedge {\bf Ur}({\bf y})\wedge(\forall {\bf z})({\bf x}\in{\bf z}\Rightarrow{\bf y}\in{\bf z})].
 \end{equation}
 
In addition to the logical axioms of the formal theory of the first order, there is a finite set of special axioms of the formal theory $NBG^-[{\cal U}]$, which are divided into axioms for sets, axioms of the existence of classes and axioms for individuals. We will note only some of them that we will use in this work.

Axiom ${\bf Ur 1}$ about individuals:
\begin{equation}
 \label{f110}
(\forall{\bf x})(Ur({\bf x})\Rightarrow (\forall{\bf y})({\bf y}\notin{\bf x})).
\end{equation} 

 The axiom of  extensionality ${\bf Ur 2}$:
\begin{equation}
\label{f1951} 
 (\forall {\bf x})({\bf Cls}({\bf x}))(\forall {\bf y})({\bf Cls}({\bf y}))(\forall {\bf z})({\bf Cls}({\bf z}))({\bf x}={\bf y}\wedge {\bf x}\in{\bf z}\Rightarrow{\bf y}\in{\bf z}).
\end{equation} 

The first (of seven) axiom of the existence of classes ${\bf Ur 5}$:

\begin{equation}
\label{f2951} 
 (\exists X)(\forall{\bf u})(\forall{\bf v})({\bf El}({\bf u})\wedge{\bf El}({\bf v})\Rightarrow[<{\bf u},{\bf v}>\in X]\Leftrightarrow{\bf u}\in{\bf v}),
\end{equation} 
where $<{\bf u},{\bf v}>$ denotes the Kuratowski ordered pair $\{\{{\bf u}\},\{{\bf u},{\bf v}\}\}$.
 
The axiom ${\bf Ur 12}$ of the existence of a class ${\bf V}[{\cal U}]$ of all sets:

\begin{equation}
 \label{f010}
(\exists X)(\forall{\bf u})({\bf u}\in X\Leftrightarrow {\bf M}({\bf u})).
\end{equation}

Axiom of the sum of sets ${\bf Ur 13}$:

\begin{equation}
\label{f3951} 
 (\forall {\bf x})({\bf M}({\bf x}))(\exists {\bf y})({\bf M}({\bf y}))(\forall{\bf u})({\bf El}({\bf u}))({\bf u}\in {\bf y}\Leftrightarrow (\exists{\bf v})({\bf M}({\bf v}))({\bf u}\in{\bf v}\wedge{\bf v}\in {\bf x})).
\end{equation}
The set ${\bf y}$ in the sum axiom is sometimes denoted for brevity by $\cup{\bf x}$, more precisely $y=\bigcup\limits_{{\bf v}\in{\bf x}}{\bf v}$ and they say that ${\bf y}$ is the union of sets ${\bf v}$ that are elements of the set ${\bf x}$. The axioms of the existence of classes allow us to define in the only way an improper class $y=\cup{\bf x}=\bigcup\limits_{{\bf v}\in{\bf x}}{\bf v}$ in the case when ${\bf x}$ is an improper class (see: \cite{l7}, p. 184), i.e. it is possible to combine sets ${\bf v}$ by  a number of indexes that themselves form their proper class, i.e. it makes sense, for example, such a union $\bigcup\limits_{\alpha\in A}X_\alpha$ of sets $X_\alpha$, that the index class $A$ is a proper class. In particular we can consider our proper class 
$\bigsqcup\limits_{\alpha\in A}X_\alpha\stackrel{def}{=}\{<x_\alpha,\alpha>\,|\,x_\alpha\in X_\alpha,\alpha\in A\}$, and $<x_\alpha,\alpha>=\{\{x_\alpha\},\{x_\alpha,\alpha\}\}$ -- an ordered pair, which is called {\it discrete sum} or {\it discrete union} of sets $X_\alpha$, $\alpha\in A$.

The axiom of the power of sets ${\bf Ur 14}$:

\begin{equation}
\label{f4951} 
 (\forall {\bf x})(\exists {\bf y})(\forall{\bf u})({\bf u}\in {\bf y}\Leftrightarrow {\bf u}\subset {\bf y}),
\end{equation} 
 where ${\bf u}\subset y$ denotes the formula ${\bf M}({\bf u})\wedge {\bf M}({\bf x})\wedge(\forall {\bf v})({\bf v}\in{\bf u}\Rightarrow{\bf v}\in{\bf x})$.

The axiom of replacement ${\bf Ur 15}$:
\begin{equation}
\label{f2951} 
 (\forall Y)(\forall {\bf x})({\bf Un}(Y)\Rightarrow(\exists y)(\forall {\bf u})[{\bf u}\in y\Leftrightarrow(\exists {\bf v})(<{\bf v},{\bf v}>\in Y\wedge{\bf v}\in{\bf x})]),
\end{equation} 
 where ${\bf Un}({\bf z})$ denotes the formula $(\forall {\bf x}_1),(\forall {\bf x}_2),(\forall {\bf x}_3)[El({\bf x}_1)\wedge El({\bf x}_2)\wedge El({\bf x}_3)\Rightarrow (<{\bf x}_1,{\bf x}_2>\in{\bf z}\wedge<{\bf x}_1,{\bf x}_3>\in{\bf z}\Rightarrow{\bf x}_2={\bf x}_3)]$.

The axiom of choice ${\bf Ur 16}$:
\begin{equation}
\label{f2851} 
 (\forall x)({\bf M}({\bf x}))(\exists{\bf Fnc}{F})\forall {\bf y}({\bf y}\subseteq {\bf x}\wedge {\bf y}\not=\emptyset)(F({\bf y})\in {\bf y}),
\end{equation} 
where ${\bf Fnc}(F)$, or $F$ is a function means the formula $F\subseteq{\bf V}[{\cal U}]^-\times{\bf V}[{\cal U}]^-\wedge{\bf Un}(F)$with ${\bf V}[{\cal U}]^-\times{\bf V}[{\cal U}]^-\stackrel{def}{=}\{{\bf z}\Leftrightarrow({\exists u})({\exists v})({\bf z}=<{\exists u},{\exists v}>\wedge{\bf u},{\bf v}\in{\bf V}[{\cal U}]^-)\}$ and $F({\bf u}={\bf v})$ if $<{\bf u},{\bf v}>\in F$.

 The remaining axioms -- the axiom of infinity, the axiom of separation, the axioms of the existence of an empty set, an unordered pair, the remaining axioms of the existence of classes -- we will not remind you here, but if we want to attach to the theory of $NBG[{\cal U}]^-$ the axiom of regularity ${\bf RA}$:

\begin{equation}
 \label{f210}
 (\forall X)(X\not=\emptyset\Rightarrow(\exists{\bf u})({\bf u}\in X\wedge\neg(\exists{\bf v}\in X\wedge{\bf v}\in{\bf u}))),
 \end{equation}
which, in the presence of the axiom of choice ${\bf Ur 16}$, is equivalent to the foundation axiom we need ${\bf FA}$, which has the form:
\begin{equation}
 \label{f1950}
 (\neg \exists x)(Fnc (x)\wedge {\cal D}(x)=\omega\wedge(\forall u\in\omega)\Rightarrow x(u')\in x(u)), 
 \end{equation}
 where $u'=u+1$, $Fnc(x)$ denotes that $x$ is a function, and ${\cal D}(x)$ is the domain of definition of the function $x$, which here coincides with the set $\omega$ of all natural numbers $n=0,1,2,...$; i.e., there is no infinitely descending $\in$ sequence ${\bf x}_1\ni{\bf x}_2\ni...\,\,$, then we get the extension $NBG[{\cal U}]^-+{\bf RA}$theory $NBG[{\cal U}]^-$, which we denote by $NBG[{\cal U}]$. Both formal theories are formal theories (systems) of the first order with equality.
\begin{equation}
 \label{f210}
  ((\forall X)(X\not=\emptyset\Rightarrow(\exists{\bf u})({\bf u}\in X\wedge\neg(\exists{\bf v}\in X\wedge{\bf v}\in{\bf u}))),
  \end{equation}
which, in the presence of the axiom of choice ${\bf Ur 16}$, is equivalent to the foundation axiom we need ${\bf FA}$, which has the form:
\begin{equation}
  \label{f1950}
  (\neg \exists x)(Fnc (x)\wedge {\cal D}(x)=\omega\wedge(\forall u\in\omega)\Rightarrow x(u')\in x(u)), 
  \end{equation}
 where $u'=u+1$, $Fnc(x)$ denotes that $x$ is a function, and ${\cal D}(x)$ is the domain of the definition of the function $x$, which here coincides with the set $\omega$ of all natural numbers $n=0,1,2,...$; i.e., there is no infinitely descending $\in$ sequence ${\bf x}_1\ni{\bf x}_2\ni...\,\,$, then we get the extension $NBG[{\cal U}]^-+{\bf RA}$theory $NBG[{\cal U}]^-$, which we will denote by $NBG[{\cal U}]$. Both formal theories are formal theories (systems) of the first order with equality.

 In simple words, in the theory of $NBG[{\cal U}]$, there are three types of objects: individuals, sets, and proper classes. Individuals do not contain any objects of this theory, but they themselves can be contained both in non-empty sets and in their proper classes; non-empty sets can contain individuals or other sets, including the empty set $\emptyset$, but do not contain their proper classes as their elements, on the other hand, any set can be contained in another a nonempty set or a proper class; proper classes contain sets or individuals, but are not contained either in sets or in proper classes. Moreover, all its objects are founded, i.e. in one of its objects ${\bf x}$ the maximum descending $\in$ sequences ${\bf x}={\bf x}_1\ni{\bf x}_2\ni...\ni{\bf x}_n$ end either with an empty set ${\bf x}_n=\emptyset$, or with an individual ${\bf x_n}$, i.e. ${\bf Ur}({\bf x}_n)$. (Note that the empty set $\emptyset$ is the only special individual that is conveniently called a set and is given by the following axiom: $(\exists{\bf x})(M({\bf x}))(\forall{\bf y})({\bf y}\notin{\bf x})$; differs from other individuals in the definition of equality.)

{\bf 2.} Here we denote the proper classes of all sets of theories $NBG[{\cal U}]^-$ and $NBG[{\cal U}]$ by ${\bf V}^-[{\cal U}]$ and ${\bf V}[{\cal U}]$ respectively, the universal classes of elements of the theories $NBG[{\cal U}]^-$ and $NBG[{\cal U}]$ via ${\bf U}^-[{\cal U}]\stackrel{def}{=}{\bf V}^-[{\cal U}]\cup {\cal U}$ and ${\bf U}[{\cal U}]\stackrel{def}{=}{\bf V}[{\cal U}]\cup {\cal U}$, respectively, and the class of all ordinal numbers (ordinals) through ${\bf On}$. 
 
 Note that ${\bf On}\subset{\bf V}$, and by the ordinal $\alpha\in{\bf On}$ we mean \grqq the pure set\grqq\, of ${\bf V}$, which is ordered by $\in$ and transitive, i.e. $y\in x\wedge z\in y\Rightarrow z\in x$. Given the presence of the axiom of regularity in $NBG$, the class ${\bf On}$ is well-ordered ordered by $\in$, i.e. any of its subclasses has the smallest element. We will designate ordinals in Greek letters. The canonical representation of ordinal numbers is as follows:
\begin{equation}
 \label{f9999}
 0=\{\},\,1=\{\{\}\}=\{0\},\,2=\{\{\},\{\{\}\}\}=\{0,1\},..., \alpha=\{0,1,...,\alpha',...\}_{\alpha'<\alpha},...
 \end{equation}
 
As usual, via $[\alpha_0,\alpha)$ we denote the interval of ordinal numbers, i.e. all such $\alpha'\in {\bf On}$ that $\alpha_0\leq\alpha'<\alpha$, and after $(\alpha_0,\alpha)$ we denote the interval of ordinal numbers, i.e. all such $\alpha'\in {\bf On}$ that $\alpha_0<\alpha'<\alpha$.

In this case, it is easy to see that the relation $<$ of the ordinal numbers $\alpha,\beta\in{\bf On}$ can be expressed by the following formula: $\alpha<\beta\Leftrightarrow[0,\alpha)\subset[0,\beta)$. This allows you to enter the symbol $\Omega$, which corresponds to the proper class of all ordinals ${\bf On}\stackrel{def}{=}[0,\Omega)$, for which $\alpha<\Omega$ for each $\alpha\in{\bf On}$, since the last class ${\bf On}$ is transitive on all ordinals, i.e. $\alpha\in\beta\in{\bf On}\Rightarrow\alpha\in{\bf On}$. 

\smallskip
 
 Cardinal numbers will be denoted by the German Gothic letters ${\frak a}$, if we exclude the axiom of choice; the Alephs $\aleph_\alpha$ are the powers of the initial ordinals $\omega_\alpha$, $\alpha\in{\bf On}$. For any set $X$, by $|X|$ we denote the cardinality of this set, or the class of sets bijective to the set $X$. The possibility of one set $X$ is less than the cardinality of another set $|Y|$, i.e. $|X|<|Y|$ if there is a bijection of the set $X$ onto some subset $Y'$ of the set $Y$. If $Y'=Y$, then $|X|=|Y|$; this happens, according to the Cantor-Bernstein theorem, if and only if $|X|<|Y|$ and $|Y|<|X|$.  If there is no axiom of choice, then two sets $X$ and $Y$ can be incomparable in power, for example, if the set $X$ is a Dedekind set, i.e. an infinite set that is not bijective to its proper subset, and $Y$ is any countable set or any well-ordered set.
 
 For simplicity (abbreviations of writing formulas) we will now denote sets or classes with capital letters $X$, omitting predicative letters, and the elements of these sets or classes $x\in X$ will be denoted with lowercase letters, and not as above in bold, also omitting predicative letters.
 
Using ${\cal P}$, we denote the operator that maps to each set $X$ of the theory $NBG[{\cal U}]$ the set ${\cal P}(X)$ of this theory, given by the formula $\forall x (x\in{\cal P}(X)\equiv x\subseteq X)$, i.e. ${\cal P}(X)$ is the class of all subsets of the set $X$. According to the axiom ${\bf W}$ of the set of all subsets, or the axiom of degree: $\forall x\exists y\forall u(u\in y\equiv u\subseteq x)$, the formula ${\bf M}({\cal P}(x))$. Often ${\cal P}(X)$ is denoted by $2^X$, $|{\cal P}(X)|$ after $2^{|X|}$.

 \smallskip
 Since our construction will be inductive, we denote the original set theory $NBG[{\cal U}]$ by $NBG[{\cal U}]^{(0)}$, a set or class of individuals ${\cal U}$ via ${\cal U}^{(0)}$, the class of all sets of this theory ${\bf V}[{\cal U}]$ via ${\bf V}[{\cal U}]^{(0)}$.

The model of the theory of $NBG[{\cal U}]$ is constructed as follows.
If  ${\cal U}$ is a set, then for $NBG[{\cal U}]$ the model is constructed as follows:
\begin{equation}
\label{f1948}
\begin{array}{cc}
\Xi(0)={\cal U}\\
\Xi(\alpha')={\cal P}(\Xi(\alpha))\\
\lim\limits(\lambda)\Rightarrow \Xi(\lambda)=\bigcup\limits_{\beta<\lambda}\Xi(\beta)\\
{\bf H}'[{\cal U}]=\bigcup\limits_{\alpha\in{\bf On}}\Xi(\alpha).
\end{array}
\end{equation}
Here $\alpha'=\alpha+1$, $\beta$ are ordinals from the class ${\bf On}$.

Then the class ${\bf H}[{\cal U}]=\bigcup\limits_{\alpha\in{\bf On}}\Xi(\alpha)$ is the model for $NBG[{\cal U}]$. Moreover, the axiom of regularity ${\bf RA}$ holds in $NBG[{\cal U}]^-$ if and only if ${\bf H}[{\cal U}]={\bf U}[{\cal U}]$.

If  ${\cal U}$ is a proper class, then for $NBG[{\cal U}]$ the model is constructed as follows:
for any set $L$ embedded in the class ${\cal U}$ and any ordinal $\gamma\in{\bf On}$ we construct the set $\Xi^\gamma_L$ in the following recursive way:

 \begin{equation}
\label{f1949}
\begin{array}{cc}
\Xi^\gamma_L(0)=L\\
\Xi^\gamma_L(\alpha')={\cal P}(\Xi^\gamma_L(\alpha)),\,\,\,\,  \alpha'<\gamma,\\
\lim\limits(\lambda)\Rightarrow \Xi^\gamma_L(\lambda)=\bigcup\limits_{\beta<\lambda}\Xi^\gamma_L(\beta), \,\,\,\,\lambda<\gamma.
\end{array}
\end{equation}
Here, as above, $\alpha'=\alpha+1$, $\beta$ are ordinals from the class ${\bf On}$.  

  \smallskip
Then the class ${\bf H}^*[{\cal U}]$, consisting of all such elements $M$ for which there is such a set of individuals $L$ and the ordinal $\gamma$ that $M\in\Xi^\gamma_{{\cal U}}$, is a model for ${\bf NBG}[{\cal U}]$. Moreover, the axiom of regularity ${\bf RA}$ holds in $NBG[{\cal U}]^-$ if and only if ${\bf H}^*[{\cal U}]={\bf U}[{\cal U}]$. 
\smallskip
 
  {\bf 4.} For technical purposes, we will need the following convention. Every set $X$ (or function represented by its graph, or the skand introduced below, as well as other similar objects of set theory) can be treated as a new individual. To do this, we need to apply to $X$ \grqq the forgetful functor\grqq\, $F$, i.e., to forget the structure in $X$ defined by the membership relation $\in$ in the object $FX$. In particular $FX$ does not contain any elements, i.e., it is an individual. For example, the natural numbers $0, 1, 2, 3,...$ can be considered as individuals by applying the oblivion functor ${\bf F}$ to the sets representing them: 
 $\{\}$,  $\{\{\}\}$, $\{\{\},\{\{\}\}\}$, $\{\{\},\{\{\}\},\{\{\},\{\{\}\}\}\}$,... or  $\{\}$, $\{0\}$, $\{0,1\},\{0,1,2\},...$.  Note that in this example, the empty set $\emptyset=\{\}$ representing the number/individual $0$ is itself an individual, which is usually called a set for convenience of terminology. The other individuals are not called sets. (It is common to write $0=\{\}$, $1=\{0\}$, etc. In our case, it would have to be written as $0=F\{\}$, $1=F\{0\}=F\{\{\}\}$, etc., but there's no confusion in the first entry either, so to speak: \grqq unit $1$ is an individual\grqq\, and \grqq $\{0\}$\grqq\, is a set consisting of one element -- individual $0$).

 If we omit the forgetful functor $F$ (more precisely, we apply the composition of the functors $F^{-1}\circ F$, counting the inverse functor $F^{-1}$ to the functor $F$ by the operator \grqq of remembering\grqq\, $R=F^{-1}$), then from the individual $FX$ we get the original set $R\circ FX=F^{-1}\circ FX=1\circ X=X$ with all the membership relations of its elements, i.e. $(\forall x)x\in X$. Similarly with functions and skands: after lowering the forgetful functor, the structure of the $\in$ membership relations in these objects is restored. $\Box$

\begin{center}
{\bf 3. Necessary knowledge from Conway numbers theory}
\end{center}

John Conway in $\cite{l22}$ (1979) presented a new and the best construction of a proper class of numbers which contains all real numbers and all ordinal numbers as well as many other   {\it great} and {\it small} numbers like magnitudes  $\omega-1$, $\frac{\omega}{2}$, $\sqrt{\omega}$, $\sqrt[3]{\omega}$, ... $\omega^{\frac{1}{\omega}}$, $\omega^{\sqrt 2}$, etc.  and $\frac{1}{\omega-1}$, $\frac{1}{\omega^2}$, $\omega^{-\sqrt{\omega}}$, $\omega^{-\frac{1}{\omega}}$, etc., respectively. 

\smallskip

{\bf 1.} Conway construction of numbers was a modification of two well-known ideas. 
One of them was the Mirimanoff's representation $\cite{l2}$ (1917), later in 1923 repeated by von Neumann $\cite{l090}$, of the ordinal numbers: 
$0=\{\}$, $1=\{\{\}\}=\{0\}$, $2=\{\{\},\{\{\}\}\}=\{0,\{0\}\}$, $3=\{\{\},\{\{\}\},\{\{\},\{\{\}\}\}\}=\{0,1,\{0,1\}\}=\{0,1,2\}$, 
..., $\omega=\{0,1,2,...,n,...\}$, $\omega+1=\{0,1,2,...,n,...,\omega\}$, ... and so on. 

Another one  was the construction of real numbers ${\bf R}$ $\cite{l4}$ (1888) via Dedekind sections $\{A\,|\,B\}$ of  rational numbers  ${\bf Q}$, i.e., $\xi=\{A\,|\,B\}\in{\bf R}$ if and only if  $A$ and $B$ are non-empty sets of rationals such that ${\bf Q}=A\cup B$, $A\cap B=\emptyset$, $A$ has no greatest number, and for every $a\in A$ and $b\in B$ one has $a<b$.

Conway definition of a number is the following.

If $L$ and $R$ are any two {\it sets} of numbers, and {\it no member of} $L$ {\it is} $\geq$ {\it any member of} $R$, then there is a number $\{L\,|\,R\}$.

If $x=\{L\,|\,R\}$ is a number, then, for short, $x=\{x^L\,|\,x^R\}$, where   $x^L$ is a typical member of $L$ and $x^R$ is a typical member of $R$, i.e., $L=\{x^L\}$ and $R=\{x^R\}$ are {\it sets}.

Further,

$x\geq y$ iff (no $x^R\leq y$ and $x\leq$ no $y^L$); 

$x\leq y$ iff $y\geq x$; as well as 

$x=y$ iff ($x\geq y$ and $y\geq x$); 

$x>y$ iff ($x\geq y$ and $y\not\geq x$); 

$x<y$ iff $y>x$.

 Different representations $\{L'\,|\,R'\}$ and $\{L\,|\,R\}$ can define the same number; that is why one must distinguish between the form $\{L\,|\,R\}$ of a number and the number itself. Conway called $L$ the lower class and $R$ the upper class of number $x=\{L\,|\,R\}$.

Operations $+$, $\cdot$ and $/$ on all Conway numbers $x,y\in{\bf No}$ are the following:

\begin{equation}
\label{f01}
x+y=\{x^L+y,x+y^L\,|\,x^R+y,x+y^R\},
\end{equation}
\begin{equation}
\label{f02}
-x=\{-x^R\,|\,-x^L\},
\end{equation}
\begin{equation}
\label{f03}
xy=\{x^Ly+xy^L-x^Ly^L,x^Ry+xy^R-x^Ry^R\,|\,x^Ly+xy^R-x^Ly^R,x^Ry+xy^L-x^Ry^L\},
\end{equation}
and when $x>0$ the inverse to it, $xy=1$, is given by the following formula:
\begin{equation}
\label{f04}
y=\{0,\frac{1+(x^R-x)y^L}{x^R},\frac{1+(x^L-x)y^R}{x^L}\,|\,\frac{1+(x^L-x)y^L}{x^L},\frac{1+(x^R-x)y^R}{x^R}\}.
\end{equation}

The formula $(\ref{f04})$ is an inductive definition of an inverse number $y=\{y^L\,|\,y^R\}$ to the number $x=\{x^L\,|\,x^R\}$, where naturally $x\not=0$ and $x$ is positive number; when $x<0$, then $y=-y'$, where $y'$ is the inverse to the number $-x$.

The class ${\bf No}$ of all Conway  is a linearly ordered field and is  constructed by transfinite induction by the \grqq birthdays\grqq\, of numbers  by ordinal numbers, which themselves are born on the same days:

 $0=\{\,|\,\}$ (born on day $0$),

$1=\{0\,|\,\}$ and $-1=\{\,|\,0\}$ (born on day $1$),

$2=\{0,1\,|\,\}$, $\frac{1}{2}=\{0\,|\,1\}$, $-\frac{1}{2}=\{-1\,|\,0\}$, $-2=\{\,|\,-1,0,\}$ (born on day $2$),

 ...

$\omega=\{0,1,2,3,...\,|\,\}$, $\pi$, $e$, $\sqrt{2}$, $1+\frac{1}{\omega}$, $1-\frac{1}{\omega}$,
$1/3=\{0,\frac{1}{4},\frac{5}{16},...\,|\,\frac{1}{2},\frac{3}{8},...\}$,... $1/\omega=\{0\,|\,1,\frac{1}{2},\frac{1}{4},\frac{1}{8},...\}$, $-\omega=\{\,|\,0,-1,-2,-3,...\}$ (born on day $\omega$),

$\omega+1=\{0,1,2,3,...\omega\,|\,\}$, $\omega-1=\{0,1,2,3,...\,|\,\omega\}$, $\sqrt{2}+\frac{1}{\omega}$, $\sqrt{2}-\frac{1}{\omega}$, 
... $-\omega-1$ (born on day $\omega+1$),

 ...

$2\omega =\{0,1,2,3,...\omega,\omega+1,...\,|\,\}$, $\frac{\omega}{2}=\{0,1,2,3,...\,|\,\omega,\omega-1,\omega-2,...\}$, $\frac{2}{\omega}=\{\frac{1}{\omega}\,|\,1,\frac{1}{2},\frac{1}{4},...\}$, $\frac{1}{2\omega}=\{0\,|\,\frac{1}{\omega}\}$,... $\frac{1}{\omega^2}=\{0\,|\,\frac{1}{\omega},\frac{1}{2\omega},\frac{1}{4\omega},...\}$,
...$-\omega 2=\{\,|\,0,-1,-2,=3,...=\omega,=\omega-1,...\}$ (born on day $2\omega$),
... and so on.

More precisely,
for each ordinal number $\alpha\in{\bf On}$ Conway defined a set $M_\alpha$ of numbers by setting $x=\{x^L\,|\,x^R\}$ in $M_\alpha$ if all the $x^L$ and $x^R$ are in the union of all the $M_\beta$ for $\beta<\alpha$. Then he putted $O_\alpha=\bigcup\limits_{\beta<\alpha}M_\beta$ and $N_\alpha=M_\alpha\setminus O_\alpha$. Then in the terminology of numbers' birth to which he  adhered:

$M_\alpha$ is the set of numbers born on or before $\alpha$ (Made numbers),

$N_\alpha$ is the set of numbers born first on day $\alpha$ (New numbers), and

$O_\alpha$ is the set of numbers born before day $\alpha$ (Old numbers).

Each $x\in N_\alpha$ defines a Dedekind section $\{L|R\}$ of $O_\alpha$, if one sets $L=\{y\in O_\alpha:y<x\}$ and $R=\{y\in O_\alpha:y>x\}$, then $x=\{L\,|\,R\}$. Moreover, $M_\alpha=O_\alpha\cup N_\alpha$ as the union of $O_\alpha$ together with all its sections, in the natural order. Notice also that for each number $x\in M_\alpha$ there are inequalities inequality $-\alpha\leq x\leq\alpha$ and only two numbers $-\alpha$ and $\alpha$ in $M_\alpha$ have forms $\{\,|O_\alpha\}$ and $\{O_\alpha\,|\,\}$, respectively, which dropped the Dedekind restriction on sets $L$ and $R$  to be nonempty.

Every number $x$ is in a unique set $N_\alpha$ (see $\cite{l22}$ p. 30). Taking into account the above description, we call a {\it birthday form} $\{x^L\,|\,x^R\}$ of $x$ when the birthdays of all $x^L$, $x^R$ are less than $\alpha$. 

The set ${\bf R}$ of real numbers in ${\bf No}$ is defined as following: $x\in{\bf R}$ if $x=\{x-(\frac{1}{n})\,|\,x+(\frac{1}{n})\}_{n>0}$ and the proper class ${\bf On}$ of ordinal numbers in ${\bf No}$ is defined as following: $\alpha\in{\bf No}$ if $\alpha=\{ordinals\,\, \beta<\alpha\,|\,\}$.

 Conway defined a natural function  $y=\omega^x$ for all numbers $x\in{\bf No}$, i.e.,  the $x$th power of $\omega$,  by the following formula:
\begin{equation}
\label{f314}
\omega^x=\{0,r\omega^{x^L}\,|\,r\omega^{x^R}\},
\end{equation}
  where $x=\{x^L\,|\,x^R\}$ and $r$ denotes a variable ranging over all positive real numbers.

 Conway proved (Theorem 21 in $\cite{l22}$, p. 33) that each number $x\in{\bf No}$ defines a unique expression (the normal form of $x$)
\begin{equation}
\label{f0201}
x=\sum\limits_{0\leq\beta<\alpha}\omega^{y_\beta}r_\beta,
\end{equation}
in which $\alpha$ denotes some ordinal, the numbers $r_\beta$ $(0\leq\beta<\alpha)$ are non-zero reals, and the numbers $y_\beta$ form a descending sequence of numbers. Moreover, normal forms for distinct $y$ are distinct, and every form satisfying these conditions occurs.

Every odd-degree polynomial $f(x)=x^n+Ax^{n-1}+Bx^{n-2}+...+K$ with coefficients $A,B,...,K$ in ${\bf No}$ has a root in ${\bf No}$.

 Every positive number $x$ in ${\bf No}$ has a unique positive $n$th root in ${\bf No}$, for each positive integer $n$.

 The  Ring ${\bf No}[{\bf i}]$ of all numbers of the form $x+{\bf i}y$ ($x,y\in{\bf No}$), ${\bf i}^2=-1$, is an algebraically closed Field.

Proofs of these facts  see in $\cite{l22}$, p. 31-33, 40-42.

\smallskip

{\bf 2.} The theory of infinite sums  enables Conway to do quite o lot of classical analysis in ${\bf No}$ as well as in ${\bf No}[{\bf i}]$ but only for all numbers $x\in{\bf No}$ whose 
normal form $(\ref{f0201})$ has $y_0\leq 0$ (so called {\it finite} numbers), e.g., sines $y=\sin x$ 
and cosines $y=\cos x$,  $y=\tan x$ and  $y=\cot x$,  if in latter cases $\cos x\not=0$ and $\sin x\not=0$, respectively,  $y=exp\, x=e^x$ and $y=\ln x=\log_e x$.

Moreover, one can extend the functions $y=a^x$ and $y=\log_ax$ with any fixed  finite number $a>1$ (or $0<a<1$ such that $\frac{1}{a}$ is a finite number) for any finite number $x\in{\bf No}$ by formulas $y=a^x=e^{x\cdot\ln a}$ and $y=\log_a x=\frac{\ln x}{\ln a}$, respectively, not only for finite numbers $x\in{\bf No}$ but also for some infinite number $x\in{\bf No}$, more precisely, for each number $x\in domain(f)$, $y=f(x)$. 

For example, Conway defined trigonometry functions $y=\sin x$, $y=\cos x$ for all finite numbers but 
one can extend them to all infinite numbers in ${\bf No}$ in the following way: for each 
 infinite number $x\in{\bf No}$ there is as above a decomposition $x=x'+x''$, where $x'$ is 
 an infinite part of $x$ and $x''$ is a finite part of $x$ and putting $\sin x=\sin x''$, $\cos x=\cos x''$, respectively. As well as for $y=\tan x=\tan x''$, and  $y=\cot x=\cot x''$,  if in latter cases $\cos x\not=0$ and $\sin x\not=0$, respectively.

Conway defined $y=exp\, x=e^x$ and $y=\ln x=\log_e x$ for all finite numbers and all positive finite numbers except infinitesimal but one can extend them to all infinite numbers in ${\bf No}$ and some positive finite numbers in the following way: for each 
 infinite number $x\in{\bf No}$ there is as above a decomposition $x=x'+x''$, where $x'$ is 
 an infinite part of $x$ and $x''$ is a finite part of $x$ and putting $e^x=e^{x'+x''}=e^{x'}\cdot e^{x''}\stackrel{def}{=}\omega^{\frac{x'}{\omega}}\cdot e^{x''}$, where $\omega^{\frac{x'}{\omega}}$ is defined by formula $(\ref{f314})$, and $\ln x=y$, where $x\in Im(e^y), y\in{\bf No}$.

Notice here that there is a natural homomorphism: $exp(x):{\bf No}\rightarrow {\bf S}$

\begin{equation}
 \label{f1507}
 exp(x)=\cos{2\pi x}+i\sin{2\pi x},\,\,\, x\in{\bf No},
 \end{equation}
 where ${\bf S}=\{z\in {\bf No}[{\bf i}]\,|\, z=a+{\bf i}b\,\&\,a^2+b^2=1\}$.

\smallskip

{\bf 3.}
Let $\zeta=\omega^{\omega^\kappa}$, $\kappa\geq0$, be the {\it main ordinal number} in the sense of Jacobstahl $\cite{l15}$, i.e., for each ordinal numbers $\alpha, \beta<\zeta$ their product $\alpha\cdot\beta<\zeta$. One can easily see that all Conway numbers $x\in O_\zeta$, i.e., born before day $\zeta$, is a subring of the Field ${\bf No}$. Consider a linearly topology on $O_\zeta$ generated by  the family ${\cal B}$ of all intervals $(a,b)\subset O_\zeta\,|\,a,b\in O_\zeta$.

By $\zeta$-sequence in $O_\zeta$ we understand a mapping $x:(0,\zeta)\rightarrow O_\zeta$, where $(0,\zeta)=\{0<\alpha <\zeta\}$ is an interval of ordinal numbers. Denote it by $x=(x_\alpha)_{0<\alpha<\zeta}$.
 
We say that $\zeta$-sequence $(x_\alpha)_{0<\alpha<\zeta}$ in  $O_\zeta$ {\it converges} to $a\in O_\zeta$, and we write $\lim\limits_{0<\alpha<\zeta}x_\alpha=a$,  if for each positive  number $\varepsilon\in O_\zeta$ there is an ordinal number $\alpha_0\in(0,\zeta)$ such that $|x_\alpha-a|<\varepsilon$ for all $\alpha_0\leq\alpha<\zeta$.  
In this case we also say that $\zeta$-sequence $(x_\alpha)_{0<\alpha<\zeta}$ is {\it convergent} in  $O_\zeta$.

One can see that if   $\lim\limits_{0<\alpha<\zeta}x_\alpha=a$ and $\lim\limits_{0<\alpha<\zeta}y_\alpha=b$, then $\lim\limits_{0<\alpha<\zeta}(x_\alpha+y_\alpha)=a+b$, $\lim\limits_{0<\alpha<\zeta}(x_\alpha\cdot y_\alpha)=a\cdot b$ and $\lim\limits_{0<\alpha<\zeta}(\frac{x_\alpha}{y_\alpha})=\frac{a}{b}$, when in latter case $b\not=0$.

 A $\zeta$-sequence $(x_\alpha)_{0\leq\alpha<\zeta}$ in  $O_\zeta$ is called  {\it fundamental}  f for each positive  number $\varepsilon\in O_\zeta$  there is an ordinal number $\alpha_0$ such that $|x_\alpha-x_{\alpha'}|<\varepsilon$, for all $\alpha_0\leq\alpha<\alpha'<\zeta$.

Two $\zeta$-fundamental sequences $(x_\alpha)_{0<\alpha<\zeta}$ and $(y_\alpha)_{0<\alpha<\zeta}$ in  $O_\zeta$ are $\zeta$-{\it equivalent}, denoted by $(x_\alpha)_{0<\alpha<\zeta}\sim(y_{\alpha})_{0<\alpha<\zeta}$,  if for each  positive number $\varepsilon\in O_\zeta$  there are ordinal numbers $\alpha_0$ and $\alpha'_0$ such that $|x_\alpha-y_{\alpha'}|<\varepsilon$, for all $\alpha_0\leq\alpha<\zeta$ and all $\alpha'_0\leq\alpha'<\zeta$. 

Since $(x_\alpha)_{0<\alpha<\zeta}\sim(y_{\alpha})_{0<\alpha<\zeta}$ is an equivalent relation we can consider the set of all classes of equivalent $\zeta$-sequences in $O_\zeta$. Denote it by ${\bf R}_\zeta$.

One can see that the remainder ${\bf R}_\zeta\setminus O_\zeta$ are Conway numbers (more precisely can be identified with Conway numbers) are some but not all numbers born on day $\zeta$.

One can see that ${\bf R}_\zeta$ is a subfield of ${\bf No}$, it is a complete field, i.e., each $\zeta$-fundamental sequence $x=(x_\alpha)_{0<\alpha<\zeta}$ in  ${\bf No}$ converges to a number $x_0\in {\bf No}$, and thus all above results are valid including the last result that ${\bf C}_\zeta={\bf R}_\zeta[{\bf i}]$ is  an algebraically closed field. 

One can distinguish a continuous function $y=f(x)$ in a domain $X\subseteq{\bf R}_\zeta$ of $f$, i.e., when for all $x_0\in X$ and for all $\zeta$-sequences $x=(x_\alpha)_{(0<\alpha<\zeta)}$ the formula $\lim\limits_{0<\alpha<\zeta}x_\alpha=x_0$ implies $\lim\limits_{0<\alpha<\zeta}f(x_\alpha)=f(x_0)$ and a discontinuous function when for at least one point $x_0\in X$ there is no such implication.

Moreover, all mentioned above functions $y=x^n+Ax^{n-1}+Bx^{n-2}+...+K$, $A,B,...,K\in {\bf R}_\zeta$, $y=\sin x$, $y=\cos x$, $y=\tan x$ are continuous in their domains. As to $y=a^x$ and $y=\log_a x$ ($a$ is a positive finite number which is not a infinitesimal and not equal 1), $y=a^x$ is continuous at each  number $x\in{\bf L}_\zeta$, where ${\bf L}_\zeta=\{x\}_{x\in{\bf R_\zeta}\,\&\,a^x\in{\bf R}_\zeta}$; 
$y=\log_a x$ is continuous at each number $x\in{\bf D}_\zeta$, where  ${\bf D}_\zeta=\{a^x\}_{x\in{\bf L}_\zeta}$. 

Note  that ${\bf L}_\zeta={\bf R}_\zeta$ if and only if  $\zeta=\omega$ or $\zeta$ is an $\varepsilon$-number $\varepsilon_\alpha$, $0\leq\alpha$, or ${\bf L}={\bf No}$. In the end,  $y=a^x$ is continuous at each number  $x\in{\bf L}={\bf No}$ and $y=\ln x$ is continuous at each number $x\in{\bf D}=\{x\}_{a^x\in{\bf No}}$.

A natural homomorphism: $exp(x):{\bf R}_\zeta\rightarrow {\bf S}_\zeta$
 \begin{equation}
 \label{f2507}
 exp(x)=\cos{2\pi x}+i\sin{2\pi x},\,\,\, x\in{\bf R}_\zeta,
 \end{equation}
 where ${\bf S}_\zeta=\{z\in {\bf R}_\zeta[{\bf i}]\,|\, z=a+{\bf i}b\,\&\,a^2+b^2=1\}$, is also continuous.
See proofs and details in \cite{l1332}.

\bigskip

\begin{center}
{\bf Part One}
\end{center}

\begin{center}
{\bf 1. Definition of skands}
\end{center}

{\bf Definition 1.}  A {\it skand} with a {\it clutch  region} $[\alpha_0,\alpha)$ is called an object $X$,  which for each ordinal $\alpha'\in[\alpha_0,\alpha)$ uniquely associates the set $X(\alpha')\in{\bf V[{\cal U}]}$   and has the following {\it internal structure} of the membership relations $\in$ of its constituents and parts: for any ordinal $\alpha'\in[\alpha_0,\alpha)$ and for any element $x\in X(\alpha')$ there is a relation $x\in X|_{[\alpha',\alpha)}$; and for each ordinal $\alpha'\in[\alpha_0,\alpha)$, such that $\alpha'+1\in[\alpha_0,\alpha)$, there is a relation $X|_{[\alpha'+1,\alpha)}\in  X|_{[\alpha',\alpha)}$, where $X|_{[\alpha',\alpha)}$ and $X|_{[\alpha'+1,\alpha)}$ are restrictions of $X$ on $[\alpha',\alpha)$ and $[\alpha'+1,\alpha)$, respectively.

Here for each ordinal $\alpha'\in[\alpha_0,\alpha)$ via $X|_{[\alpha',\alpha)}$ we denote the skand, which is a  constituent of $X$, with a {\it clutch  region}  $[\alpha',\alpha)$ and $X|_{[\alpha',\alpha)}(\alpha'')=X(\alpha'')$, for any ordinal $\alpha''\in[\alpha',\alpha)$. In particular $X|_{[\alpha',\alpha)}$ coincides with $X$ for $\alpha'=\alpha_0$ and is the {\it remainder} of the  skand $X$ on  a  clutch  region $[\alpha',\alpha)$ for each $\alpha'\in(\alpha_0,\alpha)$. We also denote sets $X(\alpha')$ by $X_{\alpha'}$, $\alpha_0\leq\alpha'<\alpha$,  calling them {\it constituents} of the skand $X$ and skands $X|_{[\alpha',\alpha)}$, $\alpha_0\leq\alpha'<\alpha$,  calling them {\it parts} of the skand $X$.
It is also clear that all remainders $X|_{[\alpha',\alpha)}$, $\alpha'\in(\alpha_0,\alpha)$, of the skand $X$ are skands indeed with  clutch  regions $[\alpha',\alpha)$ and with induced relations of their constituents and parts that are portion of the relationship and portion of the constituents of the skand $X$ itself.

It is also clear that a skand $X$ with a   clutch  region $[\alpha_0,\alpha)$ uniquely defines a mapping $f:[\alpha_0,\alpha)\rightarrow{\bf V[{\cal U}]}$ if we put $f(\alpha')=X(\alpha')$ for each ordinal $\alpha'\in[\alpha_0,\alpha)$. This mapping $f=FX$, where $F$ is the forgetful functor of the membership relation $\in$ of the constituents and parts of the skand $X$, we will call {\it the mapping associated with the skand} $X$.

We call {\it the length} of a skand $X$ with a clutch  region $[\alpha_0,\alpha)$ an ordinal number $l=\alpha-\alpha_0$ , i.e.,  a single ordinal such that $\alpha=\alpha_0+l$ (here $+$ is an ordinary sum of ordinal numbers as types of well-ordered sets and itis different from adding of Conway numners), and we will distinguish here between {\it finite} skands whose length are equal to a natural numbers, i.e., the ordinals $l\in(0,\omega)$, and {\it infinite} skands whose length are $l\geq\omega$, where $\omega$ is the smallest countable ordinal number.

To define equality of skands, we need the notions of counting functions and equality of counting functions.

{\bf Definition 2.} Let $f:C'\rightarrow{\bf V}[{\cal U}]$ and $f': C'\rightarrow{\bf V}[{\cal U}]$ be arbitrary mappings of subclasses $C$ and $C'$ of the class ${\bf On}$ of all ordinal numbers in the universal class ${\bf V}[{\cal U}]$ of all sets of the theory $NBG[{\cal U}]$, which we call them {\it counting functions}. Let's say that these counting functions $f$ and $f'$ {\it are equal} if there exists a similar mapping $\varphi:C\rightarrow C'$, i.e., a bijection $\varphi$ preserving the structure of well-ordered sets $C$ and $C'$, which we will call an isomorphism for short for which the following usual equality of mappings  $f=f'\circ\varphi$ is satisfied. 

We denote this equality of the counting functions $f$ and $f'$ also by $f=f'$, but distinguish it from the usual equality $f=f'$ of the mappings $f$ and $ f'$ when their domains of definition are equal and for each element from this domain the values of the mappings on this element are equal. 

{\bf Proposition 2} {\it The relation of equality of mappings as counting functions is an equivalence relation. }

{\bf Proof.} The reflexivity of $f=f$ is obvious, since the identity mapping $\varphi=1|_C:C\rightarrow C$ is an identity isomorphism of a well-ordered class $C$ onto itself and the usual equality of mappings $f=f\circ 1|_C$ takes place. Hence, there is equality of the mapping as a counting function of itself, i.e., $f=f$. This is the rare case when the notions of equality of mappings as counting functions and the usual equality of mappings coincide with each other. 

Let's prove symmetry. Let the mappings $f:C\rightarrow{\bf V}[{\cal U}]$ and $f':C'\rightarrow{\bf V}[{\cal U}]$ be equal as counting functions, i.e., $f=f'$. Then, by Definition 2, the usual equality of mappings $f=f'\circ \varphi$ is satisfied, where $\varphi:C\rightarrow C'$ is an isomorphism of well-ordered classes. But since $\varphi^{-1}:C'\rightarrow C$ is an inverse isomorphism of well-ordered classes, the formula $f=f'\circ \varphi$ entails the following equality: $f\circ \varphi^{-1}=(f'\circ \varphi)\circ \varphi^{-1}=f'\circ (\varphi\circ \varphi^{-1})=f'\circ 1|_{C'}=f'$ of the corresponding mappings, and hence the equality $f'=f\circ\varphi^{-1}$ of the mappings is satisfied by virtue of the symmetry of the equality relation of the mappings.  Hence, the equality $f'=f$ of mappings as counting functions is satisfied. 

If $f'':C''\rightarrow{\bf V}[{\cal U}]$ is another mapping of the class $C''\subseteq{\bf On}$ into the universal class ${\bf V}[{\cal U}]$, then the transitivity $(f=f')\&(f'=f'')\Rightarrow(f=f'')$ of the equality of mappings as counting functions is also obvious. If the premise, by Definition 2, entails the usual equalities $f=f'\circ\varphi$ and $f'=f''\circ\psi$ of the corresponding mappings, where $\psi:C'\rightarrow C'$ is an isomorphism of the well-ordered classes $C'$ and $C''$, then supposing $\chi=\psi\circ\varphi$, we obtain the isomorphism $\chi: C\rightarrow C''$ of well-ordered classes $C$ and $C''$ and the following equality: $f=f'\circ\varphi=(f''\circ\psi)\circ\varphi=f''\circ(\psi\circ\varphi)=f''\circ\chi$ mappings, and hence the equality $f=f''\circ\chi$ mappings. Hence, by Definition 2, there is equality $f=f''$ of mappings as counting functions.

\bigskip

{\bf Definition 3.} Two skands $X$ and $Y$ with clutch  regions $[\alpha_0,\alpha)$ and $[\beta_0,\beta)$, respectively, are called {\it  equal} (we write $X=Y$) if the corresponding mappings $f$ and $g$ associated with skands $X$ and $Y$ are equal as counting functions, i.e.,   there exists a similar mapping $\varphi:[{\alpha_0,\alpha})\rightarrow[{\beta_0,\beta})$ of the intervals $[{\alpha_0,\alpha})$ and $[{\beta_0,\beta})$ that preserves the structures of the well-ordered sets $[{\alpha_0,\alpha})$ and $[{\beta_0,\beta})$,  i.e., the existence of an isomorphism $\varphi$ of well-ordered sets and the following usual equality of mappings $f=g\circ\varphi$.

 \bigskip
 
The description of skands and their equalities can be given less formally, but more clearly and conveniently for their recording and working with them. Let's call this form of a skand {\it unfolded}, in the sense that in the record of the skand itself the region of its coupling will be given and all its elements, or more precisely, all the relations  of membership $\in$ of its constituents and parts will already be specified.

Let's first consider \grqq strictly decreasing\grqq\, sequences of nested curly braces $\{_{\alpha_0}\{_{\alpha_0+1}...\{_{\alpha'}...\,\,\,\,...\}_{\alpha'}...\}_{\alpha_0+1}\}_{\alpha_0}$; strictly decreasing in the sense that for any ordinals $\alpha'<\alpha''$ such that $\alpha',\alpha''\in[\alpha_0,\alpha)$, the pair of brackets $\{_{\alpha''}...\}_{\alpha''}$\grqq lies inside\grqq\, a pair of brackets $\{_{\alpha'}...\}_{\alpha'}$. 

The possibility of the existence of such a strictly decreasing sequence of nested curly braces can \grqq be modeled\grqq\, on the sequence of ends of segments $[-\frac{1}{\alpha'},\frac{1}{\alpha'}]$, $\alpha'\in[\alpha_0,\alpha)$ of Conway numbers \cite{l22}, where the Conway number $-\frac{1}{\alpha'}$ corresponds to the open  curly brace $\{_{\alpha'}$, and the Conway number $\frac{1}{\alpha'}$ corresponds to the closed  curly brace $\}_{\alpha'}$, if of course, $\alpha_0\not=0$. In this very case we put  the Conway number $-2$ corresponds to the open  curly brace $\{_{\alpha_0}$, and the Conway number $2$ corresponds to the closed  curly brace $\}_{\alpha_0}$. And since for any ordinals $\alpha_0\leq\alpha'<\alpha''<\alpha$, the closed segment $[-\frac{1}{\alpha''},\frac{1}{\alpha''}]$ of all Conway numbers lies inside the segment  $[-\frac{1}{\alpha'},\frac{1}{\alpha'}]$, then in the same sense the pair of  curly braces $\{_{\alpha''}...\}_{\alpha''}$ \grqq lies inside\grqq\, a pair of  curly brace $\{_{\alpha'}...\}_{\alpha'}$.

{\bf Remark 1}. Note only that unlike Conway numbers, which are mathematical objects, curly brackets are not at all, but play here a syntactical role, indicating only the membership relations $\in$ of  constituents and parts of the skand: everything between the corresponding pair of curly braces $\{_{\alpha'}\}_{\alpha'}$, where $\alpha'\in[\alpha_0,\alpha)$ will denote the elements of the set $X_{\alpha'}\cup\{ X|_{[\alpha'+1,\alpha)}\}$ (see below) generated by the skand $X|_{[\alpha',\alpha)}$.

We perform this replacement of Conway numbers with curly brackets as follows. Denote by $\psi$ a mapping $\psi:\{_{\alpha_0}\{_{\alpha_0+1}...\{_{\alpha'}...\,\,\,\,...\}_{\alpha'}...\}_{\alpha_0+1}\}_{\alpha_0}\rightarrow [-2,2]$ of a strictly decreasing sequence of nested curly braces to the close segment $[-2,2]$ of all Conway numbers, given by the two formulas $\psi(\{_{\alpha'})=-\frac{1}{\alpha'}$ and $\psi(\}_{\alpha'})=\frac{1}{\alpha'}$. It is evidently a linear ordered isomorphism between the set of all above nested curly braces to the set of corresponding Conway numbers in the closed segment $[-2,2]$.

{\bf Definition 1$'$.} A trivial skand ${\bf e}_{[\alpha_0,\alpha)}$ with a clutch  region $[\alpha_0,\alpha)$ is a strictly decreasing sequence of nested curly braces  $\{_{\alpha_0}\{_{\alpha_0+1}...\{_{\alpha'}...\,\,\,\,...\}_{\alpha'}...\}_{\alpha_0+1}\}_{\alpha_0}$, where $\alpha'\in[\alpha_0,\alpha)$, which we reduce for simplicity to  $\{\{...\{...\,\,\,...\}...\}\}$  by dropping the indexes. To construct a nontrivial skand $X_{[\alpha_0,\alpha)}$ with the same clutch  region $[\alpha_0,\alpha)$, we consider for each ordinal $\alpha'\in[\alpha_0,\alpha)$ an arbitrary set $X_{\alpha'}$ from the universe ${\bf V}[{\cal U}]$, i.e., $X_{\alpha'}\in{\bf V}[{\cal U}]$.  Since all Conway numbers in the segment $(-\frac{1}{\alpha'},-\frac{1}{\alpha'+1})$ is a proper class (not a set) there is an injective mapping $\psi_{\alpha'}:X_{\alpha'}\rightarrow(-\frac{1}{\alpha'},-\frac{1}{\alpha'+1})$ and thus the  injective mapping $\Psi:\bigsqcup\limits_{\alpha_0\leq\alpha'<\alpha}\psi_{\alpha'}\rightarrow(-2,2)$ of the set $\bigsqcup\limits_{\alpha_0\leq\alpha'<\alpha}\{X_{\alpha'}\}$ to the proper class $(-2,2)$ of Conway numbers. Then $X_{[\alpha_0,\alpha)}\stackrel{def}{=}(\psi\cup\Psi)^{-1}(Im(\psi\cup\Psi))$ with induce linear ordering by the linear ordering on $[-2,2]$ and the mapping $(\psi\circ\Psi)^{-1}$ is called {\it skand} with a clutch  region $[\alpha_0,\alpha).$

\parindent=0,5 cm
Sets $X_{\alpha'}=\{x^{\alpha'}_0,x^{\alpha'}_1,...,x^{\alpha'}_\lambda,...\}$, for each ordinal $\alpha'\in[\alpha_0,\alpha)$ and $0\leq\lambda<\Lambda(\alpha')$, we call $\alpha'$-{\it components} of the skand $X_{[\alpha_0,\alpha)}$ since, by Choice Axiom, we can turn into a well-ordered set $X_{\alpha'}$. By construction of the skand $X_{[\alpha_0,\alpha)}$, we understand it as a {\it set} whose elements are founded sets (individuals are not excluded when $X_{\alpha'}=\{x^{\alpha'}_\lambda\}$ is an individual $x^{\alpha'}_\lambda\in{\cal U}$)  $x^{\alpha_0}_0,x^{\alpha_0}_1,...,x^{\alpha_0}_\lambda,...$ of the universe ${\bf V}[{\cal U}]$  and one non-founded element $X_{[\alpha_0+1,\alpha)}$.

(Clearly, the order of the elements in each component plays no role and is given in Definition $1'$ only for simplicity of notation, since in the theory of $NBG[{\cal U}]$ there is a choice axiom ${\bf CA}$, and any component $X_{\alpha'}$ of the skand $X_{[\alpha_0, \alpha)}$ can be quite ordered by numbering its elements with indexes $\lambda=\lambda'\in[0,\Lambda')=[0,\Lambda_{\alpha'})$, where the well-ordered set $[0,\Lambda_{\alpha'})$ depends on $\alpha'\in[\alpha_0,\alpha)$ and on the very choice of order on $X_{\alpha'}$. )

In Definition  $1'$ we give the notation of the skand $X$ of Definition 1 with its clutch  region $[\alpha_0,\alpha)$  as an index, i.e., $X_{[\alpha_0,\alpha)}$. for the sake of completeness. Moreover, for each ordinal. $\alpha'\in(\alpha_0,\alpha)$ we will denote by $X_{[\alpha',\alpha)}$  the skand whose corresponding components $X_{\alpha''}$, $\alpha'\leq\alpha''<\alpha$, are the same as those of the skand  $X_{[\alpha_0,\alpha)}$.  Other skands different from these will be denoted by other letters, e.g., $Y_{[\beta_0,\beta)}$, where the intervals $[\beta_0,\beta)$ may or may not coincide with the intervals discussed above.

In other words, a skand $X_{[\alpha_0,\alpha)}$ is an {\it aggregate}, or a {\it  chained decreasing-tuple} of {\it founded} sets or individuals, or {\it  their absence}. A trivial skand is a chained decreasing-tuple with empty components; an arbitrary skand is a chained decreasing-tuple with random components that are elements of the universe ${\bf V}[{\cal U}]$ of the founded sets. 

(Skand in Sanskrit means \grqq random aggregate of being\grqq.)

Any skand $X_{[\alpha_0,\alpha)}$ with components $X_{\alpha'}$, $\alpha_0\leq\alpha'<\alpha$, Definitions $1'$ uniquely defines a skand $X$ with a clutch  region $[\alpha_0,\alpha)$ if we put for each ordinal $\alpha'\in[\alpha_0,\alpha)$ $X(\alpha')=X_{\alpha'}$  and  $x\in X_{\alpha'}\Rightarrow x\in X|_{[\alpha',\alpha)}$, and for any ordinal $\alpha'\in[\alpha_0,\alpha)$, such that $\alpha'+1\in[\alpha_0,\alpha)$, put $X_{[\alpha'+1,\alpha)}\in X_{[\alpha',\alpha)}\Rightarrow X|_{[\alpha'+1,\alpha)}\in X|_{[\alpha',\alpha)}$.

 The converse is also true: every skand $X$ with a clutch  region $[\alpha_0,\alpha)$ of Definition 1 uniquely defines a skand  $X_{[\alpha_0,\alpha)}$ of Definition $1'$ if for every $\alpha'\in[\alpha_0,\alpha)$ put $X_{\alpha'}=X(\alpha')$ and curly brackets  play  a syntactical role, indicating  the membership relations $\in$ of  constituents and parts of the skand $X_{[\alpha_0,\alpha)}$.

Thus, the skands $X$ and $X_{[\alpha_0,\alpha)}$ in both Definitions are one and the same object, which uniquely matches for each ordinal $\alpha'\in[\alpha_0,\alpha)$ a set-skand $X_{[\alpha',\alpha)}$ with the same  membership relations $\in$ of its constituents and parts, though represented in fact in different forms -- formal and extended (similarly to how a function can be represented in different forms: A law that maps a single value to each argument; its graph; an implicit equation; a table; a two-line matrix; etc.). We will see that skands are real and natural mathematical objects.

The general form of the skand $X$ in the unfolded form is as follows: 
\begin{equation}
\label{f023}
\begin{array}{l}
X_{[\alpha_0,\alpha)}
=\\=\{x^{\alpha_0}_0,x^{\alpha_0}_1,...,\{x^{\alpha_0+1}_0,x^{\alpha_0+1}_1,...,\{...\{x^{\alpha'}_0,x^{\alpha'}_1,...,\{...\{x^{\alpha-1}_0,x^{\alpha-1}_1,...\}...\}...\}...\}\}\},
\end{array}
\end{equation}
where the components of $X_{\alpha'}$ are either empty, i.e., $X_{\alpha'}=\emptyset$, or consist of elements $x^{\alpha'}_0,x^{\alpha'}_1,...$ of class ${\bf V}[{\cal U}]$ such that their set $\{x^{\alpha'}_0,x^{\alpha'}_1,...\}=X_{\alpha'}$ is an element of  founded set in the class ${\bf V[{\cal U}]}$, i.e., $X_{\alpha'}\in {\bf V[{\cal U}]}$, $\alpha_0\leq\alpha'<\alpha$. Thus the general formula $(\ref{f023})$ represents the case where the ordinal $\alpha$ has a predecessor. If $\alpha$ is a limit ordinal, then the last pair of brackets $\{\}$ in the formula $(\ref{f023})$, together with the elements contained there or not, is missing.

{\bf Definition 2$'$}. Two skands $X_{[\alpha_0,\alpha)}$ and $Y_{[\beta_0,\beta)}$ are called {\it  equal} when there exists a similar mapping (an isomorphism) $\varphi:[{\alpha_0,\alpha})\rightarrow[{\beta_0,\beta})$  of the intervals $[{\alpha_0,\alpha})$ and $[{\beta_0,\beta})$ that preserves the structures of the well-ordered sets $[{\alpha_0,\alpha})$ and $[{\beta_0,\beta})$, and the corresponding components are equal as sets, i.e., for any ordinals $\alpha'\in[\alpha_0,\alpha)$ and $\beta'\in[\beta_0,\beta)$ such that $\beta'=\varphi(\alpha')$, the sets are equal:
\begin{equation}
\label{f777}
X_{\alpha'}=\{x^{\alpha'}_0,x^{\alpha'}_1,x^{\alpha'}_2,...,x^{\alpha'}_\lambda,...\}=\{y^{\beta'}_0,y^{\beta'}_1,y^{\beta'}_2,...,y^{\beta'}_\mu,...\}=Y_{\beta'}.
\end{equation}

\bigskip

A real nonstandard model of the skand $X$ of Definition $1$ can be constructed as follows. Consider the non-similar class  ${\cal U}\cup{\bf V}[{\cal U}]\cup\bigcup\limits_{\alpha'<\alpha'',\alpha',\alpha''\in{\bf On}}{\bf V}[{\cal U}]^{[\alpha',\alpha'')}$,where ${\bf V}[{\cal U}]^{[\alpha',\alpha'')}$ denotes the non-similar class of all mappings $f_{\alpha'<\alpha''}$ of an interval $[\alpha',\alpha')$ in ${\bf V}[{\cal U}]$. 

Now for any fixed skand $X$ with a clutch  region $[\alpha_0,\alpha)$, denote by $M$ the set $\bigcup\limits_{\alpha'\in[\alpha_0,\alpha)}X(\alpha')\cup\bigcup\limits_{\alpha'\in[\alpha_0,\alpha)}{\bf V}[{\cal U}]^{[\alpha',\alpha)}$.

Next, on the Cartesian product $M\times M$ we define the following membership relation $R\subset M\times M$: 
\begin{equation}
\label{f0607}
(x,y)\in R\Leftrightarrow (x\in X(\alpha')\wedge y=f_{\alpha'<\alpha})\vee(x=f_{\alpha'+1<\alpha}\wedge y=f_{\alpha'<\alpha}).
\end{equation}

Comparing now to each constituent of $x^{\alpha'}_{\lambda}$, $\alpha'\in[\alpha_0,\alpha)$, $\lambda\in[0,\Lambda({\alpha'}))$, of the skand $X$  the same element $x^{\alpha'}_{\lambda}\in M$, to each constituent of $X|_{[\alpha',\alpha)}$, $\alpha'\in[\alpha_0,\alpha)$of the skand $X$ an element $f_{\alpha'<\alpha}\in M$, we obtain a correspondence of the membership relation $\in$ of the constituents, or parts of the skand $X$ with the relation $R\subset M\times M$ on the set $M$.

\bigskip

{\begin{center}
{\bf 2. Reflexive skands}
\end{center}

{\bf Definition 4} A skand $X_{[\alpha_0,\alpha)}$ is called {\it  reflexive} if the skand $X_{[\alpha_0+1,\alpha)}$ is equal to the skand $X_{[\alpha_0,\alpha)}$, i.e., $X_{[\alpha_0+1,\alpha)}=X_{[\alpha_0,\alpha)}$.

{\bf Proposition 3}. {\it A skand $X_{[\alpha_0,\alpha)}$ is  reflexive if and only if one of the following two conditions is satisfied: 

$1)$ $\alpha-\alpha_0\geq\omega$ and for any ordinal $\alpha'$ such that $\alpha_0<\alpha'<\alpha_0+\omega$, there is an equality $X_{\alpha'}=X_{\alpha_0}$, or, which is the same thing, the mapping $f|_{[\alpha_0,\alpha_0+\omega)}$, which is a restriction to the interval $[\alpha_0,\alpha_0+\omega)$ of the mapping $f: {[\alpha_0,\alpha)}\rightarrow{\bf V}[{\cal U}]$ associated with the skand $X_{[\alpha_0,\alpha)}$, is a constant mapping whose image $f([\alpha_0,\alpha_0+\omega))$ is the set $\{X_{\alpha_0}\}\subset{\bf V}[{\cal U}]$; 

$2)$ for any ordinal $\alpha'$ such that $\alpha_0<\alpha'<\alpha_0+\omega$, the equality of the skands $X_{[\alpha',\alpha)}=X_{[\alpha_0,\alpha)}$.}

{\bf Proof}. {\bf Necessity}. Indeed, let the skand $X_{[\alpha_0,\alpha)}$ be reflexive. Then, by Definition 4, $X_{[\alpha_0+1,\alpha)}=X_{[\alpha_0,\alpha)}$.

Hence, by Definition 2$'$, first, the intervals $[\alpha_0+1,\alpha)$ and $[\alpha_0,\alpha)$ are similar (isomorphic), and hence $\alpha-\alpha_0\geq\omega$. Otherwise, such an isomorphism is impossible, since any finite interval $[\alpha_0,\alpha)$ of ordinal numbers represents the ordinal type of a finite ordinal number $l=\alpha-\alpha_0$ different from the finite ordinal number $l-1$ represented by the ordinal type of the interval $[\alpha_0+1,\alpha)$.

Second, there are the following equality of sets: $X_{\alpha_0+1}=X_{\alpha_0}$, $X_{\alpha_0+2}=X_{\alpha_0+1}$,...,$X_{\alpha'+1}=X_{\alpha'}$,... for any ordinal $\alpha'\in[\alpha_0,\alpha_0+\omega)$.  By the transitivity of the set equality relation, we obtain the following formula: $X_{\alpha'}=X_{\alpha_0}$ for all ordinals $\alpha'$ such that $\alpha_0<\alpha'<\alpha_0+\omega$.

It is also obvious that the restriction $f|_{[\alpha_0,\alpha_0+\omega)}$ of the mapping $f:{[\alpha_0,\alpha)}\rightarrow{\bf V}[{\cal U}]$ associated with the reflexive skand $X_{[\alpha_0,\alpha)}$, on the set $[\alpha_0,\alpha_0+\omega)$, is a constant mapping, whose image is the set $\{X_{\alpha_0}\}\subset{\bf V}[{\cal U}]$, or equivalently, the value $f(\alpha')=X_{\alpha_0}$ for any ordinal $\alpha'\in[\alpha_0,\alpha_0+\omega)$.

The converse is also true, i.e., constancy of the restriction $f|_{[\alpha_0,\alpha_0+\omega)}$ of the mapping $f:{[\alpha_0,\alpha)}\rightarrow{\bf V}[{\cal U}]$ associated with the reflexive skand  $X_{[\alpha_0,\alpha)}$, on the set $[\alpha_0,\alpha_0+\omega)$, whose image is the set $\{X_{\alpha_0}\}\subset{\bf V}[{\cal U}]$, entails equality of the sets $X_{\alpha'}=X_{\alpha_0}$ for any $\alpha_0<\alpha'<\alpha_0+\omega$. Thus the necessity of condition $1)$ is proved.

Let's now show that condition $1)$ entails the equality $X_{[\alpha',\alpha)}=X_{[\alpha_0,\alpha)}$ for any ordinal $\alpha'$ such that $\alpha_0<\alpha'<\alpha_0+\omega$. Indeed, we associate for each ordinal $\alpha''\in[\alpha',\alpha_0+\omega)$,  an ordinal $\alpha''-1\in[\alpha_0,\alpha)$, and for each ordinal $\alpha''\in[\alpha_0+\omega,\alpha)$ we associate this ordinal $\alpha''\in[\alpha_0,\alpha)$ itself and obtain an isomorphism of well-ordered sets $[\alpha',\alpha)$ and  $[\alpha_0,\alpha)$. Thus the corresponding components of $X_{\alpha''}$ of skands $X_{[\alpha',\alpha)}$ and $X_{[\alpha_0,\alpha)}$ skands are equal to each other: for ordinals $\alpha''\in[\alpha',\alpha_0+\omega)$ and $\alpha''\in[\alpha_0,\alpha_0+\omega)$, respectively, they are equal to the set $X_{\alpha_0}$, and for each ordinal $\alpha''\in[\alpha_0+\omega,\alpha)\subset[\alpha',\alpha)$ and $\alpha''\in[\alpha_0+\omega,\alpha)\subset[\alpha_0,\alpha)$, respectively, the components $X_{\alpha''}$ and $X_{\alpha''}$ of the skands $X_{[\alpha',\alpha)}$ and $X_{[\alpha',\alpha)}$, respectively, are equal to each other, i.e., $X_{\alpha''}=X_{\alpha''}$.  Thus the necessity of the condition $2)$ is proved.

\bf Sufficiency}. If condition $1)$ of Proposition 3 is satisfied, then, first, $\alpha'-\alpha_0\geq\omega$, and hence there is an isomorphism $\varphi$ of well-ordered sets $[\alpha_0+1,\alpha)$ and $[\alpha_0,\alpha)$, we associate for each ordinal $\alpha'\in[\alpha_0+1, \alpha_0+\omega)$, the ordinal $\varphi(\alpha')=\alpha'-1\in[\alpha_0,\alpha_0+\omega)$, and for each ordinal $\alpha'\in[\alpha_0+\omega,\alpha)$ we associate the same ordinal $\varphi(\alpha')=\alpha'\in[\alpha_0+\omega,\alpha)$, what is the desired isomorphism of $\varphi: [\alpha_0+1,\alpha)\rightarrow[\alpha_0,\alpha)$ since $[\alpha_0+1,\alpha)=[\alpha_0+1, \alpha_0+\omega)\sqcup[\alpha_0+\omega,\alpha)$ and $[\alpha_0,\alpha)=[\alpha_0,\alpha_0+\omega)\sqcup[\alpha_0+\omega,\alpha)$, and the mapping $\varphi$ defined above preserves in an obvious way the structure of the well-ordered sets $[\alpha_0+1,\alpha)$ and $[\alpha_0,\alpha)$.

Second, condition $1)$ entails the equality $X_{\alpha'}=X_{\alpha_0}$ for any ordinal $\alpha'\in[\alpha_0+1,\alpha_0+\omega)$, and hence, the equality $X_{[\alpha_0+1,\alpha)}=X_{[\alpha_0,\alpha)}$, since for any component $X_{\alpha'}$, $\alpha_0+1\leq\alpha'<\alpha_0+\omega$, skand $X_{[\alpha_0+1, \alpha'}$, and any component $X_{\alpha'}$, $\alpha_0\leq\alpha'<\alpha_0+\omega$, of the skand $X_{[\alpha_0,\alpha)}$, there is an equality $X_{\alpha'}=X_{\alpha_0}$, and the components of $X_{\alpha'}$ of the skands $X_{[\alpha_0+1,\alpha)}$ and $X_{[\alpha_0,\alpha)}$ for all ordinals $\alpha'\in[\alpha_0+\omega,\alpha)$ are equal to each other,  i.e., $X_{\alpha'}=X_{\alpha'}$.

If condition $2)$ of Proposition 3 is satisfied, then the equality $X_{[\alpha_0+1,\alpha)}=X_{[\alpha_0,\alpha)}$ is a special case where $\alpha'=\alpha_0+1$. $\Box$

\bigskip
\begin{center}
{\bf 3. Self-similar skands}
\end{center}

{\bf Definition 5}.
 An infinite skand $X_{[\alpha_0,\alpha)}$ is called {\it self-similar} if for each ordinal $\alpha'$ such that $\alpha_0<\alpha'<\alpha$, the skand $X_{[\alpha',\alpha)}$ is equal to the skand $X_{[\alpha_0,\alpha)}$, i.e., $X_{[\alpha',\alpha)}=X_{[\alpha_0,\alpha)}$. 

{\bf Proposition 4}. {\it A skand $X_{[\alpha_0,\alpha)}$ is self-similar if and only if $\alpha-\alpha_0=\omega^\kappa$, $0<\kappa< \Omega$, and thus for any ordinal $\alpha'$, $\alpha_0<\alpha'<\alpha$, there is an equality of the sets $X_{\alpha'}=X_{\alpha_0}$ or, similarly, the mapping $f: {[\alpha_0,\alpha)}\rightarrow{\bf V}[{\cal U}]$ associated with the skand $X_{[\alpha_0,\alpha)}$ is constant. }

{\bf Proof}. {\bf Necessity}. By Definition 5, for any ordinal $\alpha'$ such that $\alpha_0<\alpha'<\alpha$, there is an equality of the skands $X_{[\alpha',\alpha)}=X_{[\alpha_0,\alpha)}$, and hence an isomorphism $\varphi: [\alpha',\alpha)\rightarrow[\alpha_0,\alpha)$ of the intervals $[\alpha',\alpha)$ and $[\alpha_0,\alpha)$. Since the formula $\psi(\alpha'')=\alpha_0+\alpha''$, for any ordinal $\alpha''\in[0,\alpha-\alpha_0)$, gives the isomorphism $\psi: [0,\alpha-\alpha_0)\rightarrow[\alpha_0,\alpha)$ of the intervals $[0,\alpha-\alpha_0)$ and $[\alpha_0,\alpha)$, then for any residue $\rho$ of the ordinal $\alpha-\alpha_0$ (i.e., such an ordinal $\rho$ that $\alpha-\alpha_0=\sigma+\rho$ for some ordinal $\sigma$) there exists  an ordinal $\alpha''\in[0,\alpha-\alpha_0)$ such that $\rho$ is defined by the interval $[\alpha'',\alpha-\alpha_0)$ and vice versa: any ordinal $\alpha''$ such that $0<\alpha''<\alpha-\alpha_0$ defines the residue $\rho$ of the ordinal $[0,\alpha-\alpha_0)$ by the interval $[\rho,\alpha-\alpha_0)$. But for any ordinal $\alpha'=\psi(\alpha'')$ there is an isomorphism of the intervals $\varphi:[\alpha',\alpha)\rightarrow[\alpha_0,\alpha)$ and hence the isomorphism $\psi^{-1}\circ\varphi\circ\psi: [\alpha'',\alpha-\alpha_0)\rightarrow[0,\alpha-\alpha_0)$ intervals $[\alpha'',\alpha-\alpha_0)$ and $[0,\alpha-\alpha_0)$, and hence the equality of the numbers $\rho$ and $\alpha-\alpha_0$. 

 And if any residue $\rho$ of the ordinal number $\alpha-\alpha_0$ is equal to this number $\alpha-\alpha_0$, then it happens if and only if $\alpha-\alpha_0=\omega^\kappa$, $\kappa\geq 1$, (see: \cite{l5}, Chap. VII, $\S 7$, Theorem 7).

The equality of the skands $X_{[\alpha',\alpha)}=X_{[\alpha_0,\alpha)}$, for any ordinal $\alpha'$ such that $\alpha_0<\alpha'<\alpha$, by Definition $2'$, entails equality of the first components of these skands, i.e., $X_{\alpha'}=X_{\alpha_0}$, for any ordinal $\alpha_0<\alpha'<\alpha$. And since $X_{\alpha'}=f(\alpha')$, where  $f:{[\alpha_0,\alpha)}\rightarrow{\bf V}[{\cal U}]$ is the mapping associated with the skand $X_{[\alpha_0,\alpha)}$, it is a constant mapping, whose image $f([\alpha_0,\alpha))$ is the subset $\{X_{\alpha_0}\}$ of the universal class ${\bf V}[{\cal U}]$.

{\bf Sufficiency}. Let $\alpha-\alpha_0=\omega^\kappa$, $0<\kappa\leq \Omega$. Then by the already mentioned Theorem 7 of \cite{l5}, Chap. VII, \S 7, and the above, we obtain an isomorphism of the intervals $[\alpha_0,\alpha)$ and $[\alpha',\alpha)$ for any $\alpha'$, $\alpha_0<\alpha'<\alpha$. And since the mapping $f:{[\alpha_0,\alpha)}\rightarrow{\bf V}[{\cal U}]$ associated with the skand $X_{[\alpha_0,\alpha)}$ is constant, then, by  Definition 2$'$, for any $\alpha'$, $\alpha_0<\alpha'<\alpha$, we obtain an equality of the skands $X_{[\alpha',\alpha)}=X_{[\alpha_0,\alpha)}$. $\Box$

Clearly, any self-similar skand $X_{[\alpha_0,\alpha)}$ is reflexive, but not vice versa. For a self-similar skand $X_{[\alpha_0,\alpha)}$, the ordinal $\alpha-\alpha_0\geq\omega$ must be a limit ordinal of the form $\alpha-\alpha_0=\omega^\kappa$, $\kappa\geq 1$, and for the reflexive skand $X_{[\alpha_0,\alpha)}$ the ordinal $\alpha-\alpha_0\geq\omega$ can be either a limit of any kind or a non-limit ordinal number.

\bigskip
\begin{center}
{\bf 4. Periodic skands}
\end{center}

{\bf Definition 6}. A skand $X_{[\alpha_0,\alpha)}$ is called {\it  weakly periodic} with period  $\tau=\omega^{\xi_1}\cdot n_1+\omega^{\xi_2}\cdot n_2+...+\omega^{\xi_r}\cdot n_r$, where  $\xi_1>\xi_2>...>\xi_r$, $0\leq n_i<\omega$ for each natural numbers $n_i$, $0<i\leq r$,
  if $\tau$ is the smallest ordinal  $0<\tau<\Omega$ such that for each ordinal $\alpha'\in[\alpha_0,\alpha_0+\omega^{\xi_1+1})$  the equality    $X_{[\alpha',\alpha)}=X_{[\alpha'(+)\tau,\alpha)}$ holds;  it is called {\it periodic} if for each ordinal $\alpha'\in[\alpha_0,\alpha)$  the skand    $X_{[\alpha',\alpha)}$ weakly periodic with the same period $\tau$; and it is called {\it strictly periodic} if it is periodic and for any  ordinal $\lambda=\alpha_0+\omega^{\xi_1+1}\cdot\kappa\in[\alpha_0,\alpha)$,  $0<\kappa<\kappa(\alpha)$, the equality $X_{[\alpha_0,\alpha)}=X_{[\lambda,\alpha)}$ holds. 

{\bf Remark 2}. One has distinguish here and further  so-called a {\it natural sum} of ordinal numbers $\alpha'(+)\tau$ and a usual sum $\alpha_0+\omega^{\xi_1+1}$   of ordinal numbers. In general, they are different  (a natural sum $(+)$ is commutative, at that time a usual sum $+$ is not commutative). Definitions of natural sum $(+)$ and of natural product $(\cdot)$ see in \cite{l5} Chap. VII, $\S$ and \cite{l3363}, p. 591-594. 
Notice only that for any ordinals $\nu$ and $\mu$ one hasf $\nu(+)\mu=\max\{\nu+\mu,\mu+\nu\}$, where $+$ is a usual sum of ordinal numbers, which, clearly, for finite ordinals $0<\nu<\omega$ coincides with a natural sum for every $\mu$ (finite or infinite ordinals).

{\bf Proposition 5}. {\it A skand $X_{[\alpha_0,\alpha)}$ is weakly periodic with period $\tau=\omega^{\xi_1}\cdot n_1+\omega^{\xi_2}\cdot n_2+...+\omega^{\xi_r}\cdot n_r$, where  $\xi_1>\xi_2>...>\xi_r$, $0\leq n_i<\omega$ for each natural numbers $0<i\leq r$ are ordinal numbers,
 if and only if  $\alpha\geq\alpha_0+\omega^{\xi_1+1}$ and the interval $[\alpha_0,\alpha_0+\omega^{\xi_1+1})$ is a discrete union of intervals of minimal length, i.e., 

\begin{equation}
\label{l771}
[\alpha_0,\alpha_0+\omega^{\xi_1+1})=\bigcup_{0\leq\sigma<\omega}[\alpha_0(+)\tau\cdot\sigma,\alpha_0(+)\tau\cdot(\sigma+1))
\end{equation}
  such that all restrictions $f|_{[\alpha_0(+)\tau\cdot\sigma,\alpha_0(+)\tau\cdot(\sigma+1))}$ of the mapping $f: [\alpha_0,\alpha)\rightarrow{\bf V}[{\cal U}]$, associated with the skand $X_{[\alpha_0,\alpha)}$, on the intervals $[\alpha_0(+)\tau\cdot\sigma,\alpha_0(+)\tau\cdot(\sigma+1))$, $0\leq\sigma<\omega$, are equal to each other as counting functions, or what is the same, are related as follows: 
\begin{equation}
\label{l7782}
\begin{array}{ccc}
f|_{[\alpha_0(+)\tau\cdot\sigma,\alpha_0(+)\tau\cdot(\sigma+1))}=f|_{[\alpha_0(+)\tau\cdot(\sigma+1),\alpha_0(+)\tau\cdot(\sigma+2))}
 \circ\varphi_\sigma,
\end{array}
 \end{equation}
 where 
\begin{equation}
\label{l7783}
\varphi_\sigma: [\alpha_0(+)\tau\cdot\sigma,\alpha_0(+)\tau\cdot(\sigma+1))\rightarrow[\alpha_0(+)\tau\cdot(\sigma+1),\alpha_0(+)\tau\cdot(\sigma+2))
\end{equation}
  is an isomorphism  of well-ordered intervals $[\alpha_0(+)\tau\cdot\sigma,\alpha_0(+)\tau\cdot(\sigma+1))$ and $[\alpha_0(+)\tau\cdot(\sigma+1),\alpha_0(+)\tau\cdot(\sigma+2))$, for each $\sigma\in[\alpha_0,\omega)$, given by the formula $\varphi_\sigma(\alpha')=\alpha'(+)\tau$ for each $\alpha'\in[\alpha_0(+)\tau\cdot\sigma,\alpha_0(+)\tau\cdot(\sigma+1))$. }

{\bf Proof}. {\bf Necessity}. Let $X_{[\alpha_0,\alpha)}$ be a weakly periodic skand with period $\tau=\omega^{\xi_1}\cdot n_1+\omega^{\xi_2}\cdot n_2+...+\omega^{\xi_r}\cdot n_r$. We will prove that $\alpha\geq\alpha_0+\omega^{\xi_1+1}$. If it is not and $\alpha<\alpha_0+\omega^{\xi_1+1}$, then for any ordinal  $\alpha'$ such that $\alpha\leq\alpha'<\alpha_0+\omega^{\xi_1+1}$  there is no equality $X_{[\alpha',\alpha)}=X_{[\alpha'(+)\tau,\alpha)}$   because there is no terms like $X_{[\alpha',\alpha)}$ and  $X_{[\alpha'(+)\tau,\alpha)}$ what contradicts  the Definition 6.

  Further, for each ordinal $\alpha_0\leq\alpha'<\alpha_0+\omega^{\xi_1+1}$ the equality    $X_{[\alpha',\alpha)}=X_{[\alpha'(+)\tau,\alpha)}$ holds.  Then we can represent the interval $[\alpha_0,\alpha_0+\omega^{\xi_1+1})$ as a discrete sum $\bigcup\limits_{0\leq\sigma<\omega}[\alpha_0(+)\tau\cdot\sigma,\alpha_0(+)\tau\cdot(\sigma+1))$  of isomorphic intervals 
$$[\alpha_0(+)\tau\cdot\sigma,\alpha_0(+)\tau\cdot(\sigma+1)),\,\,\, 0\leq\sigma<\omega,$$ 
where each such isomorphism is some finite composition of the following elementary isomorphisms 
$$\varphi_\sigma:[\alpha_0(+)\tau\cdot\sigma,\alpha_0(+)\tau\cdot(\sigma+1))\rightarrow[\alpha_0(+)\tau\cdot(\sigma_1),\alpha_0(+)\tau\cdot(\sigma+2)),$$
 given by the following formula:
$$\varphi_\sigma(\alpha')=\alpha'(+)\tau,\,\,\,\sigma\in[0,\omega)$$.

Now, by virtue of the equality $X_{[\alpha',\alpha)}=X_{[\alpha'(+)\tau,\alpha)}$ for any ordinal $\alpha'\in[\alpha_0,\alpha_0+\omega^{\xi_1+1})$, and by Definition $2'$, we obtain the equality $(\ref{l7782})$,  i.e., for any ordinal $\alpha'\in[\alpha_0(+)\tau\cdot\sigma,\alpha_0(+)\tau\cdot(\sigma+1))$ the equality $X(\alpha')=X(\varphi_\sigma(\alpha'))$ holds. Which means that the restrictions $f|_{[\alpha_0(+)\tau\cdot\sigma,\alpha_0(+)\tau\cdot(\sigma+1))}$ of the mapping $f: [\alpha_0,\alpha)\rightarrow{\bf V}[{\cal U}]$, associated with the skand $X_{[\alpha_0,\alpha)}$, on the intervals $[\alpha_0(+)\tau\cdot\sigma,\alpha_0(+)\tau\cdot(\sigma+1))$, $0\leq\sigma<\omega$, are equal to each other as counting functions, i.e., the equality $(\ref{l7783})$ holds.

{\bf Sufficiency}. Let $\alpha$ in a skand $X_{[\alpha_0,\alpha)}$ be an ordinal with the inequality $\alpha\geq\alpha_0+\omega^{\xi_1+1}$ and  the conditions $(\ref{l7782})$ of Proposition $5$ be satisfied, i.e., all restrictions $f|_{[\alpha_0(+)\tau\cdot\sigma,\alpha_0(+)\tau\cdot(\sigma+1))}$  (as counting functions) of the mapping $f: [\alpha_0,\alpha)\rightarrow{\bf V}[{\cal U}]$ associated with the skand $X_{[\alpha_0,\alpha)}$, on the intervals $[\alpha_0(+)\tau\cdot\sigma,\alpha_0(+)\tau\cdot(\sigma+1))$, $0\leq\sigma<\omega$, be equal as counting functions, i.e., 
$$f|_{[\alpha_0(+)\tau\cdot\sigma,\alpha_0(+)\tau\cdot(\sigma+1))}=f|_{[\alpha_0(+)\tau\cdot(\sigma+1),\alpha_0(+)\tau\cdot(\sigma+2))}\circ\varphi_\sigma,\,\,\, 0\leq\sigma<\omega,$$ where 
$$\varphi_\sigma: [\alpha_0(+)\tau\cdot\sigma,\alpha_0(+)\tau\cdot(\sigma+1))\rightarrow[\alpha_0(+)\tau\cdot(\sigma+1),\alpha_0(+)\tau\cdot(\sigma+2))$$ 
is an  isomorphism, given by the formula $\varphi_\sigma(\alpha')=\alpha'(+)\tau$, $\alpha'\in[\alpha_0(+)\tau\cdot\sigma,\alpha_0(+)\tau\cdot(\sigma+1))$, the well-ordered segments $[\omega\cdot\sigma, \alpha_0(+)\omega\cdot(\sigma+1))$ and $[\alpha_0(+)\tau\cdot(\sigma+1),\alpha_0(+)\tau\cdot(\sigma+2))$, for each $\sigma\in[0,\omega)$, respectively. 

We have to prove that for any ordinal $\alpha'$, $\alpha_0\leq\alpha'<\alpha_0+\omega^{\xi_1+1}$, the equality $X_{[\alpha',\alpha)}=X_{[\alpha'(+)\tau,\alpha)}$ of the skands holds. 

For this purpose, let's show that there exists an isomorphism $\varphi: [\alpha',\alpha)\rightarrow[\alpha'(+)\tau,\alpha)$ of intervals $[\alpha',\alpha)=[\alpha', \alpha_0+\omega^{\xi_1+1})\cup[\alpha_0+\omega^{\xi_1+1},\alpha)$ and $[\alpha'(+)\tau,\alpha)=[\alpha'(+)\tau,\alpha_0+\omega^{\xi_1+1})\cup[\alpha_0+\omega^{\xi_1+1},\alpha)$. 

Consider the mapping 
$\varphi'=\bigcup\limits_{0\leq\sigma<\omega}\varphi_\sigma\cup 1|_{[\alpha_0+\omega^{\xi_1+1},\alpha)},$ which is a gluing of the above isomorphisms $\varphi_\sigma: [\alpha_0(+)\tau\cdot\sigma,\alpha_0(+)\tau\cdot(\sigma+1))\rightarrow[\alpha_0(+)\tau\cdot(\sigma+1),\alpha_0(+)\tau\cdot(\sigma+2))$, $0\leq\sigma<\omega$, and the identical mapping $1|_{[\alpha_0+\omega^{\xi_1+1},\alpha)}$ of the interval $[\alpha_0+\omega^{\xi_1+1},\alpha)$. Then the desired isomorphism $\varphi=\varphi'|_{[\alpha',\alpha_0+\omega^{\xi_1+1})}\cup 1|_{[\alpha_0+\omega^{\xi_1+1},\alpha)}$.
 
Now consider any ordinal $\alpha''\in[\alpha',\alpha_0+\omega^{\xi_1+1})$. Then, by virtue of the formula $(\ref{l771})$, there is a smallest ordinal $\sigma_0\in[0,\omega)$ such that $\alpha''\in[\alpha_0(+)\tau\cdot\sigma_9,\alpha_0(+)\tau\cdot(\sigma_0+1))$.  By virtue of the equality between the restrictions $$f|_{[\alpha_0(+)\tau\cdot\sigma,\alpha_0(+)\tau\cdot(\sigma+1))}=f|_{[\alpha_0(+)\tau\cdot(\sigma+1),\alpha_0(+)\tau\cdot(\sigma+2))}$$ as counting functions for any ordinal $\sigma$, $0\leq\sigma<\omega$, we obtain an equality 
$$f|_{[\alpha_0(+)\tau\cdot\sigma,\alpha_0(+)\tau\cdot(\sigma+1))}(\alpha'')=f|_{[\alpha_0(+)\tau\cdot(\sigma_1),\alpha_0(+)\tau\cdot(\sigma+2))}(\varphi_\sigma(\alpha''))$$ of values of counting functions 

$f|_{[\alpha_0(+)\tau\cdot\sigma,\alpha_0(+)\tau\cdot(\sigma+1))}$ and $f|_{[\alpha_0(+)\tau\cdot(\sigma_1),\alpha_0(+)\tau\cdot(\sigma+2))}(\varphi_\sigma(\alpha''))$ for any ordinal  $\alpha'\leq\alpha''<\alpha_0+\omega^{\xi_1+1}$, i.e.,  the equality  $f(\alpha'')=f(\varphi(\alpha''))$. Note that, $f(\alpha'')=f(\varphi(\alpha''))$ for each ordinal $\alpha_0+\omega^{\xi_1+1}\leq\alpha''<\alpha$ because $\varphi(\alpha'')=\alpha''$. Thus we obtain the equality of the components $X_{\alpha''}=X_{\varphi(\alpha'')}$ of the skands $X_{[\alpha',\alpha)}$ and $X_{[\alpha'(+)\tau,\alpha)}$,  for each ordinal   $\alpha_0\leq\alpha'<\alpha_0+\omega^{\xi_1+1}$. $\Box$

{\bf Proposition 6}. {\it A skand $X_{[\alpha_0,\alpha)}$ is  periodic with period $\tau=\omega^{\xi_1}\cdot n_1+\omega^{\xi_2}\cdot n_2+...+\omega^{\xi_r}\cdot n_r$, where  $\xi_1>\xi_2>...>\xi_r$, $0\leq n_i<\omega$ for each natural numbers $0<i\leq r$ are ordinal numbers, 
 if and only if $\alpha=\alpha_0+\omega^{\xi_1+1}\cdot\kappa$, for a fixed ordinal $0<\kappa<\Omega$, and for each ordinal $\lambda=\alpha_0+\omega^{\xi_1+1}\cdot\nu$, $0\leq\nu<\kappa$, there is an equality $X_{[\lambda(+)\tau\cdot\sigma,\alpha)}=X_{\lambda(+)\tau\cdot(\sigma+1),\alpha)}$ for each ordinal $0<\sigma<\omega$.}

{\bf Proof}. {\bf Necessity}. Let a skand $X_{[\alpha_0,\alpha)}$ be periodic with period $\tau=\omega^{\xi_1}\cdot n_1+\omega^{\xi_2}\cdot n_2+...+\omega^{\xi_r}\cdot n_r$. Then, by Definition 6, for any ordinal $\alpha'$, $\alpha_0\leq\alpha'<\alpha$, the skand $X_{[\alpha',\alpha)}$ is weakly periodic with the same period $\tau$.
In particular when $\alpha'=\alpha_0$ we obtain, by Theorem 5, an inequality
$\alpha_0+\omega^{\xi_1+1}\leq\alpha$ and the interval $[\alpha_0,\alpha_0+\omega^{\xi_1+1})$ is a discrete union  $\bigcup\limits_{0\leq\sigma<\omega}[\alpha_0(+)\tau\cdot\sigma,\alpha_0(+)\tau\cdot(\sigma+1))$  of intervals $[\alpha_0(+)\tau\cdot\sigma,\alpha_0(+)\tau\cdot(\sigma+1))$, $0\leq\sigma<\omega$, of minimal length such that there are equalities $X_{[\alpha_0(+)\tau\cdot\sigma,\alpha)}=X_{[\alpha_0(+)\tau\cdot(\sigma+1),\alpha)}$   for all ordinals $0\leq\sigma<\omega$. If $\alpha=\alpha_0+\omega^{\xi_1+1}$ then the proof of necessity is complete. 

If not and $\alpha_0+\omega^{\xi_1+1}<\alpha$, then for $\lambda=\alpha_0+\omega^{\xi_1+1}$ the skand $X_{[\lambda,\alpha)}$ is weakly periodic with period $\tau$ an as above the interval $[\lambda,\lambda+\omega^{\xi_1+1})$ is a discrete union  $\bigcup\limits_{0\leq\sigma<\omega}[\lambda(+)\tau\cdot\sigma,\lambda(+)\tau\cdot(\sigma+1))$  of intervals $[\lambda(+)\tau\cdot\sigma,\lambda(+)\tau\cdot(\sigma+1))$, $0\leq\sigma<\omega$, of minimal length such that there are equalities $X_{[\lambda(+)\tau\cdot\sigma,\alpha)}=X_{\lambda(+)\tau\cdot(\sigma+1),\alpha)}$ for each ordinal    for all ordinals $0\leq\sigma<\omega$ such that there are equalities $X_{[\lambda+\omega^{\xi_1+1}(+)\tau\cdot\sigma,\alpha)}=X_{[\alpha_0+\omega^{\xi_1+1}(+)\tau\cdot(\sigma+1),\alpha)}$   for all ordinals $0\leq\sigma<\omega$. If $\alpha=\lambda+\omega^{\xi_1+1}=\alpha_0+\omega^{\xi_1+1}+\omega^{\xi_1+1}=\alpha_0+
=\alpha+\omega^{\xi_1+1}\cdot 2$, then the proof of necessity is complete. 

If not and $\lambda=\alpha_0+\omega^{\xi_1+1}\cdot 2<\alpha$, then we continue the process for all possible ordinals $\nu$ such that $\alpha_0+\omega^{\xi_1+1}\cdot\nu\leq\alpha$ and the equalities $X_{[\alpha_0+\omega^{\xi_1+1}\cdot\nu,\alpha)}=X_{[\alpha_0+\omega^{\xi_1+1}\cdot(\nu+1),\alpha)}$.

If $\alpha=\lambda+\alpha_0+\omega^{\xi_1+1}\cdot\nu$, then the proof of necessity is complete. If not, then we put $\alpha=\sup\limits_{0<\nu<\kappa}(\alpha_0+\omega^{\xi_1+1}\cdot\nu)=\alpha_0+\omega^{\xi_1+1}\cdot\kappa$ and
for each ordinal $\lambda=\alpha_0+\omega^{\xi_1+1}\cdot\nu$, $0<\nu<\kappa$, there is an equality $X_{[\alpha'(+)\tau\cdot\sigma,\alpha)}=X_{\lambda'(+)\tau\cdot(\sigma+1),\alpha)}$ for all $\lambda'=\alpha_0+\omega^{\xi_1+1}\cdot\nu$ for each ordinal $0\leq\nu<\kappa$ and each ordinal $0<\sigma<\omega$. And the proof of necessity is complete.

{\bf Sufficiency}. Let $\alpha=\alpha_0+\omega^{\xi_1+1}\cdot\kappa$ be an ordinal for a fixed ordinal $0<\kappa<\Omega$ such that for each ordinal $\lambda=\alpha_0+\omega^{\xi_1+1}\cdot\nu$, $0<\nu<\kappa$, there is an equality $X_{[\lambda(+)\tau\cdot\sigma,\alpha)}=X_{\lambda(+)\tau\cdot(\sigma+1),\alpha)}$ for each ordinal $0<\sigma<\omega$, where $\tau=\omega^{\xi_1}\cdot n_1+\omega^{\xi_2}\cdot n_2+...+\omega^{\xi_r}\cdot n_r$ is a canonical form of an ordinal $\tau$ and   $\xi_1>\xi_2>...>\xi_r$, $0\leq n_i<\omega$ for each natural numbers $0<i\leq r$ are unique ordinal numbers.  By Definition 6, we need to prove that for each ordinal $\alpha_0\leq\alpha'<\alpha$ the skand $X_{[\alpha',\alpha)}$ is weakly periodic with preiod $\tau$. But for any fixed $0\leq\alpha'<\alpha$ there is a minimal ordinal $0<\nu_0<\kappa$  and a minimal $0\leq\sigma_0<\alpha_0+\omega^{\xi_1+1}$ such that $\alpha_0+\omega^{\xi_1+1}\cdot\nu_0(+)\tau\cdot\sigma_0\leq\alpha'<
\alpha_0+\omega^{\xi_{1}+1}\cdot \nu_0(+)\tau\cdot(\sigma_0+1)<\alpha_0+\omega^{\xi_{1}+1}\cdot (\nu_0+1)$.  Since the above equality $X_{[\lambda(+)\tau\cdot\sigma,\alpha)}=X_{\lambda(+)\tau\cdot(\sigma+1),\alpha)}$ implies $X_{[\alpha_0+\omega^{\xi_1+1}\cdot\nu_0(+)\tau\cdot\sigma_0,\alpha)}=X_{[\alpha_0+\omega^{\xi_1+1}\cdot\nu_0(+)\tau\cdot(\sigma_0+1),\alpha)}$ the equality $X_{[\alpha'(+)\tau\cdot\sigma,\alpha)}=X_{\alpha'(+)\tau\cdot(\sigma+1),\alpha)}$ for each ordinal $0\leq\sigma<\omega$. Just to be clear, $\alpha'$ plays the same role here as $\alpha_0$ in Theorem 6. And the proof of sufficiency is complete and thus the proof of Proposition 6 is complete.$\Box$

{\bf Proposition 7}. {\it A skand $X_{[\alpha_0,\alpha)}$ is strictly periodic with period $\tau=\omega^{\xi_1}\cdot n_1+\omega^{\xi_2}\cdot n_2+...+\omega^{\xi_r}\cdot n_r$, where  $\xi_1>\xi_2>...>\xi_r$, $0\leq n_i<\omega$ for each natural numbers $0<i\leq r$ are ordinal numbers,
 if and only if $\alpha=\omega^\mu$ for some fixed ordinal $\xi_1+1\leq\mu<\Omega$ and the interval $[\alpha_0,\alpha_0+\omega^{\xi_1+1})$ is a discrete union of the following intervals:

\begin{equation}
\label{l779}
[\alpha_0,\alpha)=\bigcup\limits_{0\leq\nu<\alpha}\bigcup\limits_{0<\sigma<\omega^{\xi_1+1}}[\alpha_0+\omega^{\xi_1+1}\cdot\nu(+)\tau\cdot\sigma,\alpha_0+\omega^{\xi_1+1}\cdot\nu\cdot(+)\tau(\sigma+1)),
\end{equation}
such that restrictions  
 $f|_{[\alpha_0+\omega^{\xi_1+1}\cdot\nu(+)\tau\cdot\sigma,\alpha_0(+)\tau\cdot\nu\cdot(\sigma+1))}$,  $0<\sigma<\omega$, $0\leq\nu<\alpha$, of mapping $f:[\alpha_0,\alpha)\rightarrow{\bf V}[{\cal U}]$, associated with the skand  $X_{[\alpha_0,\alpha)}$, on the intervals $[\alpha_0+\omega^{\xi_1+1}\cdot\nu(+)\tau\cdot\sigma,\alpha_0+\omega^{\xi_1+1}(+)\tau\cdot\nu\cdot(\sigma+1))$,  $0\leq\nu<\alpha$, $0<\sigma<\alpha_0+\omega^{\xi_1+1}$,  are equal to each other as counting functions.}

{\bf Proof. Necessity } Let $X_{[\alpha_0,\alpha)}$ be a strictly periodic skand with period $\tau=\omega^{\xi_1}\cdot n_1+\omega^{\xi_2}\cdot n_2+...+\omega^{\xi_r}\cdot n_r$ as in Proposition 7. Then it is periodic with the  periodic $\tau$ and for each  ordinal $\lambda=\alpha_0+\omega^{\xi_1+1}\cdot\kappa\in[\alpha_0,\alpha)$,  $0<\kappa<\kappa(\alpha)$, the equality $X_{[\lambda,\alpha)}=X_{[\alpha_0,\alpha)}$ holds.  Thus, there is an isomorphism $\varphi: [\alpha_0,\alpha)\rightarrow[\lambda,\alpha)$, for each ordinal  $\lambda=\alpha_0+\omega^{\xi_1+1}\cdot\kappa\in[\alpha_0,\alpha)$, where $0<\kappa<\kappa(\alpha)$. For any fixed ordinal $0<\alpha'<\alpha$ there is a minimal $\lambda_0$ and $\kappa_0$ such that $\lambda_0=\alpha_0+\omega^{\xi_1+1}\cdot\kappa_0\leq\alpha'<
\alpha_0+\omega^{\xi_1+1}\cdot(\kappa_0+1)<\alpha$ and there is an isomorphism  $\varphi': [\alpha_0,\alpha)\rightarrow[\alpha',\alpha)$, given by the formula  $\varphi'(\alpha'')=\varphi(\alpha''(+)\lambda_0)$ for each ordinal $0\leq\alpha''<\alpha$. Since  for each $0<\alpha'<\alpha$ there is an isomorphism $\varphi': [\alpha_0,\alpha)\rightarrow[\alpha',\alpha)$ we conclude that $\alpha=\omega^\mu$ for some ordinal $0<\mu<\Omega$, (see: \cite{l5}, Chap. VII, $\S 7$, Theorem 7). Clearly, $\mu\geq\xi_1+1$.

Moreover, since by Definition 6, for each  ordinal $\lambda=\alpha_0+\omega^{\xi_1+1}\cdot\kappa\in[\alpha_0,\alpha)$,  $0<\kappa<\kappa(\alpha)$, the equality $X_{[\lambda,\alpha)}=X_{[\alpha_0,\alpha)}$ holds we conclude that the interval $[\alpha_0,\alpha_0+\omega^{\xi_1+1})$ is a discrete union of intervals $(\ref{l779})$ such that restrictions  
 $f|_{[\alpha_0+\omega^{\xi_1+1}\cdot\nu(+)\tau\cdot\sigma,\alpha_0(+)\tau\cdot\nu\cdot(\sigma+1))}$,  $0<\sigma<\omega$, $0\leq\nu<\alpha$, of mapping $f:[\alpha_0,\alpha)\rightarrow{\bf V}[{\cal U}]$, associated with the skand  $X_{[\alpha_0,\alpha)}$, on the intervals $[\alpha_0+\omega^{\xi_1+1}\cdot\nu(+)\tau\cdot\sigma,\alpha_0+\omega^{\xi_1+1}(+)\tau\cdot\nu\cdot(\sigma+1))$,  $0\leq\nu<\alpha$, $0<\sigma<\omega$,  are equal to each other as counting functions. Necessity  is completely proved.

{\bf Sufficiency}. Let the conditions of Proposition 7 be satisfied, i.e., 
$\alpha=\omega^\mu$ for some fixed ordinal $\xi_1+1\leq\mu<\Omega$ and
the interval $[\alpha_0,\alpha_0+\omega^{\xi_1+1})$ is a discrete union of the  intervals $[\alpha_0+\omega^{\xi_1+1}\cdot\nu(+)\tau\cdot\sigma,\alpha_0+\omega^{\xi_1+1}\cdot\nu\cdot(+)\tau(\sigma+1))$ $(\ref{l779})$
such that restrictions  
 $f|_{[\alpha_0+\omega^{\xi_1+1}\cdot\cdot(+)\tau\cdot\sigma,\alpha_0(+)\tau\cdot\nu\cdot(\sigma+1))}$,  $0<\sigma<\omega$, $0\leq\nu<\alpha$, of mapping $f:[\alpha_0,\alpha)\rightarrow{\bf V}[{\cal U}]$, associated with the skand  $X_{[\alpha_0,\alpha)}$, on the intervals $[\alpha_0+\omega^{\xi_1+1}\cdot\nu(+)\tau\cdot\sigma,\alpha_0+\omega^{\xi_1+1}(+)\tau\cdot\nu\cdot(\sigma+1))$,  $0\leq\nu<\alpha$, $0<\sigma<\omega$,  are equal to each other as counting functions.

Then for each ordinal $\lambda=\alpha_0+\omega^{\xi_1+1}\cdot\kappa\in[\alpha_0,\alpha)$,  $0<\kappa<\kappa(\alpha)$, the equality $X_{[\lambda,\alpha)}=X_{[\alpha_0,\alpha)}$ holds. Indeed, since $\alpha=\omega^\mu$ for some ordinal $\xi_1+1\leq\mu<\Omega$ we obtain an isomorphism $\varphi: [\lambda,\alpha)\rightarrow[\alpha_0,\alpha)$ (see: \cite{l5}, Chap. VII, $\S 7$, Theorem 7).

And by above condition on restrictions $f|_{[\alpha_0+\omega^{\xi_1+1}\cdot\nu(+)\tau\cdot\sigma,\alpha_0(+)\tau\cdot\nu\cdot(\sigma+1))}$,  which are equal to each other as counting functions, where $0<\sigma<\omega$, $0\leq\nu<\alpha$, of mapping $f:[\alpha_0,\alpha)\rightarrow{\bf V}[{\cal U}]$, associated with the skand  $X_{[\alpha_0,\alpha)}$, on the intervals $[\alpha_0+\omega^{\xi_1+1}\cdot\nu(+)\tau\cdot\sigma,\alpha_0+\omega^{\xi_1+1}(+)\tau\cdot\nu\cdot(\sigma+1))$,  $0\leq\nu<\alpha$, $0<\sigma<\omega$, we obtain a require isomorphism $\varphi': [\alpha_0,\alpha)\rightarrow[\lambda,\alpha)$, given by the formula $\varphi'(\alpha')=\alpha'(+)\alpha_0+\omega^{\xi_1+1}$, $0\leq\alpha'<\alpha$,  the equality $f|_{[\alpha_0,\alpha_0(+)\omega^{\xi_1+1})}=f|_{[\lambda,\lambda(+)+(\alpha_0+\omega^{\xi_1+1}))}$ and thus, by Definition $2'$ the equality $X_{[\alpha_0,\alpha)}=X_{[\lambda,\alpha)}$. Sufficiency is completely proved. $\Box$

{\bf Remark 3}. If $\alpha-\alpha_0=\omega$, then all three concepts -- weak periodicity, periodicity, and strict periodicity -- coincide. 
If $\alpha-\alpha_0>\omega$, then all three concepts -- weak periodicity, periodicity and strict periodicity are different, generally speaking. Clearly, that a strictly periodic skand is periodic and a periodic skand is weakly periodic but in both case non vica virsa. Indeed, a weakly periodic skand $X_{[\alpha_0,\alpha)}$ such that for some $\alpha_0+\omega\leq\lambda'<\alpha$ a skand $X_{[\alpha',\alpha)}\not= X_{[\alpha'',\alpha)}$ for every ordinal $\alpha_0\leq\alpha''<\alpha_0+\omega$, is not periodic. Similarly, a  periodic skand $X_{[\alpha_0,\alpha)}$ such that for some $\alpha_0+\omega\leq\lambda'<\alpha$ a skand $X_{[\alpha',\alpha)}\not= X_{[\alpha'',\alpha)}$ for every ordinal $\alpha_0\leq\alpha''<\alpha_0+\omega$, is not strictly periodic. A strictly periodic skand with period $\tau=1$ is a reflexive skand and a strictly periodic skand with period $\tau=1$ is a self-similar skand, and  in both cases vice versa. The notion of strict periodicity so defined here was applied \cite{l332} to extend the definition of the trigonometric functions $\sin x$ and $\cos x$ to {\it of all} Conway numbers $x\in{\bf No}$, which were originally defined by Conway only for such $x\in{\bf No}$ for which there is a natural number $n$ such that $-n<x<n$.

{\bf Example 1}. Consider a skand $X_{[0,\omega^2+1)}=$
$$\{_0 1,2,3\{_1 1,2,3\{_2 1,2,3\{_3...\{_\omega\omega\{_{\omega+1}\omega+1\{_{\omega+2}\omega+2...\{_{\omega^2}\omega^2\}_{\omega^2}...\}_{\omega21}\}_{\omega+1}\}_{\omega}...\}_3\}_2\}_1\}_0$$.

This skand is weakly periodic with period $\tau=3$ but it is not weakly periodic. Moreover, all its restrictions $X_{[\alpha',\omega^2+1}$ are not periodic with any period $\tau'$.

{\bf Example 2}. Consider a skand $X_{[0,\omega^2)}=$
$$\{_0 1,2,3\{_1 1,2,3\{_2 1,2,3\{_3...\{_\omega 4,5,6\{_{\omega+1}4,5,6\{_{\omega+2}4,5,6... ...\}_{\omega+2}\}_{\omega+1}\}_{\omega}...\}_3\}_2\}_1\}_0$$.

This skand is  periodic with period $\tau=3$ but it is nor strictly periodic. Nevertheless,  its restriction $X_{[\lambda,\omega^2)}$ are  periodic with the same period $\tau$ for each ordinal $\lambda=\omega\cdot\nu$, $0<\nu<\omega^2$ but not equal to each other.

{\bf Example 3}. Consider a skand $X_{[0,\omega^2)}=$
$$\{_0 1,2,3\{_1 1,2,3\{_2 1,2,3\{_3...\{_\omega 1,2,3\{_{\omega+1}1,2,3\{_{\omega+2}1,2,3... ...\}_{\omega+2}\}_{\omega+1}\}_{\omega}...\}_3\}_2\}_1\}_0$$.

This skand is  strictly periodic with period $\tau=3$ because it is weakly periodic and  its restrictions $X_{[\lambda,\omega^2}=X_{[0,\omega^2}$  for each ordinal $\lambda=\omega\cdot\nu$, $0<\nu<\omega^2$.

\bigskip
\parindent=0 cm
\begin{center}
{\bf 5. The elementary formal theory $NBG[{\cal U}]^{(1)}$ of non-founded sets and its model in the universe ${\bf U}[{\cal U}]^{(1)}$}
\end{center}

\smallskip

\parindent=0,5 cm
 
 Let's now define an extension $NBG[{\cal U}]^{(1)}$ of the formal theory $NBG[{\cal U}]^-$, where ${\cal U}$ can be either proper class of individuals or a set of individuals in particular the empty set, and then we obtain the theory of \grqq pure non-founded sets\grqq\, $NBG^{(1)}$mentioned above in the Abstract and Introduction.
 
By $NBG[{\cal U}]^{(1)}$, we mean the formal theory $NBG[{\cal U}]^{-}$ with two more axioms attached to its own axioms. 

\smallskip

{\bf Skand existence axiom} $({\bf SEA})$. For any ordinals $\alpha_0,\alpha\in{\bf On}$ such that $\alpha-\alpha_0\geq\omega$, and any mapping $f: [\alpha_0,\alpha)\rightarrow{\bf V}[{\cal U}]$ there exists a {\it skand} $X_{[\alpha_0,\alpha)}$, for which the equality $X_{\alpha'}=f(\alpha')$ holds for every ordinal $\alpha'\in[\alpha_0,\alpha)$, and it belongs to the universe $X_{[\alpha_0,\alpha)}\in{\bf V}[{\cal U}]^-$, i.e.,  the set $X_{\alpha_0}\cup\{X_{[\alpha_0+1,\alpha)}\}\in{\bf V}[{\cal U}]^-$.

{\bf Remark 4}. The class of all such sets $X_{\alpha_0}\cup\{X_{[\alpha_0+1,\alpha)}\}$ we denote by ${\bf V}[{\cal U}]^{(1)}$ if we interpret ${\bf V}[{\cal U}]={\bf V}[{\cal U}]^{(0)}$. If the length $l=\alpha-\alpha_0$ of the skand $X_{[\alpha_0,\alpha)}$ is finite, then the set $X_{[\alpha_0,\alpha)}=X_{\alpha_0}\cup\{X_{[\alpha_0+1,\alpha)}\}$ is founded; if the length is infinite, then the set $X_{[\alpha_0,\alpha)}=X_{\alpha_0}\cup\{X_{[\alpha_0+1,\alpha)}\}$ is non-founded but pseudo-founded what arises from the following axiom. Notice only that any set $X\in{\bf V}[{\cal U}]$ can be represented as a finite skand $X_{[\alpha_0,\alpha)}$, with the length $\alpha-\alpha_0=1$. Indeed, if $X=\{x_0,x_1,...,x_\lambda,...\}$, $0\leq\lambda<\Lambda$, is a founded set, then the skand $X_{[\alpha_0,\alpha_0+1)}$ with the length $\alpha-\alpha_0=1$ and $X_{\alpha_0}=X$ is a desired skand.

{\bf Pseudo-founding Axiom $({\bf PFA})$}.  For each set $X\in{\bf V}[{\cal U}]^{(1)}$, any maximal descending $\in$-sequence $X\ni x_1\ni x_2\ni x_3\ni...\ni x_{n-1}\ni x_n\ni.... $ has the following property: if the sequence is {\it finite}, then its last element $x_n$ is either the empty set $\emptyset$ or an individual $u\in{\cal U}$; if it is {\it infinite}, then there exists the smallest number $n$ such that $x_n$ is a skand of infinite length.

The latter means that $x_n=X_{[\alpha_0,\alpha)}$ for some ordinals $\alpha_0,\alpha\in{\bf On}$ such that $l=\alpha-\alpha_0\geq\omega$. 

 In the case when all maximal descending $\in$-sequences in $X$ are finite, then the set $X$ is obviously founded, and otherwise, when at least one $\in$-sequence is infinite, then the set $X$ is non-founded, but obviously {\it is pseudo-founded}, i.e., the $\in$-sequence is finite in the sense that its last term is a skand of infinite length.

Let's now show that the collection of all skands is a proper class. Since any skand $X_{[\alpha_0,\alpha)}$ is equal to the skand $Y_{[0,\bar\alpha)}$, where $\bar\alpha=\alpha-\alpha_0$ is ordinal and $Y_{\bar\alpha'}=X_{\alpha'}$, where $\bar\alpha'=\alpha_0$, $\alpha_0\leq\alpha'<\alpha$, then we restrict ourselves to the collection of all skands of the form $X_{[0,\alpha)}$, $0<\alpha<\Omega$, is evidently a proper class.

{\bf Proposition 8} {\it There exists a class ${\bf S}[{\cal U}]$ of all skands $X_{[0,\alpha)}$ with arbitrary components and lengths $0<\alpha<\Omega$ and it is a proper class.}

{\bf Proof.}
Any skand $X_{[0,\alpha)}\in{\bf S}[{\cal U}]$ with a clutch  region $[0,\alpha)$ and  components $X_{\alpha'}$, $0\leq\alpha'<\alpha$, we map the founded set $X^{{\bf v}_\alpha}_\alpha=\bigcup\limits_{0\leq\alpha'<\alpha}\{<v_{\alpha'},\alpha'>\}$, where ${\bf v}_\alpha=\bigsqcup\limits_{0\leq\alpha'<\alpha}\{v_{\alpha'}\}$  is the discrete sum of the single-element sets $\{v_{\alpha'}\}$ of all indexes $0\leq\alpha'<\alpha$, $v_{\alpha'}=X_{\alpha'}$ and $<v_{\alpha'},\alpha'>=\{\{v_{\alpha'}\},\{v_{\alpha'},\alpha'\}\}$ is an ordered pair, which is  a founded set since the ordinals of$\alpha'$, by the formula $(\ref{f9999})$, are \grqq pure\grqq\, founded sets, and every  $v_{\alpha'}$ for any ordinal  $0\leq\alpha'<\alpha$ is a founded set.

Clearly, for two skands $X_{[0,\alpha)}$ and $Y_{[0,\beta)}$ of ${\bf S}[{\cal U}]$ such that $X_{[0,\alpha)}\not=Y_{[0,\beta)}$,  the corresponding sets $X_\alpha^{{\bf v}_\alpha}$ and $Y_\beta^{{\bf w}_\beta}$ form the class ${\bf V}[{\cal U}]$ are not equal, i.e.,  $X_\alpha^{{\bf v}_\alpha}\not=Y_\beta^{{\bf w}_\beta}$.

Consider now the class ${\bf T}=\bigcup\limits_{{\bf v}_\alpha, \alpha}\{X^{{\bf v}_\alpha}_\alpha\}$ of all one-element founded sets $\{X^{{\bf v}_\alpha}_\alpha\}$, where the sum is taken over all indexes of ${\bf v}_\alpha$, represented as a discrete sum of $\bigsqcup\limits_{0\leq\alpha'<\alpha}\{v_{\alpha'}\}$ of different or no-founded sets $v_{\alpha'}$, $0\leq\alpha'<\alpha$ and over all ordinals $\alpha\in(0<\alpha<\Omega)$. It is clear that ${\bf T}\subset{\bf V}[{\cal U}]$ and ${\bf T}$ is a proper class, and the mapping discussed above defines an injective mapping $i:{\bf S}[{\cal U}]\rightarrow{\bf T}$.

The converse is also true, to any element $X^{{\bf v}_\alpha}_\alpha\in{\bf T}$  we map uniquely to the skand $X_{[0,\alpha)}\in{\bf S}[{\cal U}]$, assuming $X_{\alpha'}=v_{\alpha'}$ for each ordinal $0\leq\alpha'<\alpha$, where $v_{\alpha'}$ are the constituents of the set $X^{{\bf v}_\alpha}_\alpha=\bigcup\limits_{0\leq\alpha'<\alpha}\{<v_{\alpha'},\alpha'>\}$. Moreover, we can see that different elements $X^{{\bf v}_\alpha}_\alpha$ and $Y^{{\bf w}_\beta}_\beta$ of the class ${\bf T}$  correspond to different skands $X_{[0,\alpha)}$ and $Y_{[0,\beta)}$ from ${\bf S}[{\cal U}]$. Thereby there is an injective mapping $j:{\bf T}\rightarrow{\bf S}[{\cal U}]$.

The equality $j\circ i=1_{{\bf S}[{\cal U}]}$ and $i\circ j=1_{{\bf T}}$ follows directly from the construction. Thus these classes ${\bf S}[{\cal U}]$ and ${\bf T}$ are bijective. 
 Hence, the class ${\bf S}[{\cal U}]$ exists and is a proper class.
 
 {\bf Corollary 1}. The collection ${\bf S}'[{\cal U}]$ of all skands $X_{[0,\alpha)}$ of length $\alpha\geq\omega$ is a proper class.

 {\bf Corollary 2}. Proposition 8 is a model of a skand in the universe ${\bf V}[{\cal U}]^-$.
 
 {\bf Proof} exactly repeats the proof of Proposition 8, only we need to consider all infinite skands.
 $\Box$

 {\smallskip}

{\bf Exact extensionality axiom} $({\bf EEA})$. $(\forall X)(\forall Y)[X=Y\Rightarrow(\forall Z)( X\in Z)\Rightarrow Y\in Z)$, where  sets $X,Y,Z\in{\bf V}[{\cal U}]^{(1)}$ (the elements of the sets can be individuals, founded sets, and unfounded sets) and $X$ and $Y$  {\it are equal} , if for each element $x\in X$, there exists an element $y\in Y$ such that $x=y$, and for each element $y\in Y$, there exists an element $x\in X$ such that $y=x$, where the equality \grqq$=$\grqq\ means the following:.

 $1)$ the equality of individuals $x,y\in {\cal U}$, i.e., $(x=y)\Leftrightarrow(\forall Z)(x\in Z\Leftrightarrow y\in Z)$; 
 
 $2$) the equality of founded sets  $Z,Z'\in{\bf V}[{\cal U}]^{(0)}$, i.e., by equality in Extensionality axiom ${\bf EA}$ in $NBG[{\cal U}]^{(0)}$ (i.e., two founded sets are equal if they contain the same elements);
 
 $3)$ the equality of skands when $x,y\in {\bf S}'[{\cal U}]$, i.e., by Definition  $2'$; 
 
 $4)$ the equality of sets when $x,y\notin{\cal U}$ and $x,y\notin{\bf V}[{\cal U}]^{(0)}$, i.e., the subsequent iteration using $1)$, $2)$, $3)$ in ${\bf EES}$.
 
 By virtue of the axiom ${\bf PFA}$, for every element $x,y$ in $4)$, such an iteration is finite and ends up with either an empty set, an individual, or an infinite skand to which the rules $2)$, $1)$, and $3)$ apply, respectively.

{\bf Remark 5}. The last condition means, for example, that if a set $x$ contains an element such that it contains as its elements individuals, founded sets and infinite skands, then a set $y$ must have the same element containing equal individuals, founded sets and infinite skands, respectively, and vice versa. This brings us back to the usual definition of an equality of sets: sets $X$ are equal to $Y$ if they consist of the same elements. A priori, if in the universe ${\bf V}[{\cal U}]$ there are two non-founded sets, e.g., $X\in X$ and $Y\in Y$, then their equality $X=Y$ without a specific definition of reflexive sets is impossible by the definition of equality in the extensionality axiom ${\bf EA}$, viz. i.e. by the formula $(\ref{f1947})$, since it reduces to the tautology $X=Y\Leftrightarrow X=Y$, but in $({\bf EEA})$ this equality agrees with the fact that equal sets do consist of the same elements and vice versa.

\bigskip

The last class of sets ${\bf V}[{\cal U}]^{(1)}$ is defined by the following transfinite recursion. We first extend the proper class or the  set ${\cal U}={\cal U}^{(0)}$ of individuals of the theory $NBG[{\cal U}]^{(0)}$ by adding to it elements from the proper class ${\cal U}^{(1)}$ of new individuals, i.e., we assume ${\cal U}'={\cal U}^{(0)}\cup{\cal U}^{(1)}$. The elements $u_\lambda^{(1)}\in {\cal  U}^{(1)}$, $\lambda\in{\bf On}$, are arbitrary infinite skands $X_{[\alpha_0,\alpha)}\in {\bf S}'[{\cal U}]$ together with the operator \grqq forgetting functor\grqq\, $F$, which \grqq forgets\grqq\, the internal structure of the membership relation $\in$ of the constituents, or parts of the skand  $X_{[\alpha_0,\alpha)}$ (cf. {\bf 2} in {\bf 2}); thus we can briefly denote ${\cal U}^{(1)}\stackrel{def}{=}F{\bf S}'[{\cal U}]$. More precisely,  $u_\lambda^{(1)}=FX_{[\alpha_0,\alpha)}$ and no element or object is an element in $u_\lambda^{(1)}$. Moreover, individuals are equal to $FX_{[\alpha_0,\alpha)}=FX_{[\beta_0,\beta)}$ if and only if $X_{[\alpha_0,\alpha)}=X_{[\beta_0,\beta)}$. 

\bigskip

(To avoid unnecessary complications, we assume that ${\cal U}\cap{\cal U}^{(1)}=\emptyset$; this is true, for example, when ${\cal U}=\emptyset$ and we are constructing a theory of \grqq pure\grqq\ sets; if ${\cal U}\cap{\cal U}^{(1)}\not=\emptyset$, then it's technically possible to separate the individuals of the intersection of ${\cal U}\cap{\cal U}^{(1)}$, denoting them by different letters in ${\cal U}$ and ${\cal U}^{(1)}$, or by considering the class ${\cal U}'$ as a discrete sum of ${\cal U}\sqcup{\cal U}^{(1)}$, especially since we have the class of individuals ${\cal U}^{(1)}$ playing only an auxiliary role. Below we show that these complications can be avoided by not applying the oblivion functor $F$).

 To build a model of $NBG[{\cal U}]^{(1)}$ theory, we need to make sure that the subclass ${\cal U}^{(1)}$ of the class ${\cal U}'$ really consists of individuals and satisfies all axioms of set theory with individuals. We only know that for any individual $u'\in{\cal U}^{(1)}$, i. e., $u'=FX_{[\alpha_0,\alpha)}$, $(\forall {\bf x})({\bf x}\notin u')$ and the equality $u'=v'$ is the equality of their corresponding skands $X_{[\alpha_0,\alpha]}$ and $Y_{[\beta_0,\beta]}$. Now we must show the equality of individuals in set theory with individuals, i.e., the formula $(\ref{f1947})$: $(\forall{\bf x})(u'=v'\wedge u'\in{\bf x}\Rightarrow v'\in{\bf x}).$

  Indeed, by the pair axiom, there exist one-element sets $\{u'\}$ and $\{v'\}$, and the equality $u'=v'$ entails the equality $\{u'\}=\{v'\}$. By the sum axiom, for any set ${\bf x}$ there exist $\cup\{\{u'\},{\bf x}\}$ and $\cup\{\{v'\},{\bf x}\}$, and the latter equality entails the equality $\cup\{\{u'\},{\bf x}\}=\cup\{\{v'\},{\bf x}\}$. By virtue of the premise $u'\in{\bf x}$, we get $u'\in\cup\{\{u'\},{\bf x}\}$. Since  $\cup\{\{u'\},{\bf x}\}\subset{\bf x}$ and ${\bf x}\subset\cup\{\{u'\},{\bf x}\}$, then $\cup\{\{u'\},{\bf x}\}={\bf x}$ and therefore $\cup\{\{v'\},{\bf x}\}={\bf x}$. But then so is $v'\in{\bf x}$.
  
  The converse is also true. Let $(\forall{\bf x})(u'\in{\bf x}\Rightarrow v'\in{\bf x})$ in particular for ${\bf x}=\{u\}$ and hence $v'\in\{u'\}$. Then $v'=u'$. $\Box$

The model of the theory $NBG[{\cal U}]^{(1)}$ is constructed as follows: by formulas (\ref{f1948}) and (\ref{f1949}) we define a class ${\bf H}[{\cal U}']$ which is a model for the theory $NBG[{\cal U}']$.

We now apply the operator of \grqq remembering\grqq\, $R=F^{-1}$ to all individuals from ${\cal U}^{(1)}$, which are constituents of the sets in ${\bf V}[{\cal U}']={\bf H}[{\cal U'}]\setminus{\cal U}^{(0)}$ or to the individuals themselves from ${\bf H}[{\cal U'}]\setminus{\cal U}^{(0)}$. In other words, in all objects $X\in{\bf H}[{\cal U'}]\setminus{\cal U'}^{(0)}$ we change all individuals from ${\cal U}^{(1)}$ to their corresponding skands, then treated as sets, which satisfy all axioms of the theory of $NBG[{\cal U}]^-$ except the foundation axiom, and by the axiom ${\bf SEA}$ contained in the universe ${\bf V}[{\cal U}]^-$ and are non-founded sets. Indeed, each individual $FX_{[\alpha_0,\alpha)}$ becomes a set $X_{[\alpha_0,\alpha)}$ containing a set $X_{[\alpha_0+1,\alpha)}$, which in turn contains a set $X_{[\alpha_0+2,\alpha)}$, and so on: $X_{(\alpha_0,\alpha)}\ni X_{[\alpha_0+1,\alpha)}\ni X_{[\alpha_0+2,\alpha)}\ni. ..$\,\, since the skand $X_{[\alpha_0,\alpha)}$ has infinite length $l=\alpha-\alpha_0\geq\omega$.

It is easy to see that after such changes, the class ${\bf V}[{\cal U}']$ of the theory $NBG[{\cal U}']$ turns into the class ${\bf V}[{\cal U}]^{(1)}$ of sets of both the $NBG[{\cal U}]^{(1)}$-founded and non-founded theory $NBG[{\cal U}]^{(1)}$. Moreover, it follows from Paragraph 1 that ${\bf U}[{\cal U}]^{(1)}={\bf V}[{\cal U}]^{(1)}\cup{\cal U}^{(0)}$ is a model for $NBG[{\cal U}]^{(1)}$. (For details, see \cite{l7}, p. 282-283, 302).

 These definitions and the description of the theory $NBG[{\cal U}]^{(1)}$ lead to the following theorem.
 
 {\bf Theorem 1.} {\it  The theory $NBG[{\cal U}]^{(1)}$ is consistent, provided of course that the theory $NBG^-$ is consistent.}

{\bf Proof.}
It is well known that if $NBG^-$ is consistent, then the extension $NBG$ of this theory is also consistent, obtained by attaching to it the founding axiom ${\bf FA}$ as the only new axiom. (See \cite{l7} Proposition 4.45 and Corollary 4.86, p. 283.) Bernays proved in \cite{l3210} (using a more complicated model than the one considered here, and proving the compatibility of the ${\bf FA}$-founding axiom) that this axiom is {\it independent} of all other axioms of the $NBG^-$ theory. So both the axiom ${\bf FA}$ and its negation of the axioms ${\bf SEA}$ and ${\bf PFA}$ can be attached to the $NBG^-$-theory without contradiction.

It follows from Proposition 1 that $NBG^-$ and $NBG$ are consistent if and only if $NBG[{\cal U}]^-$ and $NBG[{\cal U}]$ are respectively consistent. Given the above, the consistency of $NBG[{\cal U}]^-$ follows from the consistency of $NBG[{\cal U}]$. Similarly, by Bernays, the regularity axiom ${\bf RA}$, which in the presence of ${\bf CA}$ is equivalent to the founding axiom ${\bf FA}$, is independent of all other axioms of the theory $NBG[{\cal U}]^-$. Therefore, joining the axioms ${\bf SEA}$ and ${\bf PFA}$ that negate the axiom ${\bf FA}$ to $NBG[{\cal U}]^-$ leads to a consistent theory $NBG[{\cal U}]^{(1)}$ if of course the theory $NBG[{\cal U}]^-$ itself is consistent. The latter is consistent if the theory $NBG^-$ is consistent (See details in \cite{l56} or Chap. 4, Proposition 4.50 in \cite{l7}.) $\Box$

Since $\emptyset\subset{\cal U}$, from the axiom  ${\bf PFA}$ and the formulas $(\ref{f1948})$ and $(\ref{f1949})$, we get the following inclusions:
 \begin{equation}
 \label{f0990}
 {\bf V}\subset{\bf V}[{\cal U}]={\bf V}[{\cal U}]^{(0)}\subset{\bf V}[{\cal U}]^{(1)}\subset {\bf V}[{\cal U}]^-.
 \end{equation}
Simple examples (e.g., any skand $X_{[\alpha_0,\alpha)}$ of length $l=\alpha-\alpha_0\geq\omega$, seen as a set that is obviously non-founded and does not belong to ${\bf V}[{\cal U}]$) show that the first inclusion is proper. Next we will see that the last inclusion is also proper value.
$\Box$

\bigskip

\parindent=0 cm
\begin{center}
{\bf 6. $\nu$-generalized skands, the formal theory $NBG[{\cal U}]^{(\nu)}$ and its model in the universal  ${\bf U}[{\cal U}]^{(\nu)}$}
\end{center}
\bigskip

\parindent=0.5cm
In the previous paragraph, we extended the formal set theory of $NBG[{\cal U}]^-$ by attaching two axioms ${\bf SEA}$ and ${\bf PFA}$ which are less restrictive, than the founding axiom ${\bf FA}$ in $NBG[{\cal U}]$, and we obtained a new formal theory of already non-founded sets $NBG[{\cal U}]^{(1)}$ in contrast to the theory of $NBG[{\cal U}]$ where all sets are founded. We also extended this class of all founded sets ${\bf V}[{\cal U}]={\bf V}[{\cal U}]^{(0)}$ by using infinite length skands, which are essentially non-founded sets, to the universal class of all sets ${\bf V}[{\cal U}]^{(1)}$ of $NBG[{\cal U}]^{(1)}$ theory, with ${\bf V}[{\cal U}]^{(1)}\subset{\bf V}[{\cal U}]^-$, and increased the class of elements ${\bf U}[{\cal U}]^{(1)}={\bf V}[{\cal U}]^{(1)}\cup{\cal U}$ of this theory, which is a model for it. We can naturally continue such a process of extending the theory $NBG[{\cal U}]^-$ for every ordinal number $\nu\in{\bf On}$, $\nu>1$, defining inductively the theories $NBG[{\cal U}]^{(\nu')}$, $1\leq\nu'\leq\nu$, extending ${\bf V}[{\cal U}]^{(0)}$ with new {\it $\nu'$-generalized skands} of infinite length, $1\leq\nu'\leq\nu$. Moreover, there is the following sequence of inclusions: 

\begin{equation}
\label{f0911}
{\bf V}[{\cal U}]^{(0)}\subset{\bf V}[{\cal U}]^{(1)}\subset{\bf V}[{\cal U}]^{(2)}\subset...\subset{\bf V}[{\cal U}]^{(\nu)}\subset{\bf V[{\cal U}]^-}.
 \end{equation}

We do this by the following transfinite induction. 
 
 Let $\nu$ be an arbitrary ordinal such that $\nu>1$. Suppose we have already constructed formal set theories $NBG[{\cal U}]^{(\nu')}$ for every $\nu'$, $1\leq\nu'<\nu$, universal classes ${\bf V}[{\cal U}]^{(\nu')}$ and ${\bf U}[{\cal U}]^{(\nu')}$ of all sets and all elements of these theories, respectively, by means of new objects $\nu''$- generalized skands $X^{(\nu'')}\in{\bf S}'[{\cal U}]^{(\nu'')}$ of infinite length, $1\leq\nu''\leq\nu'$, and ${\bf U}[{\cal U}]^{(\nu'')}$ are models of the theories $NBG[{\cal U}]^{(\nu'')}$, $1\leq\nu''\leq\nu'$.

Now, to construct a formal set theory $NBG[{\cal U}]^{(\nu)}$ and the class ${\bf V}[{\cal U}]^{(\nu)}$ of all sets of this theory with the universal class of elements ${\bf U}[{\cal U}]^{(\nu)}={\bf V}[{\cal U}]^{(\nu)}\cup{\cal U}$, which will represent a model of the theory $NBG[{\cal U}]^{(\nu)}$, we need the following definitions of the new objects.
 
{\bf Definition 7.}
A $\nu$-{\it generalized skand}, $\nu>1$, with region {\it coupling} $[\alpha_0,\alpha)$ is an object $X^{(\nu)}$, which uniquely maps to each ordinal $\alpha'\in[\alpha_0,\alpha)$ the set $X^{(\nu)}(\alpha')\in\bigcup\limits_{0\leq\nu'<\nu}{\bf V}[{\cal U}]^{(\nu')}$ such that for each ordinal $0\leq\nu_0<\nu$ there exists an ordinal $\alpha'_0$, $0\leq\alpha'_0<\alpha$, and $X^{(\nu)}(\alpha'_0)\notin\bigcup\limits_{0\leq\nu'<\nu_0}{\bf V}[{\cal U}]^{(\nu')}$, i.e., $X^{(\nu)}$ is a new object,  and it has the following internal structure of the membership relation  $\in$ of its constituents and parts: For any ordinal $\alpha'\in[\alpha_0,\alpha)$ and for any element $x\in X^{(\nu)}(\alpha')$ there is a relation $x\in X'=X^{(\nu)}|_{[\alpha',\alpha)}$; and any ordinal $\alpha'\in[\alpha_0,\alpha)$ such that $\alpha'+1\in[\alpha_0,\alpha)$, there is a relation $X''=X^{(\nu)}|_{[\alpha'+1,\alpha)}\in X'$.

Here, for any ordinal $\alpha'\in[\alpha_0,\alpha)$, we denote by  $X^{(\nu)}|_{[\alpha',\alpha)}$ the restriction of the skand $X^{(\nu)}$ with a clutch  region  $[\alpha_0,\alpha)$ to the interval $[\alpha',\alpha)$, in particular, $X^{(\nu)}|_{[\alpha',\alpha)}$ coincides with  $X^{(\nu)}$ for $\alpha'=\alpha_0$ and is the {\it residue} of the  $\nu$-generalized skand $X^{(\nu)}$ on the clutch  region $[\alpha',\alpha)$ for each $\alpha'\in(\alpha_0,\alpha)$ Clearly, all residues of $X|_{[\alpha',\alpha)}$, $\alpha'\in(\alpha_0,\alpha)$, $\nu$-generalized skand $X^{(\nu)}$ are themselves $\nu'$-generalized skands, where $1\leq\nu'\leq\nu$ and for $\nu'=1$ are usual skands as in Definition 1, with clutch  regions $[\alpha',\alpha)$ are usual skands, with clutch  regions $[\alpha',\alpha)$ with induced relations of their constituents and parts that are portions of the relations and portions of the constituents of the skand  $X^{(\nu)}$  itself.

As before, a mapping  $f: [\alpha_0,\alpha)\rightarrow\bigcup\limits_{0\leq\nu'<\nu}{\bf V}[{\cal U}]^{(\nu')}$ given by the formula $f(\alpha')=X^{(\nu)}(\alpha')$ for each ordinal $\alpha'\in[\alpha_0, \alpha)$, we  call {\it the mapping associated with the skand} $X^{(\nu)}$.

Similarly to Definition $1'$, an equivalent extended definition of $\nu$-generalized skand is given below.

{\bf Definition 7$'$.}  
 By a $\nu$-{\it generalized skand} $X^{(\nu)}_{[\alpha_0,\alpha)}$  with a clutch  region $[\alpha_0,\alpha)$, we mean, as usual, a system of nested curly braces $\{_{\alpha_0}\{_{\alpha_0+1}...\{_{\alpha'}..._{\alpha'}\}..._{\alpha_0+1}\}_{\alpha_0}\}$, with ordinal indexes $\alpha'\in[\alpha_0,\alpha)$ (which we will omit for simplicity), between some or all of the neighboring brackets $\{\{$ and $\{\}$ (if the last pair of brackets exists), and for each ordinal $\alpha'\in[\alpha_0,\alpha)$ an arbitrary set $X^{(\nu)}_{\alpha'}\in\bigcup\limits_{0\leq\nu'<\nu}{\bf V}[{\cal U}]^{(\nu')}$ such that for each ordinal $0\leq\nu_0<\nu$ there exists an ordinal $\alpha'_0$, $0\leq\alpha'_0<\alpha$, and $X^{(\nu)}_{\alpha'_0}\notin\bigcup\limits_{0\leq\nu'<\nu_0}{\bf V}[{\cal U}]^{(\nu')}$, i.e., $X^{(\nu)}_{[\alpha_0,\alpha)}$ is a new object.  Since all Conway numbers in the segment $(-\frac{1}{\alpha'},-\frac{1}{\alpha'+1})$ is a proper class (not a set) there is an injective mapping $\psi_{\alpha'}:X^{(\nu')}_{\alpha'}\rightarrow(-\frac{1}{\alpha'},-\frac{1}{\alpha'+1})$ and thus the  injective mapping $\Psi:\bigsqcup\limits_{\alpha_0\leq\alpha'<\alpha}\psi_{\alpha'}\rightarrow(-2,2)$ of the set $\bigsqcup\limits_{\alpha_0\leq\alpha'<\alpha}\{X^{(\nu')}_{\alpha'}\}$ to the proper class $(-2,2)$ of Conway numbers. Then $X^{(\nu')}_{[\alpha_0,\alpha)}\stackrel{def}{=}(\psi\cup\Psi)^{-1}(Im(\psi\cup\Psi))$ with induce linear ordering by the linear ordering on $[-2,2]$ and the mapping $(\psi\circ\Psi)^{-1}$ is called $\nu$-{\it generalized skand} with a clutch  region $[\alpha_0,\alpha).$  
   The sets $X^{(\nu')}_{\alpha'}$ we call  {\it components} of the $\nu$--generalized skand $X^{(\nu)}_{[\alpha_0,\alpha)}$ and denote by $X^{(\nu)}_{\alpha'}=\{x^{\alpha'}_0,x^{\alpha'}_1,x^{\alpha'}_2,...,x^{\alpha'}_\lambda,...\}$ for each index $\alpha'\in[\alpha_0,\alpha)$.

In other words, the $\nu$-generalized skand $X^{(\nu)}_{[\alpha_0,\alpha)}$ is an aggregate, or chained tuple of {\it sets previously constructed}, consisting of elements in the universe ${\bf U}[{\cal U}]^{(\nu)}=\bigcup\limits_{0\leq\nu'<\nu}{\bf V}[{\cal U}]^{(\nu')}\cup{\cal U}$ or containing no elements, and at least one of these previously constructed sets does not belong to the class $\bigcup\limits_{0\leq\nu'<\nu''}{\bf V}[{\cal U}]^{(\nu')}$ for all ordinals $\nu''<\nu$.

The general form of the $\nu$-generalized skand is as follows: \begin{equation}
\label{f0233}
\begin{array}{l}

X^{(\nu)}_{[\alpha_0,\alpha)}
=\\=\{x^{\alpha_0}_0,x^{\alpha_0}_1,...,\{x^{\alpha_0+1}_0,x^{\alpha_0+1}_1,...,\{...\{x^{\alpha'}_0,x^{\alpha'}_1,...,\{...\{x^{\alpha-1}_0,x^{\alpha-1}_1,...\}...\}...\}...\}\}\},
\end{array}
\end{equation}
where components $X^{(\nu)}_{\alpha'}=\{x^{\alpha'}_0,x^{\alpha'}_1,...\}$ are either empty or are elements of the class $\bigcup\limits_{0\leq\nu'<\nu}{\bf V}[{\cal U}]^{(\nu')}$, and at least one of these components does not belong to the class $\bigcup\limits_{0\leq\nu'<\nu''}{\bf V}[{\cal U}]^{(\nu')}$ for all ordinals $\nu''<\nu$.

{\bf Proposition 9} {\it There exists a class ${\bf S}[{\cal U}]^{(\nu)}$ of all $\nu$-generalized skands $X^{(\nu)}_{[0,\alpha)}$ with arbitrary components and lengths $0<\alpha<\Omega$ and it is a proper value.}

{\bf Proof} is similar to the proof of Proposition 8.

{\bf Definition 8.} Two $\nu$-generalized skands  $X^{(\nu)}$ and $Y^{(\nu)}$ with clutch  regions $[\alpha_0,\alpha)$ and $[\beta_0,\beta)$ respectively are called {\it equal}, $X^{(\nu)}=Y^{(\nu)}$ if the corresponding mappings $f$ and $g$ associated to the skandas $X^{(\nu)}$ and $Y^{(\nu)}$are equal as counting functions, i.e.   there exists a similar mapping (isomorphism) $\varphi:[{\alpha_0,\alpha})\rightarrow[{\beta_0,\beta})$  of the intervals $[{\alpha_0,\alpha})$ and $[{\beta_0,\beta})$ that preserves (preserves) the structures of well-ordered sets $[{\alpha_0,\alpha})$ and $[{\beta_0,\beta})$ and the following usual equality of mappings is satisfied: $f=g\circ\varphi$. 

{\bf Definition 8$'$.}
Two $\nu$-generalized skands $X^{[\nu)}_{[\alpha_0,\alpha)}$ and $Y^{[\nu)}_{[\beta_0,\beta)}$ let's call {\it equal} if the segments $[\alpha_0,\alpha)$ and $[\beta_0, \beta)$ are isomorphic as well-ordered sets and the corresponding components $X^{(\nu)}_{\alpha'}=Y^{(\nu)}_{\beta'}$, where $\varphi: [\alpha_0,\alpha)\rightarrow[\beta_0,\beta)$ the corresponding isomorphism and $\beta'=\varphi(\alpha')$, $\alpha_0\leq\alpha'<\alpha$, $\beta_0\leq\beta'<\beta$, are equal as sets in $\bigcup\limits_{0\leq\nu'<\nu}{\bf V}[{\cal U}]^{(\nu')}$.

\smallskip
{\bf $\nu$-Generalized skand existence axiom  ({\bf $\nu$-GSEA}, $1<\nu<\Omega$)}.
For any ordinals $\alpha_0,\alpha\in{\bf On}$ such that $\alpha-\alpha_0\geq\omega$, and any mapping $f:[\alpha_0,\alpha)\rightarrow\bigcup\limits_{0\leq\nu'<\nu}{\bf V}[{\cal U}]^{(\nu')}$ such that $f([\alpha_0,\alpha))$ is not a subclass of $\bigcup\limits_{0\leq\nu'<\nu_0}{\bf V}[{\cal U}]^{(\nu')}$, for any ordinal $1\leq\nu_0<\nu$, there exists {\it $\nu$-generalized skand} $X^{(\nu)}$ with a clutch  region $[\alpha_0,\alpha)$, for which $X^{(\nu)}(\alpha')=f(\alpha')$ holds for any ordinal $\alpha'\in[\alpha_0,\alpha)$ and $X^{(\nu)}=X^{(\nu)}_{\alpha_0}\cup\{X^{(\nu)}|_{[\alpha_0+\alpha)}\}\in{\bf V}[{\cal U}]^-$.

\smallskip

 {\bf $\nu$-Generalized pseudo-founded axiom ({\bf $\nu$-GPFA})}.
For any set $X$, any maximal descending $\in$-sequence $X\ni x_1\ni x_2\ni x_3\ni...\ni x_{n-1}\ni x_n\ni.... $ has the following property: if the sequence is finite, then its last element $x_n$ is either an empty set $\emptyset$ or an individual $u\in{\cal U}$; if it is infinite, then there exists such a smallest number $n$ that $x_n$ is an infinite skand or $\nu'$-generalized skand of infinite length for some ordinal $1<\nu'\leq\nu$.

  \smallskip
 Let's now denote by ${\bf S}'[{\cal U}]^{(\nu)}$ the class of all $\nu'$-generalized skands whose length is $l\geq\omega$, $1\leq\nu'\leq\nu$, and by ${\bf V}[{\cal U}]^{(\nu)}$ -- the class of all sets of theory $NBG[{\cal U}]^{(\nu)}$.

  Since we postulate the existence of new objects and sets in $NBG[{\cal U}]^{(\nu)}$, we must refine the equality of objects in the extensionality axiom ${\bf EA}$ of the theory $NBG[{\cal U}]^{(0)}$ by replacing it with a more precise axiom:

{\bf $\nu$-Generalized exact extensionality axiom ($\nu$-GEEA}). $(\forall X)(\forall Y)(\forall Z)(X)=Y\wedge X\in Z\Rightarrow Y\in Z)$, where the two sets $X$ and $Y$, whose elements can be individuals, founded sets and non-founded sets, as $\nu$-generalized skands, are equal, if for every element $x\in X$ there exists an element $y\in Y$ such that $x=y$ and for every element $y\in Y$ there exists an element $x\in X$ such that $y=x$, where by the equality \grqq$=$\grqq\, we mean the following: 

 $1)$ the equality of individuals when $x,y\in {\cal U}^{(0)}$, i. e.,
$(x=y)\Leftrightarrow(\forall z)(x\in z\Leftrightarrow y\in z)$; 
 
 $2$) the equality of the founded sets when $x,y\in{\bf V}[{\cal U}]^{(0)}$, i. e., by the equality of sets in the extensionality axiom ${\bf EA}$ in $NBG[{\cal U}]$;
 
 $3)$ the equality of $\nu'$-generalized skands when  $x,y\in {\bf S}'[{\cal U}]^{(\nu)}$, $1<\nu'\leq\nu$,  i.e., by Definition 7; 
 
 $4)$ the equality of sets in $x,y\in{\bf V}[{\cal U}]^{(\nu)}$ when  $x,y\notin{\cal U}$, $x,y\notin{\bf V}[{\cal U}]^{(0)}$ and $x,y\notin{\bf S}'[{\cal U}]^{(\nu)}$, i.e., by Definition 7, i.e., the next iteration using $1)$, $2)$, $3)$ in $\nu$-{\bf EEA}.

\smallskip

Now we construct a model of this theory $NBG[{\cal U}]^{(\nu)}$ by formula (\ref{f1949}), assuming instead of $L$ any subclass of the class ${\cal U}'={\bf S}'[{\cal U}]^{(\nu)}\cup\,{\cal U}$, which is a set. Then the class ${\bf H}[{\cal U}]^{(\nu)}$ consisting of all elements $M$ for which there exists a set $L$ and an ordinal $\gamma$ such that $M\in\Xi^\gamma_{{\cal U}'}$ is a model for the theory $NBG[{\cal U}]^{(\nu)}$.

It is easy to see that after such changes, the class ${\bf V}[{\cal U}']$ of the theory $NBG[{\cal U}']$ becomes the class ${\bf V}[{\cal U}]^{(\nu)}$ of sets of both $NBG[{\cal U}]^{(\nu)}$-founded and non-founded sets of the theory $NBG[{\cal U}]^{(\nu)}$. Moreover, it follows from Paragraph 1 that ${\bf U}[{\cal U}]^{(\nu)}={\bf V}[{\cal U}]^{(\nu)}\cup{\cal U}$ is a model for $NBG[{\cal U}]^{(\nu)}$. (For details, see \cite{l7}, p. 282-283, 302).

 These definitions and the description of the theory $NBG[{\cal U}]^{(\nu)}$ lead to the following theorem.
 
 {\bf Theorem 2} {\it The theory $NBG[{\cal U}]^{(\nu)}$  is consistent, provided of course that the theory $NBG^-$ is consistent.}

{\bf Proof} is similar to the proof of Theorem 1. 
 $\Box$.

The following inclusion follows from the definitions and construction:
\begin{equation}
\label{f8888}
{\cal U}={\cal U}^{(0)}\subset{\cal U}\cup{\bf S}'[{\cal U}]^{(1)}\subset{\cal U}\cup{\bf S}'[{\cal U}]^{(2)}\subset...\subset{\cal U}\cup{\bf S}'[{\cal U}]^{(\nu)}={\cal U}'.
\end{equation}
Clearly, by the formula (\ref{f8888}) and the fact that the sets in theory $NBG[{\cal U}]^{(\nu)}$ satisfy all axioms of theory $NBG[{\cal U}]^-$, we get the following inclusions:

 \begin{equation}
 \label{f0990}
 {\bf V}[{\cal U}]={\bf V}[{\cal U}]^{(0)}\subset{\bf V}[{\cal U}]^{(1)}\subset{\bf V}[{\cal U}]^{(2)}\subset...\subset{\bf V}[{\cal U}]^{(\nu)}\subset {\bf V}[{\cal U}]^-.
 \end{equation}

It is also clear that these inclusions are indeed proper states, since the next step of $\nu+1$ transfinite induction is possible. $\Box$

\smallskip

Since we have constructed $NBG[{\cal U}]^{(\nu)}$ and ${\bf V}[{\cal U}]^{(\nu)}$  for all  $\nu\in {\bf On}$, then we can combine all these classes  ${\bf V}[{\cal U}]^{(\nu)}$  and get the class $\bigcup\limits_{\nu\in{\bf On}}{\bf V}[{\cal U}]^{(\nu)}$as the following proper inclusions: 
\begin{equation}
\label{f8889}
{\bf V}[{\cal U}]={\bf V}[{\cal U}]^{(0)}\subset{\bf V}[{\cal U}]^{(1)}\subset{\bf V}[{\cal U}]^{(2)}\subset...\subset{\bf V}[{\cal U}]^{(\nu)}\subset...\subset\bigcup\limits_{\nu\in{\bf On}}{\bf V}[{\cal U}]^{(\nu)}\subset {\bf V}[{\cal U}]^-.
\end{equation}

 \smallskip
 
 Let's now construct a formal set theory for the symbol $\Omega$ mentioned in paragraph 1, point 2, for which the inequalities  $0\leq\nu<\Omega$ are satisfied for any ordinal $\nu\in{\bf On}$. Thus we can consider, within the framework of the $NBG[{\cal U}]^-$theory, intervals $[\alpha_0,\Omega)$ of ordinal numbers $\alpha'$ such that $\alpha_0\leq\alpha'<\Omega$. By $\Omega$-skand we mean any $\nu$-skand where $1\leq\nu<\Omega$. By ${\bf V}[{\cal U}]^{(\Omega)}$ we denote the sum of $\bigcup\limits_{\nu\in{\bf On}}{\bf V}[{\cal U}]^{(\nu)}$ classes of  ${\bf V}[{\cal U}]^{(\nu)}$, by ${\bf S}[{\cal U}]^{(\Omega)}$ the sum $\bigcup\limits_{\nu\in{\bf On}}{\bf S}[{\cal U}]^{(\nu)}$ of classes ${\bf S}[{\cal U}]^{(\nu)}$, and by ${\bf S}'[{\cal U}]^{(\Omega)}$ the sum $\bigcup\limits_{\nu\in{\bf On}}{\bf S}'[{\cal U}]^{(\nu)}$ of classes ${\bf S}'[{\cal U}]^{(\nu)}$.

 \smallskip
 
The formal theory of non-founded sets $NBG[{\cal U}]^{(\Omega)}$ is a first-order theory with its own axioms $NBG[{\cal U}]^-$ and attached to them two axioms ${\bf \Omega}$-${\bf SAE}$ and ${\bf \Omega}$-${\bf PFA}$ with a refinement of the exact extensional axiom ${\bf \Omega}$-${\bf EEA}$. 

The axiom ${\bf \Omega}$-${\bf SAE}$ is that for any ordinals $\alpha_0,\alpha\in{\bf On}$ such that $\alpha-\alpha_0\geq\omega$, and for any mapping $f:[\alpha_0,\alpha)\rightarrow{\bf V}[{\cal U}]^{(\Omega)}$ there exists  an ordinal $\nu$, $1\leq\nu<\Omega$,  such that a $\nu$-generalized skand $X^{(\nu)}$ has the following property: $X^{(\nu)}(\alpha')=f(\alpha')$, for each ordinal $\alpha_0\leq\alpha'<\alpha$.

The axiom ${\bf \Omega}$-${\bf PFA}$ is that for any set $X$, any maximal descending  $\in$-sequence $X\ni x_1\ni x_2\ni x_3\ni...\ni x_{n-1}\ni x_n\ni...$ has the following property: if the sequence is finite, then its last element $x_n$ is either an empty set $\emptyset$ or an individual $u\in{\cal U}$; if it is infinite, then there exists such a smallest number $n$ that $x_n$ is a $\nu$-generalized skand of infinite length for some ordinal $1\leq\nu\leq\Omega$.

 The axiom ${\bf \Omega}$-${\bf EEA}$ is the following: $(\forall X)(\forall Y)(\forall Z)(X)=Y\wedge X\in Z\Rightarrow Y\in Z)$, where the two sets $X$ and $Y$, whose elements can be individuals, founded sets and non-founded sets, as $\nu$-generalized skands, are equal, if for every element $x\in X$ there exists an element $y\in Y$ such that $x=y$ and for every element $y\in Y$ there exists an element $x\in X$ such that $y=x$, where by the equality \grqq$=$\grqq\, we mean the following:

$1)$ the equality of individuals when $x,y\in {\cal U}^{(0)}$, i.e., $(x=y)\Leftrightarrow(\forall z)(x\in z\Leftrightarrow y\in z)$; 
 
 $2$) the equality of the founded sets when $x,y\in{\bf V}[{\cal U}]^{(0)}$, i.e., by the equality of sets in the extensionality axiom ${\bf EA}$ in $NBG[{\cal U}]$;
 
 $3)$ the equality of $\Omega$-generalized skands when $x,y\in {\bf S}'[{\cal U}]^{(\Omega)}$; 
 
 $4)$ the equality of sets in $x,y\in{\bf V}[{\cal U}]^{(\Omega)}$ when $x,y\notin{\cal U}$, $x,y\notin{\bf V}[{\cal U}]^{(0)}$ and $x,y\notin{\bf S}'[{\cal U}]^{(\Omega)}$, i.e., the next iteration using $1)$, $2)$, $3)$ in ${\bf \Omega}$-{\bf EEA}.

Assuming ${\cal U}'=\bigcup\limits_{1\leq\nu<\Omega}{\bf S}'{\cal U}^{(\nu)}\cup{\cal U}$ and extending the formal theory $NBG^-$ to the formal theory $NBG[{\cal U}]^{(\Omega)}$, whose model ${\bf U}[{\cal U}]^{(\Omega)}$ s constructed similarly to the previous constructions by formulas (\ref{f1948}) and  (\ref{f1949}). We omit the details, except that the following proper value embeddings take place: \begin{equation}
 \label{f09909}
 {\bf V}[{\cal U}]={\bf V}[{\cal U}]^{(0)}\subset{\bf V}[{\cal U}]^{(1)}\subset{\bf V}[{\cal U}]^{(2)}\subset...\subset{\bf V}[{\cal U}]^{(\nu)}\subset...\subset{\bf V}[{\cal U}]^{(\Omega)}\subset {\bf V}[{\cal U}]^-.
 \end{equation}

 These definitions and the description of the theory $NBG[{\cal U}]^{(\Omega)}$ lead to the following theorem.
 
 {\bf Theorem 3.} {\it The theory $NBG[{\cal U}]^{(\Omega)}$ is consistent, provided of course that the theory $NBG^-$ is consistent.}

{\bf Proof} is similar to the proof of Theorem 1. $\Box$

In the formula $(\ref{f09909})$, the last inclusion is proper since, for example, a Mirimanoff's set of the form $X=\{a,b,c,...,\{X\},X\}$, where the skands $X=X_{[0,\omega)}$, do not belong to the class ${\bf V}^{(\Omega)}$ (consider for simplicity that ${\cal U}=\emptyset$), but belong to the class ${\bf V}[{\cal U}]$ by virtue of some less rigid restriction axioms than ${\bf \Omega}$-${\bf SEA}\&{\bf PFA}$, e. g.,  such as the Aczel's anti-foundation axiom ${\bf AFA}$ 
in \cite{l93}. Thus Axel's non-founded set theory does not coincide with any of our non-founded set theories. On the other hand, it does not even include $NBG^{(1)}$, since in $NBG^{(2)}$ the skand ${\bf E}_{[0,\omega)}$ whose components are all equal to the skand ${\bf e}_{[0,\omega)}$, by the definition of equality of generalized $2$-skands, ${\bf E}_{[0,\omega)}=\{{\bf e}_{[0,\omega)},\{{\bf e}_{[0,\omega)},\{...\,\,\,...\}\}\}\not={\bf e}_{[0,\omega)}$, and in Aczel's Theorem   ${\bf E}_{[0,\omega)}={\bf e}_{[0,\omega)}$, since the equation $X=\{{\bf e}_{[0,\omega)},X\}$ has a single solution $X={\bf e}_{[0,\omega)}$ and since the obvious equality of the sets ${\bf e}_{[0,\omega)}=\{{\bf e}_{[0,\omega)}\}=\{{\bf e}_{[0,\omega)},{\bf e}_{[0,\omega)}\}$  is satisfied. Hence, $X={\bf e}_{[0,\omega)}$, (see: \cite{l93}, p. 8). $\Box$

\bigskip

\begin{center}
{\bf  Part two}
\end{center}

\begin{center}
{\bf 1. Definition of coskands}
\end{center}

The notion of a skand discussed above allows us to dualize it and obtain a new mathematical object, which in a special case was not previously considered {\it as definable}, e.g., see: \grqq property of an expression having the form  $...\{...\{\{0\}\}...\}...$ is not  \grqq definable\grqq, (see the translator's comment in \cite{l144}, footnote $3)$, p. 109). 

Here it will be a well-defined object as an individual when system of pairs of curly braces $...\{...\{\{0\}\}...\}...$ has exactly $\omega$ pieces and when it has $\omega+1$ pieces, i.e., $\{...\{...\{\{0\}\}...\}...\}$, it is a one-element set whose element is a trivial coskand $...\{...\{\{0\}\}...\}...$.

The axioms of set theory $NBG[{\cal U}]$ do not, in general, specify whether the class ${\cal U}$ of individuals is a set or a proper class, except that it follows from the axioms that there exists a proper  class ${\bf U}[{\cal U}]$ of all elements; the existence of the class ${\cal U}$ or, equivalently, the existence of the universal class ${\bf V}[{\cal U}]$, is given by a separate axiom.
 
A sufficient number (abundance) of individuals is provided, for example, by the axiom  of plenitude $({\bf AP})$.

 {\bf Axiom of Plenitude} for the class ${\cal U}$. For every set $S$ there exists an injection $f:S\rightarrow{\cal U}$ whose image $f(S)$ is disjoint from $S$, \cite{l98}, p. 23.

Axiom of Plenitude as it stands does not guarantee that the collections we will need do exist. For this purpose there is the following stronger form.

 {\bf Strong Axiom of Plenitude} for the class ${\cal U}$. There is an operation {\bf new}$(a,b)$ so that

1. For all sets $a$ and all subsets, $b\subseteq{\cal U}$,  {\bf new}$(a,b)\in{\cal U}\setminus b$. 

2. For all $a\not= a'$ and all $b\subseteq{\cal U}$, {\bf new}$(a,b)\not=${\bf new}$(a',b)$, see: \cite{l98}, p. 23.

\smallskip

The {\bf existence coskand axiom} bellow   will provide a huge and diverse number (abundance) of well-distinguished individuals, for which the usual axiom of choice will be sufficient. The importance of the existence of a sufficient number of individuals was necessary to prove the independence of the axiom of choice from the other axioms.

For a symmetric statement of the theory of coskands, we will below fix some class ${\cal U}$ of individuals, which may be an empty, finite and infinite set, or a proper class and  it does not matter for the initial definitions but without individuals we will intruduce now, i.e., without coskands.

{\bf Definition 9}   {\it Coskand} with a  {\it clutch  region} $[\alpha_0,\alpha)$  is an object $X$ such that it associates uniquely for each ordinal $\alpha'\in [\alpha_0,\alpha)$  a set $X(\alpha')\in{\bf V}[{\cal U}]$ and has the following internal structure of the membership relation $\in$ of its constituents and parts: for each ordinal $\alpha'\in[\alpha_0,\alpha)$ and each element $x\in X(\alpha')$ there is a relation $x\in X|_{[\alpha_0,\alpha'+1)}$ and  for each ordinal $\alpha'\in(\alpha_0,\alpha)$ there is a relation $X|_{[\alpha_0,\alpha')}\in X|_{[\alpha_0,\alpha'+1)}$, where $X|_{[\alpha_0,\alpha')}$ is the restriction of $X$ to the interval $[\alpha_0,\alpha')$  and $X|_{[\alpha_0,\alpha'+1)}$ is the restriction of $X$ to the interval $[\alpha_0,\alpha'+1)$.

For each ordinal $\alpha'\in(\alpha_0,\alpha)$  the coskand $X|_{[\alpha_0,\alpha')}$ will be  called a {\it  partial  coskand} of $X$. We also denote sets $X(\alpha')$ by $X_{\alpha'}$, $\alpha_0\leq\alpha'<\alpha$,  calling them {\it constituents} of the coskand $X$ and coskands $X|_{[\alpha_0,\alpha')}$, $\alpha_0<\alpha'<\alpha$,  we call {\it parts} of the coskand $X$.
It is also clear that all partial  coskands $X|_{[\alpha_0,\alpha')}$, $\alpha'\in(\alpha_0,\alpha)$, of the coskand $X$ are real coskands  with  clutch  regions $[\alpha_0,\alpha')$ and with induced relations of their constituents and parts that are portions of the relationship and portions of the constituents of the coskand $X$ itself.

It is also clear that a coskand $X$ with a clutch  region $[\alpha_0,\alpha)$ uniquely defines a mapping $f:[\alpha_0,\alpha)\rightarrow {\bf V}[{\cal U}]$ if we put $f(\alpha')=X(\alpha')$ for each ordinal $\alpha'\in[\alpha_0,\alpha)$. We will call this mapping $f$ the {\it mapping associated with the coskand} $X$.

Let's call {\it the length} of a caskand $X$ with a clutch  region $[\alpha_0,\alpha)$ the ordinal number $l=\alpha-\alpha_0$ and distinguish here between {\it finite} caskands, which have length equal to an ordinal $l\in(0,\omega)$, clearly, being usual sets in ${\bf V}[{\cal U}]$ and {\it infinite} coskands which have length $l\geq\omega$, where $\omega$ is the smallest countable ordinal, being new objects.

{\bf Definition 10.} Two coskands $X$ and $Y$ with clutch  regions $[\alpha_0,\alpha)$ and $[\beta_0,\beta)$, respectively, are called {\it  equal}, we write $X=Y$, if their corresponding mappings $f$ and $g$ associated to the coskands $X$ and $Y$ are equal as counting functions,  i.e., $f=g\circ\varphi$, where $\varphi:[\alpha_0,\alpha)\rightarrow[\beta_0,\beta)$ is an isomorphism of well-ordered sets.

 \bigskip
 
The description of coskands and their equalities can be given in a less formal way, but in a more visual and convenient way for their recording and working with them. Let's call this form of a coskand as {\it  unfolded}, in the sense that in the record of the caskand itself the  clutch  region will be given and all its elements, or more precisely, all the relations $\in$ of  belonging to its constituents and parts will be specified.

Consider first a nested \grqq strictly increasing\grqq$\,$  system of pairs of curly braces $...\{_{\alpha'}...\{_{\alpha_0+1}\{_{\alpha_0}\}_{\alpha_0}\}_{\alpha_0+1}...\}_{\alpha'}...$ strictly increasing in the sense that for any ordinals $\alpha'<\alpha''$, $\alpha',\alpha''\in[\alpha_0,\alpha)$, the pair of brackets $\{_{\alpha'}...\}_{\alpha'}$ \grqq lies inside\grqq\, the pair of brackets $\{_{\alpha''}...\}_{\alpha''}$. The possibility of such a system of nested brackets can be modeled on the ends of segments of Conway numbers \cite{l22} $[-\alpha',\alpha']$, $\alpha'\in[\alpha_0,\alpha)$ when $\alpha_0\not=0$ and $[-\frac{1}{2},\frac{1}{2}]$ when $\alpha_0=0$, where a Conway number $-\alpha'$ corresponds to the open bracket $\{_{\alpha'}$, and a Conway number $\alpha'$ corresponds to the closed bracket $\}_{\alpha'}$, and if $\alpha_0=0$, then we take the ends of the segment $[-\frac{1}{2},\frac{1}{2}]$, i.e., Conway numbers $-\frac{1}{2},\frac{1}{2}$ will correspond to a pair of curly brackets $\{_0\}_0$, respectively. Note only that here too, unlike in Conway numbers, curly brackets play a syntactic role.

{\bf Remark 6}. Note only that unlike Conway numbers, which are mathematical objects, curly brackets are not at all, but play here a syntactical role, indicating only the membership relations $\in$ of  constituents and parts of the coskand: everything between adjacent curly brackets $\{_{\alpha'+1}\{_{\alpha'}$, where $\alpha'\in[\alpha_0,\alpha)$ such that $\alpha'+1\in[\alpha_0,\alpha)$,  will denote the elements of the set $X_{\alpha'+1}\cup\{ X|_{[\alpha_0,\alpha')}\}$, generated by the coskand $X$ with a clutch  region $[\alpha_0,\alpha)$ and everything between  the last pair of curly braces $\{_{\alpha_0}\}_{\alpha_0}$ will denote the elements of the set $X_{\alpha_0}$. 

We perform this replacement of Conway numbers with curly brackets as follows. Denote by $\psi$ a mapping $\psi:\{_{\alpha_0}\{_{\alpha_0+1}...\{_{\alpha'}...\,\,\,\,...\}_{\alpha'}...\}_{\alpha_0+1}\}_{\alpha_0}\rightarrow [-\alpha,\alpha]$ of a strictly increasing sequence of nested curly braces to the close segment $[-\alpha,\alpha]$ of all Conway numbers, given by the two formulas $\psi(\{_{\alpha'})=-\alpha'$ and $\psi(\}_{\alpha'})=\alpha'$ when $\alpha'\not=0$, and $\psi(\{_0)=-\frac{1}{2}$ and $\psi(\}_0)=\frac{1}{2}$, when $\alpha'=0$. It is evidently a linear ordered isomorphism between the set of all above nested curly braces to the set of corresponding Conway numbers in the closed segment $[-\alpha,\alpha]$.

{\bf Definition 9$'$.} A trivial coskand ${\bf e}_{[\alpha_0,\alpha)}$ with a clutch  region $[\alpha_0,\alpha)$ is a system of nested \grqq strictly increasing\grqq\, pairs of curly braces $...\{_{\alpha'}...\{_{\alpha_0+1}\{_{\alpha_0}\}_{\alpha_0}\}_{\alpha_0+1}...\}_{\alpha'}...$,
$\alpha'\in[\alpha_0,\alpha)$, which for simplicity we reduce to $...\{...\{\{\}\}...\}...\,$ . To construct a nontrivial coskand $X_{[\alpha_0,\alpha)}$ with the same clutch  region $[\alpha_0,\alpha)$ we consider for each ordinal $\alpha'\in[\alpha_0,\alpha)$ an arbitrary set $X_{\alpha'}$ from the universe ${\bf V}[{\cal U}]$, i.e., $X_{\alpha'}\in{\bf V}[{\cal U}]$.  Since all Conway numbers in the segment $(-\alpha'-1,-\alpha')$, for each ordinal $\alpha_0\leq\alpha'<\alpha$, when $\alpha_0\not=0$, and $(-\frac{1}{2},\frac{1}{2})$, when $\alpha_0=0$, is a proper class (not a set) there is an injective mapping $\psi_{\alpha'}:X_{\alpha'}\rightarrow(-\alpha'-1,-\alpha')$, for each ordinal $\alpha_0<\alpha'<\alpha$, and $\psi_{\alpha_0}:X_{\alpha_0}\rightarrow(-\alpha_0,\alpha_0)$ and in particular case $\alpha_0=0$ an injective mapping  $\psi_0:X_0\rightarrow(-\frac{1}{2},\frac{1}{2})$, thus the  injective mapping $\Psi:\bigsqcup\limits_{\alpha_0\leq\alpha'<\alpha}\psi_{\alpha'}\rightarrow(-\alpha,\alpha)$ of the set $\bigsqcup\limits_{\alpha_0\leq\alpha'<\alpha}\{X_{\alpha'}\}$ to the proper class $(-\alpha,\alpha)$ of Conway numbers. Then $X_{[\alpha_0,\alpha)}\stackrel{def}{=}(\psi\cup\Psi)^{-1}(Im(\psi\cup\Psi))$ with induce linear ordering by the linear ordering on $[-\alpha,\alpha]$ and the image Im$((\psi\circ\Psi)^{-1})$ of the mapping $(\psi\circ\Psi)^{-1}$ is called a {\it coskand} with a clutch  region $[\alpha_0,\alpha).$

\smallskip

One can see that if  $\alpha$ is a limit ordinal, then a coskand $X_{[\alpha_0,\alpha)}$ is an individual because there is no relation $x\in X_{[\alpha_0,\alpha)}$ for any element $x\in{\bf U}[{\cal U}]$; if  $\alpha$ is not a limit ordinal, then a coskand $X_{[\alpha_0,\alpha)}$ is a founded set $X_{\alpha-1}\in{\bf V}[{\cal U}]$ when $l=\alpha-\alpha_0<\omega$, and if   $l=\alpha-\alpha_0>\omega$, then a coskand $X_{[\alpha_0,\alpha)}$ is  a founded set with a new individual  $X_{[\alpha_0,\alpha-1)}$ or, possibly,  a set, which has been constructed by using coskands with limit ordinal as above,  because there is  relation $X_{[\alpha_0,\alpha'})\in X_{[\alpha_0,\alpha'-1)}$.

In other words, a coskand $X_{[\alpha_0,\alpha)}$ is an {\it  aggregate} or a {\it  chained increasing tuple} of {\it is the founded sets} or {\it is their absence} in that tuple. A trivial coskand is a tuple with empty components; an arbitrary coskand is a tuple with random components that are elements of the universe ${\bf V}[{\cal U}]$ of the founded sets.

Any coskand $X_{[\alpha_0,\alpha)}$ with components $X_{\alpha'}$, $\alpha_0\leq\alpha'<\alpha$, in Definition $9'$ uniquely defines a coskand $X$ with a clutch  region $[\alpha_0,\alpha)$ if we put $X(\alpha')=X_{\alpha'}$ and  $x\in X|_{[\alpha_0,\alpha'+1)}\equiv x\in X_{\alpha'}$, for each ordinal $\alpha'\in[\alpha_0,\alpha)$, and $X|_{[\alpha_0,\alpha')}\in X|_{[\alpha_0,\alpha'+1)}\equiv X_{[\alpha_0,\alpha')}\in X_{[\alpha_0,\alpha'+1)}$ for any ordinal $\alpha'\in(\alpha_0,\alpha)$ such that $\alpha'+1\in[\alpha_0,\alpha)$.   The reverse is also true: every coskand $X$ with a clutch  region  $[\alpha_0,\alpha)$ of Definition 9 uniquely defines a coskand $X_{[\alpha_0,\alpha)}$ of Definition $9'$ if we put $X_{\alpha'}=X(\alpha')$ for every $\alpha'\in[\alpha_0,\alpha)$ and insert elements of these sets between neighboring open or last pair of brackets, if any, of the trivial skand 
$\bf e_{[\alpha_0,\alpha)}$.

 Thus, the coskands $X$ and $X_{[\alpha_0,\alpha)}$ in both definitions are one and the same object, which uniquely maps for each ordinal $\alpha'\in[\alpha_0,\alpha)$ a founded set, with certain relations of $\in$ membership of its constituents, or parts, although represented in fact in different forms (similarly to how a function can be represented in different forms: A law that maps a single value to each argument; its graph; an implicit equation; a table; a two-line matrix; etc.).

The general view of the $X$ coskand in unfolded form is as follows: 
\begin{equation}
\label{f123}
X_{[\alpha_0,\alpha)}
=...\{x^{\alpha'}_0,x^{\alpha'}_1,...,\{...\{x^{\alpha_0+1}_0,x^{\alpha_0+1}_1,...,\{x^{\alpha_0}_0,x^{\alpha_0}_1,...\}\}...\}\}...,
\end{equation}
where the components of $X_{\alpha'}$ are either empty, i.e., $X_{\alpha'}=\emptyset$,  or consist of elements  $x^{\alpha'}_0,x^{\alpha'}_1,...$ of the class ${\bf U}[{\cal U}]$ such that their set $\{x^{\alpha'}_0,x^{\alpha'}_1,...\}=X_{\alpha'}$ is an element of the class ${\bf V[{\cal U}]}$, i.e., $X_{\alpha'}\in {\bf V[{\cal U}]}$, $\alpha_0\leq\alpha'<\alpha$. Thus the general formula $(\ref{f123})$ highlights the case where the ordinal $\alpha'$ has an antecedent; if $\alpha'$ has no predecessor, the notation of the formula $(\ref{f123})$looks a little different and is clear as.

\bigskip

{\bf Definition 10$'$.} Two coskands $X_{[\alpha_0,\alpha)}$ and $Y_{[\beta_0,\beta)}$ are called {\it is equal} when there exists a similar mapping (isomorphism): $\varphi:[{\alpha_0,\alpha})\rightarrow[{\beta_0,\beta})$ of the intervals $[{\alpha_0,\alpha})$ and $[{\beta_0,\beta})$ that preserves the structures of the well-ordered sets $[{\alpha_0,\alpha})$ and $[{\beta_0,\beta})$, and the corresponding components are equal as sets,  i.e., for any ordinals $\alpha'\in[\alpha_0,\alpha)$ and $\beta'\in[\beta_0,\beta)$ such that $\beta'=\varphi(\alpha')$,  the following sets are equal:
first of all  $X_{\alpha_0}=Y_{\beta_0}$ as sets and then inductively for each $(\alpha_0<\alpha'<\alpha)$ and each $(\beta_0<\beta'<\beta)$ the following sets are equal, i.e.,
\begin{equation}
\label{f777}
X_{\alpha'+1}\cup\{ X|_{[\alpha_0,\alpha')}\}=Y_{\beta'+1}\cup\{ X|_{[\beta_0.\beta')}\},
\end{equation}

\bigskip

It is easy to see that if the coskand $X_{[\alpha_0,\alpha)}$ is of finite length, then it is a usual founded set of the theory $NBG[{\cal U}]$ and contributes nothing new. Therefore, we will only consider coskands of infinite length. We denote by ${\cal U}'$ the class of all coskands $X_{[\alpha_0,\alpha)}$ such that $l=\alpha-\alpha_0\geq\omega$, and $\alpha$ is a limit ordinal. We will interpret  elements of this class as individuals. To avoid the complication that such an individual could be contained a priori in a fixed class ${\cal U}$, we assume ${\cal U}\cap{\cal U'}=\emptyset$ and denote by ${\cal U}^{(0)}$ the class ${\cal U}$ and by  ${\cal U}^{(1)}$ the class ${\cal U}'$ defined above.

We can now extend the formal theory $NBG=NBG[{\cal U}^{(0)}]$ to the formal theory $NBG[{\cal U}^{(1)}]$ if we introduce the following additional axiom:

 {\bf Coskand existence axiom  and its interpretation as an individual or a set with  individuals} $({\bf CSEA}\&{\bf ICSU}$). 

For each ordinals $\alpha_0,\alpha\in{\bf On}$ such that $\alpha-\alpha_0\geq\omega$, and any mapping $f: [\alpha_0,\alpha)\rightarrow{\bf V}[{\cal U}^{(0)}]$ there exists a coskand $X$ with a clutch  region $[\alpha_0,\alpha)$ such that $X(\alpha')=f(\alpha')$ holds for any ordinal $\alpha'\in[\alpha_0,\alpha)$, and for each  $\alpha_0\leq\alpha'<\alpha$  when $\alpha'$ a limit ordinal a coskand $X_{[\alpha_0,\alpha')}$ is an individual in ${\cal U}^{(1)}$, and when $\alpha'$ is not a limit ordinal, then $X_{[\alpha_0,\alpha')}$ is a set from the universe ${\bf V}[{\cal U}^{(0)}\cup{\cal U}^{(1)}]$.

The axiom ${\bf CSEA}\&{\bf ICSU}$ is quite strong and saves us from proving the axioms of $NBG[{\cal U}^{(0)}\,\,\cup\,\,{\cal U}^{(1)}]$-theory, concerning in addition the class of individuals ${\cal U}^{(1)}$, some of which are obvious, e.g., that these individuals do not contain elements. This axiom is a consistent $NBG[{\cal U}^{(0)}\cup{\cal U}^{(1)}]$-theory if the $NBG$-theory is consistent and there is a model for it that builds on the formulas $(\ref{f1948})$ and $(\ref{f1949})$.

Since any coskand $X_{[\alpha_0,\alpha)}$ is equal to the coskand $X_{[0,\alpha_1)}$ for the ordinal $\alpha_1=\alpha-\alpha_0$ in particular for a limit ordinal $\alpha$ the ordinal $\alpha-\alpha_0$ is also limit, then, obviously, the class ${\cal U}^{(1)}$ of such coskands is a proper class. However, we cannot prove the axioms of plenitude ({\bf AP}) and strong axiom of plenitude ({\bf SAP}) for all coskands ${\cal U}^{(1)}$.

Nevertheless, if by transfinite induction we extend this class to the class of individuals ${\cal U}^{\Omega}=\bigcup\limits_{0<\nu<\Omega}{\cal U}^{(\nu)}$ of  the corresponding formal  $NBG[{\cal U}\cup{\cal U}^{(1)}\cup...\cup{\cal U}^{(\nu)}]$-theories for any ordinal $1<\nu<\Omega$, then the axioms of plenitude ({\bf AP}) and strong axiom of plenitude ({\bf SAP}) are fulfilled. This means that they are compatible with the other axioms of the theory $NBG[{\cal U}]$ and joining each of them to $NBG[{\cal U}]$ gives a consistent extension under the condition that the theory $NBG$ is consistent.

Suppose that for any ordinals $0<\nu'<\nu$ we construct a formal  $NBG[{\cal U}\cup{\cal U}^{(1)}\cup...\cup{\cal U}^{(\nu')}]$-theory, the class  ${\bf V}[{\cal U}\cup{\cal U}^{(1)}\cup...\cup{\cal U}^{(\nu')}]$ of all sets of this theory with the axioms of set theory $NBG[{\cal U}]$ plus the $\nu'$-coskand existence axiom and its interpellations as individuals or sets with  individuals, where the class of individuals ${\cal U}^{(\nu')}$ consists of all $\nu'$-coskands $X^{(\nu')}_{[\alpha_0,\alpha)}$ such that length $l=\alpha-\alpha_0\geq\omega$ is a limit ordinal and the components $X^{(\nu')}_{\alpha'}$, $\alpha'\in[\alpha_0,\alpha)$, belong to the universal class of sets 
$\bigcup\limits_{0<\nu''<\nu'}{\bf V}[{\cal U}\cup{\cal U}^{(1)}\cup...\cup{\cal U}^{(\nu'')}]$.  
Then we can define the $\nu$-skand $X^{(\nu)}_{[\alpha_0,\alpha)}$ as follows: its $\alpha'$-components $X^{(\nu)}_{\alpha'}$, $\alpha'\in[\alpha_0,\alpha)$ are elements of the class $\bigcup\limits_{0\leq\nu'\leq\nu}{\bf V}[{\cal U}\cup{\cal U}^{(1)}\cup...\cup{\cal U}^{(\nu')}]$, the membership relations of its constituents and parts, are analogous to the definition of such relations in a usual coskand (cf. above). An additional the coskand existence axiom  and its interpretation as an individual or a set with  individuals $({\bf CSEA}\&{\bf ICSU})^{(\nu)}$, where the class of individuals ${\cal U}^{(\nu)}$ is axiomatically given by the existence of $\nu$-caskads, and the class of all elements of ${\bf U}^{(\nu)}$ of a theory follows from its axioms, and the class of all sets of this theory is ${\bf V}[{\cal U}\cup{\cal U}^{(1)}\cup...\cup{\cal U}^{(\nu)}]={\bf U}[{\cal U}\cup{\cal U}^{(1)}\cup...\cup{\cal U}^{(\nu)}]\setminus({\cal U}\cup{\cal U}^{(1)}\cup...\cup{\cal U}^{(\nu')})$.

The limit  $NBG[{\cal U}]^{(\Omega)}$-theory is an extension of the $NBG$-theory with the addition of the axiom $({\bf CSEA}\&{\bf ICSU})^{(\Omega)}$, i.e., the coskand existence axiom  and its interpretation as an individual or a set with  individuals $({\bf CSEA}\&{\bf ICSU})^{(\nu)}$ for all ordinals $0<\nu<\Omega$  from the class ${\cal U}^{(\Omega)}=\bigcup\limits_{0<\nu<\Omega}{\cal U}^{(\nu)}$ and as sets with these individuals from the class ${\bf V}[{\cal U}\cup{\cal U}^{(\Omega)}]=\bigcup\limits_{\nu\in{\bf On}}{\bf V}[{\cal U}\cup{\cal U}^{(\nu)}]$, where the sums are taken over all limit ordinals $\nu\in{\bf On}$. The class of individuals ${\cal U}^{(\Omega)}$ is axiomatically given by the existence of $\nu$-caskads for any limit ordinal $0<\nu<\Omega$, and the class of all elements ${\bf U}^{(\Omega)}$ follows from the axioms of the theory $NBG[{\bf U}^{(\Omega)}]$, and the class of all sets of this theory is ${\bf V}[{\cal U}\cup{\cal U}^{(\Omega)}]={\bf U}^{(\Omega)}\setminus({\cal U}\cup{\cal U}^{(\Omega)})$.

{\bf Proposition 8}. {\it The axiom of plenitude ${\bf AP}$ is fulfilled for the class ${\cal U}^{(\Omega)}$ of the theory $NBG[{\cal U}^{\Omega}]$.}

{\bf Proof}. Let $S$ be an arbitrary set from the class ${\bf V}[{\cal U}\cup{\cal U}^{(\Omega)}]$.
 The {\it support of} $S$, {\bf support}$(S)$ is $TC(S)\cap({\cal U}\cup{\cal U}^{(\Omega)})$, where $TC(S)=\bigcup\{S,\bigcup S,\bigcup\bigcup S,...\}$ the so-called transitive closure of  $S$, i.e., the smallest transitive set including $S$. The elements of {\bf support}$(S)$ are individuals which are \grqq somehow involved\grqq in $X$. A set $X$ is called {\it transitive} if every set $Y$ which is an element of $X$ has the property that all of {\it its }  elements also belong to $X$. Definitions of the transivity, transitive closure and support of a set, see \cite{l98}, p. 14, 22.) 
 
 Actually, the support of {\bf support} $(S)$ of a set $S$ is the set of those individuals who participate in the formation of set $S$.  Note that {\bf support}$(S)$ is always a set, not a proper class, since the class $S$ formed by them is a set. We consider only subset of it, i.e., {\bf support}$(S)\cap{\cal U}^{(\Omega)}$.
 
Since any element from the set {\bf support}$(S)\cap{\cal U}^{(\Omega)}$ is a coskand $X_{[0,\lambda)}$ for some limiting ordinal $0<\lambda<\Omega$ and {\bf support}$(S)\cap{\cal U}^{(\Omega)}$ is not a proper class of the theory of $NBG[{\cal U}\cup{\cal U}^{(\Omega)}]$, then there are such limiting ordinals $\mu\in M$ such that the following inequalities are satisfied $\mu>\lambda$ for each $\lambda\in \Lambda$, where $$\Lambda=\{\lambda (\exists X_{[0,\lambda)}), X_{[0,\lambda)}\in\,{\bf support}(S)\cap{\cal U}^{(\Omega)}\}.$$ Otherwise, {\bf support}$(S)\cap{\cal U}^{(\Omega)}$ would be a proper class, not a set.
 
  Since ${\bf On}$ is a well-ordered class, any subclass of it has a minimal element. Let's denote this limit ordinal by $\nu=\min M$. Then it is obvious from Definition 10 that any coskand of length $\nu$ does not belong to the set {\bf support} $(S)\cap{\cal U}^{(\Omega)}$. We now assume for any element $s\in S$ a coskand $X^s_{[0,\nu)}$ for which $X^s_0=s$ if $s$ is a set, and $X^s_0=\{s\}$ if $s$ is an individual, and for any ordinal $\nu'\in(0,\nu)$ we assume $X^s_{\nu'}=\emptyset$.
 
 Thus we define the mapping $f:S\rightarrow{\cal U}^{(\Omega)}$ by assuming $f(s)=X^s_{[0,\nu)}$. The mapping $f$ is injective. Indeed, if $s\not=s'$, then $X^s_{[0,\nu)}\not=X^{s'}_{[0,\nu)}$, since their first components are different.

It is also clear that $f(S)\cap S=\emptyset$. Otherwise, $f(S)\cap${\bf support}$(S)\not=\emptyset$. But for any individual $X_{[0,\lambda)}$ of {\bf support}$(S)$, the equality $X_{[0,\lambda)}=X^s_{[0,\nu)}$ is not satisfied, since their clutch  regions are not isomorphic due to the fact that $\nu>\lambda$, for any $\lambda\in \Lambda$. As for any individual $u\in{\cal U}$ such that $u\in$ {\bf support}$(S)$ we see that $u\notin f(S)$ because each element of $f(S)$ is a coskand with a clutch  region $[0,\nu)$, where $\nu$ is a limit ordinal, and evidently $u\not=X^s_{[0,\nu)}$, for any $s\in S$. Therefore, the formula $f(S)\cap${\bf support}$S\not=\emptyset$ is incorrect. So the formula is fulfilled by $f(S)\cap S=\emptyset$. $\Box$

{\bf Proposition 9}. {\it The  strong axiom of plenitude ${\bf SAP}$ is satisfied for the class ${\cal U}^{(\Omega)}$ of the theory $NBG[{\cal U}\cap{\cal U}^{\Omega}]$.}

{\bf Proof}. Let $X$ be an arbitrary set of ${\bf V}[{\cal U}\cup{\cal U}^{(\Omega)}]$
 $Y$ be any set of individuals and $Y\subset {\cal U}^{(\Omega)}$.  Since any element from set $Y$ is a coskand of $X_{[0,\lambda)}$ for some limit ordinal $0<\lambda<\Omega$ and $Y$ is not a proper class of the theory $NBG[{\cal U}\cap{\cal U}^{(\Omega)}]$, then there exist such limit ordinals $\mu$ such that $\mu>\lambda$, where $\lambda\in \Lambda$ and $\Lambda=\{\lambda (\exists X_{[0,\lambda)}), X_{[0,\lambda)}\in Y\}$. Otherwise, $Y$ would be a proper class, not a set. Since ${\bf On}$ is a well-ordered class, any subclass of it has a minimal element. Let's denote this limit ordinal by $\nu=\min\Lambda$.
  Then it is obvious from Definition 10 that any coskand of length $\nu$ does not belong to the set $Y$. We now assume {\bf new}$(X,Y)=X_{[0,\nu)}$, where $X_0=X$ and $X_{\nu'}=\emptyset$ for any ordinal $0<\nu'<\nu$. As we have already observed, {\bf new}$(X,Y)\in{\cal U}^{(\Omega)}\setminus Y$, then the first condition of the strong axiom of plenitude ${\bf SAP}$ is satisfied. If $X'$ is any other set from the class ${\bf V}[{\cal U}\cup{\cal U}^{(\Omega)}]$ such that $X'\not=X$, then {\bf new}$(X,Y)=X_{[0,\nu)}\not=X'_{[0,\nu)}=${\bf new}$(X',Y)$ since $X_0\not=X'_0$. $\Box$

\bigskip
\begin{center}
{\bf 2. Reflexive sets and Mirimanoff extraordinary sets in $NBG[{\cal U}]^{(\Omega)}$.}
\end{center}

\bigskip

The axioms of set theory $NBG[{\cal U}]^-$ are compatible with (see \cite{l144}, p. 117) non-founded sets, which include the so-called reflexive sets and the Mirimanoff extraordinary sets \cite{l2}.

{\bf Definition 11}. A set $X\in{\bf V}[{\cal U}]^-$ is called {\it is reflexive} if one of its elements is $X$ itself, i.e., $X\in X$. Otherwise, the set $X$ is called {\it is non-reflexive}, i.e., it is not its  element, i.e., $X\notin X$.

{\bf Definition 12}. A set $X\in{\bf V}[{\cal U}]^-$ is called {\it is extraordinary} in the sense of Mirimanoff \cite{l2} if there exists at least one of its elements $x_0$ with infinite descending $\in$-sequence $X\ni x_0\ni x_1\ni x_2\ni...\,$. Otherwise, the set $X$ is called {\it  ordinal}, i.e., any maximal descending $\in$-sequence $X\ni x_0\ni x_1\ni x_2\ni x_n$ ends  either the empty set or an individual, i.e., $x_n=\emptyset\vee u$, $u\in{\cal U}$.

 Already the theory of $NBG[{\cal U}]^{(1)}$ is enough to explain the indeterminacy of the admissible relation $X\in X$ and the indeterminacy of the Mirimanoff definition of an extraordinary set, if by default, one does not assume that the Mirimanoff set is representable as a skand $X=X_{[\alpha_0,\alpha)}$ of length $l=\alpha-\alpha_0=\omega$.  

A reflexive set $X$ of the universe of sets ${\bf V}[{\cal U}]^-$ can be represented as a solution of the following \grqq implicit equation\grqq\, Mirimanoff \cite{l4}, p. 34:

\begin{equation}
\label{f1}
X=\{x_0,x_1,x_2,...,x_\lambda,...,X\},
\end{equation}
where $x_0,x_1,x_2,...,x_\lambda,...$ are elements from ${\bf U}[\cal U]^-$ such  that $\{x_0,x_1,x_2,...,x_\lambda,.... \}\in{\bf V}[{\cal U}]^-$, or their absence, in which $\{\}=\emptyset\in{\bf V}[{\cal U}]^-$ is trivial, and $X$ is an independent variable taking values in ${\bf V}[{\cal U}]^-$.

The following natural question arises: \grqq Is the solution of the equation $(\ref{f1})$ unambiguous, or is it ambiguous and we need to find all solutions of this equation? In other words, do the elements $x_0,x_1, x_2, ..., x_\lambda,...$ in the equation $(\ref{f1})$ uniquely determine the set $X$ or not? The answer to this question depends essentially on one or another theory of non-founded sets. For example, from the anti-founded axiom {\bf AFA} of Aczel's \cite{l93}, p. 10, follows from the uniqueness of the solution of the equation $(\ref{f1})$, although Aczel did not exclude that, in general, the equation $(\ref{f1})$ can have many different solutions in particular for the simplest equation he considered, $X=\{X\}$ (in our case the set represented by the empty skand ${\bf e}_{[\alpha_0,\alpha)}$ of length $l=\alpha-\alpha_0\geq\omega$), cf. : \cite{l93}, p. 6-7.

Let's now consider the simplest case of the definition of a reflexive set in formal set theory $NBG[{\cal U}]^{(1)}$, i.e., the case where in the equation (\ref{f1}) the elements $x_0,x_1,x_2,...,x_\lambda,...$ belong to the universe ${\bf U}[\cal U]$ so that $\{x_0,x_1,x_2,...,x_\lambda,...\}\in{\bf V}[{\cal U}]$, or they are absent, then trivially $\{\}=\emptyset\in{\bf V}[{\cal U}]$, and $X$ is an independent variable taking values in ${\bf V}[{\cal U}]^{(1)}$. In this case, we will briefly say that there is an equation $(\ref{f1})$ in theory  $NBG[{\cal U}]^{(1)}$whose solutions are reflexive sets $X$ of the universe ${\bf V}[{\cal U}]^{(1)}$. 

Similarly, we can talk about the equation $(\ref{f1})$ in the theory   $NBG[{\cal U}]^{(\nu)}$, $2\leq\nu\leq\Omega$ and its solutions if the elements $x_0,x_1,x_2,...,x_\lambda,...$ belong to the universals ${\bf U}[\cal U]^{(\nu')}$, $0\leq\nu'<\nu$, such that $\{x_0,x_1,x_2,...,x_\lambda,...\}\in\bigcup\limits_{0\leq\nu'<\nu}{\bf V}[{\cal U}]^{(\nu')}$, or they are absent, then trivially $\{\}=\emptyset\in\bigcup\limits_{0\leq\nu'<\nu}{\bf V}[{\cal U}]^{(\nu')}$, and  $X$ is the independent variable taking values in ${\bf V}[{\cal U}]^{(\nu)}$. Then the equation $(\ref{f1})$ in the theories $NBG[{\cal U}]^{(\nu)}$, $1\leq\nu\leq\Omega$ has a huge number of different solutions, which is a parametric family of sets whose parameters constitute a proper class, as shown by the following 
Proposition:

{\bf Proposition 10}. {\it All solutions of the equation $(\ref{f1})$ in the theory $NBG[{\cal U}]^{(1)}$ are a parametrized proper class of skands $X_{[\alpha_0,\alpha)}$, $\omega\leq\alpha<\Omega$, such that}

\begin{equation}
\label{f111}
X_{\alpha'}=\{x_0,x_1,x_2,...,x_\lambda,...\},\,0\leq\alpha'<\omega.
\end{equation}

{\bf Proof.}   For any skand $X_{[\alpha_0,\alpha)}$, $\omega\leq\alpha<\Omega$, for which the conditions $(\ref{f111})$ are satisfied, by the axiom ${\bf SEA\& RSS}$, there is a mapping $\sigma$, mapping to each skand  $X_{[\alpha',\alpha)}$, $\alpha_0\leq\alpha'<\omega$, the set $\sigma(X_{[\alpha',\alpha)})=\{x_0,x_1,x_2,...,x_\lambda,...,X_{[\alpha'+1,\alpha)}\}=X_{[\alpha',\alpha)}\in {\bf V}[{\cal U}]^{(1)}$. And by the definition of equality of skands and the condition $(\ref{f111})$ and equality of sets in the theory of $NBG[{\cal U}^{(1)}]$, we obtain the following equalities: $X_{[\alpha_1,\alpha)}=X_{[\alpha_0,\alpha)}$ and $X=\sigma(X_{[\alpha_0,\alpha)})=\{x_0,x_1,x_2,...,x_\lambda,...,\sigma(X_{[\alpha_1,\omega)})\}=\{x_0,x_1,x_2,...,x_\lambda,...,\sigma(X_{[\alpha_0,\omega)})\}=\{x_0,x_1,x_2,...,x_\lambda,...,X\}$. Thus the set $X=\sigma(X_{[\alpha_0,\alpha)})$ is a solution to the equation $(\ref{f1})$. $\Box$

Since any solution $X$ of equation $(\ref{f1})$ is a reflexive set, and ${\bf V}[{\cal U}]^{(0)}$ is the class of all founded sets of theory $NBG[{\cal U}]^{(1)}$, or, what is the same, the class of all sets of theory $NBG[{\cal U}]^{(0)}$, then $X\notin {\bf V}[{\cal U}]^{(0)}$, hence $X\in{\bf V}[{\cal U}]^{(1)}\setminus{\bf V}[{\cal U}]^{(0)}$. It is known that the class ${\bf V}[{\cal U}]^{(0)}$ is not a set, but is a proper class (we will give a new proof of this fact, which does not use the Mirimanoff or Russell antinomies, in the next paragraph), then the set of all solutions of the equation $(\ref{f1})$ in the theory $NBG[{\cal U}]^{(1)}$ is not a proper class, since the parameters of this family contain arbitrary elements from ${\bf V}[{\cal U}]^{(0)}$.

Recall here the definition of a parametric family of sets.

A mapping $f$ of a set $E$ {\it onto} a set $F$ is called a {\it  parametric representation of a set $F$ by means of a set $E$ whose elements are called parameters}. Moreover, a parametric representation of any subset of set $F$ is equivalent to specifying a mapping of some set $E$ {\it into} $F$. (See \cite{l88}, p. 364-365). Let's add only one clarification: If the mapping $f$ of the set $E$ \grqq into\grqq\, (respectively \grqq onto\grqq) $F$ is injective, then the subset $f(E)$ (respectively the set $F$ itself) is parametrized by the set of parameters $E$ and is representable as $\bigcup\limits_{{e\in E}}\{f(e)\}$; if the mapping $f$ is not injective, then the parameterized set $\bigsqcup\limits_{{e\in E}}\{f(e)\}$ is actually the set $\bigsqcup\limits_{{e\in E}}\{f(e)\}$. In our special case, we extend this definition of a parametrized family of sets to arbitrary classes $E$ and $F$, which can be (and in our case are) proper classes. Indeed, as $F$ we suppose a subclass of all reflexive sets $X\in{\bf V}[{\cal U}]^{(1)}$ with the same subsets $X\setminus\{X\}$, and as $E$ we suppose all reflexive skands $X_{[\alpha_0, \alpha)}$ for which $X_{\alpha'}=X\setminus\{X\}$, $\alpha_0\leq\alpha'<\omega$, and we map to each such skand $X_{[\alpha_0,\alpha)}\in E$ a reflexive set $\{x_0,x_1,. ..x_\lambda,...,X_{[\alpha_0+1,\alpha)}\}$.

{\bf Proposition 11}. {\it All solutions to the equation $(\ref{f1})$ in the theory $NBG[{\cal U}]^{(\nu)}$, $2\leq\nu<\Omega$, represent a parameterized proper class of sets, with parameters in the form of a reflexive skand 
$$X^{(\nu)}_{[0,\omega)}=\{x_0,x_1,x_2,...,x_\lambda,...\{x_0,x_1,x_2,...,x_\lambda,...\{...\,\,\,...\}\}\}$$ and skands $X^{(\nu)}_{[0,\alpha)}$,   $\alpha>\omega$, with $X^{(\nu)}_{\alpha'}=\{x_1,x_2,...,x_\lambda,...\}$, $0\leq\alpha'<\omega$, and with arbitrary $X^{(\nu)}_{\alpha'}\in\bigcup\limits_{0\leq\nu'<\nu}{\bf V}^{(\nu')}$, $\omega\leq\alpha'<\alpha$. 

Moreover, if in the equation $(\ref{f1})$ there is such a subsequence of elements $x_{\lambda_{\nu_\sigma}}\in{\bf V}[{\cal U}]^{(\nu_\sigma)}$ and $x_{\lambda_{\nu_\sigma}}\not\in{\bf V}[{\cal U}]^{(\nu')}$, $0\leq\nu'<\nu_\sigma$, and the set of all such ordinals $\{\nu_\sigma\}$ is confinal in the set $[0,\nu)$, then no solution of $X$ equation $(\ref{f1})$ is not an element of $\bigcup\limits_{0\leq\nu'<\nu}{\bf U}[\cal U]^{(\nu')}$, and belongs to proper class ${\bf V}[{\cal U}]^{(\nu)}\setminus\bigcup\limits_{0\leq\nu'<\nu}{\bf V}[{\cal U}]^{(\nu')}$.}

{\bf Proof} is similar to the proof of Proposition 10.

{\bf Proposition 12}. {\it All solutions of the equation $(\ref{f1})$ in the theory $NBG[{\cal U}]^{(\Omega)}$ are a parameterized proper class of sets, with parameters in the form of a reflexive skand 
$$X^{(\Omega)}_{[0,\omega)}=\{x_0,x_1,x_2,...,x_\lambda,...\{x_0,x_1,x_2,...,x_\lambda,...\{...\,\,\,...\}\}\}$$ and skands $X^{(\Omega)}_{[0,\alpha)}$,   $\alpha>\omega$, with $X^{(\nu)}_{\alpha'}=\{x_1,x_2,...,x_\lambda,...\}$, $0\leq\alpha'<\omega$, and with arbitrary $X^{(\Omega)}_{\alpha'}\in\bigcup\limits_{0\leq\nu<\Omega}{\bf V}^{(\nu)}$, $\omega\leq\alpha'<\alpha$. 

Moreover, in this limiting case, any solution $X$ of the equation $(\ref{f1})$ belongs to the universe ${\bf V}[{\cal U}]^\Omega$ and the set of all solutions constitutes a proper class of reflexive sets}.

{\bf Proof} is similar to the proof of Proposition 10.
 \bigskip

Consider another \grqq implicit equation\grqq\, of Mirimanoff \cite{l4}, p. 34:

\begin{equation}
\label{f121}
X=\{x^{0}_0,x^{0}_1,...,x^{0}_\lambda,...,\{x^{1}_0,x^{1}_1,...,x^{1}_\lambda,...,\{...\{x^{n-1}_0,x^{n-1}_1,...,x^{n-1}_\lambda,...,X\}\}\},
\end{equation}
where $x^{i}_0,x^{i}_1,...,x^{i}_\lambda,...$ are such elements from ${\bf U}[\cal U]$, $0\leq i\leq n-1$ such that $\{x^{i}_1,x^{i}_2,...,x^{i}_\lambda,...\}\in{\bf V}[{\cal U}]$, or the absence for some or even all $i$, $0\leq i\leq n-1$, in which it is trivial $\{\}_i=\emptyset\in{\bf V}[{\cal U}]^-$, and $X$ is an independent variable taking values in  ${\bf V}[{\cal U}]^{(1)}$. 

Sets of this kind $(\ref{f121})$ will be called periodic sets with period $n\geq 1$, or $n$-periodic in the theory $NBG[{\cal U}]^{(1)}$.  

It is clear that for $n=1$ such sets will be reflexive. Note, however, that periodic skands with period $1\leq n<\omega$ considered as sets will be periodic sets with the same period, but not vice versa.
 
 Similarly, we can talk about the equation $(\ref{f121})$ in the theory  $NBG[{\cal U}]^{(\nu)}$, $2\leq\nu\leq\Omega$ and its solutions if the elements  $x^{i}_0,x^{i}_1,...,x^{i}_\lambda,...$, $0\leq i\leq n-1$, belong to the universals ${\bf U}[\cal U]^{(\nu')}$, $0\leq\nu'<\nu$,  so that $\{x_0,x_1,x_2,...,x_\lambda,...\}\in\bigcup\limits_{0\leq\nu'<\nu}{\bf V}[{\cal U}]^{(\nu')}$ or there are none, then trivially $\{\}=\emptyset\in\bigcup\limits_{0\leq\nu'<\nu}{\bf V}[{\cal U}]^{(\nu')}$ and $X$  is the independent variable taking values in ${\bf V}[{\cal U}]^{(\nu)}$. Then the equation $(\ref{f121})$ in the theory $NBG[{\cal U}]^{(\nu)}$, $1\leq\nu\leq\Omega$, has a huge number of different solutions, which is a parametric family of sets whose parameters constitute a proper class, as shown by the following Propositions.

{\bf Proposition 13}. {\it Solutions to equation $(\ref{f121})$ in theory $NBG[{\cal U}]^{(\nu)}$, $1\leq\nu\leq\Omega$, represent a parameterized proper class of $n$-periodic sets, with parameters in the form of the $n$-periodic skand $X^{(\nu)}_{[0,\omega)}=$
 $$
\{x^{0}_0,x^{0}_1,...,x^{0}_\lambda,...,\{x^{1}_0,x^{1}_1,...,x^{1}_\lambda,...,\{...,\{x^{n-1}_0,x^{n-1}_1,...,x^{n-1}_\lambda,...,\{x^{0}_0,x^{0}_1,...,x^{0}_\lambda,...,\{...\}\}\}\}\}\}$$ and skands $X^{(\nu)}_{[0,\alpha)}$},
    
{\bf Proof.}  $n$-periodic skand $X^{(\nu)}_{[0,\omega)}=$
$$
\{x^{0}_0,x^{0}_1,...,x^{0}_\lambda,...,\{x^{1}_0,x^{1}_1,...,x^{1}_\lambda,...,\{...,\{x^{n-1}_0,x^{n-1}_1,...,x^{n-1}_\lambda,...,\{x^{0}_0,x^{0}_1,...,x^{0}_\lambda,...,\{...\,\,\,...\}\}\}\}\}\},
$$ with period $n$ and components $X^{(\nu)}_i=\{x^{i}_0,x^{i}_1,...,x^{i}_\lambda,...\}\in{\bf V}[{\cal U}]$, $0\leq i\leq n-1$, considered as the set $X=\{x^{0}_0,x^{0}_1,...,x^{0}_\lambda,...,X^{(\nu)}_{[1,\omega)}\}\in {\bf V}[{\cal U}]^{(\nu)}$, is indeed a solution to the equation $(\ref{f121})$, which follows from Definitions 1$'$ and 2$'$, and from the axiom  ${\bf SEA\& RSS}^{(\nu)}$. But, obviously, any skand $X^{(\nu)}_{[0,\alpha)}$ whose components of $X^{(\nu)}_{\alpha'}$ coincide with the components of the skand $X^{(\nu)}_{(0,\omega)}$, for all such $\alpha'$ that $0\leq\alpha'<\omega$, and the components of $X^{(\nu)}_{\alpha'}$ -- arbitrary founded sets for all such $\alpha'$ such that $\omega\leq\alpha'<\alpha'$, are solutions of the equation $(\ref{f121})$, where $\alpha$ is any ordinal such that $\alpha>\omega$.

Since the parameters of $n$-periodic solutions of the equation $(\ref{f121})$ in the theory of $NBG[{\cal U}]^{(\nu)}$ are collectively a proper class, then the set of all solutions of the equation $(\ref{f121})$ in the theory of $NBG[{\cal U}]^{(\nu)}$ is an proper class for any ordinal $\nu$, $1\leq\nu\leq\Omega$.

 \bigskip
Consider the following \grqq explicit equation\grqq\, of Mirimanoff \cite{l4}, p. 34:

\begin{equation}
\label{f122}
X=\{x^{0}_0,x^{0}_1,...,x^{0}_\lambda,...,\{x^{1}_0,x^{1}_1,...,x^{1}_\lambda,...,\{...\{x^{i}_0,x^{i}_1,...,x^{i}_\lambda,...,\{...\,\,\,...\}\}\}\}\},
\end{equation}
where $x^{i}_0,x^{i}_1,...,x^{i}_\lambda,...$, $0\leq i<\omega$, such elements from ${\bf U}[\cal U]^{(\nu)}$, $0\leq\nu\leq\Omega$ that $\{x^{i}_0,x^{i}_1,...,x^{i}_\lambda,... \in\bigcup\limits_{0\leq\nu'<\nu}{\bf V}[{\cal U}]^{(\nu')}$, or the absence for some and even all $i$, $0\leq i\leq \omega$, where $\{\{\}_i=\emptyset\in\bigcup\limits_{0\leq\nu'<\nu}{\bf V}[{\cal U}]^{(\nu')}$ is trivial, and $X$ is an independent variable taking values in ${\bf V}[{\cal U}]^{(\nu')}$. Sets of this kind $(\ref{f122})$ will be called Mirimanoff extraordinary sets.

{\bf Proposition 14}. {\it Solutions of the equation $(\ref{f122})$ in the theory $NBG[{\cal U}]^{(\nu)}$, $1\leq\nu\leq\Omega$, represent a parameterized proper class of Mirimanoff sets, with parameters in the form of a skand 
$X^{(\nu)}_{[0,\omega)}=$
 $$
\{x^{0}_0,x^{0}_1,...,x^{0}_\lambda,...,\{x^{1}_0,x^{1}_1,...,x^{1}_\lambda,...,\{...,\{x^{i}_0,x^{i}_1,...,x^{i}_\lambda,...,\{...\,\,\,...\}\}\}\}\}$$ and skands $X^{(\nu)}_{[0,\alpha)}$,
    $\alpha>\omega$, with the same components
 $X_{\alpha'}$, $1\leq\alpha'<\omega$ as in the skand $X_{[1,\omega)}$ and with arbitrary components $X^{(\nu)}_{\alpha'}\in\bigcup\limits_{0\leq\nu'<\nu}{\bf V}[{\cal U}]^{(\nu')}$, $\omega\leq\alpha'<\alpha$.}

{\bf Proof}.  The proof is similar to the proof of Proposition 10.

\begin{center}
{\bf 3. Naive and formal notions of reflexive and Mirimanoff extraordinary sets}
\end{center}

Reflexive sets in an implicit and unarticulated form first appeared in 1899 by Georg Cantor, the founder of set theory, when he distinguished a {\it  consistent multiplicity}, or such a definable multiplicity, which has the property that the assumption of {\grqq coexistence of all its elements {\it does not lead to a contradiction} and it can be regarded as a unity, as some \grqq completed or final thing\grqq, and he also distinguished  {\it a inconsistent plurality}, i.e., such plurality, which has the property that the assumption of \grqq coexistence\grqq\, {\it of all} its elements {\it leads to a contradiction} and it cannot be regarded as a unity as some \grqq completed or final thing\grqq\, \cite{l555} (Cantor's letter to Dedekind from Halle,  July 28, 1899). To paraphrase Cantor, a inconsistent multiplicity is an incomplete thing that cannot be fixed or defined in any way. It is a inconsistent multiplicity that formally contained itself as its element, but since it were contradictory and, consequently, non-existent, this reflexivity of theirs turned out to be merely fictitious.

Examples of such inconsistent sets are, according to Cantor, the set (collection) ${\bf U}$ {\it of all sets} $M$, the system $\Omega$ {\it of all ordinal numbers} $\alpha$ and the system $S$ {\it of all cardinal numbers} $\frak a$. It is easy to see that all these multiplicities, assuming their existence, contain themselves as their own element, which after simple reasoning leads to a contradiction, although their contradictory nature was proved differently by Cantor. This contradiction indicates only that such entities simply do not exist. 
Therefore, they cannot be considered as examples of reflexive sets, contrary to the standard example of the set of all sets, \grqq For other 
sets, one would as little hesitate to regard them as being members of  
themselves; the set of all sets is an obvious example\grqq, according to the authors of the famous monograph \cite{l14}, p. 5.

Indeed, the set of all sets contains itself as an element by the definition of ${\bf U}\in{\bf U}$, or because, for example, it contains as its subset the set of all its subsets, i.e., $2^{{\bf U}}\subset{\bf U}$, which leads to a contradiction by Cantor's famous theorem about a cardinality of all subsets of a set, and thus makes the multiplicity ${\bf U}$ inconsistent.

 If $\Omega$ were a consistent multiplicity, Cantor showed that it would be a well-ordered multiplicity to which corresponds some single ordinal number $\delta$, and then by virtue of the maximality of $\Omega$, the ordinal $\delta$ would be contained in $\Omega$ and would be greater than all ordinal numbers, and thus greater than itself, i.e., $\delta<\delta$,  \cite{l555} (Cantor's letter to Dedekind from Halle,  July 28, 1899), but following von Neumann the ordinal numbers are founded sets ordered by the $\in$ membership relation, we obtain the implicit presence in Cantor's proof of the relation $\delta\in\delta$, since $\delta<\delta$, the system $\Omega$ is a reflexive set, without the relative fact that it does not exist in this particular case. 

Similarly, Cantor showed that if $M_{\frak a}$ is any definite set of the class $\frak a$, $\frak a\in S$, and $S$ were a set, then the discrete union $T=\bigsqcup\limits_{\frak a\in S} M_{\frak a}$ of all sets $M_{\frak a}$, each of which corresponds in power to one and only one cardinal number $\frak a\in S$, and hence $T$ would belong to some class, say, ${\frak a}_0$. But by Cantor's famous theorem ${\frak a}_0'=2^{{\frak a}_0}>{\frak a}_0$, which means that the set $T$ would contain as a part a set, for example, $M_{{\frak a}_0'}=2^T$ of greater power ${\frak a}_0$, which is a contradiction according to Cantor. Hence, the multiplicity $T$ as well as the multiplicity $S$ are inconsistent. What is important for us here is that in Cantor's theorem $M_{{\frak a}_0'}=2^T\ni T$; and since $M_{{\bf a}_0'}\subset T$ and $T\in M_{{\frak a}_0'}$, then $T\in T$, which, however, is unarticulated by Cantor himself. But this unarticulated formula implies that $T$ is a reflexive set. True, the contradiction obtained by Cantor says only that such a system $T$ is not a consistent multiplicity and is not a completed or final thing or unity that is a mathematical object. See \cite{l555} (Cantor's letter to Dedekind from Hahnenklee,  August 31, 1899).

Later, reflexive sets appeared in the same fictitious way in the well-known set-theoretic paradoxes -- Zermelo, Russell, Mirimanoff, Hilbert and others. 

Let's first note the original definition of the paradoxical paradox of Russell's class ${\bf w}$, which caused in 1903 a painful reaction of almost all mathematicians and especially tragic for Gottlieb Frege. It's Russell's famous phrase:  \grqq We consider a contradiction arising from the obvious fact that if ${\bf w}$ is the class of all such classes which as individual members are not elements of themselves, then it can be proved that ${\bf w}$ as some object is and is not an element of itself\grqq. \cite{l47}, Chap. X, p. 107. Thus. 
the following paradoxical formula of Russell is already explicitly articulated here:
\begin{equation}
\label{f00001}
{\bf w}\in {\bf w}\Leftrightarrow {\bf w}\notin {\bf w},
\end{equation}
which says that the class ${\bf w}$ is a reflexive class with respect to the membership relation $\in$ if and only if it is not reflexive. Thus even here the reflexive multiplicity is fictitious, since Russell's class ${\bf w}$ is an inconsistent multiplicity and therefore does not exist, and this can be shown (see below) without reference to Russell's paradox.
Moreover, not knowing how to resolve his paradox, Russell writes later directly about any set: \grqq The relation $X\in X$ must always be meaningless\grqq, \cite{l55}, p. 81, and later clarifies: 
\grqq Thus $X\in X$ must be considered meaningless, since $\in$ requires that the inclusive be a class composed of objects of the same type as itself\grqq, (\cite{l47}, Chap. X, p. 107). 
In this case, elementary logical reasoning compels Russell to state that: \grqq $X\not\in X$ must always also be meaningless\grqq (\cite{l55}, p. 81). He wrote: \grqq If $\alpha$ is a class, then the statement \grq $\alpha$ is not an element of the class $\alpha$' is always meaningless, so there is no meaning in the phrase \grqq the class of those classes which are not members of themselves\grqq$\,\,$(\cite{l55}, p. 66.) Specifically: \grqq A class consisting of only one element cannot be identified with that element\grqq. And then instantly adds: \grqq the class $X=\{X\}$ must be absolutely meaningless, not simply false\grqq, \cite{l55}, p. 81. 

The reason for such \grqq fictitiousness\grqq\, of reflexive and more generally non-founded sets was and often still is the naive idea that it is possible, by some condition or premise like \grqq the existence of the totality of all sets\grqq, to draw a false conclusion that this totality is an element of itself. Even at the present time, not all mathematicians clearly understand how to correctly construct a reflexive or non-founded set, since in classical monographs on set theory, such sets are almost ignored. A good example is the formal construction of non-founded sets in such a serious monograph as \grqq Set theory and the continuum hypothesis\grqq\, by Paul J. Cohen with a characteristic comment: \grqq If non-founded sets exist, then they must be objects of a very artificial nature. Let's show how you can build a world containing such sets. Let $x_1,x_2,...$ be formal symbols; let's formally put $x_{n+1}\in x_n$ and $x_i\notin x_j$, at $i\not=j+1$. Let $R(0)=\{x_1,x_2,...\}$. We define $R(\alpha)$ by recursion as the set-degree for $\bigcup\limits_{\beta<\alpha}R(\beta)$.  The objects in $R(\alpha)$ for $\alpha>0$ are sets, whereas $R(0)$ consist only of $x_i$. Now, for the elements of sets $R(\alpha)$, we define some $\in$-relation as follows: if $u$ is a set, then $v\in u$, if $v$ is an element of $u$; if $u=x_i$, then $v\in u$, if $v=x_{i+1}$. The axioms are checked quite simply\grqq. (See \cite{l222}, Chap. II, $\S\, 5$.)

There is a minor inaccuracy in this construction of non-founded sets since the axiom of extensionality is not fulfilled, because the sets $\{x_{j+1}\}$ and $x_j$, $0<j<\omega$, are equal, but in Cohen's construction they are different elements. In fact, Cohen's construction can be imagined in a more natural way. Since, by the definition of a set, $x_j=\{x_{j+1}\}$ for each $0<j<\omega$,  then $x_1=\{x_2\}=\{\{x_3\}\}=...\,\,$. That's in the language of the skand the set $x_1$,  by default, is the empty skand ${\bf e}_{[1,\omega)}=\{_1\{_2\{...\}_2\}_1$. The only difference is that the equality of Cohen's symbols $x_j$, $0<j<\omega$, and the equality of the skands are different because the skand ${\bf e}_{[1,\omega)}$  is reflexive, while Cohen's set is not reflexive. If we consider a \grqq strong equality\grqq\, of skands ${\bf e}_{[j,\omega)}=\{_j\{_{j+1}\{...\}\}_{j+1}\}_j$, $0<j<\omega$ as equalities of ordered sequences $(j,j+1,j+2,...)$, then ${\bf e}_{[j,\omega)}=\{_j\{_{j+1}\{...\}\}_{j+1}\}_j={\bf e}_{[i,\omega)}=\{_i\{_{i+1}\{...\}\}_{i+1}\}_i$ if and only if $j=i$ and hence $x_j=x_i$ if and only if $j=i$ like in Cohen's construction. Nevertheless, since Cohen does not specify, by default, what the  symbols $x_i$, $0<i<\omega$ themselves mean, then in his construction they are not at all unambiguous, i.e., a variety of non-founded sets can correspond to the same characters, and the totality of such sets makes up a proper class like ${\bf e}_{[1,\alpha)}$, $\omega<\alpha<\Omega$.

{\bf Remark 7.} All of the above is done in the logic of the mathematicians of the time mentioned. In fact, many statements do not correspond to reality. For example, an indication that the totality of all sets, all ordinals, all cardinals, etc., i.e., inconsistent multiplicities in Cantor's terminology, are contradictory. This is how Cantor and other mathematicians of that time thought. In fact, this is not the case. These multiplicative formations are simply not definable, i.e., they are not mathematical objects: they do not exist, but {\it subsist} (bestehen, {\it Ger}). They {\it don't exist}, they {\it mean}. The difference in this case can be described as the difference between the non-existence of a number that is not equal to itself, i.e., there is no such number! And the totality of {\it all Conway numbers}, which does not exist in the sense of a complete process of their formation, but there is any number individually or there is a plurality of numbers (sets of numbers). This is the same difference between an object and a process, between a constant and a proportional form, between a number and a sequence that converges in some metric to this number, in particular, a constant sequence $(a_\alpha)_{0\leq\alpha<\lambda}$, such that for each ordinal $\alpha\in[0,\lambda)$ there is an equality $a_\alpha=a$,  and the number $a$ itself, etc. 

\begin{center}
{\bf 4. Known \grqq paradoxes\grqq\, cannot be derived in formal systems but they are redundant}
\end{center}

Due to the considerable uncertainty of the naive set-theoretic concepts of that time (absence of the axiom of extensionality, theories of non-founded sets, notions of compatibility and independence of set-theoretic axioms, etc.) and the absence of formal first-order logical theory (predicate calculus, models, and interpretations), Russell's correct conjectures and intuitions can only be marveled at.

Consider now Russell's paradoxical class in formal set theory $NBG^-$, in which $X\in X$ and $X\notin X$ are properly formed formulas (pff), and the class ${\bf R}\stackrel{def}{=}\{X|\,\,X\notin X\,\,\&\,\,X\in{\bf V}^-\}$ exists because the formula $X\notin X$ is predicative, and, by Proposition 4. 4 (Theorem on the existence of classes) in \cite{l07}, p. 238, we obtain that the class of sets ${\bf R}$ does exist; obviously, $ {\bf R}\subset {\bf V}^-$. But Russell's paradox $ {\bf R}\in {\bf R}\Leftrightarrow {\bf R}\notin {\bf R}$ in the theory of $NBG^-$ is not derivable \cite{l07}, p. 246, but proves only that ${\bf R}$ is not a set but a proper class, in contrast to Russell's general paradox $ {\bf w}\in {\bf w}\Leftrightarrow {\bf w}\notin {\bf w}$, which showed in the terminology of his time that ${\bf w}$ does not exist because it is inconsistent, and in modern terminology that ${\bf w}$ is not a universal set,  i.e., containing {\it all} sets $X$ such that $X\notin X$. It is important for us to note here only that even in axiomatic set theory Russell's paradox formally involves the formula ${\bf R}\in {\bf R}$ of the reflexive set (in the premises and conclusions), which turns out to be a proper class and, hence, the formula ${\bf R}\in {\bf R}$ is not legal and hence the relation ${\bf R}\in {\bf R}$  is simply a fiction since it is not provided for in $NBG^-$ theory, i.e., it is not a pff,  as well as the formula ${\bf R}\notin {\bf R}$.

In all the conclusions of these paradoxes, in addition to the problem of the existence of these universal objects, there is the same truth-functional operation of succession (implication), which is false, for example, in Russell's paradox: If ${\bf R}\stackrel{def}{=}\{X|\,\,X\notin X\}$ is a set, then from the premise ${\bf R}\notin {\bf R}$ follows, by definition of ${\bf R} $, the statement ${\bf R}\in {\bf R}$. But Proposition 8 shows us that no predicate, no definition (in this case the universality of the set ${\bf R}$) entails the statement ${\bf R}\in {\bf R}$, since the latter relation is, in general, indefinable and multivalued with respect to the reflexive set ${\bf R}$. Moreover, even if such unambiguity can be postulated, for example by Axcel's anti-founded axiom, the above implication continues to be false, since if we consider the set ${\bf M}$ to the set of all founded sets (Mirimanoff's paradox), then the premise ${\bf M}\notin {\bf M}$ is true, and the conclusion as the formula ${\bf M}\in {\bf M}\subseteq {\bf V}$ is false. This is what Proposition 8 shows: any solution $X$ of the equation $(\ref{f1})$ $X=\{X_0,X_1,...,X_\lambda,...,X\}$ where $\{X_0,X_1,...,X_\lambda,...,\}$ the set consisting of the founded sets is a non-founded set, so it cannot belong to the system of all founded sets.  In fact, a {\it reflexive set $X$ is never the result of putting together non-reflexive sets}.

Note also that reflexivity (formal, not real) in Russell's paradox, as in other set-theoretic paradoxes, arises also in the formal theory of the first order itself.  In a formal theory (a logistic system, or uninterpreted system, see, e.g., \cite{l333}, Introduction, $\S\, 07$) the atomic formula $X\in X$ is a \grqq properly formed formula\grqq\, (pff) (well-formed formula). One of the inference rules in formal logistic systems, namely {\it is the substitution rule}, leads to a \grqq well-formed reflexive formula\grqq. 

In order to formulate the rules for the derivation of the calculus
$P_1$ we will introduce the sign $S\, |$ to indicate the substitution operation,
so that $S^{\bf b}_{\bf B}{\bf A}\, |$ is the formula resulting from
the substitution of formula $B$ instead of each occurrence of the variable ${\bf b}$ in ${\bf A}$. The rule of inference is the following:

If {\bf b} is a variable, then ${\bf A}$ implies $S^{\bf b}_{\bf B}{\bf A} |$.
(Substitution rule.)  See \cite{l333}, Chap. 1, $\S\, 10$ and Chap. 3, $\S\, 30$. Note also that every freestanding constant and freestanding variable are properly formed formulas. See \cite{l333}, Chap. 1, $\S\, 10$. 

(The substitution rule follows from the second rule of inference of first-order theory -- {\it Generalization rule}: if ${\bf a}$ is an individual variable, then $\forall{\bf a}{\bf A}$ follows from ${\bf A}$, namely, $\forall{\bf a}{\bf A}\Rightarrow{\bf S}^{{\bf a}}_{{\bf b}}{\bf A}|$, where ${\bf b}$ is an individual variable or an individual constant and where there is no free occurrence ${\bf a}$ in ${\bf A}$ is not in such a correctly formed part of the formula ${\bf A}$, which has the form $\forall {\bf b}{\bf C}$, and the sign \grqq${\bf S}\,|$\grqq\, means a substitution operation so that ${\bf S}^{{\bf a}}_{{\bf b}}{\bf A}|$ is a formula that results from substituting the formula ${\bf b}$ instead of each occurrence of the variable ${\bf a}$ in {\bf A}.  Note also that every free-standing constant and free-standing variable are correctly formed formulas. See \cite{l333}, p. 66.) Indeed,  let ${\bf A}$ be a properly formed formula (pff) and ${\bf b}$ be an
individual variable or an individual constant, and no
free occurrence of {\bf a} in ${\bf A}$ is included in the pff-part of formula ${\bf A}$
of the form $({\bf b}){\bf C}$. If $\vdash{\bf A}$, then $\vdash S^{\bf a}_{\bf b}{\bf A}$. See \cite{l333}, Chap. $\S$ 34, Prop. *351. 
(The rule of substitution of individual variables.)

Let's verify now that the usual argument for Russell's paradox does not
hold in NBG. By the class existence theorem, there is a class $Y = \{x\, |\, x\notin x\}$.
Then $(\forall x)(x \in Y \Leftrightarrow x \notin x)$. In unabbreviated notation this becomes $(\forall X)(M(X)\Rightarrow(X \in Y \Leftrightarrow X \notin X))$. Assume $M(Y)$. Then $Y \in Y \Leftrightarrow Y \notin Y$, which, by the tautology $(A \Leftrightarrow \neg A) \Rightarrow (A\, \&\, \neg A)$, yields $Y \in Y\,\&\, Y \notin Y$. Hence, by the derived
rule of proof by contradiction, we obtain $\neg M(Y)$. Thus, in $NBG^-$, the argument for Russell's paradox merely shows that Russell's class $Y$ is a proper
class, not a set. $NBG^-$ will avoid the paradoxes of Cantor and Burali-Forti in
a similar way. (See  this literal proof in \cite{l07}, p. 246.)

However, due to the lack of interpretations, i.e., semantics, in the logical system, a correctly formed sub-formula ${\bf R}\in {\bf R}$ is  {\it  indefinable} or, similarly, is {\it  non-well-defined}, and the implication ${\bf R}\notin {\bf R}\Rightarrow {\bf R}\in {\bf R}$ is false, as opposed to the true implication ${\bf R}\in {\bf R}\Rightarrow {\bf R}\notin {\bf R}$. In spite of all this, the proof of \grqq multiplicity of the class ${\bf R}$ by means of \grqq of Russell's paradox is correct but redundant: only one true implication ${\bf R}\in {\bf R}\Rightarrow  {\bf R}\notin {\bf R}$ is sufficient to circumvent this paradox. This workaround is applicable to all set-theoretic paradoxes and will be shown below. Moreover, it gives other, stronger possibilities for proving other, different results (\grqq is not a set, but a proper class\grqq).

This is justified also by the fact that the formal attraction of Russell's paradox to the proof of \grqq multiplicity\grqq\ of some classes is {\it  superfluous}, and sometimes simply ridiculous. For example, this can be seen in the so-called relative (or parameterized) Russell paradox. If for any subclass ${\bf a}\subseteq{\bf V}^-$ we denote the class ${\bf a}\cap {\bf R}$ by ${\bf R}_{\bf a}$, then the correctly formed formula $(x\in {\bf R}_{\bf a}\Leftrightarrow x\in{\bf a}\,\,\& \,\, x\notin x)$ after substituting into it instead of $x$ the set ${\bf R}_{\bf a}$, assuming of course that ${\bf R}_{\bf a}\in{\bf a}$, turns into Russell's paradoxical relative formula: $({\bf R}_{\bf a}\in {\bf R}_{\bf a}\Leftrightarrow {\bf R}_{\bf a}\notin {\bf R}_{\bf a})$, which proves such an assumption wrong. 

(Note that the assumption $M({\bf R})$ leading to Russell's paradox is nothing but the assumption ${\bf R}\in{\bf V}^-$. Since, assuming ${\bf a}={\bf V}^-$, we get ${\bf R}={\bf R}_{\bf a}$. Indeed, ${\bf R}_{{\bf a}}={\bf R}_{{\bf V}^-}\stackrel{def}{=}{\bf V}^-\cap {\bf R}={\bf R}$ and the assumption ${\bf R}_{{\bf V}^-}\in{\bf a}$ is nothing but $R\in {\bf V}^-$).

In particular, in the theory  $ZF^-$, which considers only sets, the non-existence of the universal set is deduced by this method; Zermelo showed by the same reasoning that his set $M$, containing as its elements each of its subsets, does not exist; Mirimanoff showed that the set of all ordinary sets (i.e., founded sets) is inconsistent and thus does not exist. 

{\bf Remark 8.} For the first time, the idea that \grqq Russell's paradox is not  such a paradox\grqq\, was expressed by the authors of the book \cite{l8}, who considered Russellian propositions and accounts (other than, say, Austinian propositions and accounts) to the Liar Paradox, and proved that a subset of all non-reflexive sets ${\bf R}_X$ in a fixed set $X$ does not belong to $X$ using the relative Russell's paradox ${\bf R}_X\in X\Leftrightarrow{\bf R}_X\notin X$ via proof by contradiction supposing that ${\bf R}_X\in X$. This inconsistency is \grqq Russell's Paradox is not a paradox and proving it by Russell's paradox\grqq\, was noted in a paper \grqq The paradox of Russell's paradox\grqq\, \cite{l99} that first demonstrated the maximality principle, which circumvents Russell's paradox and similar ones, i.e., if we suppose that ${\bf R}_X\in X$, then it is easy to see that ${\bf R}_X\notin {\bf R}_X$ and hence ${\bf R}_X\subset{\bf R}_X\cup\{{\bf R}_X\}\subset X$, i.e., the set ${\bf R}_X$ is a proper subset of ${\bf R}_X\cup\{{\bf R}_X\}\subset X$, what contradicts to the maximality of ${\bf R}_X$.

The authors of the monograph \cite{l98}, p. 59, cite the following formal theorem of first-order logic: 
\begin{equation}
\label{f1100}
\neg\exists x\forall y[yEx\Leftrightarrow\neg(yEy)] 
\end{equation}
without any assumptions about the notion of set, other than that the displayed predicate determines a unique set, independent of what $E$ happens to mean. The authors of the mention monograph mean of course the substitution rule: if under the assumption that such a set $x$ exists, i.e., the formula $(\exists x\forall y[yEx\Leftrightarrow\neg(yEy)])$  is true, where $y$ is a variable and $x$ is a constant of above proportional form $(yEx\Leftrightarrow\neg(yEy))$, then we substitution $x$ instead of $y$ into it of course by the substitution rule  and obtain the paradoxical formula \begin{equation}
\label{f101}
xEx\Leftrightarrow\neg(xEx) 
\end{equation}
which implies, by formal methods, that there is no such $x$ or in $NBG^-$ theory $x$ is a proper class, i.e., $Pr(x)$. Indeed, if we suppose that  $x$ is a set, i.e., $M(x)$ and that it exists, then by the logical rule $B\Rightarrow A\& \neg A\Rightarrow \neg B$, we obtain  $\neg M(x)$, i.e., $x$ is a proper class. Due to a contradiction it is easily to deduce  that there is no such set, i.e.,  $x$ does not exist. In particular, if the Russel's class ${\bf R}\stackrel{def}{=}\{y\,|\,y\notin y\}$ is understood as $x$, and the membership relation $\in$ is understood as $E$, then the assumption $(\exists {\bf R}\forall y[yE{\bf R}\Leftrightarrow\neg(yEy)])$ implies $({\bf R}\in {\bf R}\Leftrightarrow\neg({\bf R}\in {\bf R}))$, and thus the Russell's paradox is obtained. 

(Note that if the same Russell class ${\bf R}$ is understood as $x$, and the non-membership relation $\notin$ is understood as $E$, then the assumption  $(\exists {\bf R}\forall y[y\notin {\bf R}\Leftrightarrow(y\in y)])$ implies formally the same Russell paradox  ${\bf R}\notin {\bf R}\Leftrightarrow\neg({\bf R}\notin {\bf R}))={\bf R}\in {\bf R}$ is obtained, although, as we will see below, there is a qualitative difference between them.)

Formal proofs in logistic systems were invented precisely in order to follow the condition of rigor, i.e., to ensure that the proof of a theorem does not use any interpretation, but proceeds only according to the rules of the logistic system (see \cite{l333}, $\S$ 07). \grqq The purpose of this formalism is, -- wrote Church, -- that a proof satisfying our condition of rigor should remain valid under any interpretation of the logistic system, so that we obtain in the end an economy by proving different things by a single reasoning. The size of the economy is determined by the fact that there is no need to repeat an indefinable number of times proofs which coincide in form but differ in content,  so 
how can they be held all at once, once and for all\grqq\, (see \cite{l333}, Introduction, $\S$ 07). Perhaps the most illustrative  of this formalism are the following examples.

Let's consider the cases when in $(\ref{f1100})$ $E=(\sim)$ or $E=(\neg\sim)$, in particular, $E=(=)$ or $E=(\not=)$, and  $E=(\subset)$ or $E=(\neg\subset)$, $E=(\supset)$ or $E=(\neg\supset)$ as well as $E=(\subseteq)$ or $E=(\neg\subseteq)$, $E=(\supseteq)$ or $E=(\neg\supseteq)$ and apply a formal method like the substitution rule to prove the formula $(\ref{f1100})$. We can do it but with great care because in mathematics, a {\it variable} is a symbol containing 
the name of which coincides with the content of the proper name, or
constant, except that the single denotation
of the constant is replaced here by the possibility of different
variable {\it values}. 
Since it is usually necessary to limit the values that
a variable can take we will assume that with a variable $y$ 
a {\it nonempty range} of its possible {\it values} is associated with it. Therefore, the content of a variable includes, in a sense, the content
of the proper name of its meaning area. And if we consider the cases when $E=(\sim)$, in particular, $E=(=)$,  $E=(\subset)$, $E=(\supset)$ as well as $E=(\subseteq)$, $E=(\supseteq)$ in $NBG$, then a region of the forms $\exists x\forall y[yEx\Leftrightarrow\neg(yEy)]$ with a variable $y$ is empty.  Thus Substitution rule does not work but formally  under the assumption $\exists x\forall y[yEx\Leftrightarrow\neg(yEy)]$ we can obtain $[xEx\Leftrightarrow\neg(xEx)]$. In particular, when $E=(=)$, then we obtain a fake paradox $[x=x\Leftrightarrow\neg(x=x)]$ proving that such $x$ does not exist what was clear from the beginning in $(\ref{f1100})$ when $E$ is $=$. Note also that in the cases when  $E=(\neg\sim)$, in particular, when  $E=(\not=)$,   $E=(\neg\subset)$,  $E=(\neg\supset)$ as well as  $E=(\neg\subseteq)$,  $E=(\neg\supseteq)$ we obtain a fake paradox $[\neg(x=x)\Leftrightarrow(x=x)]$, which proves in $NBG^-$ theory that $\neg M(x)$, i.e., $x$ is a proper class, clear from the beginning in $(\ref{f1100})$ when $E$ is $\not=$ because every set $y\in{\bf V}^-$ is  equal to $y$ and hence $\neg\exists x$, $M(x)$, such that $\forall y [x\not=y]$. It is true only when  $x={\bf V}^-$.

The following cases are even more unexpected. Let $E$ in (\ref{f1100}) be a predicate $<_{Card}$, i.e., an inequality of cardinality (power) of sets: $y<_{Card} x\stackrel{def}{\Leftrightarrow}Card(y)<Card(x)$, where $Card (x)$ is a cardinal number corresponding to the  cardinality (power) of a set $x$.

Then (\ref{f1100}) turns to the following formula
\begin{equation}
\label{f102}
\neg\exists x\forall y[y<_{Card}x\Leftrightarrow\neg(y<_{Card}y)] 
\end{equation}
what is the same as
\begin{equation}
\label{f103}
\neg\exists x\forall y[y<_{Card}x\Leftrightarrow(Card(y)\geq Card(y))].
\end{equation}
Suppose now  the opposite that  $\exists x\forall y[y<_{Card}x\Leftrightarrow(Card(y)\geq Card(y))]$. Then the propositional form $[y<_{Card}x\Leftrightarrow(Card(y)\geq Card(y))]$ has a region the class of all sets since the relation $(Card(y)\geq Card(y))$ fulfills for every set $y$, in particular, for a supposed set $x$, and thus, by the substitution rule, we can substitute $x$ instead of $y$ and obtain a paradoxical formula $x<_{Card}x\Leftrightarrow x\geq_{Card} x$ what of $NBG$ as above proves that $x$ is not a set but a proper class. Moreover, in $NBG$ it is a true conclusion since for any set $y$ and the universal proper class $x$ of all sets which exists by axioms in $NBG$, because always $Card(y)<Card(x)$.

But everything is not so simple in the case when we consider a predicate $>_{Card}$, where $y>_{Card} x\stackrel{def}{\Leftrightarrow}Card(y)>Card(x)$, and suppose that such $x$ exists. Then the propositional form $[y>_{Card}x\Leftrightarrow(Card(y)\leq Card(y))]$ has a region the class of all sets since the relation $Card(y)\leq Card(y)$  is executed for every set $y$. After this the similar paradoxical formula $x>_{Card}x\Leftrightarrow x\leq_{Card} x$ implies  that in $NBG$ $x$ is not only set but it is not also a proper class because for supposed $x$ the true formula is the following $\forall y (Card(y)>Card(x))$. Thus such $x$ does not exist in $NBG$ even for the theory with individuals for which such notion exists, i.e., an individual ${\bf Ur}(x)=\neg{\bf Cls}(x)$ is not a class  at all. That is because at present we do not know and cannot even imagine a mathematical object $x$ such as a set, family, class, etc. whose cardinality is strictly less than zero, i.e., $Card(x)<0$.

If we consider $E$ as $\leq_{Card}$ or $\geq_{Card}$, then the corresponding propositional forms $[y\leq_{Card}x\Leftrightarrow(Card(y)>Card(y))]$ and $[y\geq_{Card}x\Leftrightarrow(Card(y)<Card(y))]$ with a variable $y$ and a constant $x$ have an empty regions. Thus in these cases the substitution rule does not work at all.

We can speak about serious attitude and correct understanding of the reflexive and extraordinary only in their axiomatic postulation: \cite{l91}, \cite{l31}, \cite{l92}, \cite{l94}, \cite{l93}, \cite{l8}, \cite{l98}, and others. Note only that in all monographs devoted to set theory there is not a single example of a reflexive set, except for fictitious and naive examples, e.g., in \cite{l14}, p. 5, \grqq Other sets can just as naturally be considered as proper elements without hesitation: an obvious example is the set of all sets\grqq.  In the monograph \cite{l144}, Chap. II, p. 137-138, an example of a non-founded set is constructed, but the construction is accompanied by the following words: \grqq If non-founded sets exist, they must be objects of artificial nature\grqq. 

 And only in special monographs, e.g., in \cite{l93}, \cite{l8}, \cite{l98}, devoted directly to non-founded sets, there are these and other natural examples.

\bigskip
\bigskip
\bigskip

\begin{center}
{\bf 5. Bypassing the \grqq  paradoxical\grqq\, sets  in resolving set-theoretic antinomies:

 $(\exists {\bf x})(\forall{\bf y})[{\bf y}\in {\bf x}\Leftrightarrow {\bf y}\notin {\bf y}]\Rightarrow(\neg\exists {\bf x})(\forall {\bf y})[{\bf y}\in {\bf x}\Leftrightarrow {\bf y}\notin y]$ {\bf VS} ${\bf x}\in {\bf x}\Leftrightarrow {\bf x}\notin {\bf x}$}
\end{center}

In the previous section, we showed how to prove the  formula  $(\ref{f1100})$, i.e., $$\neg\exists x\forall y[yEx\Leftrightarrow\neg(yEy)]$$
 by the method of reduction to absurdity, i.e., assuming that it  is false and thus 
$$\exists x\forall y[yEx\Leftrightarrow\neg(yEy)]$$
is true, and establish the formula $xEx\Leftrightarrow xEx$  what is absurd. After that, we conclude that the formula $(\ref{f1100})$ is true.

As we noted above, the formal substitution rule, although acceptable, leads to \grqq half-empty\grqq\, formula. We will now show the proof, which also leads to absurdity and is more meaningful and preferable, since the poof is based on an  assumption of the formula $\exists x\forall y[yEx\Leftrightarrow\neg(yEy)]$, which leads to contradiction:  $(\exists x)(\forall y) [yEx\Leftrightarrow\neg(yEy)]\Rightarrow(\neg\exists x)(\forall y) [yEx\Leftrightarrow\neg(yEy)]$, which entails saying that  the maximal set $x=\{y\,|\, \forall y [yEx\Leftrightarrow\neg(yEy)]\}$ is not maximal because $x\notin x$, or what is the same, that $x\subset x\cup\{x\}$, i.e., $x$ is a proper subset of $ x\cup\{x\}$ which contradicts the maximality of $x$. Thus the $x$ does not exisis or $x$ is not a set, but a proper class, in some cases, e.g, in $NBG$.

Moreover, this method can be applied in the case 
\begin{equation}
\label{f4114}
\neg\exists x\forall y[yEx\Leftrightarrow(yEy)]
\end{equation} that is not covered by the Russell's method, as it leads not to a contradiction, but to a tautology 
\begin{equation}
\label{f1415}
xEx\Leftrightarrow xEx.
\end{equation}

Indeed, let's  consider  the case  $E=\in$ or $E=\notin$,   when $x=\{y\,\,|\,\, y\in{\bf V}[{\cal U}]^-\,\,\,\&\,\,\, y\in y\}$, i.e.,  {\it the collection of all reflexive sets} $x$ in the theory  $NBG[{\cal U}]^-$, and try to prove that $x$ is not a set but a proper class, by an analogous formal substitution method, and we would yield the tautology $x\in x\Leftrightarrow x\in x$ or $x\notin x\Leftrightarrow x\notin x$, which in this case does not solve the problem at hand, unlike Russell's class ${\bf R}$ above. This shows that the way of formal substitutions and \grqq paradoxical sets\grqq\, is not obligatory or necessary in solving such problems. In general, there are no universal methods for solving problems and the method, say, of the \grqq maximality principle\grqq\, is not universal either: there are problems to which it is not applicable and it is necessary to resort to completely new methods.

The redundancy of the formal proof that ${\bf R}_{\bf a}\notin{\bf a}$ by means of Russell's relative paradox appears in particular cases, e.g., if we consider a subclass ${\bf a}\subset{\bf V}^-$ such that {\bf a} is founded set $X$, i.e., ${\bf a}=X\in {\bf V}\subset{\bf V}^-$  (note that the formula $X\subset{\bf V}^-$ is also true since $X$ is a set), then, first, ${\bf R}_X=X$ and, second, assuming that ${\bf R}_X\in X$, or what is the same: $X\in X$ (which is itself incorrect by virtue of the definition of a founded set) Russell's relative paradox ${\bf R}_X\in {\bf R}_X\Leftrightarrow{\bf R}_X\notin {\bf R}_X$, or which is the same: $X\in X\Leftrightarrow X\notin X$, formally proves that ${\bf R}_X\notin X$, or which is the same: $X\notin X$. In other words, by assuming the obviously false formula $X\in X$ for the founded set $X$, we show by means of Russell's relative paradox that this formula is false, and that the true formula is the formula $X\notin X$. Moreover, we see that one of the implication ${\bf R}_X\in {\bf R}_X\Rightarrow {\bf R}_X\notin {\bf R}_X$ of Russell's relative paradox is true because its premise is false, and the second implication ${\bf R}_X\notin {\bf R}_X\Rightarrow {\bf R}_X\in {\bf R}_X$ of Russell's relative paradox is false because its premise is true and its conclusion is false (Modus ponens).

The excess is especially noteworthy in the simplest case, when $X$ is equal to the empty set $\emptyset$.
 In this case the above proposition is the following: ${\bf R}_\emptyset\notin\emptyset$. Indeed, truth of this proposition is evident by the definition of the empty set. Nevertheless, following the ideology of formal systems and formal methods to give one proof for all possible cases in one step, we repeat it with a help of Russell's paradox. Suppose the opposite ${\bf R}_\emptyset\in\emptyset$. Then by definitions ${\bf R}_\emptyset\in\emptyset\Leftrightarrow {\bf R}_\emptyset\notin\emptyset$, or more simple, if we suppose that $\emptyset\in\emptyset$, then, by Russell's paradox, obtain $\emptyset\in\emptyset\Leftrightarrow \emptyset\notin\emptyset$. Contradiction. Thus $\emptyset\notin\emptyset$. 

Another extreme case is no less significant, when we consider the theory of well-founded sets $NBG$ and want to prove formally with the help of Russell's paradox that the universal set ${\bf V}$ is not a set, but a proper class. The premise (suppose that ${\bf V}$ is a set) in this proof is that this set gives the only correct reasoning that ${\bf V}\notin {\bf V}$ and thus the implication ${\bf V}\in {\bf V}\Rightarrow {\bf V}\notin {\bf V}$ is true and ${\bf V}\notin {\bf V}\Rightarrow {\bf V}\in {\bf V}$ is false. So a formal Russell's paradox ${\bf V}\notin {\bf V}\Leftrightarrow {\bf V}\in {\bf V}$ which gives us that  the hypothesis that  ${\bf V}$ is a set is false and thus it is a proper class.
A more logical proof without fictitious incorrect assumptions is the following. If the class ${\bf V}$  of all founded sets is a founded set, then it follows from the obvious relation that ${\bf V}\notin{\bf V}$ and ${\bf V}\subset{\bf V}\cup\{{\bf V}\}$, i.e., it  is a proper subset of the set ${\bf V}\cup\{{\bf V}\}$, which contradicts universality of ${\bf V}$.

Note here that {\it all proofs with a help of Russell's paradox have the same intension}, e.g., when $NBG$ or $ZF$ includes Regular axiom and $X$ is always founded and thus there is no formulas $X\in X$ in $NBG$ or $ZF$ and thus ${\bf R}_X\in {\bf R}_X\Rightarrow {\bf R}_X\notin{\bf R}_X$ is true implication on the other hand ${\bf R}_X\notin {\bf R}_X\Rightarrow {\bf R}_X\in {\bf R}_X$ is always false because of Regular axiom as well as a nature of reflective sets and in initial Russell's paradox because of false premise of an existence of universal set (class or family).

\bigskip

The main purpose of this section is to exclude, although logically permissible, but from our point of view a meaningfully empty implication $x\notin x\Rightarrow x\in x$ in the proof of formula $(\ref{f1100})$, which leads either to the non-existence of the set $x$, or to the proof that $x$ is a proper class. We will conditionally called this method above the \grqq maximality principle\grqq.

{\bf 1.} {\it The Russell-type paradoxes}.

Here we will show how to bypass all Russell's type paradoxes what is more and simpler way to understand the state of affairs in such situations.

Let's consider the Russell class ${\bf R}\stackrel{def}{=}\{y\,|\,y\notin y\}$ in $NBG^-$ theory, i.e., in $(\ref{f1100})$ we put $x={\bf R}$, and the membership relation $\in$ is understood as $E$. Then the assumption $(\exists {\bf R}\forall y[yE{\bf R}\Leftrightarrow\neg(yEy)])$ implies $({\bf R}\in {\bf R}\Leftrightarrow\neg({\bf R}\in {\bf R}))={\bf R}\notin{\bf R}$, and thus the Russell's paradox is obtained. 
Now we miss a formal and, in our opinion, dubious or far from reality implication ${\bf R}\notin {\bf R}\Rightarrow{\bf R}\in {\bf R}$ and  obtain the same result. 

Indeed, by the definitions of reflexive and non-reflexive sets and under assumption  that ${\bf R}$ is a set, i.e., $M({\bf R})$, the relation ${\bf R}\in {\bf R}$ is false because all elements of ${\bf R}$ are non-reflexive and hence the formula ${\bf R}\notin{\bf R}$ is true. And this is no longer a hypothesis, but a real fact. Then obviously ${\bf R}$ is the proper subset of the set ${\bf R}\cup\{{\bf R}\}$, i.e., ${\bf R}\subset{\bf R}\cup\{{\bf R}\}$, that contradicts the maximality of the set ${\bf R}$. Thus the assumption  that ${\bf R}$ is a set, i.e., $M({\bf R})$, is false in $NBG^-$ theory and we conclude that ${\bf R}$ is a proper class, i.e., $Pr({\bf R})$.

The same arguments  are suitable for the  Russell-type class ${\bf Q}$ which is understood in $(\ref{f1100})$ as $x$, and the non-membership relation $\notin$ is understood as $E$, then the assumption  $\exists {\bf {\bf Q}}\forall y[y\notin {\bf Q}\Leftrightarrow(y\in y)]$ implies formally the  Russell-type paradox  ${\bf Q}\notin {\bf Q}\Leftrightarrow\neg({\bf Q}\notin {\bf Q})={\bf Q}\in {\bf Q}$.

Now we miss a formal and, in our opinion, dubious or far from reality implication ${\bf Q}\notin {\bf Q}\Rightarrow{\bf Q}\in {\bf Q}$ and obtain  the same result. Indeed, by the definition of ${\bf Q}$, one can see that ${\bf}$ consists of {\it all} non-reflexive sets, and by definitions of reflexive and non-reflexive sets and under assumption that ${\bf Q}$ is a set, i.e., $M({\bf Q})$, the relation ${\bf Q}\in {\bf Q}$ is false because all elements of ${\bf Q}$ are non-reflexive  hence the formula ${\bf Q}\notin{\bf Q}$ is true. And this is no longer a hypothesis, but a real fact. Then obviously ${\bf Q}$ is the proper subset of the set ${\bf Q}\cup\{{\bf Q}\}$, i.e., ${\bf Q}\subset{\bf Q}\cup\{{\bf Q}\}$, that contradicts the maximality of the set ${\bf Q}$. Thus the assumption  that ${\bf Q}$ is a set, i.e., $M({\bf Q})$, is false in $NBG^-$ theory and we conclude that ${\bf Q}$ is a proper class, i.e., $Pr({\bf Q})$. 

(Notice that the proofs in these two cases are correct in $NBG$-theory but here we have bypassed a substitution in $NBG$-theory with Regular Axiom ${\bf RA}$ when in formal methods the {\it empty region} of the propositional form $y\notin {\bf Q}\Leftrightarrow (y\in y)$ of the variable $y$.)

The same arguments are valid for the relative Russell's paradox in the following way. 
Let's consider the Russell class ${\bf R}_X\stackrel{def}{=}\{y\,|\,y\notin y\,\&\,y\in X\}$ in $NBG^-$ theory, i.e., in $(\ref{f1100})$ we put $x={\bf R_X}$, and the membership relation $\in$ is understood as $E$. Then the assumption $(\exists {\bf R}_X\forall y[yE{\bf R}_X\Leftrightarrow\neg(yEy)])$ is nothing more than $({\bf R}_X\in {\bf R}_X\Leftrightarrow\neg({\bf R}_X\in {\bf R}_X))={\bf R}_X\notin{\bf R}_X$, and thus the relative Russell's paradox is obtained.

Suppose now that ${\bf R}_X\in X$. We see that ${\bf R}_X\notin {\bf R}_X$ because all elements of ${\bf R}_X$ are non-reflexive sets and hence the formula ${\bf R}_X\notin{\bf R}_X$ is true. Thus ${\bf R}_X\cup\{{\bf R}_X\}\subseteq X$ and ${\bf R}_X$ is a proper subset of ${\bf R}_X\cup\{{\bf R}_X\}$, i.e.,  ${\bf R}_X\subset {\bf R}_X\cup\{{\bf R}_X\}$ what contradicts the maximality of ${\bf R}_X$. Thus ${\bf R}_X\notin X$.

The same arguments are valid for the relative Russell-type paradox in the following similar example. 
Let's consider the Russell class ${\bf Q}_X\stackrel{def}{=}\{y\,|\,y\notin y\,\&\,y\in X\}$ in $NBG^-$ theory, i.e., in $(\ref{f1100})$ we put $x={\bf Q}_X$, and the membership relation $\in$ is understood as $E$. Then the assumption $(\exists {\bf Q}_X\forall y[yE{\bf Q}_X\Leftrightarrow\neg(yEy)])$ is nothing but $({\bf Q}_X\in {\bf Q}_X\Leftrightarrow\neg({\bf Q}_X\in {\bf Q}_X)={\bf Q}_X\notin{\bf Q}_X$, and thus the relative Russell paradox is obtained.

We  bypass this paradox in the following way. Suppose  that ${\bf Q}_X\in X$. We see that ${\bf Q}_X\notin {\bf Q}_X$ because all elements of ${\bf Q}_X$ are non-reflexive sets. Thus ${\bf Q}_X\cup\{{\bf Q}_X\}\subseteq X$ and ${\bf Q}_X$ is a proper subset of ${\bf Q}_X\cup\{{\bf Q}_X\}$, i.e.,  ${\bf Q}_X\subset {\bf Q}_X\cup\{{\bf Q}_X\}$ what contradicts the maximality of ${\bf Q}_X$. Thus ${\bf Q}_X\notin X$.

(Notice that $\{{\bf R}_X\}$ always exists either it is a set or a proper class. Only fake formal Russell's paradox deducts that proper class is not an element. But it is very artificial to identify \grqq property\grqq\, with a \grqq proper class\grqq. The model of $NBG$, the existence of an inaccessible cardinal number contradicts it since there is a one-element set with this inaccessible cardinal number.)

Moreover, this circumvention of the Russell-type paradox makes it possible to prove propositions that are impossible by such the above formal methods. Indeed, denote by ${\bf S}=\{y\,|\, y\in y\}$, i.e., a collection of {\it all reflexive sets}. Then in $NBG$ the following Proposition is true: 
\begin{equation}
\label{f1416}
\neg\exists x\forall y[y\in x\Leftrightarrow(y\in y)]
\end{equation}
 and hence $x$ is not a set but a proper class. The proof by the above formal method is useless because, by opposite and subsistence rule, we obtain a tautology $x\in x\Leftrightarrow x\in x$.

But by an assumption that such set ${\bf S}$ exists we consider a reflexive set ${\bf S}'$ such that ${\bf S}'\setminus\{{\bf S}'\}={\bf S}$, say, $S_{[0,\omega)}$ such that $S_\alpha={\bf S}$ for each $0\leq\alpha<\omega$, which exists by Definition 2$'$ above. Now we prove that ${\bf S}'=S_{[0,\omega)}\notin {\bf S}$. Really, for each element $y\in {\bf S}$ the set $y\setminus\{y\}$ does not have $y$ as its element. But $y$ is an element of ${\bf S}'$. Consequently, ${\bf S}\subset {\bf S}\cup\{{\bf S}'\}$, i.e., ${\bf S}$ is a proper subset of ${\bf S}\cup\{{\bf S}'\}$ what contradicts the maximality of ${\bf S}$. Thus ${\bf S}$ is a proper class and there is no a singleton $\{{\bf S}\}$. Note of course that we mentioned here a theory of reflexive sets (see, e.g., \cite{l1333}) for which an equality $y={\bf S}'$ of two reflexive sets $y$ and ${\bf S}'$ implies the identity of $y\setminus\{y\}={\bf S}'\setminus\{{\bf S}'\}$. If it is not so (see, e.g., \cite{l93}), then one must give another proof that ${\bf S}$ is a proper class.

(Notice that if we consider instead formula $(\ref{f1416})$ the following formula:
\begin{equation}
\label{f1417}
\neg\exists x\forall y[y\notin x\Leftrightarrow(y\notin y)],
\end{equation}
then the solutuion of this statement depends on the formal set theory, i.e., if we consider $NBG^-$, then $x$ is not a set but a proper calss of all reflexive sets; but in $NBG$, where ${\bf V}$ is a proper class of all non-founded sets, then this statement is fauls because there exists the empty set $\emptyset=x$, which satisfies the negation of formula $(\ref{f1417})$: $\exists x\forall y[y\notin x\Leftrightarrow(y\notin y)]$, i.e., $\forall y[y\notin \emptyset\Leftrightarrow(y\notin y)]$.)

So choosing the axiom of the existence of a singleton set $\{X\}$ and getting the result of the non-existence of a universal class ${\bf V}$ for sets is preferable from the point of view of the state of affairs or subsistence than postulating the axiom $\forall X Pr(X)\neg\{X\}$ of denying the existence of a singleton set $\{X\}$ for proper classes and getting a universal class ${\bf V}$ with questionable \grqq geography\grqq\, (say clumsily in the case of all Conway numbers, \grqq numbers geography\grqq\, or more precisely \grqq geometry of all Conway numbers\grqq).

\bigskip

{\bf 2.} {\it Zermelo-type paradoxes}

Zermelo posed the following question: \grqq Does there exist a set $X$ such that each subset of it is an element of the set $X$?\grqq  His own answer was \grqq No\grqq\, \cite{l922}. Moreover, he proved it by using a paradoxical set like in above examples. He supposed that such set $X$ existed and consider the following subspace ${\bf Z}\stackrel{def}{=}\{x\,\,|\,\, x\subseteq X\,\,\&\,\,x\notin x\}$ of the set $x$, i.e., ${\bf R}_X=X\cap{\bf R}$ and showed that ${\bf Z}\in{\bf Z}\Leftrightarrow{\bf Z}\notin{\bf Z}$, i.e., the existence of paradoxical set. Hence the above set $X$ does not exist.

Here we also exclude, although logically permissible, but from our point of view a meaningfully empty implication ${\bf Z}\notin {\bf Z}\Rightarrow {\bf Z}\in {\bf Z}$ in the Zermelo's proof, which leads either to the non-existence of the set $X$, or to the proof that $X$ is a proper class in $NBG^-$-theory and prove the same result by the maximality principal.

We see that, by definition, ${\bf Z}\notin {\bf Z}$ and by supposition that the set $X$ exists and clearly ${\bf Z}\subseteq X\,\,\&\,\,{\bf Z}\in X$, we conclude that ${\bf Z}\subset{\bf Z}\cup\{{\bf Z}\}$, i.e., ${\bf Z}$ is a proper subset of ${\bf Z}\cup\{{\bf Z}\}$ what in contradiction with maximality of ${\bf Z}$. Thus, in $ZF$-theory the set $X$ does not exist and in $NBG^-$-theory
$X={\bf V}^-$ is a proper class, i.e., $Pr({\bf X})$.

{\bf 3.} {\it Mirimanoff-type paradoxes}.

This is also the basis of Mirimanoff's proof (the paradox of his name) that the set of all founded sets (note -- outside any formal theory) does not exist: if it exists, it is a founded set, since otherwise it would have an non-founded element (true statement); if it is a founded set, it is an element of itself (false statement, since the hypothesis of the existence of a universal class of founded sets has not yet been proved, but is put forward as a hypothesis, and an incorrect one at that), hence it is an non-founded set. See \cite{l2}, p. 43. In the formal theory of $NBG$ or $NBG^-$, where such a universal class of founded sets exists as a consequence of the existence axioms for classes, Mirimanoff's paradox proves only that ${\bf V}$ is not a set, but is a proper class. Although even here the logistic formalism is redundant.

\bigskip

Mirimanoff's paradox is solved in the same way: the family {\it of all} of ordinal (in modern terminology -- founded) sets {\it does not exist} \cite{l2}, p. 43. Its proof reduces to the proof of the non-existence of the family {\it all} of sets of the first kind (in modern terminology -- non-reflexive sets) \cite{l2}, p. 40, since the latter set is a subset of the former. A curious proof of Mirimanoff's paradox is given later in \cite{l131} and \cite{l2}, p. 15: let ${\bf M}$ be a family {\it of all} founded sets. If the family ${\bf M}$ exists, then it is founded (true statement); indeed, otherwise there would be an infinite sequence of
 $\in$-relations ${\bf M}\ni X_1\ni X_2\ni X_3\ni...$and hence there would be an infinite sequence of $\in$-relations $ X_1\ni X_2\ni X_3\ni.... $, which would mean that the set $X_1\in{\bf M}$ is non-founded, which contradicts the fact that all elements of the family ${\bf M}$ are founded sets; now, if it is proved that ${\bf M}$ is a founded set, and it is itself a family of all founded sets, then it contains itself, i.e. ${\bf M}\in{\bf M}$, hence it is non-founded. By Mirimanoff's non-existence, for subsequent authors the assumption that ${\bf M}$ is a set is incorrect, so this set ${\bf M}$ is a proper class. There is an incorrect implication in this proof: from the universality of the family ${\bf M}$ and the correct proof that ${\bf M}$ is a founded family or presumably a set, it does not follow that ${\bf M}\in{\bf M}$, since reflexive sets (classes, hyper-classes, etc., if they are admissible in certain axiom systems) are not obtained by simply attaching a set to itself, as can be seen from the experience of non-founded set theory. Like the case of Russell's paradox, Mirimanoff's paradox is solved as follows.

Let's denote by ${\bf M}$ the class {\it of all} founded sets in $NBG^-$, which we will call the class {\it of all ordinal Mirimanoff sets}, and by ${\bf N}$ the class {\it of all} non-founded sets in $NBG^-$, which we will call the class {\it of all extraordinary Mirimanoff sets}. Both classes exist since the class ${\bf M}$ coincides with the class ${\bf V}$ of the theory $NBG$, which in exactly this theory $NBG$ is defined by the predicative formula $X\notin X$ and obviously ${\bf M}\subset{\bf V}^-$. The class ${\bf N}$ also exists by the aforementioned Axiom of Additivity in $NBG^-$ theory.

 Let ${\bf a}\subseteq{\bf V}^-$ be a fixed subclass (set or a proper class, indifferent) of the class ${\bf V}^-$. We denote by ${\bf M}_{\bf a}$ the class {\it of all} sets $X$ in ${\bf a}$ which are founded sets. Let's call ${\bf M}_{\bf a}$ {\it is a relative class of ordinal Mirimanoff sets}, since the ordinariness property for set $X$ does not apply to all sets, but only to elements of the class ${\bf a}$.

It is obvious that
 ${\bf M}_{\bf a}={\bf a}\cap{\bf M}$ and ${\bf M}_{\bf a}={\bf M}$ if ${\bf a}={\bf V}^-$ or ${\bf a}={\bf M}$.

 Let ${\bf a}\subset{\bf V}$ be a fixed subclass of the class ${\bf V}^-$. Let's denote by ${\bf N}_{\bf a}$ the class {\it of all} sets $X$ in ${\bf a}$ which are non-founded sets, i.e., are extraordinarily Mirimanoff sets.
We call ${\bf N}_{\bf a}$ {\it is a relative class of Mirimanoff extraordinary sets}, since the extraordinarity property for set $X$ does not apply to all sets, but only to elements of the class ${\bf a}$.

And this the most general case can be easily proved in the following way. It is clear that ${\bf M}_{\bf a}\notin{\bf M}_{\bf a}$ because all elements of ${\bf M}_{\bf a}$ are founded sets. Then ${\bf M}_{\bf a}$ can not be an element of ${\bf a}$. Indeed, if we suppose that ${\bf M}_{\bf a}\in X$, then ${\bf M}_{\bf a}\subset{\bf M}_{\bf a}\cup\{{\bf M}_{\bf a}\}$, i.e., ${\bf M}_{\bf a}$ is a proper set of ${\bf M}_{\bf a}\cup\{{\bf M}_{\bf a}\}$ what contradicts the maximality of ${\bf M}_{\bf a}$.

It is obvious tha ${\bf N}_{\bf a}={\bf a}\cap{\bf N}$ and ${\bf N}_{\bf a}={\bf N}$ if ${\bf a}={\bf V}^-$ or ${\bf a}={\bf N}$.

{\bf 4.} {\it Cantor-type paradoxes}

Note here that Cantor's classical \grqq diagonal method\grqq\, (\cite{l112}, p. 75-78) later acquired a formal form  similar in form to the  Russell's paradox. Indeed, if $X$ is an arbitrary set, and ${\cal P}(X)$ is the set of all its subsets, then Cantor's theorem states that there is no bijection between the sets $X$ and ${\cal P}(X)$. Indeed, if such a bijection $f:X\rightarrow {\cal P}(X)$ existed, then the subset $K\stackrel{def}{=}\{x\in X|\,\,x\notin f(x)\}$ would be contradictory: the element $x_0=f^{-1}(K)$ would satisfy the following formula: $x_0\in K\Leftrightarrow x_0\notin K$ (indeed, the formula is similar to Russell's paradoxical formula: $R\in R\Leftrightarrow R\notin R$). Therefore, there is no such element $x_0$ in $X$, but then there is no correspondence to the subset $K$ of the set $X$ in $X$, which means there is no bijection of $f$.

The same proof is available in a more general formulation of the theorem (see, for example, \cite{l5}, Chap. V, \S 3.), from which Cantor's theorem follows.: 

{\it If the domain of definition $T$ of an arbitrary function $f$ is contained in $X$, the values of which are subsets of the set $X$, then the set $K\stackrel{def}{=}\{x\in T|\,\, x\notin f(x)\}$ is not the value of the function $f$, or in other words, $K\notin f(T)$}.

Indeed, if for some $x_0\in X$ we have the formula $f(x_0)=K$, or in another way, $K\in f(T)$, then we get the same contradiction: $x_0\in K\Leftrightarrow x_0\notin K$, which proves the statement of the theorem. As in the case of the non-elimination of Russell's paradox in the formal theory of $NBG^-$, the evidence in various specific cases is redundant and ridiculous, obscuring the essence of the matter. For example, if $T=X$, and $f$ is an injective mapping that maps to each element of $x\in X$ a one-element subset of $\{x\}\subset X$, or, equivalently, $\{x\}\in{\cal P}(X)$, then it is obvious that in this particular case $K=\emptyset$ and it is obvious that $\emptyset\notin\{\{x\}|\,\,x\in X\}=f(X)$. Although from the general formal proof we get a contradictory formula: $x_0\in\emptyset\Leftrightarrow x_o\notin\emptyset$, which is really contradictory, but trivially ridiculous, since the formula $x_0\in\emptyset$ is incorrect, and the formula $x_0\notin\emptyset$ is correct.

The second extreme. If $T=X$, and $f$ is a constant map whose value is an empty subset of the set $X$, then it is obvious that in this case $K=X$ and that the formula $f(x_0)=X$ does not work out for any $x_0\in X$, or, which is the same as $X\notin f(X)$. Although the same is obtained from a formal proof through a contradictory formula: $x_0\in X\Leftrightarrow x_o\notin X$, which is really contradictory, but trivially ridiculous, since the formula $x_0\in X$ is correct, and the formula $x_0\notin X$ is incorrect.

 Meanwhile, the above Cantor's Theorem can be proved by bypassing the above-mentioned contradictory formulas. Indeed, the set is $K\in f(T)$. Then there is an element of $t\in T$ such that $f(t)=K$. By defining the set $K$, we obtain the correct relation $K\notin K$. Then, according to the Pair Axiom, there exists a one-element set $\{K\}$, according to the Sum Axiom, there exists $K'=K\sqcup\{K\}$, which contains the set $K$ as its proper subset, which contradicts the maximality of the set $K$. Therefore, $K\notin f(T)$. $\Box$

{\bf 5.} {\it Hilbert-type paradoxes}

Hilbert's paradox  is solved in a similar way, which he
never published, because, according to him, this paradox was known to all
experts in set theory, especially Cantor (\cite{l212}, p. 204).
Let's present it according to the source \cite{l292}, p. 505-506. The paradox is based on the special
concept of a set, given by Hilbert on the basis of two principles.

The first principle of additivity, which says that any two
sets can be combined into one set as well as combine infinitely
many sets into one set that contains elements of all these
sets. The second principle is the mapping principle, which says,
 that for any set ${\cal M}$ there exists a set ${\cal M}^{\cal M}$ of all mappings
of the same into itself. Then he constructs the universe ${\cal U}$ of all sets obtained
by applying these two principles indefinitely to
the set ${\bf N}$ of all natural numbers. Finally, he gets the last
set ${\cal F}={\cal U}^{\cal U}$ and concludes \grqq that since ${\cal F}$ is formed by only two
principles, it is obtained only from natural numbers\grqq, he concludes,
 what is ${\cal F}\in{\cal U}$. (By the way, the first principle instantly implies inclusion
 ${\cal F}\subseteq{\cal U}$, since the union of two sets ${\cal F}\cup{\cal U}$ contains elements
of each of them.) Considering \grqq that there are no more elements in the set ${\cal F}$ than in the set ${\cal U}$\grqq,
 Hilbert supposed that \grqq there is a bijection between ${\cal U}$ and ${\cal F}$
that corresponds to each $u_i\in{\cal U}$ a map $f_i\in {\cal F}$ and, by Cantor's diagonal method,
  constructed such a mapping $g:{\cal U}\rightarrow{\cal U}$ such that $g\notin{\cal F}$, what is in contradiction with the fact that ${\cal F}$ is the set of all maps from ${\cal U}$ to ${\cal U}$\grqq. 

(The assumption
of a bijection between ${\cal F}$ and ${\cal U}$ follows from the Cantor-Berstein theorem, since
there is already one bijection ${\cal F}\rightarrow{\cal U}$, by consequence of ${\cal F}\subseteq{\cal F}$ from the previous
predicate ${\cal F}\in{\cal U}$ and the second ${\cal U}\rightarrow{\cal F}$, which  can easily be solved by comparing each
element $u\in{\cal U}$ to a constant map ${\cal U}\rightarrow{\cal U}$, the image of which coincides
with $u$.) 

Hilbert concluded from this \grqq that ${\cal U}^{\cal U}$ is always greater than ${\cal U}$, but
according to the construction ${\cal U}^{\cal U}\subset{\cal U}$ what is a contradiction . 

But it could
have been simpler: to show the same diagonal Cantor's method.,
 that for any set ${\cal M}$, the set ${\cal M}^{\cal M}$ is not a subset
of the set ${\cal M}$, and therefore ${\cal M}$ is a proper subset of ${\cal M}\cup\{{\cal M}\}$. This instantly
implies that there is no universal set in the Hilbert sense. At
the same time, Hilbert's uncertainty does not matter, as \grqq combining
sets an infinite number of times into one set\grqq: no matter how many times
you combine all constructed sets into a universal set U,
the set ${\cal U}^{{\cal U}}$ always expands the set ${\cal U}$. 

Without any antinomy, this
Hilbert method shows that there are no universal sets, there are no
universal sets in the understanding of classical formal theories
$ZFC$ and $NBG^-$ , there are no universal classes for those formal theories
that allow their own classes to be elements of hyper classes,
etc.

{\bf 6.} {\it Propositions' paradox}

Let's recall the following \grqq Paradox of Propositions\grqq\,
 \cite{l212}, 2.3). Let ${\bf m}$ be a set of statements, and $\prod{\bf m}$ be a
statement such that \grqq every statement $m\in{\bf m}$ is true\grqq\, (considered as
a possible infinite conjunction of statements). Then, if ${\bf m}$ and ${\bf n}$ are
distinct sets, then the statements $\prod{\bf m}$ and $\prod{\bf n}$ are also distinct, and the mapping which for 
 every statement ${\bf m}$ associates the  statement   $\prod{\bf m}$ is
an injection. Therefore, if we consider the set ${\bf P=\{p\,|\,\exists{\bf m}(\prod{\bf m}=p\,\,\&\,\,p\notin{\bf m})\}}$, by virtue of the above-mentioned injectivity, we obtain a
contradiction. Indeed, $\prod{\bf P}\in{\bf P}\Leftrightarrow\prod{\bf P}\notin{\bf P}$.

{\it In the propositional paradox, the family of propositions
${\bf P}$ is not a set}.

Indeed, let's assume the opposite, and the family ${\bf P}$ of all such
statements $\{p\,|\,\exists{\bf m}(\prod{\bf m}=p\,\,\&\,\,p\notin{\bf m})\}$ is a set. Let $\prod{\bf P}$
be a one-element set consisting of a single statement
${\bf P}$. Then, by definition 
$\prod\{{\bf P}\}$, this is a statement that is true,
 since the only utterance of the set $\{{\bf P}\}$ consists of true
statements of $p$, which, by definition, are equal to $\prod{\bf m}$ for some set
${\bf m}$ of true statements of $ m$. But, obviously, that 
$\{{\bf P}\}$ does not combine
with ${\bf P}$, because of the same injection, since ${\bf P}\not=\{\prod{\bf P}\}$. 
Let's now consider the union of sets of statements ${\bf P}\cup\{\prod{\bf P}\}$, which is
a set and contains its proper subset ${\bf P}$, since we
have already shown that ${\bf P}$ is ${\bf P}$, which contradicts the maximality of the set
${\bf P}$.  Therefore, ${\bf P}$ is not a set, but is a proper class,
for which the concept of class conjunction $\prod{\bf P}$ is not applicable

\begin{center}
{\bf 6. There are essentially no universal sets for most predicates: $\neg\exists {\bf X}(\forall {\bf x}) ({\bf x}\in{\bf X})$ vs $\exists{\bf X}\vdash\neg\{{\bf X}\}$}
 \end{center}

In the previous two paragraphs, we have shown that the well-known set-theoretic paradoxes are imaginary, i.e., they  are really  not contradictions,  and have circumvented them using the maximality principle method. An important consequence of this approach is that the reality, generally speaking, is the following: \grqq There is no universal set containing all sets with some property\grqq, e.g, $X=\{x\,|\,x\notin x\}$ . More generally, \grqq For any class of elements, there is an element that does not belong to this class\grqq, in other words, any class of elements is incomplete. This of course contradicts the first axiom of the existence of the class of all pairs of sets $(a,b)$ for which $a\in b$ in $NBG$-theory. This axiom prohibits proper classes from being elements of other classes, for example, $\{{\bf V}\}$, but this does not mean that in fact there is no such pair as $\{\{{\bf V}\},\{{\bf V},\{{\bf V}\}\}$, when $a={\bf V}$ and $b=\{{\bf V}\}$. Formal theories allowing for this exist and are consistent: 7.5 Classes Taken Seriously -- the System of Quine and Morse; 7.6. Classes not Taken Seriously -- Systems of Bernays and Quine; 7.7. The System of Ackermann. See \cite{l14}, p. 138-153; see also  \cite{l222}, Chap. II, \S 7.

In this case, we conclude that the very statement, for example, by Russell: \grqq The class of all non-reflexive classes is contradictory\grqq, is incorrect, but it is true that \grqq the class of all non-reflexive sets does not exist, more precisely, it is indefinable\grqq, and according to Cantor's intuition, \grqq it is not a completed or  final thing\grqq. One can find another names for such thing as \grqq trans-definite sets without boundary\grqq\, see \cite{l315}, p. 311. Paul Cohen expressed the intuition of this indeterminacy of the eternal process of forming new  sets or numbers (see below about Conway numbers): \grqq Perhaps the main reason for accepting the axiom of infinity is the absurdity of the idea that adding only one set at a time could exhaust the entire world\grqq. (See \cite{l222}, Chap. IV, \S 13. Conclusion.)

This does not mean at all that there are no universal classes in mathematics, in particular, in set theory. On the contrary, they exist and there are quite a lot of them. For example, the set of all natural numbers. In a strict sense, the existence of this set is determined by the axiom of infinity. The \grqq incompleteness or non-final things\grqq, mentioned above, differ from the set of all natural numbers because they cannot be defined even axiomatically, since any such definition will make them complete or final, which is easily proved by the Russell functor: \grqq The subclass of all non-reflexive objects of this universal set is not contained in the universe\grqq, which contradicts a possible axiom: \grqq axiomatic existence of the universe\grqq. 

Such universe is strikingly different from the second example, namely, the existence of the set of all countable ordinals. Otherwise, if the set of all countable ordinals is not countable, then it itself, as an ordered set, defines the cardinal that is larger than all  cardinals included in it, and the set consisting of these ordinals and the last one will obviously be countable, which contradicts the universality of the assumed set of all countable ordinals. Thus, our assumption about its countability is false.

These \grqq incompleteness or non-final things\grqq\, are not mathematical objects because they can be compared and called as unfinished processes (\grqq eternal processes\grqq, perpetuum mobile processes sequences, etc.) that are continuing and cannot end, by definition. Indeed, for any set $X$ we imply Russell's functor $F_{\bf R}(X)\stackrel{def}{=}X\cup\{\{x\,|\,x\in X\,\&\,x\notin x\}\}$ and obtain a new set $X_1=F_{\bf R}(X)\supset X$. Then continue this process for each non-limit ordinal and for a limit ordinal $\alpha$ we consider $X_\alpha=\bigcup\limits_{\alpha'<\alpha}X_{\alpha'}$ and put $X_{\alpha+1}=F_{\bf R}(X_\alpha)$, and so on, and forever. Of course it is only an idea because there is no final object only tendency.

The Russell functor is nothing more than a counting functor, or a construction functor for a class of ordinal numbers. Indeed, let's set $X=\emptyset$ and obtain the class of all ordinals according to Mirimanoff and von Neumann:
$$
0=\emptyset, 1=\{\emptyset\}, 2=\{\emptyset,\{\emptyset\}\},...,\omega=\{0,1,2,...\},  \omega=\{0,1,2,...,\omega+1\},...
$$
for $0\leq\alpha<\Omega$ but actually we can continue this process $\Omega, \Omega+1,...,\Omega+\omega,...$  and go on forever, obtaining an indefinable object of arbitrary ordinal numbers $\alpha\geq 0$, or by Ulrich Blau, \grqq trans-definite ordinals without boundary\grqq \cite{l315}.

 Since each ordinal number is determined by a set of previously defined ordinal numbers, the statement about the possibility of the existence of a set of all ordinal numbers is incorrect. The Burali-Forti paradox consisted precisely in this incorrect assumption. And the following reality or situation is true: \grqq The ordinal number formation functor through all the previous ones (Cantor, Mirimanoff, von Neumann, Conway, Blau, etc.) starts the following well-ordered counting process: $F_{\bf Ord}(\alpha+1)=[0,\alpha)\cup\{\alpha\}$\grqq. And this process continues without completion or finality, let's say that it is a \grqq  limitlessness\grqq.

Considering that too extensive sets lead to contradictions, von Neumann in 1925 proposed his formal set theory, which excluded such extensive sets, requiring that they should not be elements of either the sets themselves or other equally extensive classes, i.e., in the terminology of Bernays and G{\" o}del, proper classes.

However, like the previous formal system of Zermelo Fraenkel  ($ZFC$) with Choice Axiom, which excluded such vast sets altogether, the $NBG$ system proceeded from the same intuition and the premise of the presence (not formalized in any way) of some well-defined nonempty class (objects), the so-called universe of reasoning, or the world under consideration (universe of discourse) is in the $ZFC$ system, and some class in which subsets are identified with their characteristic functions is in $NBG$, and in the latter case, functions that are element sets, i.e. arguments, and functions that are not arguments and not elements, but proper classes. And both assumed classes are static and do not contain any elements, sets, or classes outside of themselves. But in reality this is not the case, and such restrictions have no advantage other than some formal convenience and the imaginary elimination of set-theoretic antinomies.

Note that in an unfinished or non-final thing, in the sense of Cantor, there is no characteristic function for this thing, since it itself is incomplete or non-final.

A set is a collection of definite,  distinguishable objects of perception or thought which is perceived as a whole. The objects are known as elements or members of a set.

Following Cantor: \grqq Namely, multiplicity may have the property that the perception of the 
\grq co-existence\grq\, of all its elements leads to a contradiction, so that this multiplicity cannot be considered as a unity, as \grq some kind of complete thing\grq. I call such multiplicities absolutely infinite or inconsistent multiplicities. 
As it is easy to see, \grq the totality of everything imaginable\grq, for example, 
it is a similar multiplicity; other examples will appear later. 
On the contrary, if a set of elements of a certain multiplicity
can be thought of without contradiction as a set of 
\grq co-existing\grq\, elements, so that they can be combined into a \grq single
thing\grq, then I call it a consistent set or \grq set\grq\, (in French and Italian, this concept is appropriately
expressed by the words \grq ensemble\grq\, and \grq insieme\grq\,\grqq) \cite{l1}.

\bigskip

So, in  $NBG^-$-theory, there are objects that are not elements of any sets or classes, but are proper classes of ${\bf V}^-$. We will show that these objects can be subjected, for example, Conway numbers ${\bf No}$ to algebraic operations similar to operations on proper classes: intersection of $X\subset{\bf V}^-$ and $Y\subset{\bf V}^-$, i.e., $X\cap Y\subset{\bf V}^-$, union of $X$ and $Y$, i.e., $X\cup Y\subset{\bf V}^-$, difference of $X$ and $Y$, i.e., $X\setminus Y\subset{\bf V}^-$. Although $X$ and $Y$ cannot be elements, they as do other similar proper classes of ${\bf V}^-$ can be consider as an algebraic structure of a linearly ordered field similar to Boolean algebra for subsets. Thus, the formula $\vdash {\bf Pr}( X)\Rightarrow\neg\exists\{X\}$, i.e., non-existence of the {\it singleton} of $X$, that prohibits proper classes from being elements does not prevent these proper classes $X\subseteq{\bf No}$ from having an algebraic and topological structure. Thus, by introducing the concept of super-classes, we obtain the usual form of a topological linearly ordered field, but for the sake of argument, this can be done in the language of proper classes, although it is not as convenient as introducing super-classes. Let us do this carefully.

For this purpose we consider a linear ordering topology on ${\bf No}$, i.e., a class $X\in{\bf No}$ is {\it open} in a linear ordering topology on ${\bf No}$, if for every $x\in X$ there exists a positive number $\varepsilon\in{\bf No}$ such that the interval $(x-\varepsilon,x+\varepsilon)\subset X$.

Now we elucidate what transfinite $\Omega$-fundamental sequences of numbers in ${\bf No}$ and how they can be called {\it convergent} and what transfinite $\Omega$-fundamental  sequences of numbers in ${\bf No}$ why they can be called {\it non-convergent} in the ordering topology on ${\bf No}$.

{\bf Definition 13}. A mapping $x:{\bf On}\rightarrow {\bf No}$ is called a {\it transfinite sequence} of  $\Omega$-type of Conway numbers in ${\bf No}$, or $\Omega$-{\it transfinite sequence} in ${\bf No}$, or shortly a $\Omega$-{\it sequence} in ${\bf No}$. We denote it as above by $(x_\alpha)_{0\leq\alpha<\Omega}$, i.e., $x(\alpha)-x_\alpha$, ${0\leq\alpha<\Omega}$.

{\bf Definition 14}.
We say that $\Omega$-sequence $(x_\alpha)_{0\leq\alpha<\Omega}$ in  ${\bf No}$, {\it converges} to $a\in {\bf No}$, and we write $\lim\limits_{0\leq\alpha<\Omega}x_\alpha=a$,  if for each positive  number $\varepsilon\in {\bf No}$ there is an ordinal number $\alpha_0\in[0,\Omega)$ such that $|x_\alpha-a|<\varepsilon$ for all $\alpha_0\leq\alpha<\Omega$.  
In this case we also say that $\Omega$-sequence $(x_\alpha)_{0<\alpha<\Omega}$ is {\it convergent} in  ${\bf No}$. If $\lim\limits_{0\leq\alpha<\Omega}x_\alpha=0$, then $(x_\alpha)_{0\leq\alpha<\Omega}$ is called an $\Omega$-{\it infinitesimal sequence} in ${\bf No}$.

{\bf Definition 15.}
 An $\Omega$-sequence $(x_\alpha)_{0\leq\alpha<\Omega}$ in ${\bf No}$,  is called  {\it fundamental}  or a {\it Cauchy} $\Omega$-sequence,  if for each positive  number $\varepsilon\in {\bf No}$  there is an ordinal number $\alpha_0$ such that $|x_\alpha-x_{\alpha'}|<\varepsilon$, for all $\alpha_0\leq\alpha\leq\alpha'<\Omega$. 

{\bf Definition 16.}
Two $\Omega$-fundamental sequences $(x_\alpha)_{0\leq\alpha<\Omega}$ and $(y_\alpha)_{0\leq\alpha<\Omega}$ in ${\bf No}$ are $\Omega$-{\it equivalent}, denoted by as  $(x_\alpha)_{0\leq\alpha<\Omega}\sim(y_{\alpha})_{0\leq\alpha<\Omega}$,  if for each  positive number $\varepsilon\in {\bf No}$  there are ordinal numbers $\alpha_0$ and $\alpha'_0$ such that $|x_\alpha-y_{\alpha'}|<\varepsilon$, for all $\alpha_0\leq\alpha<\Omega$ and all $\alpha'_0\leq\alpha'<\Omega$. 

One can see that $(x_\alpha)_{0\leq\alpha<\Omega}\sim(y_{\alpha})_{0\leq\alpha<\Omega}$ and $\lim\limits_{0\leq\alpha<\Omega}x_\alpha=a$, then $\lim\limits_{0\leq\alpha<\Omega}y_\alpha=a$. Moreover, if $\lim\limits_{0\leq\alpha<\Omega}x_\alpha=a$ and $\lim\limits_{0\leq\alpha<\Omega}y_\alpha=b$, then $\lim\limits_{0<\alpha<\Omega}(x_\alpha+y_\alpha)=a+b$, $\lim\limits_{0\leq\alpha<\Omega}x_\alpha\cdot y_\alpha=a\cdot b$ and $\lim\limits_{0\leq\alpha<\Omega}\frac{x_\alpha}{y_\alpha}=\frac{a}{b}$ (the latter when $b\not= 0)$.

Recall here the definition of a Dedekind section $(L,R)$, induced by an $\Omega$-fundamental sequence $(x_\alpha)_{0\leq\alpha<\Omega}$. We put $L$ as  the subclass of ${\bf No}$ of all $x^L\in {\bf No}$ such that there exists an $\alpha_0$ and inequalities $ x^L<x_\alpha$ for all $\alpha_0\leq\alpha<\Omega$ and put $R= {\bf No}\setminus L$, i.e., two disjoint classes $L$ and $R$ whose union $L\cup R={\bf No}$, with  every member of $R$ exceeding every member of $L$; and in our special case (a section $(L,R)$ generated by {\it $\Omega$-fundamental sequence}), $L$ has no greatest and $R$ may or may not have a least number.

Clearlly, such Dedekind section $(L,R)$ cannot define a number $a=\{L\,|\,R\}$ because $L$ and $R$
are proper classes being illegal objects in $NBG$ (Conway called it a {\it gap} and speak about a collection of all gaps which us not even a proper class but an {\it improper class}, and speak of there being an impropriety of gaps, see \cite{l22}, 37-38).

We will now show that these non-existent objects in $NBG$-theory behave exactly like the usual Conway numbers, and they can also be expressed in the language of non-proper classes in theory $NBG$. Indeed, we identify non-existent number $a=\{L\,|\,R\}$ with an existent in $NBG$-theory a proper calss $L$. The equality relation $a=a'=\{L'\,|\,R'\}$ is given by an equality relation of proper class $A=A'$ in $NBG$-theory. A linear ordering $a<a'$ is given by a subclass relation  of proper class $A\subset A'$ in $NBG$-theory.

To define arithmetic operations on numbers $a=\{L\,|\,R\}$ and $b=\{L'\,|\,R'\}$, we must choose an $\Omega$-sequences that defines  Dedekind sections $(L,R)$ and $(L',R')$, respectively. Let's do it for $(L,R)$ (for $(L',R')$ it is clearlly similar). If $R$ has a least number say also $a$ (it is really a Conway number $a$ see below), then  the desired $\Omega$-sequence is $(a-\frac{1}{\alpha})_{0<\alpha<\Omega}$. One can see that generated by it Dedekind section is $(L',R')=(L,R)$. We call this $\Omega$-sequence is $(a-\frac{1}{\alpha})_{0<\alpha<\Omega}$ canonical.

If a Dedekind sections $(L,R)$ is such that $R$ has no least number, then, e.g., strong Choice Axiom or von Neumann Choice Axiom let for simplicity use the last: there is a bijection $\psi:{\bf On}\rightarrow{\bf No}$, i.e., each class $X\subseteq{\bf No}$ is well-ordered and thus has a least element $\psi(\min\psi^{-1}(X))$. Then we choose an $\Omega$-sequence is $(x_\alpha)_{0<\alpha<\Omega}$ which generates a Dedekind section $(L,R)$ by the following inductive construction. Put $x_0=\psi(\min\psi^{-1}(L))$ and  $x'_0=\psi(\min\psi^{-1}(R))$, respectively. Then consider a number $x''_0=\frac{x_0+x'_0}{2}$ and put $x_1=x''_0$ if $x''_0\in L$ and otherwise, put $x'_1=x_0$. For each $0<\alpha<\omega$ consider a number $x''_{\alpha-1}=\frac{x_{\alpha-1}+x'_{\alpha-1}}{2}$ and put $x_\alpha-x''_{\alpha-1}$ if $x''\in L$ and otherwise, put $x'_\alpha=x''_{\alpha-1}$. We also call this $\Omega$-sequence canonical.

Put $x_\omega=\psi(\min\psi^{-1}(\{x\,|\, x\in L\,\&\, x>x_\alpha, 0\leq\alpha<\omega\}))$ and $x'_\omega=\psi(\min\psi^{-1}(\{x'\,|\, x'\in R\,\&\, x'<x'_\alpha, 0\leq\alpha<\omega\}))$. And then consider number $x''_\omega=\frac{x_\omega+x'_\omega}{2}$ and so on like above.

Now we can define arithmetic operations on numbers $a=\{L\,|\,R\}$ and $b=\{L'\,|\,R'\}$ as follows: $a+ b=\{L\oplus L'\,|\,{\bf No}\setminus L\oplus L'\}$, $a\cdot b=\{L\odot L'\,|\,{\bf No}\setminus L\odot L'\}$ and $\frac{a}{b}=\{\frac{L}{L'}\,|\,{\bf No}\setminus\frac{L}{L'}\}$, where $L\oplus L'$, $L\odot L'$ and $\frac{L}{L'}$ are proper classes defined by canonical $\Omega$-sequences $(x_\alpha+y_\alpha)_{0\leq\alpha<\Omega}$, $(x_\alpha \cdot y_\alpha)_{0\leq\alpha<\Omega}$ and $(\frac{x_\alpha}{y_\alpha})_{0\leq\alpha<\Omega}$, respectively, $L'\not={\bf No}^-$ in the last case.

Note that  $\Omega$-sequences $(x_\alpha+y_\alpha)_{0\leq\alpha<\Omega}$, $(x_\alpha \cdot y_\alpha)_{0\leq\alpha<\Omega}$ and $(\frac{x_\alpha}{y_\alpha})_{0\leq\alpha<\Omega}$ are not canonical $\Omega$-sequences but are $\Omega$-equivalent to them.
 
Thus we can express non-legal objects in terms of legal proper classs in ${\bf V}$ of $NBG$-theory. In other words, we can reject  axiom $\exists X\vdash\neg\exists\{X\}$ and  abandon this axiom, by introducing the axiom of the existence of super-classes containing proper classes in ${\bf V}$ as their elements, and then express the previous statement in a more convenient form: define the collection  of all fundamental Dedekind sections in ${\bf No}$, denoted by ${\bf R}_\Omega$, as the completion of the field ${\bf No}$. In this case, the following statement is true: \grqq For every Dedekind section $(L,R)$ in ${\bf R}_\Omega$, defined by a fundamental $\Omega$-sequence $(x_\alpha)_{0\leq\alpha<\Omega}$, $x_\alpha\in{\bf R}_\Omega$,  has a smallest element  in $R$\grqq.

The second step is to define not a symbol $\Omega$, as above, but a real ordinal number $\Omega$ as well as others like $\Omega+1$, $\Omega+2$,..., $\Omega+\omega$,... in the language of proper classes in ${\bf No}$ of  $NBG$-theory, which in the form of a Conway number would be much better expressed by Conway's form: $\Omega=\{{\bf On}\,|\,\}$,  $\Omega+1=\{{\bf On},\Omega\,|\,\}$, $\Omega+2=\{{\bf On},\Omega,\Omega+1\,|\,\}$, ..., $\Omega+\omega=\{{\bf On},\Omega,...,\Omega+\alpha\,|\,\}_{0<\alpha<\omega}$, ...\,.

Let's recall the original Cantor's definition of the first derivative  set $P'$ (another form $P^{(1)}$) of a fixed set $P$ of real numbers ($P\subset{\bf R}$) as the set  of all its limiting points in linearly ordered topology (actually in our notation $\omega$-topology), then further derivatives of this set may in some cases either {\it end in an empty set}, or become a {\it perfect set} whose derivative set coincides with it. In the first case, the sequence of derived sets determines the ordinal number $0\leq\alpha<\omega$ when $P^{(\alpha)}\not=\emptyset$  the $(\alpha+1)$-derivative $P^{(\alpha+1)}=\emptyset$.

In our case ${\bf No}$ with above $\Omega$-topology we can construct an example of a proper class $P=(\bigcup\limits_{0\leq\alpha<\Omega}P_\alpha)\cup\{0\}\subset[0,2]\subset{\bf No}$ such that $P_\alpha^{(\alpha)}={\frac{1}{\alpha+1}}$, $0\leq\alpha<\Omega$,  and $P^{(0)}_0=P$ and $P^{(\Omega)}=\{0\}$ with the following property: $P\supset P'\supset ... \supset P^{(\alpha)}\supset P^{(\alpha)}\supset...\supset P^{(\Omega)}=\{0\}$, $0\leq\alpha\leq\Omega$.
Thus proper classes $P^{(\alpha)}\subset{\bf No}$, $0\leq\alpha<\Omega$  together with $\emptyset$ determine the ordinal  number $\Omega$.

Here is a construction of this example, preceding this construction by defining   a {\it sum} of of two classes $S\subset{\bf No}$ and $T\subset{\bf No}$ as $S+T\stackrel{def}{=}\{s+t\,\,|\,\,s\in S\,\&\,t\in T\}$, and using a {\it symbol of infinity} $+\infty_{\frac{1}{\omega}}$, induced by a sequence $(\frac{1}{\lambda})_{0<\alpha<\lambda}$, where $\lambda$ is a limit ordinal, which will be defined below in Part Three.

Put $Q_0=\{0,1,2,\frac{1}{2},\frac{1}{3},...,\frac{1}{\alpha},...\}_{0<\alpha<\Omega}$; put $P_0=Q_0\cap[2,3]=\{2\}$; clearly, $P'_0=\{2\}'=\emptyset$ and thus $P_0$ presents the ordinal $0$.

Put $Q_1=(Q_0+\{0,1,\frac{1}{2},\frac{1}{3},...,\frac{1}{\alpha},...\}_{0<\alpha<\Omega})\cap[0,2]$; put $P_1=Q_1\cap[1,2]$; clearly,  $(P_1')'=\{1\}'=\emptyset$ and thus $P_1$ presents the ordinal $1$.

Put $Q_2=(Q_1+\{0,1,\frac{1}{2},\frac{1}{3},...,\frac{1}{\alpha},...\}_{0<\alpha<\Omega})\cap[0,1]$; put $P_2=Q_2\cap[\frac{1}{2},1]$; clearly, $(P_2'')'=\{\frac{1}{2}\}'=\emptyset$ and thus $P_2$ presents the ordinal $2$.

.............................

For each $1<n<\omega$ we put $Q_n=(Q_{n-1}+\{0,1,\frac{1}{2},\frac{1}{3},...,\frac{1}{\alpha},...\}_{0<\alpha<\Omega})\cap[0,\frac{1}{n-1}]$; put $P_n=Q_n\cap[\frac{1}{n},\frac{1}{n-1}]$; clearly, $(P_n^{(n)})'=\{\frac{1}{n}\}'=\emptyset$ and thus $P_n$ presents the ordinal $n$.

.............................

Put $Q_\omega=(\bigcup\limits_{0\leq n<\omega}Q_n+\{0,1,\frac{1}{2},\frac{1}{3},...,\frac{1}{\alpha},...\}_{0<\alpha<\Omega})\cap[0,+\infty_{\frac{1}{\omega}}]$; put $P_\omega=Q_\omega\cap[\frac{1}{\omega},+\infty_{\frac{1}{\omega}})$; clearly, $(P_\omega^{(\omega)})'=\{\frac{1}{\omega}\}'=\emptyset$ and thus $P_\omega$ presents the ordinal $\omega$. 

Put $Q_{\omega+1}=Q_\omega+\{0,1,\frac{1}{2},\frac{1}{3},...,\frac{1}{\alpha},...\}_{0<\alpha<\Omega})\cap[0,\frac{1}{\omega}]$; put $P_{\omega+1}=Q_{\omega+1}\cap[\frac{1}{\omega+1},\frac{1}{\omega}]$; clearly, $(P_{\omega+1}^{(\omega+1)})'=\{\frac{1}{\omega+1}\}'=\emptyset$ and thus $P_{\omega+1}$ presents the ordinal $\omega+1$.

.............................

If $\beta$ is a limit ordinal, then put $Q_\beta=(\bigcup\limits_{0< \beta'<\beta}Q_{\beta'}+\{0,1,\frac{1}{2},\frac{1}{3},...,\frac{1}{\alpha},...\}_{0<\alpha<\Omega})\cap[0,+\infty_{\frac{1}{\beta}})$; put $P_\beta=Q_\beta\cap[\frac{1}{\beta},+\infty_{\frac{1}{\beta}})$; clearly, $(P_{\beta}^{(\beta)})'=\{\frac{1}{\beta}\}'=\emptyset$ and thus $P_{\beta}$ presents the ordinal $\beta$. 
 
If $\beta$ is not a limit ordinal, then put $Q_\beta=(Q_{\beta-1}+\{0,1,\frac{1}{2},\frac{1}{3},...,\frac{1}{\alpha},...\}_{0<\alpha<\Omega})\cap[0,\frac{1}{\beta-1}]$; put $P_\beta=Q_\beta\cap[\frac{1}{\beta},\frac{1}{\beta-1}]$; clearly, $(P_\beta^{(\beta)})'=\{\frac{1}{\beta}\}'=\emptyset$ and thus $P_\beta$ presents the ordinal $\beta$.

And so on
.............................

At last, obtain the proper class
$P=(\bigcup\limits_{ \alpha\in{\bf On}}P_\alpha)\cup\{0\}$; clearly, $(P^{(\Omega)})'=\{0\}'=\emptyset$ and thus $P$ presents the ordinal $\Omega$.

 The property $P\supset P'\supset ... \supset P^{(\alpha)}\supset P^{(\alpha)}\supset...\supset P^{(\Omega)}=\{0\}$, $0\leq\alpha\leq\Omega$ is clear by the construction because the number $0$ is a limit point of all proper classes $P^{(\alpha)}$, $0<\alpha<\Omega$ since in each proper set there is a coinitial subsequence of $\Omega$-sequence $(\frac{1}{\alpha})_{)<\alpha<\omega}$.

Notice that $U_1=P+Q_0$ has the following property: $U_1^{(\Omega)}=Q_0$ and thus $({U_1^{(\Omega)}})''=(Q_0)''=\{0\}'=\emptyset$. Consequently, the class $U_1$ and all its limit points determine the ordinal number $\Omega+1$; similar $U_2=U_1+Q_0$ for which $({U_2^{(\Omega+1)}})''=(Q_0)''=\{0\}'=\emptyset$. Consequently, the class $U_2$ and all its limit points determine the ordinal number $\Omega+2$;  and such similar process has no boundary.

(This example was first presented in \cite{l1332}, p. 35-36, unfortunately with some typos that we have corrected here.)

{\bf Remark 9.} The above examples of proper classes $P$, $U_1$, $U_2$,...\,, which in $\Omega$-topology define a well-ordered proper class of limit points, raise an open question: \grqq Is it possible to use $\Omega$-topology to construct any trans-ordinal number $\gamma\geq\Omega$? Or do we need $\Gamma$-topologies of higher ordinals $\Gamma>\Omega$ for this?" A related question is: "Can the $\omega$-topology of real numbers ${\bf R}$ define all ordinal numbers $0\leq\alpha<\Omega$?" The last question is connected with the Continuum Hypothesis $CH$ because we can, using Choice Axiom, well-ordered all irrational numbers $\xi\in[0.1]\subset{\bf R}$ and in above example in the $\omega$th step, we put 
$P^1_\omega=(\bigcup\limits_{0< \alpha<\omega}P_\alpha)\cup\{0\}$, where
for each $1<n<\omega$ we put $Q_\beta=(Q_{\beta-1}+\{0,1,\frac{1}{2},\frac{1}{3},...,\frac{1}{\alpha},...\}_{0<\alpha<\omega})\cap[0,\frac{1}{\beta-1}]$, $P_\beta=Q_\beta\cap[\frac{1}{\beta},\frac{1}{\beta-1}]$, $0<\beta<\alpha$, and $Q_0=\{0,1,\frac{1}{2},\frac{1}{3},...,\frac{1}{\alpha},...\}_{0<\alpha<\omega})$, $Q_1=Q_0\cap[0,2]$ and $P_1=Q_1\cap[1,2]$.

In the $(\omega+1)$th step, we put  $Q^1_{\omega+1}=P^1_\omega+\{0,\xi_1,\frac{\xi_1}{2},\frac{\xi_1}{3},...,\frac{\xi_1}{\alpha},...\}_{0<\alpha<\omega}\}\cap[\xi_1,2]$ and for each $1<n<\omega$ we put $Q^1_{\omega+\beta}=(Q^1_{\omega+\beta-1}+\{0,\xi_1,\frac{\xi_1}{2},\frac{\xi_1}{3},...,\frac{\xi_1}{\alpha},...\}_{0<\alpha<\omega})\cap[0,\frac{\xi_1}{\beta-1}]$, $P^1_{\omega+\beta}=Q_{\omega+\beta}\cap[\frac{\xi_1}{\beta},\frac{\xi_1}{\beta-1}]$, $0<\beta<\alpha$. 

Now we put 
$P^2_{2\cdot\omega}=\bigcup\limits_{0<\beta<\omega}P^1_{0<\alpha<2\cdot\omega}\cup\{0\}$
and after the $2\cdot\omega$th step we assume 
$Q^2_{2\cdot\omega+1}=P^2_{2\cdot\omega}+\{0,\xi_2,\frac{\xi_2}{2},\frac{\xi_2}{3},...,\frac{\xi_2}{\alpha},...\}_{0<\alpha<\omega}\}$ and so on, i.e., using the choice function $\psi:X\rightarrow X$ on subsets of all irrational numbers $X\subset[0,1]\subset{\bf R}$ and putting $\psi(X)=\xi_1\in X$, $\psi(X_1)=\xi_2\in X_1$, ..., $\psi(X_\alpha)=\xi_{\alpha+1}\in X_\alpha$, ..., $0\leq\alpha<\Omega$, and so on, where $X_1=X\setminus\{a\cdot\xi_1+b\}_{a,b\in{\bf Q}\cup[0,1]}$, ... , $X_\alpha=X\setminus\bigcup\limits_{0<\beta<\alpha} X_\beta$ and ${\bf Q}$ are all rational numbers,  and  $\psi(X_\alpha)=\xi_{\alpha+1}\in X_\alpha$ and so on. 

Then in $ZFC$-theory with addition Axiom of Non-existence of an Inaccessible cardinal conclude that, by $\neg CH$, $|{\bf R}|=2^{\bf \aleph_0}$ can be greater than or equal to any ${\bf \aleph}_\alpha$, $1<\alpha<\Omega$ (see \cite{l222}, Chap. IV, \S 13. Conclusion) and thus irrational numbers enough to define ordinals via limit points of some set $P\subset[0.1]\subset{\bf R}$. Notice also that $2^{\bf \aleph_0}$ can be naturaly equal to $\aleph_2$  (see \cite{l939}, Theorem 12).

We are going to explain this state of affairs.

In spite of the fact that the formal logical system $NBG^-$ is convenient in many cases (e.g., Conway theory of numbers and games, as well as  Category Theory and many other mathematical theories and constructions) by the finiteness of its axioms, by existence of universal objects (proper classes) identified with \grqq properties\grqq\, or \grqq singular proposition functions\grqq, etc., there is a lack in it which, by thought of P. J. Cohen, the \grqq theory $NBG^-$ is a less intuitive system than $ZF$\grqq (\cite{l222}, Cpt. 3, $\S 5$). Indeed, such formal systems are very formal and they were inputed to avoid so-called set-theoretic paradoxes of Russell's type. Often these \grqq tricks\grqq\, (proofs by contradiction via paradoxes) were out of contents and concepts of mathematical notions and definitions like \grqq reflexive sets\grqq, \grqq universal objects\grqq, \grqq premises of false propositions\grqq, \grqq undetectable  objects\grqq, etc.

The only correct and meaningful conclusion from Russell's paradox is the following: \grqq Any family of non-reflexive elements is always non-reflexive and is not contained in the original family\grqq. As a consequence of it is the following: \grqq There is no universal families like all sets, all ordinals, all groups, all rings, all fields, etc.\grqq, but there are internal processes, perpetuum mobile processes sequences, ordinal numbers without border, endlessness or something like that.

 \begin{center}
{\bf Part Three}
\end{center}

\begin{center}
{\bf 1. Conway numbers as a platform for various mathematical and other models}
\end{center}

As we saw in the previous section, there is no such mathematical object as the collection of {\it all ordinals} although each ordinal exists individually, since it is determined by all the preceding ordinals  the totality of which is a mathematical object because it is a set or a proper class.

Moreover, if for each ordinal we can enter a symbol, as, for instance, above $\alpha$, then for a non-existent \grqq collection\grqq\, of all ordinals we cannot formally enter a symbol, for instance, ${\cal O}$ and write the formula $\alpha\in{\cal O}$ because the right part of it ${\cal O}$ does not exist. 

On the other hand, since we can consider a similar {\it eternal process} of ordinal numbers of the form $\xi=\omega^{\omega^\nu}$, $\nu\geq 0$, we obtain the corresponding process of increasing linearly ordered fields of Conway numbers ${\bf R}_\zeta\subset {\bf R}_{\zeta'}$, where $\zeta'=\omega^{\omega^{\nu+1}}$, $\nu\geq 0$, and then for any ordinal $\alpha$ we can write the correct formula $\alpha\in{\bf R}_\zeta$ for some fixed and even smallest ordinal number $\nu\geq 0$, since both $\alpha$ and ${\bf R}_\zeta$ are determined by their existing as mathematical objects.

Thus, a linearly ordered field of Conway numbers ${\bf R}_\zeta$ of the highest power can be a platform for models of mathematical formal systems.

However, now we will consider strictly increasing and strictly decreasing sequences of Conway numbers $(x_\alpha)_{0<\alpha<\lambda}$, $x_\alpha<x_{\alpha+1}$, and $(x_\alpha)_{0<\alpha<\lambda}$, $x_\alpha>x_{\alpha+1}$, for limit ordinals $\lambda\geq\omega$, including eternal processes $(x_\alpha)_{0<\alpha}$, $x_\alpha<x_{\alpha+1}$, and $(x_\alpha)_{0<\alpha}$, $x_\alpha>x_{\alpha+1}$ for every ordinal $\alpha$. 

The simplest such eternal sequences are the following: $(\alpha)_{0<\alpha}$ and $(-\alpha)_{0<\alpha}$, $(\frac{1}{-\alpha})_{0<\alpha}$ and $(\frac{1}{\alpha})_{0<\alpha}$, $(\frac{1}{\omega^{\frac{1}{\alpha}}})_{0<\alpha}$ and $(\omega^{\frac{1}{\alpha}})_{0<\alpha}$, respectively.

The structure of a field of the Conway numbers  allows us to construct similar sequences, for example, around each Conway number $a$ as $(a-\frac{1}{\alpha})_{0<\alpha}$
and $(a+\frac{1}{\alpha})_{0<\alpha}$, i.e., say, two infinite processes  tending to each other. 

Then the natural question arises: \grqq  What tends towards a strictly increasing sequence of $(\alpha)_{0<\alpha<\omega}$ or a strictly decreasing sequence of $(\frac{1}{\alpha})_{0<\alpha<\omega}$, respectively?\grqq. The answer is very simple they are $(\omega^{\frac{1}{\alpha}})_{0<\alpha}$ and $(\frac{1}{\omega^{\frac{1}{\alpha}}})_{0<\alpha}$, respectively. But in this case there is no Conway number $a$ between these sequences  tending to each other. In $NBG$ theory between them there is a gap which Conway called a \grqq gap of the second kind form\grqq\,  \cite{l22}, p. 37.

The structure of a field of the Conway numbers  allows us to construct similar consecutive pairs: 
$$(\alpha)_{0<\alpha<2\cdot\omega}\,\,\,\&\,\,\, (\omega +\omega^{\frac{1}{\alpha}})_{0<\alpha};$$
$$(\alpha)_{0<\alpha<3\cdot\omega}\,\,\,\&\,\,\,(2\cdot\omega +\omega^{\frac{1}{\alpha}})_{0<\alpha},...\,;$$
$$(\alpha)_ {0<\alpha<\omega^2}\,\,\,\&\,\,\,(\omega\cdot\omega^{\frac{1}{\alpha}})_{0<\alpha},...\,;$$
$$(\alpha)_ {0<\alpha<\omega^3}\,\,\,\&\,\,\,(\omega^2\cdot\omega^{\frac{1}{\alpha}})_{0<\alpha},...\,;$$
$$(\alpha)_{0<\alpha<\omega^\omega}\,\,\,\&\,\,\,(\omega^{\omega^{\frac{1}{\alpha}}})_{0<\alpha},  ...\,;$$ 
$$ (\alpha)_{0<\alpha<\omega^{\omega^\omega}}\,\&\,(\omega^{\omega^{\omega^{\frac{1}{\alpha}}}})_{0<\alpha},  ...\,;$$
 $$(\alpha)_{0<\alpha<\varepsilon_0}\,\&\, (\varepsilon_{-\alpha})_{0<\alpha},...,$$
where 
  $$\varepsilon_{-\alpha}=\{ordinals<\varepsilon_0\,|\,\varepsilon_{-\alpha+1}-1,\omega^{\varepsilon_{-\alpha+1}-1},\omega^{\omega^{\varepsilon_{-\alpha+1}-1}},\omega^{\omega^{\omega^{\varepsilon_{-\alpha+1}-1}}},...\},$$
if $0<\alpha$ is not a limiting ordinal; otherwise,
$$\varepsilon_{-\alpha}=\{ordinals<\varepsilon_0\,|\,\varepsilon_0,\varepsilon_{-1},\varepsilon_{-2},...,\varepsilon_{-\alpha'},...\},$$
where $0\leq\alpha'<\alpha$, in particular, the case $\alpha=0$ is also limiting since there is no the previous ordinal because $a'<0$ is a Conway number but it is not an ordinal number. Nevertheless, we obtain $\varepsilon_0=\{ordinals<\varepsilon_0\,|\,\}=\varepsilon_0\stackrel{def}{=}\{\omega,\omega^{\omega}.\omega^{\omega^{\omega}},...\,|\,\}.$ Note also that for all limit ordinals $0\leq\lambda<\alpha$ in series $\varepsilon_0,\varepsilon_{-1},\varepsilon_{-2},...,\varepsilon_{-\alpha'},...$, $0\leq\alpha'<\alpha$,  by induction, we have already defined numbers $\varepsilon_{-\lambda}$.
$$(\alpha)_{0<\alpha<\varepsilon_{1}}\,\&\,(\varepsilon_{\frac{1}{2^\alpha}})_{0<\alpha},...\,;$$
where $$\varepsilon_{\frac{1}{2}}=\{ordinals<\varepsilon_1\,|\,\varepsilon_1-1,\omega^{\varepsilon_1-1},\omega^{\omega^{\varepsilon_1-1}},...\}\,,$$
$$\varepsilon_{\frac{1}{2^2}}=\{ordinals<\varepsilon_1\,|\,\varepsilon_{\frac{1}{2}}-1,\omega^{\varepsilon_{\frac{1}{2}}-1},\omega^{\omega^{\varepsilon_{\frac{1}{2}}-1}},...\,,\}$$
$$...................$$
$$\varepsilon_{\frac{1}{2^\alpha}}=\{ordinals<\varepsilon_1\,|\,\varepsilon_{\frac{1}{2^{\alpha-1}}}-1,\omega^{\varepsilon_{\frac{1}{2^{\alpha-1}}}-1},\omega^{\omega^{\varepsilon_{\frac{1}{2^{\alpha-1}}}-1}},...\}\,,$$
if $\alpha$ is not a limit ordinal; otherwise,
$$\varepsilon_{\frac{1}{2^\alpha}}=\{ordinals<\varepsilon_1\,|\,\varepsilon_1,\varepsilon_{\frac{1}{2^2}},...,\varepsilon_{\frac{1}{2^{\alpha'}}},...\}\,,$$
for all $0\leq\alpha'<\alpha$, and so on.

More general,
$$(\alpha)_{0<\alpha<\varepsilon_\nu}\,\&\,(\varepsilon_{\nu-1+\frac{1}{2^\alpha}})_{0<\alpha}\,,$$
where $$\varepsilon_{\nu-1+\frac{1}{2^\alpha}}=\{ordinals<\varepsilon_{\nu}\,|\,\varepsilon_{\nu-1+\frac{1}{2^{\alpha-1}}}-1,\omega^{\varepsilon_{\nu-1+\frac{1}{2^{\alpha-1}}}-1},\omega^{\omega^{\varepsilon_{\nu-1+\frac{1}{2^{\alpha-1}}}-1}},...\}\,,$$
if $\alpha$ and $\nu$ are not limit ordinals; otherwise,
$$\varepsilon_{\nu-1+\frac{1}{2^\alpha}}=\{ordinals<\varepsilon_{\nu}\,|\,\varepsilon_\nu,\varepsilon_{\nu-\frac{1}{2}},\varepsilon_{\nu-\frac{3}{2^2}},...,\varepsilon_{\nu-1+\frac{1}{2^{\alpha'}}},...\},$$ for $0\leq\alpha'<\alpha$, and
$$\varepsilon_{\nu+\frac{1}{2^\alpha}}=\{ordinals<\varepsilon_{\nu}\,|\,\varepsilon_\nu-\varepsilon_{\frac{1}{2^\alpha}},\varepsilon_{\nu}-\varepsilon_{1+\frac{1}{^\alpha2}},...,\varepsilon_\nu-\varepsilon_{\nu'+\frac{1}{2^\alpha}},...\}\,,$$ for $0\leq\nu'<\nu$,  respectively.

In particular of the latter case, considering that the initial ordinals are $\varepsilon$-numbers, i.e., $\omega_\nu=\{\varepsilon_0,\varepsilon_1,...,\varepsilon_{\nu'},...\,|\,\},$ where $0\leq\nu'<\omega_\nu$, and thus $\omega^{\omega_\nu}=\omega_\nu$, $\nu>0$, we obtain
$$(\alpha)_{0<\alpha<\omega_1}\,\&\,(\delta_\alpha)_{0<\alpha},...\,;$$
as well as in general
$$ (\alpha)_{0<\alpha<\omega_\nu}\,\&\,(\delta_\alpha)_{0<\alpha},...\,,$$
 for any ordinal $\nu> 0$ and so on, we correspond the following decreasing sequences:
 $$\delta_1=\{ordinals<\omega_\nu\,|\,\omega_\nu-1,\omega^{\omega_\nu-1},\omega^{\omega^{\omega_\nu-1}},...\,;$$
  $$\delta_2=\{ordinals<\omega_\nu\,|\,\delta_1-1,\omega^{\delta_1-1},\omega^{\omega^{\delta_1-1}},...\}\,;$$
$$....................$$
  $$\delta_\alpha=\{ordinals<\omega_\nu\,|\,\delta_{\alpha-1}-1,\omega^{\delta_{\alpha-1}-1},\omega^{\omega^{\delta_{\alpha-1}-1}},...\}\,;$$
if $\alpha$ is not a limiting ordinal; otherwise,
$$\delta_\alpha=\{ordinals<\omega_\nu\,|\,\omega_\nu,\delta_1,\delta_2,...,\delta_{\alpha'},...\}\,,$$
where $0<\alpha'<\alpha$.

Thus, the set of all ordinals that define a limit ordinal forms a sequence that defines a second-kind gap in $NBG$ theory.
\bigskip

{\bf Remark 10.} It is worth noting here that  Cantor's general form of the $\varepsilon$-numbers $\varepsilon_\nu$, $0\leq\nu<\omega_1$, (\cite{l1}, $\S 20$) was given by the  formulas   
$E(\gamma)=\lim\limits_{\nu}\gamma_\nu$, where $$\gamma_1=\omega^{\gamma},\,\,\,\, \gamma_2=\omega^{\gamma_1},\,\,\,\,... ,\,\,\,\,\gamma_\nu=\omega^{\gamma_{\nu-1}}, \,\,\,\,... \,, \,\,\,\,
0<\nu<\omega,\,\,\,\, \nu_0=\gamma,$$ i.e., $E(\gamma)$ is the smallest ordinal which is greater than $\gamma_\nu$ for all $0\leq\nu<\omega$. One can see that $\omega^{E(\gamma)}=E(\gamma)$ since the function $y=\omega^\alpha$, $0\leq\alpha<\omega_1$ is continuous, i.e., $\omega^{E(\gamma)}=\lim\limits_{0\leq\nu<\omega}\omega^{\gamma_\nu}=\lim\limits_{0\leq\nu<\omega}\gamma_\nu=
E(\gamma)$. Then Cantor put $\varepsilon_0=E(1)$ and
$\varepsilon_\nu=E(\varepsilon_{\nu-1}+1)$,  when $\nu$ is not a limit ordinal and $\varepsilon_\nu=\lim\limits_{0\leq\nu'<\nu}\varepsilon_{\alpha'}$ when $\alpha$ is a limit ordinal. Clearly, $\omega^{\varepsilon_\nu}=\varepsilon_\nu$ for all $0\leq\nu<\omega_1$ (in fact, also for all ordinals $(\omega_1\leq\nu<\Omega)$ and for every super-ordinal $\nu\geq\Omega$, which Cantor did not consider).
Conway (\cite{l22}, p. 33-36) defined the irreducible numbers that generalize the concept of ordinal $\varepsilon$-numbers, moreover, the birthday of any  irreducible number is an $\varepsilon$-number, and  $\varepsilon_\alpha$-numbers for all $0\leq\alpha<\Omega$ are given by the following formulas: $\varepsilon_\alpha=\{\varepsilon_{\alpha-1}+1, \omega^{\varepsilon_{\alpha-1}+1},\omega^{\omega^{\varepsilon_{\alpha-1}+1}}, ...\,|\,\}$ when $\alpha$ is not a limit ordinal and $\varepsilon_\alpha=\{\varepsilon_0, \varepsilon_1, ..., \varepsilon_{\alpha'}, ...\,|\,\}_{0\leq\alpha'<\alpha}$ when $\alpha$ is not a limit ordinal, where $\varepsilon_0$ is the first ordinal $\varepsilon$-number is greater than $\omega$, namely the number $\varepsilon_0=\{\omega,\omega^\omega,\omega^{\omega^\omega}, ...\,|\,\}$.  Then using the exponential function $y=\omega^x$ defined (here formula $(\ref{f314})$),  he determined the first to the left and the first to the right of $\varepsilon_0$ epsilon-numbers $\varepsilon_{-1}$ and $\varepsilon_1$, respectively, and thens $\varepsilon_2$ the first to the right of $\varepsilon_1$ and $\varepsilon_{-\frac{1}{2}}$ the fist between $\varepsilon_0$ and $\varepsilon_1$ and so on, i.e., all  epsilon-numbers $\varepsilon_a$, for $a\in{\bf No}$. Notice also that the function $y=e^x$ defines the only epsilon-number $\varepsilon=\omega$, since by definition $e^\omega=\omega$ and $e^x\not=x$ for each $x\not=\omega$.  (In \cite{l1333} $\omega$ was included in a similar sequence of different epsilon-numbers $\tilde\varepsilon_\gamma$ such that  $\tilde\varepsilon_{\gamma+1}= {\tilde\varepsilon_\gamma}^{ \tilde\varepsilon_{\gamma+1}}$, namely in the form of a sequence $\tilde\varepsilon_\gamma=\{\tilde\varepsilon_{\gamma-1}, {\tilde\varepsilon_{\gamma-1}}^{\tilde\varepsilon_{\gamma-1}}, {\tilde\varepsilon_{\gamma-1}}^{{\tilde\varepsilon_{\gamma-1}}^{\tilde\varepsilon_{\gamma-1}}},\,|\,\}$ when $\gamma$ is not a limit ordinal and $\tilde\varepsilon_\gamma=\{\tilde\varepsilon_0, \tilde\varepsilon_1, ..., \tilde\varepsilon_{\alpha'}, ...\,|\,\}_{0\leq\gamma'<\gamma}$ when $\gamma$ is not a limit ordinal, where $\tilde\varepsilon_0=\omega$ is the first ordinal $\tilde\varepsilon$-number is greater than $\alpha<\omega$, namely the number $\tilde\varepsilon_0=\{2,2^2,2^{2^2}, ...\,|\,\}$; clearly,  $\tilde\varepsilon_1=\varepsilon_0$.) It is also interesting that for each $\varepsilon_\alpha$-number in the sense of Conway there is $\ln \varepsilon_\alpha=\omega\cdot\varepsilon_\alpha$ and the logarithmic function $y=\ln x$ behaves in an unexpected and atypical way, e.g., $\ln \varepsilon_\alpha=\omega\cdot\varepsilon_\alpha>\varepsilon_\alpha$ although it behaves, say  \grqq normally\grqq\, later on, e.g., $\ln (\varepsilon_\alpha\cdot\omega^\omega)=\ln \varepsilon_\alpha+\ln \omega^\omega=\omega\cdot\varepsilon_\alpha+\omega^2<\varepsilon_\alpha\cdot\omega^\omega$ but say \grqq periodically\grqq, i.e., $\ln \varepsilon_{\alpha+1}=\omega\cdot\varepsilon_{\alpha+1}>\varepsilon_{\alpha+1}$, $\ln (\varepsilon_{\alpha+1}\cdot\omega^\omega)=\ln \varepsilon_{\alpha+1}+\ln \omega^\omega=\omega\cdot\varepsilon_{\alpha+1}+\omega^2<\varepsilon_{\alpha+1}\cdot\omega^\omega$  and so on, for each $0\leq\alpha<\Omega$. (Note at last that $\varepsilon_\alpha$-numbers have an important application in solving the problems of consistency in formal logical systems and in the theory of artificial intelligence $AI$.)
\bigskip

Now we ask the last question: \grqq Can two eternal processes of sequences of Conway numbers $(x_\alpha)_{0<\alpha}$, $x_\alpha<x_{\alpha+1}$, and $(y_\alpha)_{0<\alpha}$, $y_\alpha>y_{\alpha+1}$, respectively, tend to each other so that there would be no Conway number $a$ between them?\grqq. We will show now that the answer is positive.

We present an algorithm,  called \grqq left-right choice\grqq\, for short (clearly like a possible 
 \grqq right-left choice algorithm\grqq). Consider the \grqq segment\grqq\, $[0,1]$ and let's call its end $0$ left and end $1$ right. We divide the \grqq segment\grqq\, $[0,1]$ in half and let's call the \grqq segment\grqq\, $[0,\frac{1}{2}]$ left and \grqq segment\grqq\, $[\frac{1}{2},1]$ right. 

Note that by a \grqq segment\grqq\, $[a,b]$ in quotation marks we mean so far only two Conway numbers -- the beginning $a$ and the end $b$ of this \grqq segment\grqq, where $a<b$, and the possible Conway numbers inside it, i.e., all numbers $x\in{\bf R}_\zeta$, $\zeta=\omega^{\omega^\nu}$, for some $0\leq\nu$, such that $a\leq x\leq b$, moreover, we can choose the smallest such ordinal $\nu$. Since we cannot say that the \grqq segment\grqq\, $[0,1]$ is the totality of all or even the set of all Conway numbers  since in the first case there is no such totality and in the second case it is a proper class in $NBG$ theory. Formally, this can be expressed as an eternal process of an increasing sequence of embedded \grqq segments\grqq\, $[a,b]=[a,b]_{\nu_0}\subset[a,b]_{\nu_0+1},...,\subset[a,b]_{\nu_0+\alpha}\subset,...$\, where $\nu_0=\min\limits_{0\leq\nu}\{\nu\}$ such that $a,b\in{\bf R}_{\zeta}$, $\zeta=\omega^{\omega^\nu}$ and $0\leq\alpha$ or a completed process in $NBG$ when $0\leq\alpha<\Omega$. 

 The construction of \grqq dividing a segment in half\grqq\, does not give anything significant. But for simplicity, we omit the quotation marks further. Using this method, we are going to construct two sequences of left and right ends of the nested segments, which will tend to each other without an intermediate number. If we conditionally call the segment $[0,1]$ the left one, then after dividing it in half $[0,1]=[0,\frac{1}{2}]\cup[\frac{1}{2},1]$, we choose the right segment $[\frac{1}{2},1]$, which we divide in half $[\frac{1}{2},1]=[\frac{1}{2},\frac{3}{2^2}]\cup[\frac{3}{2^2},1]$, and choose the left segment $[\frac{1}{2},\frac{3}{2^2}]$, which we divide in half $[\frac{1}{2},\frac{3}{2^2}]=[\frac{1}{2},\frac{5}{2^3}]\cup[\frac{5}{2^3},\frac{3}{2^2}]$, similar choosing the right one $[\frac{5}{2^3},\frac{3}{2^2}]$ and dividing $[\frac{5}{2^3},\frac{3}{2^2}]=[\frac{5}{2^3},\frac{11}{2^4}]\cup[\frac{11}{2^4},\frac{3}{2^2}]$, further $[\frac{5}{2^3},\frac{11}{2^5}]=[\frac{5}{2^3},\frac{21}{2^5}]\cup[\frac{21}{2^5},\frac{11}{2^4}]$ and so on, i.e., \grqq $\omega$ times\grqq\, or more precisely for all $0<\alpha<\omega$. Thus we obtain
$$
[0,1]\supset[1/2,1]\supset[1/2,3/2^2]\supset[5/2^3,3/2^2]\supset[5/2^3,11/2^4]\supset[21/2^5,11/2^4]\supset...
$$
 Then we consider the Conway number $a_1=<x^L\,|\,x^R>$, where $\{x^L\}$ is the set of all left ends and $\{x^R\}$ is the set of all right ends of the above embedded segments, i.e., $a_1=\{0,1/2,5/2^3,21/2^5,...\,|\,1,3/2^2,11/2^4,11/2^6...\}$ and call $a_1$ a left end and of the following segment $[a_1,b_1]$, where $b_1=<x^L,a_1\,|\,x^R>$, i.e.,

 $a_1=\{0,1/2,5/2^3,21/2^5,...,\,a_1\,|\,1,3/2^2,11/2^4,43/2^6...\}$ and call it a left segment. One can see that $a_1=\frac{2}{3}$ but for general formu, we will leave the notation $a_1$.

Note that we can easily find recursive formulas for the left and right components of number $a_1=<x^L\,|\,x^R>$. Indeed, $x^L=\frac{y^L_\alpha}{2^{2\cdot\alpha+1}}$ and $x^R=\frac{y^R_\alpha}{2^{2\cdot\alpha}}$, $0<\alpha<\omega$,  $y^L_\alpha=2\cdot y^L_{\alpha-1}+y^R_{\alpha-1}$ and $y^R_\alpha= y^L_{\alpha-1}+2\cdot y^R_{\alpha-1}$, where $y^L_0=\frac{1}{3}=y^R_0$.

Now we divide $[a_1,b_1]$ in half and choose the right segment, which we divide in half and choose the left one, and so on, i.e., \grqq$\omega$ times\grqq\, again, and we define the Conway number $a_2=a_1+\frac{a_1}{\omega}=$
$$=\{x^L_1,...,x^L_\alpha,...,a_1+x^L_1/\omega,...,
a_1+x^L_\alpha/\omega,...\,|\,x^R_1,...,x^R_\alpha,...,a_1+x^R_1/\omega,...,
a_1+x^R_\alpha/\omega,...\}$$
 through the set of left ends of the segment and the right ends of the segments above, which make up an increasing $2\cdot\omega$-sequence and a decreasing $2\cdot\omega$-sequence, respectively, and similarly we define the Conway number $b_2=a_1+\frac{a_1}{\omega}+\frac{1}{\omega^2}$ and the segment $[a_2,b_2]=[a_1+\frac{a_1}{\omega},a_1+\frac{a_1}{\omega}+\frac{1}{\omega^2}]$, which is embedded in all above segments. And so on
$$
[0,1]\supset[1/2,1]\supset...\supset[a_1,a_1+\frac{1}{\omega}]\supset[a_1+\frac{a_1}{\omega},a_1+\frac{a_1}{\omega}+\frac{1}{\omega^2}]\supset...\supset[a_1+\frac{a_1}{\omega},a_1+\frac{a_1}{\omega}+\frac{1}{\omega^\alpha}]\supset...
$$
i.e., for all $0<\alpha<\Omega$ in $NBG$ theory and for ever, i.e., for any $\alpha>0$ in the case of trans-definite Conway numbers. 

One can see that we obtain an increasing sequence $(a_\alpha)_{0<\alpha}$ of left ends and a decreasing sequence $(b_\alpha)_{0<\alpha}$ of right ends such that for any positive Conway number $\varepsilon$ there is an ordinal number $\alpha_0$ such that $b_\alpha -a_\alpha<\varepsilon$ for all $\alpha>\alpha_0$ and, by left-right construction, there is no Conway numbers $c\in{\bf No}$ between them, i.e., $a_\alpha<c<b_\alpha$, $0<\alpha<\Omega$ because the gap of the first kind due to Conway (see \cite{l22}, p.38) has in our case the normal form:
$$
a_1+\frac{a_1}{\omega} +\frac{a_1}{\omega^2}+\frac{a_1}{\omega^\omega}+...+\frac{a_1}{\omega^\alpha}+...\,,
$$
summed over {\it all} ordinals $\alpha$.

Thus we can conclude that there is no even trans-definable numbers $c\notin{\bf No}$ between these increasing   and  decreasing sequences $(a_\alpha)_{0<\alpha}$ and $(b_\alpha)_{0<\alpha}$.

\bigskip
$$
$$

\begin{center}
{\bf 2. The difference between infinity and an infinite number  in the structure of the Conway numbers Field}
\end{center}

The concept of {\it  infinity} is well known to everyone from the school curriculum, and the symbol of infinity $\infty$ has firmly entered the universal vocabulary.

Recall that in the field of real numbers ${\bf R}$, a sequence $(x_n)_{0\leq n<\omega}$ is called {\it infinitely large} if for any positive number $E\in{\bf R}$ there exists a natural number $n_0$ such that $|x_n|>E$ for all $n>n_0$. 

At the same time, the so-called infinities with a sign are distinguished, i.e., $+\infty$ and $-\infty$, i.e.,  a sequence $(x_n)_{0\leq n<\omega}$ is called infinitely large with a sign $+$ if for any positive number $E\in{\bf R}$ there exists a natural number $n_0$ such that $x_n>E$ for all $n>n_0$ and a sequence $(x_n)_{0\leq n<\omega}$ is called infinitely large with a sign $-$ if for any positive number $E\in{\bf R}$ there exists a natural number $n_0$ such that $x_n<-E$ for all $n>n_0$ . (Let's mention also dual notions of sequences  $(x_n)_{0\leq n<\omega}$ which are called {\it infinitely small}, i.e., if for any positive number $\varepsilon\in{\bf R}$ there exists a natural number $n_0$ such that $|x_n|<\varepsilon$ for all $n>n_0$ and also with sign as for sequences $(\frac{1}{n})_{0<n<\omega}$ and $(-\frac{1}{n})_{0<n<\omega}$ when $\lim\limits_{n\rightarrow \infty}\frac{1}{n}=+0$ and $\lim\limits_{n\rightarrow \infty}-\frac{1}{n}=-0$.)

One writes these definitions like this $\lim\limits_{n\rightarrow\infty}x_n=\infty$, $\lim\limits_{n\rightarrow\infty}x_n=+\infty$ and $\lim\limits_{n\rightarrow\infty}x_n=-\infty$, e.g., $\lim\limits_{n\rightarrow\infty}(-1)^nn=\infty$, $\lim\limits_{n\rightarrow\infty}2^n=+\infty$ and $\lim\limits_{n\rightarrow\infty}(-1^{2n+1})n=-\infty$.

Symbols of infinity $+\infty$ and $-\infty$ became convenient abbreviations: $(-\infty<x<+\infty)$, $(-\infty<x<b)$ and $(a<x<+\infty)$ as infinite intervals of all real numbers $x$ between these symbols and numbers  and studying the behavior of functions at infinity (asymptotic behavior, classification or order of infinitely large quantities, etc.). However, there was no understanding of the fundamental difference between symbols of infinity and the infinite numbers that appeared in Cantor's set theory (ordinal and cardinal numbers). Nevertheless, there is a fundamental difference between them, there are heterogeneous concepts -- numbers and symbols of infinity.

Conway kept these symbols for his numbers, explaining the difference between a number and a symbol of infinity. A number is a section in sets of numbers, and symbols of infinity are Dedekind sections (gaps) in a proper class ${\bf No}$ of all Conway numbers, which are illegal in $NBG$ theory and therefore cannot be numbers, but only symbols:

 ${\bf On}$ for the gap $({\bf No},\emptyset)$ at the end of the number line (in our notation $+\infty_\Omega$),

 $\frac{1}{{\bf On}}$ for the gap between $0$ and all positive numbers, (in our notation $+\infty_{\frac{1}{\Omega}}$),

$\infty$ for the gap between reals and positive infinite numbers (in our notation 
$+\infty_\omega$), 

$\frac{1}{\infty}$ for that between infinitesimals and the positive reals (in our notation $+\infty{\frac{1}{\omega}}$. 

For these gaps, we have the following formulas-symbol: 

${\bf On}=\omega^{\bf On}$, $\frac{1}{{\bf On}}=\omega^{-{\bf On}}$, $\infty=
\omega^{\frac{1}{{\bf On}}}$, $\frac{1}{\infty}=\omega^{-\frac{1}{{\bf On}}}$ (in our notation more adequately $\Omega=\omega^\Omega$, $\frac{1}{\Omega}=\omega^{-\Omega}$, $\infty=+\infty_\omega$, $\frac{1}{\infty}=+\infty_{\frac{1}{\omega}}$, respectively). 

Moreover, \grqq these gaps definable as upper or lower bounds are particularly important in the theory of games\grqq\, according to Conway's theory (see \cite{l22}, p. 37-38).

In fact, the Conway's symbols above are by breed nature transdefinite numbers and they are no different from definitite numbers -- they are mathematical entities that can be added, multiplied, subtracted, and divided (not by zero, of course) like any Conway number, using identical formulas, what cannot be done with symbols of infinity $+\infty$ and $-\infty$, except  of course primitive operations  with them sometimes (in Mathematical Calculus or in Measure Theory), e.g., $\infty+\infty=\infty$, $x\cdot\infty=\infty$, for any positive number $x$, but $\infty-\infty$ and $0\cdot\infty$ are undefined.

We also recall from the course of mathematical analysis that every upper bounded subset $X$ of real numbers ${\bf R}$ has an exact upper bound $\sup\,X\in{\bf R}$, i.e., the smallest upper bound; and every lower bounded subset $X\subseteq {\bf R}$ has an exact lower bound $\inf\,X\in {\bf R}$, i.e., the greatest lower bound (a fact equivalent to Dedekind's Theorem). For a set $X$ that is unbounded from above, it is assumed $\sup\,X=+\infty$, and for a set $X$ that is unbounded from below, it is assumed $\sup\,X=-\infty$. In particular, for the empty subset $\emptyset$ of the set of real numbers, we obtain $\sup\,\emptyset=-\infty$ and $\inf\,\emptyset=+\infty$.

There is a completely different nature of subsets and proper subclasses $X\subseteq{\bf No}$ of the collection  ${\bf No}$ of all Conway numbers.  If a set  $X\subset{\bf No}$ does not have a maximum element, then there is no exact upper bound $\sup X$ in ${\bf No}$; and if a set  $X\subset{\bf No}$ does not have a minimum element in ${\bf No}$, then there is no exact lower bound $\inf\, X$ in ${\bf No}$. Conway posited, as we saw above $\sup{\bf No}=\{{\bf No}\,|\,\}=+\Omega$ as the end of the number line (or the last term $\Omega$ in our notation here) as the Dedekind section $({\bf No},\emptyset)$, which is a gap and an illegal object and called it {\it improper class} and speak of there being an {\it impropriety of gap} (see \cite{l22}, p. 38) (we can add naturally the Dedekind section $(\emptyset,{\bf No})$, which is a gap and also an illegal object and put $\{\,|\,{\bf No}\}=-\Omega$), but Conway did not define a natural logically possible particular degenerate example $\sup\emptyset$ and $\inf\emptyset$, respectively. The latter can also be defined in $NBG$ theory as $\sup{\emptyset}=-\infty_\Omega$ and $\inf{\emptyset}=+\infty_\Omega$, but for (say not quite correctly) \grqq all\grqq\, trans-definite Conway numbers, the empty set has  the exact upper bound $\sup\emptyset=\inf\limits_{0\leq\alpha}-\alpha=-\infty_0$ and the  exact lower bound $\inf\emptyset=\sup\limits_{0\leq\alpha}\alpha=+\infty_0$ and they are not numbers bu symbols (see below).

Indeed, if a set  $X\subset{\bf No}$ does not have a maximum element in ${\bf No}$, then consider a proper class $Y=\{y\,|\, y>x\,(\forall x)\, x\in X\}$, then $Y$ does not have a smallest element $y_0$. Otherwise, $y_1=\{X\,|\, y_0\}\in{\bf No}$ has the following property: $y_1>x \, (\forall x)\,x\in X$ and $y_1<y_0$ and hence $y_1\notin Y$. Contradiction. And if a set  $X\subset{\bf No}$ does not have a minimal element, then consider a proper class $Y=\{y\,|\, y<x\,(\forall x)\, x\in X\}$, then $Y$ does not have a smallest element $y_0$. Otherwise, $y_1=\{y\,|\,X\}\in{\bf No}$ has the following property: $y_1<x \, (\forall x)\,x\in X$ and $y_1>y_0$ and hence $y_1\notin Y$. Contradiction.

The existence of trans-definite Conway  numbers (e.g., formal systems with super-classes) allows us to define the exact upper and exact lower bounds of sets or proper classes of trans-definite numbers such that do not have maximum and minimum numbers, respectively.

Indeed, if there is any set  $X$ or proper class of Conway numbers, or any collection of trans-definite numbers  obtained by Conway's formulas and extending beyond the class of all Conway numbers such that it does not have maximum or minimum element, it has the exact upper bound and exact lower bound, respectively, but these bounds are not Conway numbers at all in either the $NBG$ theory or the theory that extends beyond $NBG$ but {\it symbols}, i.e., the symbols that separate $X$ and any numerical upper bounds and any numerical lower bounds, respectively. These symbols are like \grqq immaterial\grqq\, boundaries, if we assume that Conway numbers are ideal material points or \grqq material entities\grqq, figuratively we can compare them our earthly understanding of boundaries on the sea or ocean, in contrast to the visible boundaries on the earth. Or, in another way, if we consider the geometric interpretation of Conway numbers and trans-definite Conway  numbers as points on a line, then the symbols of infinity are not points at all, but rather conditional separators of points on this line. These symbols can also be considered as punctuation marks that are part of the syntax of formal systems.

Below, we will try to depict these symbols  with signs and indexes (labels) as the earliest (their birthdays) or the simplest Conway numbers $a=\{X\,|\,Y\}$ and $a'=\{Y\,|\,X\}$, which are determined by $X$ for some $Y$, convenient for a unique index (label) of this symbol of infinity, defined by the set $X$, i.e., $\sup X=\pm\infty_a$ and  $\inf X=\pm\infty_{a'}$, when $a,a'>0$ and $a,a'<0$, respectively. In the case of an empty set $X=\emptyset$, the indexes are defined as follows:  $\sup\emptyset=-\infty_0$ and $\sup\emptyset=+\infty_0$, respectively, because the simplest Conway number, defined by $X=\emptyset$ is when $Y=\emptyset$, i.e., $a=\{X\,|\,Y\}=\{\,|\,\}=0$. In other form  $\sup\emptyset=\inf\limits_{0\leq\alpha}-\alpha=-\infty_0$ and $\inf\emptyset=\sup\limits_{0\leq\alpha}\alpha=+\infty_0$. These are different from above case of the real numbers ${\bf R}$, for which $\sup\emptyset=-\infty_\omega$ and $\sup\emptyset=+\infty_\omega$, respectively.

Generally speaking, every strictly increasing or strictly decreasing $\lambda$-sequence $(x_\alpha)_{0\leq\alpha<\lambda}$ of Conway numbers with a limit ordinal $\lambda$ defines the symbol of infinity $\infty_c\stackrel{def}{=}+\infty_c$ if $c$ is a positive number and $\infty_c\stackrel{def}{=}-\infty_{c}$ if $c$ is a negative number, respectively, where $c$ is uniquely defined by the corresponding sequences $(x_\alpha)_{0\leq\alpha<\lambda}$, i.e., $c=\{members\,\, x_\alpha,0\leq \alpha<\lambda\,|\,Y\}$ and $c=\{Y'\,|\,members\,\, x_\alpha,0\leq \alpha<\lambda\}$ for some $Y$ and $Y'$, respectively. E.g.  $(\pm\alpha)_{0<\alpha<\omega}$ and $(\pm\frac{1}{\alpha})_{0<\alpha<\omega}$ considered above, then $c=\pm\omega$ and $c=\pm\frac{1}{\omega}$, respectively, since $+\omega=\{members\,\, \alpha,0\leq \alpha<\lambda\,|\,\}$, $-\omega=\{\,|\,members\,\, -\alpha,0\leq \alpha<\lambda\}$ and $+\frac{1}{\omega}=\{0\,|\,members\,\, \frac{1}{\alpha},0\leq \alpha<\lambda\,|\,\}$, $-\frac{1}{\omega}=\{members\,\, -\frac{1}{\alpha},0\leq \alpha<\lambda\,|\,0\}$, respectively, and $\pm\omega$, $\pm\frac{1}{\omega}$ are Conway numbers born first on day $\omega$ and these numbers are defined by corresponding sequences and sets $Y$ of Conway numbers, i.e., $Y=\emptyset$ and $Y'=\{0\}$, respectively.

Now we will use canonical examples to demonstrate a more precise place and meaning of these symbols of infinity. Let's start with the simplest classical example, namely, a sequence of increasing natural numbers $(n)_{0<n<\omega}$ of the field of real numbers ${\bf R}_\xi$, $\xi=\omega^{\omega^\nu}$, in the case $\nu=0$ and clear $\xi=\omega$. 

We denote $+\infty_\omega=\sup\limits_{0<\alpha<\omega}\{1,2,3,...,\alpha,...\}=\inf\limits_{0<\beta}\{\omega,\omega^{\frac{1}{2}},\omega^{\frac{1}{3}},...,\omega^{\frac{1}{\beta}},...\}$ in the following sense: for every Conway number $a$ such that $a>\alpha$, for all ordinals $0<\alpha<\omega$, there is an ordinal $\beta_0$ such that $a>\omega^{\frac{1}{\beta}}$ for every ordinal $\beta>\beta_0$, and for other Conway number $b$ there is an ordinal $\alpha_0$ such that $b<\alpha_0<\omega$. 

	The second simplest classical example, namely, a sequence of decreasing fractions $(\frac{1}{\alpha})_{0<\alpha<\omega}$ of the field of real numbers ${\bf R}_\omega$.

 We denote $+\infty_{\frac{1}{\omega}}=\inf\limits_{0<\alpha<\omega}\{1,\frac{1}{2},\frac{1}{3},...,\frac{1}{\alpha},...\}=\sup\limits_{0<\beta}\{\frac{1}{\omega},\frac{1}{\omega^{\frac{1}{2}}},\frac{1}{\omega^{\frac{1}{3}}},...,\frac{1}{\omega^{\frac{1}{\beta}}},...\}$ in the following sense: for every Conway number $\varepsilon>0$ such that $\varepsilon<\frac{1}{\alpha}$, for all ordinals $0<\alpha<\omega$, there is an ordinal $\beta_0$ such that $\varepsilon<\frac{1}{\omega^{\frac{1}{\beta}}}$ for every ordinal $\beta>\beta_0$, and for other Conway number $b$ there is an ordinal $\alpha_1$ such that $b>\frac{1}{{\alpha_1}}$. 

Formally, we expand these examples by \grqq multiplying them by minus one\grqq\, swapping places $\sup$ and $\inf$, i.e., $$-\infty_{-\omega}=\inf\limits_{0<\alpha<\omega}\{-1-,2,-3,...,-\alpha,...\}=\sup\limits_{0<\beta}\{-\omega,-\omega^{\frac{1}{2}},-\omega^{\frac{1}{3}},...,-\omega^{\frac{1}{\beta}},...\}$$
$$-\infty_{-{\frac{1}{\omega}}}=\sup\limits_{0<\alpha<\omega}\{-1,-\frac{1}{2},-\frac{1}{3},...,-\frac{1}{\alpha},...\}=\inf\limits_{0<\beta}\{-\frac{1}{\omega},-\frac{1}{\omega^{\frac{1}{2}}},-\frac{1}{\omega^{\frac{1}{3}}},...,-\frac{1}{\omega^{\frac{1}{\beta}}},...\}$$

There may be other designations  symbols of infinity such as $+\infty_\omega=\infty_\omega$, $-\infty_{\omega}=\infty_{-\omega}$, $+\infty_{\frac{1}{\omega}}=\infty_{\frac{1}{\omega}}$ and $-\infty_{\frac{1}{\omega}}=\infty_{-\frac{1}{\omega}}$ or like in Conway's paper $+\infty_\omega=\infty_\omega$, $-\infty_{\omega}=\infty_{-\omega}$, $+\infty_{\frac{1}{\omega}}=\frac{1}{\infty_\omega}$ and $-\infty_{\frac{1}{\omega}}=-\frac{1}{\infty_\omega}$.

\bigskip 

Although all strictly increasing and strictly decreasing sequences of Conway numbers cannot  be indexed by all Conway numbers (in $NBG$ formal theory), as they are not even a proper class but, by the well-known Cantor theorem have a greater power (note that in the imaginary  formal theory here it is possible, since the sequences are either sets or proper classes but not endless sequences), we will conventionally indicate the initial numbers (say numbers born on or before day $\omega_1$, i.e., the first non-countable ordinal number) in order to understand the essence of both the symbols of infinity and the objective nature of Conway numbers, particularly the ordinal numbers known to everyone after Cantor (see the following paragraph) that will shed light on the \grqq geography\grqq\, of all Conway numbers or, more precisely, the Conway's line.

 We  start with the simplest case when  strictly increasing or strictly decreasing $\omega$-sequences $(x_\alpha)_{0\leq\alpha<\omega}$  are sequences of positive dyadic rational numbers of the form $\frac{m}{2^n}$, $m$ and $n$ are integers ($n\geq 0$), i.e., by Conway's construction, members of these sequences  are $x_\alpha\in O_\omega$, for all $0\leq\alpha<\omega$, i.e., born before day $\omega$.

The four strictly increasing/decreasing and strictly decreasing/increasing $\omega$-sequences $(\pm\alpha)_{0<\alpha<\omega}$ and $(\pm\frac{1}{\alpha})_{0<\alpha<\omega}$  here are primordial and have already been designated by us above. The indexes of symbols of infinity $\pm\infty_{\pm\omega}$ and $\pm\infty_{\pm\frac{1}{\omega}}$ are uniquely determined by the numbers $\pm\omega$ and $\pm\frac{1}{\omega}$ born on  day $\omega$, respectively.

For the convenience of distinguishing these and other $\omega$-sequences, we will introduce the following terminology: $(\pm\alpha)_{0<\alpha<\omega}$ and $(\pm\frac{1}{\alpha})_{0<\alpha<\omega}$ are $\omega$-sequences of  the $\omega$-type and $\omega$-sequences of $\frac{1}{\omega}$-type, respectively. 

Let's first consider the set $O_\omega$ of all Conway numbers born before day $\omega$. There are only two (modulo of cofinal and coinitial $\omega$-sequences, respectively) $\omega$-sequences of the $\omega$-type, i.e., $(\alpha)_{0<\alpha<\omega}$ and $(-\alpha)_{0<\alpha<\omega}$. We see that $\{1,2,...,\alpha,...\,|\,\}=\omega$ and $\{\,|\,-1-,2,...,-\alpha,...\}=-\omega$ and thus symbols of infinity for $(\alpha)_{0<\alpha<\omega}$ and $(-\alpha)_{0<\alpha<\omega}$ are $+\infty_\omega$ and $-\infty_{-\omega}$, respectively.

 Each other   strictly increasing  $\omega$-sequence $(x_\alpha)_{0\leq\alpha<\omega}$ such that $x_\alpha\in O_\omega$ and non-cofinal  to $(\alpha)_{0<\alpha<\omega}$ defines a single number $a\in{\bf R}_\omega$ by the following rule. Consider a set $L$ of all numbers $a^L\in O_\omega$ such that there is an ordinal $0<\alpha_0<\omega$ such that $a^L<x_{\alpha_0}$ and a set $R$ of all numbers $a^R\in O_\omega$ such that $a^R>x_\alpha$ for all $0<\alpha<\omega$. If $R$ has the least number $a\in R$, then we put  $c=a-\frac{1}{\omega}$. If $R$ has no least number, then we consider $a=\{a^L\,|\,a^R\}$ and similar put $c=a-\frac{1}{\omega}$. Thus we conditionally (formally) define for  a strictly increasing  $\omega$-sequence $(x_\alpha)_{0\leq\alpha<\omega}$
the symbol of infinity $\pm\infty_{a-\frac{1}{\omega}}=\sup\limits_{0<\alpha<\omega}x_\alpha=\inf\{y\,|\,y>x_\alpha \forall\alpha\in(0,\omega)\}$ when $a>/<0$, respectively.

Each other   strictly decreasing  $\omega$-sequence $(x'_\alpha)_{0\leq\alpha<\omega}$ such that $x'_\alpha\in O_\omega$ and non-coinitial  to $(-\alpha)_{0<\alpha<\omega}$ defines a single number $a\in{\bf R}_\omega$ by the following rule. Consider a set $R$ of all numbers $a^R\in O_\omega$ such that there is an ordinal $0<\alpha_0<\omega$ such that $a^R>x'_{\alpha_0}$ and a set $L$ of all numbers $a^L\in O_\omega$ such that $a^L<x'_\alpha$ for all $0<\alpha<\omega$. If $L$ has the greatest number $a\in R$, then we put  $c=a+\frac{1}{\omega}$. If $L$ has no greatest number, then we consider $a=\{a^L\,|\,a^R\}$ and similar put $c=a+\frac{1}{\omega}$. Thus we conditionally (formally) define for  a strictly decreasing  $\omega$-sequence $(x'_\alpha)_{0\leq\alpha<\omega}$
the symbol of infinity $\pm\infty_{a+\frac{1}{\omega}}=\inf\limits_{0<\alpha<\omega}x'_\alpha=\sup\{y'\,|\,y'<x'_\alpha \forall\alpha\in(0,\omega)\}$  when $a>/<0$, respectively.

 {\bf Remark 11.} Our formalism is that  collections $Y=\{y\,|\,y>x_\alpha \forall\alpha\in(0,\omega)\}$ $Y'=\{y'\,|\,y'<x_\alpha \forall\in(0,\omega)\}$ are not mathematical objects (they do not exist), but these formulas formally explain that there is no smallest number $\hat y$ such that $\hat y>x_\alpha$ for all $0<\alpha<\omega$ and there is no greatest number $\hat y'$ such that $\hat y'<x'_\alpha$ for all $0<\alpha<\omega$, respectively, and thus only symbols of infinity $\pm\infty_{a-\frac{1}{\omega}}$ and $\pm\infty_{a+\frac{1}{\omega}}$ have the corresponding properties.

Note only that number $a\mp\frac{1}{\omega}$ born on day $\omega$ when $a\in O_\omega$  and born on day $\omega+1$ when $a\notin O_\omega$.

 {\bf Remark 12.} Notice also that each strictly increasing $\omega$-sequence $(x_\alpha)_{0<\alpha<\omega}$, $x_\alpha\in {\bf R}_\omega$, of $\frac{1}{\omega}$-type is cofinal to   $(a-\frac{1}{\alpha})_{0<\alpha<\omega}$ for some $a\in{\bf R}_\omega$ and hence  $\pm\infty_{a-\frac{1}{\omega}}=\sup\limits_{0\leq\alpha<\omega}x_\alpha=
\sup\limits_{0<\alpha<\omega}(a-\frac{1}{\alpha})$ for $a>0$ and $a<0$, respectively, and for $a=0$ we obtain $-\infty_{-\frac{1}{\omega}}=\sup\limits_{0<\alpha<\omega}x_\alpha=
\sup\limits_{0<\alpha<\omega}(-\frac{1}{\alpha})_{0<\alpha<\omega}$, and a single number $-\frac{1}{\omega}$ is defined by $(x_\alpha)_{0<\alpha<\omega}$. Similar each strictly decreasing $\omega$-sequence $(x'_\alpha)_{0<\alpha<\omega}$, $x_\alpha\in {\bf R}_\omega$, of $\frac{1}{\omega}$-type is coinitial to   $(a+\frac{1}{\alpha})_{0<\alpha<\omega}$ for some $a\in{\bf R}_\omega$ and hence  $\pm\infty_{a+\frac{1}{\omega}}=\inf\limits_{0\leq\alpha<\omega}x'_\alpha=
\inf\limits_{0<\alpha<\omega}(a+\frac{1}{\alpha})$ for $a>0$ and $a<0$, respectively, and for $0$ we obtain $+\infty_{\frac{1}{\omega}}=\inf\limits_{0<\alpha<\omega}x'_\alpha=
\inf\limits_{0<\alpha<\omega}(\frac{1}{\alpha})_{0<\alpha<\omega}$, and a single number $\frac{1}{\omega}$ is defined by $(x'_\alpha)_{0<\alpha<\omega}$.

These $\omega$-sequences of $\frac{1}{\omega}$-type can be fully described as follows. (By \grqq fully described\grqq\, we mean \grqq for all $\omega$-sequences of $\frac{1}{\omega}$-type in ${\bf No}$ of $NBG$ theory and for each $\omega$-sequence of $\frac{1}{\omega}$-type beyond this theory\grqq\,.)
Let
$b=\sum\limits_{0\leq\beta<\gamma}\omega^{y_\beta}r_\beta$ be a normal form of Conway number $b$ $(\ref{f0201})$ in which $\gamma$ denotes some ordinal, the numbers $r_\beta$ ($\beta<\gamma$) are non-zero reals, and the numbers $y_\beta$ form a strictly decreasing sequence of numbers such that $y_\beta>-1$, for all $0\leq\beta<\gamma$, and add one more number $b=0$. Then each possible strictly increasing and strictly decreasing $\omega$-sequences of $\frac{1}{\omega}$-type (of course up to cofinality and coinitiality, respectively) are of the following form: $(b-\frac{1}{2^\alpha})_{0\leq\alpha<\omega}$ and $(b+\frac{1}{2^\alpha})_{0\leq\alpha<\omega}$, respectively, and their  symbols of infinity are $\pm\infty_{b-\frac{1}{\omega}}$ when $b>0$ or $b<0$ as well as $\mp\infty_{b+\frac{1}{\omega}}$ when $b>0$ or $b<0$, respectively.

\bigskip

New $\omega$-sequences of other type than $\frac{1}{\omega}$-type arise for Conway numbers whose birthday is less than or equal to $2\cdot\omega$  and this case immediately becomes more complex and diverse, which further complicates the situation when $\lambda>2\cdot\omega$ with limit ordinal numbers $\lambda$. Therefore, we will briefly discuss only this case $\lambda=2\cdot\omega$.

Let's now consider the set $O_{2\cdot\omega}$ of all Conway numbers born before day $2\cdot\omega$. There are only four (modulo  cofinal and coinitial $\omega$-sequences, respectively) new strictly increasing and strictly decreasing $\omega$-sequences of the $\omega$-type, i.e., $(\pm\omega+\alpha)_{0<\alpha<\omega}$ and $(\mp\omega-\alpha)_{0<\alpha<\omega}$, respectively. 

We see that they are $\{\omega+1,\omega+2,...,\omega+\alpha,...\,|\,\}=2\cdot\omega$, $\{\,-\omega+1,-\omega+2,...,-\omega+\alpha,...\,|\,-1-,2,...,-\alpha,...\}=-\frac{\omega}{2}$ and $\{\,|\,-\omega-1,-\omega-2,...,-\omega-\alpha,...\}=-2\cdot\omega$, $\{1,2,...,\alpha,...\,|\,\omega-1,\omega-2,...,\omega-\alpha,...\}=\frac{\omega}{2}$,  respectively. Thus symbols of infinity for $(\pm\omega+\alpha)_{0<\alpha<\omega}$ and $(\mp\omega-\alpha)_{0<\alpha<\omega}$ are $+\infty_{2\cdot\omega}$, $-\infty_{-\frac{\omega}{2}}$ and $-\infty_{-2\cdot\omega}$, $+\infty_{\frac{\omega}{2}}$, respectively. Remember that we have already two old symbols of infinity $+\infty_\omega$ and $-\infty_{-\omega}$ of old $\omega$-sequences of the $\omega$-type $(\alpha)_{0<\alpha<\omega}$ and $(-\alpha)_{0<\alpha<\omega}$, respectively. 

As to other $\omega$-sequences in $O_{2\cdot\omega}$ for the sake of brevity, we will now describe only sequences of positive numbers in $O_{2\cdot\omega}$. Two primordial $\omega$-sequences of $\frac{1}{\omega^2}$-type and $\frac{1}{\omega^{\frac{1}{2}}}$-type are the following:
$$
\frac{1}{\omega}, \frac{1}{2\cdot\omega}=\{0\,|\,\frac{1}{\omega}\}, \frac{1}{2^2\cdot\omega}=\{0\,|\,\frac{1}{2\cdot\omega}\},..., \frac{1}{2^\alpha\cdot\omega}=\{0\,|\,\frac{1}{2^{\alpha+1}\cdot\omega}\},...,
$$
and
$$
\frac{1}{\omega}, \frac{2}{\omega}=\{\frac{1}{\omega}\,|\,1,\frac{1}{2}, \frac{1}{2^2},..., \frac{1}{2^{\alpha-1}},...\},..., \frac{2^\alpha}{\omega}=\{\frac{1}{\omega}\,|\,1,\frac{1}{2}, \frac{1}{2^2},..., \frac{1}{2^{\alpha-1}},...\},...\,\,,
$$
where $0<\alpha<\omega$, and with the following symbols of infinity $+\infty_{\frac{1}{\omega^2}}$ and $+\infty_{\frac{1}{\omega^{\frac{1}{2}}}}$,  respectively. It is really so because the first numbers born in day $2\cdot\omega$ are $\frac{1}{\omega^2}=\{0\,|\,\frac{1}{\omega}, \frac{1}{2\cdot\omega}, \frac{1}{2^2\cdot\omega},..., \frac{1}{2^\alpha\cdot\omega},...,
\}$ and $\frac{1}{\omega^{\frac{1}{2}}}=\{\frac{1}{\omega}, \frac{2}{\omega},... \frac{2^\alpha}{\omega},...\,|\,1,\frac{1}{2}, \frac{1}{2^2},..., \frac{1}{2^\alpha},...\}$, respectively, where $0\leq\alpha<\omega$.

Thus all strictly increasing and strictly decreasing $\omega$-sequences in $B=M_{2\cdot\omega}\setminus(\{a\pm\frac{1}{\omega^2}\}_{a\in O_{2\cdot\omega}}\cup\{2\cdot\omega,-2\cdot\omega\})$ of the $\frac{1}{\omega^2}$-type are defined by algebraic structure on Conway numbers and have the following form:  $(b\mp\frac{1}{2^{\alpha}\cdot\omega})_{0<\alpha<\omega}$, where $b\in B$, and their symbols of infinity are $\pm\infty_{b-\frac{1}{\omega^2}}$ for $b>0$ or $b<0$ and $\pm\infty_{b+\frac{1}{\omega^2}}$ for $b>0$ or $b<0$, respectively.

More general if $b=\sum\limits_{0\leq\beta<\gamma}\omega^{y_\beta}r_\beta$ be a normal form of Conway number $b$ in which $\gamma$ denotes some ordinal, the numbers $r_\beta$ ($\beta<\gamma$) are non-zero reals, and the numbers $y_\beta$ form a strictly decreasing sequence of numbers such that $y_\beta>-2$, for all $0\leq\beta<\gamma$. We add also number $b=0$. Then each possible strictly increasing and strictly decreasing $\omega$-sequences of $\frac{1}{\omega^2}$-type (of course up to cofinality and coinitiality, respectively) are of the following form:  $(b\mp\frac{1}{2^{\alpha}\cdot\omega})_{0<\alpha<\omega}$ and their symbols of infinity are $\pm\infty_{b-\frac{1}{\omega^2}}$ for $b>0$ or $b<0$ and $\pm\infty_{b+\frac{1}{\omega^2}}$ for $b>0$ or $b<0$, respectively.

As to the increasing  and strictly decreasing $\omega$-sequences $(\frac{2^\alpha}{\omega})_{0\leq\alpha<\omega}$ and $(-\frac{2^\alpha}{\omega})_{0\leq\alpha<\omega}$ of $\frac{1}{\omega^{\frac{1}{2}}}$-type
 in $B=M_{2\cdot\omega}\setminus(\{a\pm c\}_{a\in O_{2\cdot\omega},c\in(-\infty_{\frac{-1}{\omega^{\frac{1}{2}}}},+\infty_{\frac{1}{\omega^{\frac{1}{2}}}})}\cup\{2\cdot\omega,-2\cdot\omega\})$ are also defined by algebraic structure on Conway numbers and have the following form:  $(b\pm\frac{2^{\alpha}}{\omega})_{0<\alpha<\omega}$, where $b\in B$, and their symbols of infinity are $\pm\infty_{b+\frac{1}{\omega^{\frac{1}{2}}}}$ for $b>0$ or $b<0$  and $\pm\infty_{b-\frac{1}{\omega^{\frac{1}{2}}}}$ for $b>0$ or $b<0$, respectively.

More general if $b=\sum\limits_{0\leq\beta<\gamma}\omega^{y_\beta}r_\beta$ be a normal form of Conway number $b$ in which $\gamma$ denotes some ordinal, the numbers $r_\beta$ ($\beta<\gamma$) are non-zero reals, and the numbers $y_\beta$ form a strictly decreasing sequence of numbers such that $y_\beta\geq-\frac{1}{2}$, for all $0\leq\beta<\gamma$. And we add to them the number $b=0$. Then each possible strictly increasing and strictly decreasing $\omega$-sequences of $\frac{1}{\omega^{\frac{1}{2}}}$-type (of course up to cofinality and coinitiality, respectively) are of the following form:  $(b\pm\frac{2^{\alpha}}{\omega})_{0<\alpha<\omega}$ and their symbols of infinity are $\pm\infty_{b+\frac{1}{\omega^2}}$ for $b>0$ or $b<0$ and   $\pm\infty_{b-\frac{1}{\omega^{\frac{1}{2}}}}$ for $b>0$ or $b<0$, respectively.

In order to complete the description of all strictly increasing and strictly decreasing $\omega$-sequences of $\frac{1}{\omega}$-type in the set $O_{2\cdot\omega}$, we must consider the sequences $(b\mp\frac{1}{2^\alpha})_{0\leq\alpha<\omega}$ for all $$b\in M_{2\cdot\omega}\cap((-\infty_{-2\cdot\omega},-\infty_{-\frac{\omega}{2}})\cup(+\infty_{\frac{\omega}{2}},+\infty_{2\cdot\omega})),$$
with symbols of infinity $\pm\infty_{b-\frac{1}{\omega}}$ when $b>0$ or $b<0$ as well as $\pm\infty_{b+\frac{1}{\omega}}$ when $b>0$ or $b<0$, respectively,
which we have considered a;ready above.

{\bf Remark 13.} We will not continue to consider $\omega$-sequences as well as other $\lambda$-sequences which only begin to occur essentially (modulo confinality and initiality) from the first uncountable limit ordinal $\lambda=\omega_1$ of different types for numbers born after day $2\cdot\omega$ because each $\lambda$-sequence, where $\lambda<\omega_1$ is a limit ordinal, has an $\omega$-subsequence which is cofinal or coinitial to it.  It doesn't make much sense although the idea of how this is done is clear. It is important to see from the previous one interesting total laws of infinite numbers, symbols of infinity and intervals like $(+\infty_\omega,\omega]$, $(+\infty_{2\cdot\omega},2\cdot\omega]$, $(+\infty_{\frac{\omega}{2}},\omega]$,  and $[\frac{1}{\omega},+\infty_{\frac{1}{\omega}})$, $[\frac{1}{\omega^2},+\infty_{\frac{1}{\omega^2}})$, $(+\infty_{\frac{1}{\omega^{\frac{1}{2}}}},\frac{1}{{\omega^{\frac{1}{2}}}}]$, etc. which are not homogeneous, but with constant and increasing cycles and periods. 
\bigskip

So we can specify all  $\omega$-sequences of $\alpha\cdot\omega$- and $\frac{1}{\alpha\omega}$- and $\frac{1}{\omega^\alpha}$-types, where $0<\alpha<\omega$, and the  symbols of infinity they define in the four regions $(+\infty_{\omega^{\frac{1}{2}}},+\infty_{\omega^2})$, $(-\infty_{\omega^2},-\infty_{-\omega^{\frac{1}{2}}})$ and $(+\infty_{\frac{1}{\omega^2}},+\infty_{\frac{1}{\omega^{\frac{1}{2}}}})$, $(-\infty_{-\frac{1}{\omega^{\frac{1}{2}}}},-\infty_{-\frac{1}{\omega^2}})$ of Conway numbers, which are commensurate with $\omega$, $-\omega$, $\frac{1}{\omega}$ and $-\frac{1}{\omega}$, respectively. That means that for every real number $r\in{\bf R}$ one can uniquely define symbols of infinty for $\omega$-sequences $(r\cdot\omega\pm\alpha)_{0<\alpha<\omega}$ and $(r\cdot\frac{1}{\omega}\pm\frac{1}{\alpha\cdot\omega})_{0<\alpha<\omega}$.

Indeed, if we limit ourselves to the interval $(+\infty_{\omega^{\frac{1}{2}}},+\infty_{\omega^2})$ for simplicity, then the increasing $\omega$-sequence $(\omega+\alpha)_{0<\alpha<\omega}$ of $\omega$-type uniquely generates the number $2\cdot\omega=\{\omega,\omega+1,\omega+2,...,\omega+\alpha,...\,|\,\}$, where $0\leq\alpha<\omega$, and its symbol of  infinity is $+\infty_{2\cdot\omega}$ because for $X=\{\omega,\omega+1,\omega+2,...,\omega+\alpha,...\}$, $0\leq\alpha<\omega$, and $Y=\emptyset$ and $2\cdot\omega$ is the simplest number defined by $X$ and $Y$, i.e., $2\cdot\omega=\{X\,|\,Y\}$.

In a similar way the decreasing sequence $(\omega-\alpha)_{0<\alpha<\omega}$  of $\omega$-type uniquely generates the number $\frac{\omega}{2}=\{0,1,2,...,\alpha,...\,|\,\omega,\omega-1,...,\omega-\alpha,...\}$, where $0\leq\alpha<\omega$, and its   symbol of infinity is $+\infty_{\frac{\omega}{2}}$ because for $X=\{0,1,2,...,\alpha,...\}$, $0\leq\alpha<\omega$, $Y=\{\omega,\omega-1,...,\omega-\alpha,...\}$, $0\leq\alpha<\omega$, and $\frac{\omega}{2}$ is the simplest number defined by $X$ and $Y$, i.e., $\frac{\omega}{2}=\{X\,|\,Y\}$.

Moreover, by induction on ordinal numbers $0<\beta<\omega$, the increasing $\omega$-sequence $(\beta\cdot\omega+\alpha)_{0\leq\alpha<\omega}$ of $\omega$-type uniquely generates the number $(\beta+1)\cdot\omega=\{\beta\cdot\omega,\beta\cdot\omega+1,\beta\cdot\omega+2,...,\beta\cdot\omega+\alpha,...\,|\,\}$, where $0\leq\alpha<\omega$, and its symbol of infinity is $+\infty_{(\beta+1)\cdot\omega}$, $0\leq\beta<\omega$ because for $X=\{\beta\cdot\omega,\beta\cdot\omega+1,\beta\cdot\omega+2,...,\beta\cdot\omega+\alpha,...\}$, $0\leq\alpha<\omega$,  $Y=\emptyset$ and $(\beta+1)\cdot\omega$ is the simplest number defined by $X$ and $Y$, i.e., $(\beta+1)\cdot\omega=\{X\,|\,Y\}$.

In a similar way, by induction on ordinal numbers $0<\beta<\omega$, the decreasing $\omega$-sequence $(\frac{\omega}{2^\beta}-\alpha)_{0\leq\alpha<\omega}$ of $\omega$-type uniquely generates the number
 
$\frac{\omega}{2^{\beta+1}}=\{0,1,2,...,\alpha,...\,|\,\frac{\omega}{2^\beta},\frac{\omega}{2^\beta}-1,\frac{\omega}{2^\beta}-2,...,\frac{\omega}{2^\beta}-\alpha,...\}$, 
where $0\leq\alpha<\omega$, and its  symbol of infinity is $+\infty_{\frac{\omega}{2^{\beta+1}}}$, $0\leq\beta<\omega$, because for $X=\{0,1,2,...,\alpha,...\}$, $0\leq\alpha<\omega$,  $Y=\{\frac{\omega}{2^\beta},\frac{\omega}{2^\beta}-1,\frac{\omega}{2^\beta}-2,...,\frac{\omega}{2^\beta}-\alpha,...\}$, $0\leq\alpha<\omega$, and $\frac{\omega}{2^{\beta+1}}$ is the simplest number define by $X$ and $Y$, i.e., $\frac{\omega}{2^{\beta+1}}=\{X\,|\,Y\}$.

Thus we obtain new $\omega$-sequences $((\beta +1)\cdot\omega)_{0\leq\beta<\omega}$ and  $(\frac{\omega}{2^\beta})_{0\leq\beta<\omega}$ of $\omega^2$-type and $\omega^{\frac{1}{2}}$-type and their symbols of infinity are $+\infty_{\omega^2}$ and
 $+\infty_{\omega^{\frac{1}{2}}}$, respectively.

Now for each fixed ordinals $0\leq\beta<\omega$ and $0\leq\gamma<\omega$  the increasing $\omega$-sequence $(\beta\cdot\omega\pm\frac{\omega}{2^\gamma}+\alpha)_{0\leq\alpha<\omega}$ uniquely generates the number 
$\beta\cdot\omega\pm\frac{\omega}{2{\gamma-1}}=\{\beta\cdot\omega\pm\frac{\omega}{2^{\gamma-1}},\beta\cdot\omega\pm\frac{\omega}{2^{\gamma-1}}+1,\beta\cdot\omega\pm\frac{\omega}{2^{\gamma-1}}+2,...,\beta\cdot\omega\pm\frac{\omega}{2^{\gamma-1}}+\alpha,...\,|\,\beta\cdot\omega\pm\frac{\omega}{2^{\gamma-1}}-1,\beta\cdot\omega\pm\frac{\omega}{2^{\gamma-1}}-2,...,\beta\cdot\omega\pm\frac{\omega}{2^{\gamma-1}}-\alpha,...\}$, 
where $0\leq\alpha<\omega$, and the infinite symbol 
$+\infty_{(\beta-1)\cdot\omega+\frac{\omega}{2}}$.

$(\frac{3\cdot\omega}{2^2}+\alpha)_{0\leq\alpha<\omega}$ uniquely generates the number 
$\frac{7\cdot\omega}{2^3}=\{\frac{3\cdot\omega}{2^2},\frac{3\cdot\omega}{2^2}+1,\frac{3\cdot\omega}{2^2}+2,...,\frac{3\cdot\omega}{2^2}+\alpha,...\,|\,\omega,\omega-1,\omega-2,...,\omega-\alpha,...\}$, 
where $0\leq\alpha<\omega$, and the infinite symbol 
$+\infty_{\frac{7\cdot\omega}{2^3}}$. And so on, i.e., 
$+\infty_{\frac{15\cdot\omega}{2^4}}$,..., 
$+\infty_{\frac{(2^\gamma-1)\cdot\omega}{2^\gamma}}$,..., ($0\leq\gamma<\omega$), and the next is
$+\infty_{(1-\frac{1}{\omega})\cdot\omega}$,....

{\bf Remark 14.} It is clear that we can stop here, understanding how we can define the indexes of  symbols of infinity for numbers born on or before $\omega^2$ day. Here, we only note that there is some asymmetry in the definition of indexes and  symbols of infinity, which is already noticeable for the $\omega$-sequences of $\omega^2$-type and $\frac{1}{\omega^\alpha}$-type, $0<\alpha\leq\omega$, of numbers born before on and before day $\omega^2$. This means that there are no ordinals $\omega^\alpha$ yet, $0\leq\alpha\leq\omega$, but there are the inverse numbers $\frac{1}{\omega^\alpha}$, $0\leq\alpha\leq\omega$, and the corresponding $\omega$-sequences of $\frac{1}{\omega^\alpha}$-types , $0\leq\alpha\leq\omega$. However, on day $\omega^\omega$ we will have $\omega$-sequences of $\frac{1}{\omega^\alpha}$-types, $0\leq\alpha\leq\omega^{\omega^\omega}$, which will of course catch up with the lagging numbers and sequences on day $\varepsilon_0$. This can be visualized by the sequence of intervals (only for positive numbers here):
$$
(+\infty_{\frac{1}{\omega}},\omega), (+\infty_{\frac{1}{\omega^2}},2\cdot\omega),...,(+\infty_{\frac{1}{\omega^\alpha}},\alpha\cdot\omega),...,(+\infty_{\frac{1}{\omega^\omega}},\omega^2),(+\infty_{\frac{1}{\omega^{\omega+1}}},\omega^2+\omega), 
$$
$$
(+\infty_{\frac{1}{\omega^{\omega+2}}},\omega^2+2\cdot\omega),...,
(+\infty_{\frac{1}{\omega^{\omega+\alpha}}},\omega^2+\alpha\cdot\omega),...,
(+\infty_{\frac{1}{{\omega^{\omega^2}}}},\omega^3),(+\infty_{\frac{1}{{\omega^{\omega^3}}}},\omega^4),...,
$$
$$
(+\infty_{\frac{1}{{\omega^{\omega^{\alpha-1}}}}},\omega^\alpha),...,(+\infty_{\frac{1}{{\omega^{\omega^{\omega}}}}},\omega^\omega), ...,(+\infty_{\frac{1}{{\omega^{\omega^{\omega^\omega}}}}},\omega^{\omega^\omega}),...,(+\infty_{\frac{1}{\varepsilon_0}},\varepsilon_0),...
$$
$$
(+\infty_{\frac{1}{\omega_1}},\omega_1),..., (+\infty_{\frac{1}{\omega_2}},\omega_2),..., (+\infty_{\frac{1}{\omega_\alpha}},\omega_\alpha),..., (+\infty_{\frac{1}{\omega_\omega}},\omega_\omega),...
$$
and so on,
where $0<\alpha<\omega$, $1<\alpha<\omega$, $2<\alpha<\omega$  and $0<\alpha<\omega$, respectively.

{\bf Remark 15.} We also note that there are many intermediate, e.g.,  $\omega$-sequences of $\alpha\cdot\omega^\beta$-types, $0<\alpha,\beta<\omega$, but they appear much later. Using the  formula $(\ref{f04})$, for example, we determine the number $\frac{1}{\omega+1}=$
$$\{0,\frac{1}{\omega}-\frac{1}{\omega^2},\frac{1}{\omega}-\frac{1}{\omega^2}+\frac{1}{\omega^3}-\frac{1}{\omega^4},...\,|\,\frac{1}{\omega},\frac{1}{\omega}-\frac{1}{\omega^2}+\frac{1}{\omega^3},\frac{1}{\omega}-\frac{1}{\omega^2}+\frac{1}{\omega^3}-\frac{1}{\omega^4}+\frac{1}{\omega^5},...\}.$$
Notice only that these $\omega$-sequences define a symbol of infinity $+\infty_{\frac{1}{\omega+1}\mp\frac{1}{\omega^\omega}}$, respectively.

And in the same way we can define numbers $\frac{1}{\omega+\alpha}$, $0<\alpha<\omega$, and obtain an $\omega$-sequence $(\frac{1}{\omega+\alpha})_{0<\alpha<\omega}$ of $\frac{1}{2\cdot\omega}$-type and a symbol of infinity $+\infty_{\frac{1}{2\cdot\omega}}$. Similar for each $0<\beta<\omega$ $\omega$-sequence $(\frac{1}{\beta\cdot\omega+\alpha})_{0<\alpha<\omega}$ of $\frac{1}{\beta\cdot\omega}$-type defines a symbol of infinity $+\infty_{\frac{1}{\beta\cdot\omega}}$.

In the same way we can define numbers $\frac{1}{\omega-1}$, ..., $\frac{1}{\omega-\alpha}$, $0<\alpha<\omega$, and symbols of infinity $+\infty_{\frac{2}{\omega}}$.

The extend them to appropriate numbers $b=\sum\limits_{0\leq\beta<\gamma}\omega^{y_\beta}r_\beta$ (a canonical form $(\ref{f0201})$ of $b$) with restrictions on $y_\beta$. We omit here evident  details.

{\bf Remark 16.} We note here the remarkable effect of {\it lagging and catching up}, say, of arguments and values, i.e., the birthday of a number $2\cdot\omega$ and the $\frac{1}{\omega^2}$-type of a decreasing sequence $(\frac{1}{\alpha})_{0<\alpha<2\cdot\omega}$, whose birthday on day $\omega^2$  will still be, here on day $2\cdot\omega$  is behind day $\omega^2$, and this day $\omega^2$ is ahead of day $2\cdot\omega$. But from day $\omega^\omega$ and day $\omega^{\omega^\omega}$, respectively, the laggard will catch up with the leader and coincide with it on day $\varepsilon_0$. 
The examples of exponential functions $y=\omega^x$ and $y=e^x$, $x\in{\bf No}$,  demonstrate this effect even more clearly. One can see that these functions are strictly increasing and for the first of which $y=\omega^x$ the argument $0\leq x<\varepsilon_0$ lags behind the function value, i.e., $x<\omega^x$ but for $x=\varepsilon_0$ the argument coincides with the value of the function, i.e., $\omega^{\varepsilon_0}=\varepsilon_0$; and for the second function $y=e^x$, the argument can be greater than the value of the function, i.e., $e^{\varepsilon_0}=\omega^{\frac{\varepsilon_0}{\omega}}<\omega^\varepsilon_0=\varepsilon_0$, a kind of anti-Zeno paradox: the tortoise (input) catches up with Achilles (output) in the first case and the tortoise overtakes him in the second case, where the tortoise is an argument and Achilles is a value of functions. The latter function as we saw above (Remark 10) has a fixed point only for $x=\omega$ and then, constantly increasing and being continuous in the linearly ordered topology in ${\bf No}$, it periodically moves above, say $e^{\varepsilon_0 +\omega^2}=e^{\varepsilon_0}\cdot e^{\omega^2}=\omega^{\frac{\varepsilon_0}{\omega}}\cdot\omega^\omega>\varepsilon_0 +\omega^2$ and below $e^{\varepsilon_1}=\omega^{\frac{\varepsilon_1}{\omega}}<\varepsilon_1$ the diagonal of its graph, and its graph never crosses the diagonal. This is an unexpected model for researchers.
\bigskip

\bigskip
\begin{center}
{\bf 3. The beginning of measure theory for closed and open intervals of Conway numbers}
\end{center}

Now we will look at the natural geometry of Conway numbers, i.e., consider these numbers in an obvious metric $\rho(a,b)=|b-a|$ for every pair of Conway numbers $a$ and $b$. It is clear that $\rho(a,a)=0$, $\rho(a,b)=\rho(b,a)$ and $\rho(a,b)+\rho(b,c)\geq\rho(a,c)$, for any Conway numbers $a,b,c$. Thus we obtain the corresponding metric topology or lineary topology with a base $\{(a,b)\}$, where $a,b$ run over \grqq all\grqq\, Conway numbers (more precisely, for all $a,b\in{\bf No}$ in $NBG$ theory  and for each numbers $a,b$ in super-class theories or in the \grqq eternal process theory\grqq). The only difference of usual metric and topology in the case ${\bf R}_\omega$ and the other cases ${\bf R}_\zeta$, $\zeta=\omega^{\omega^\nu}$, $\nu>0$, is that in the case $\nu=0$ we have Archimedean field, i.e., for all numbers $x,y\in{\bf R}_\omega$ such that $0<x<y$ there is a natural number $n>0$ such that $n\cdot x>y$, and for $\nu>0$ the fields ${\bf R}_\zeta$ are not Archimedean in above sense, i.e., there is no $n$ such that $n\cdot x>y$ for all numbers $x,y\in{\bf R}_\zeta$, $\nu>0$, such that $0<x<y$, but $\zeta$-Archimedean fields, i.e., for all numbers $x,y\in{\bf R}_\zeta$, $\nu>0$, such that $0<x<y$ there is an ordinal $\alpha\in(0,\zeta)$ such that $\alpha\cdot x>y$ (see details in \cite{l1332}). Of course in eternal processes the field of the corresponding case is not $\zeta$-Archimedean at all for any $\zeta=\omega^{\omega^\nu}$, $\nu\geq 0$, but nevertheless, for any pair of numbers  such that $0<x<y$ there is an ordinal $\alpha$ such that $\alpha\cdot x>y$.

Despite the fact that the latter case excludes the existence of such a field as non-mathematical, it can be mentally contemplated or viewed as Conway did (see \cite{l22}, p. 37, 42) on the straight \grqq number line\grqq\, consisting of Conway numbers or the \grqq number plane\grqq\, in the case of imaginary Conway numbers ${\bf C}_\xi={\bf R}_\xi+{\bf i}{\bf R}_\xi$, $\xi=\omega^{\omega^\nu}$, $\nu\geq 0$.

Next, for any two different Conway numbers $a$ and $b$, we can determine the length of the segments $[a,b]$, $[a,b)$, $(a,b]$ and  $(a,b)$  as a Conway number 
$\rho(a,b)$ in any fixed linear orded field ${\bf R}_\zeta$, $\zeta=\omega^{\omega^\nu}$, $0\leq\nu\leq\Omega$ in $NBG$ theory as well as in eternal processes, i.e., in cases of super-classes.

In the same sense, we cannot determine the length of the intervals, say $(-\infty_{-\omega},+\infty_\omega)$, $(-\infty_{-\omega},b)$ and $(a,+\infty_\omega)$, etc. because as we saw above there is a fundamental difference between numbers and symbols of infinity and thus there are no numbers  that express the \grqq length\grqq\, of these intervals  because there is no numbers (like $+\infty_\omega$ and $ -\infty_{-\omega}$ and thus $+\infty_\omega -\infty_{-\omega}$, $b -\infty_{-\omega}$ and $+\infty_\omega -a$) that correspond to what we would call the \grqq length\grqq\, of an interval, say $(-\infty_{-\omega},+\infty_\omega)$, $(-\infty_{-\omega},b)$ and $(a,+\infty_\omega)$, etc. that contains all numbers $x\in{\bf No}$ (or $x\in{\bf R}_\zeta$) such that  $-\infty_{-\omega}<x<+\infty_\omega$, $-\infty_{-\omega}<x<b$ and $a<x<+\infty_\omega$, respectively, etc. as well as in  eternal processes, i.e., say \grqq each super Conway number $x$\grqq.

Nevertheless, we can compare such intervals as geometric objects by isometric embeddings one into another, i.e., for example, for intervals $(+\infty_a,b]$ and $(+\infty_{a'},b']$, if the is an isometric mapping $i:(+\infty_a,b]\rightarrow(+\infty_{a'},b']$, i.e., for all numbers $x,y\in(+\infty_a,b]$ we have $\rho(i(x),i(y))=\rho(x,y)$, then we say that the  \grqq length\grqq\, of $(+\infty_a,b]$ (more exactly, the  \grqq space\grqq\, of   all Conway numbers in it in $NBG$ theory, i.e., the proper class $(+\infty_a,b]\cap {\bf No}$)  is less   than or equal to  the \grqq length\grqq\, of $(+\infty_{a'},b']$; if $i$ is an isometric mapping onto $(+\infty_{a'},b']$, then the \grqq length\grqq\, of $(+\infty_a,b]$ and $(+\infty_{a'},b']$ are equal; if  $(+\infty_{a'},b']$ is not equal to $(+\infty_{a},b]$ then we say sorter that  $(+\infty_{a'},b']$ is greater than $(+\infty_{a},b]$ or what is the same $(+\infty_{a},b]$ is less than $(+\infty_{a'},b']$.

Note that for every ${\bf R}_\xi$, $\xi=\omega^{\omega^\nu}$, $\nu\geq 0$, the definitions are similar.

For example, for any Conway number $a\in[0,+\infty_\omega)$ the \grqq length\grqq\, of $[0,+\infty_\omega)$ and $[a,+\infty_\omega)$ are equal because of the isometric mapping $i(x)=x+a$. 

On the other hand, we will see that the \grqq length\grqq\, of an interval $(+\infty_\omega,b]$, for any $+\infty_\omega<b<+\infty_{\frac{\omega}{2}}$, is less that the \grqq length\grqq\, of $(+\infty_\omega,\omega]$ but not equal to it as opposed to intervals $(+\infty_\omega,b]$, for any $+\infty_{\frac{\omega}{2}}<b<\omega$ which are equal to $(+\infty_\omega,\omega]$.

We now consider previously unnoticed patterns of the infinite sequences $(\alpha)_{0<\alpha<\lambda}$  of ordinal numbers $0<\alpha<\lambda$ ($\lambda$ is a limit ordinal) and their reciprocals $(\frac{1}{\alpha})_{0<\alpha<\lambda}$ as \grqq pulsars\grqq\, or pulsating \grqq jumps\grqq\, and their corresponding periodic \grqq compressions\grqq.

{\bf Definition 17.} For any strictly  increasing and  strictly decreasing  $\lambda$-sequences of Conway numbers $(\alpha)_{0\leq\alpha<\lambda}$ and $(\frac{1}{\alpha})_{0<\alpha<\lambda}$, where $\lambda$ is a limit ordinal, we call the  intervals $(+\infty_\lambda,\lambda]=\bigcap\limits_{0<\alpha<\lambda}[\alpha,\lambda]$ and $[\frac{1}{\lambda},+\infty_{\frac{1}{\lambda}})=\bigcap\limits_{0<\alpha<\lambda}[\frac{1}{\lambda},\frac{1}{\alpha}]$ \grqq{\it jumps}\grqq\, and \grqq {\it compressions}\grqq\, of the corresponding sequences, compared to regular \grqq{jumps}\grqq, or rather, \grqq steps\grqq\, and regular \grqq compressions\grqq, or rather, \grqq size reductions\grqq, i.e., intervals $(\alpha,\alpha+1]$ and $[\frac{1}{\alpha+1},\frac{1}{\alpha})$, respectively.

We will show now that the  $\lambda$-sequences $(\alpha)_{0\leq\alpha<\lambda}$  of ordinal numbers $0\leq\alpha<\lambda$ for different limit ordinals $\lambda \not=\lambda'$ have usually different or equal \grqq jumps\grqq,  but especially periodic with their own laws which we are going to study. We also add that these laws cannot be revealed from the classical representation of ordinal numbers, as types of well-ordered sets. For clarity, we will set some simple questions and problems, concerning \grqq jumps\grqq.

{\bf Question 1.} The \grqq length\grqq\, of $[0,+\infty_\lambda)$ is less/greater than the \grqq length\grqq\, of $(+\infty_\lambda,\lambda]$ or it is equal to the \grqq length\grqq\, of $(+\infty_\lambda,\lambda]$?

We will see that the second part of the question can be solved by one word: \grqq Never\grqq; the first part of the question can be solved like this: \grqq Can be less and can be greater\grqq, e.g., the \grqq length\grqq\, of $[0,+\infty_\omega)$ is less than the \grqq length\grqq\, of $(+\infty_\omega,\omega]$; the \grqq length\grqq\, of $[0,+\infty_{2\cdot\omega})$ is greater than the \grqq length\grqq\, of $(+\infty_{2\cdot\omega},2\cdot\omega]$.

{\bf Question 2.} Can \grqq lengths\grqq\, of $(+\infty_\lambda,\lambda]$ and $(+\infty_{\lambda'},\lambda']$ for different $\lambda,\lambda'$  be equal?

We will see that the second  question can be solved like this: \grqq Yes and no\grqq, e.g., the \grqq length\grqq\, of $(+\infty_{\omega^2},\omega^2]$ is greater than the \grqq length\grqq\, of $(+\infty_\omega,\omega]$; and the  \grqq length\grqq\, of $(+\infty_{2\cdot\omega},2\cdot\omega]$ is equal to the \grqq length\grqq\, of $(+\infty_\omega,\omega]$.

{\bf Question 3.} What are the necessary and sufficient conditions for a positive answer to the first part of the first question?

{\bf Question 4.} Can the \grqq jumps\grqq\, of canonical sequences be arbitrarily large?

The fourth question can be solved like this: "Yes. And it can be \grqq unimaginably huge\grqq.

{\bf Question 5.} Can  the \grqq length\grqq\, of $[0,+\infty_\lambda)$ be equal to the \grqq length\grqq\, of $[\alpha',+\infty_\lambda)$ for each ordinal $0<\alpha'<\lambda$?

The answer is this: \grqq Can be and can be not\grqq, e.g., the \grqq length\grqq\, of $[0,+\infty_\omega)$ is equal to the \grqq length\grqq\, of $[n,+\infty_\omega)$, $0<n<\omega$; the \grqq length\grqq\, of $[0,+\infty_{2\cdot\omega})$ is not equal to the \grqq length\grqq\, of $[\omega,+\infty_{2\cdot\omega})$.

{\bf Proposition 15.} {\it Let   $(\alpha)_{0\leq\alpha<\lambda}$ be a  $\lambda$-sequence with limit ordinal $\lambda$. Then the isometric mapping $i:[0,+\infty_\lambda)\rightarrow{\bf No}$, given by formula $i(x)=\lambda-x$, maps $[0,+\infty_\lambda)$ into $(+\infty_\lambda,\lambda]\subset{\bf No}$ if and only if  for each pair of ordinals  $0\leq\alpha,\alpha'<\lambda$, there is an inequality $\alpha+\alpha'<\lambda$.}

{\bf Proof of \grqq if\grqq.} If  for each pair of ordinals  $0\leq\alpha,\alpha'<\lambda$, there is an inequality $\alpha+\alpha'<\lambda$, then for each $x\in[0,+\infty_\lambda)$ there is an ordinal $0\leq\alpha'<\lambda$ such that $x\leq\alpha'$ and we have $i(x)\geq i(\alpha')$, i.e., $\lambda-x\geq\lambda-\alpha'$. But $\alpha<\lambda-\alpha'$ for all $0\leq\alpha<\lambda$ because of condition $\alpha+\alpha'<\lambda$ for all ordinals $0\leq\alpha,\alpha'<\lambda$. Thus, by the above inequality $\lambda-x\geq\lambda-\alpha'$, we obtain $\alpha<\lambda-x$ for all ordinals $0\leq\alpha,\alpha'<\lambda$. And we conclude that $i(x)\in (+\infty_\lambda,\lambda]$ for every number $x\in[0,+\infty_\lambda)$.

{\bf Proof of \grqq only if\grqq.} If the isometric mapping $i:[0,+\infty_\lambda)\rightarrow{\bf No}$, given by formula $i(x)=\lambda-x$, maps $[0,+\infty_\lambda)$ into $(+\infty_\lambda,\lambda]\subset{\bf No}$, then for arbitrary ordinal $\alpha'\in[0,+\infty_\lambda)$ we have an inequality $\lambda-\alpha'>\alpha$ for all ordinals $0\leq\alpha<\lambda$ and what means that $\alpha+\alpha'<\lambda$ for each pair of ordinals  $0\leq\alpha,\alpha'<\lambda$.
 $\Box$

{\bf Proposition 16.} {\it Let   $(\alpha)_{0\leq\alpha<\lambda}$ be a  $\lambda$-sequence with limit ordinal $\lambda$. Then the isometric mapping $i:[0,+\infty_\lambda)\rightarrow{\bf No}$, given by formula $i(x)=\lambda-x$, maps $[0,+\infty_\lambda)$ into $(+\infty_\lambda,\lambda]\subset{\bf No}$ if and only if   $\lambda=\omega^\mu$ for some ordinal number $\mu\geq 1$.}

{\bf Proof of \grqq if\grqq.} If $\lambda=\omega^\mu$ for some ordinal $\mu\geq 1$. We will show that for each pair of ordinals $0\leq\alpha,\alpha'<\omega^\mu$ there is an inequality $\alpha+\alpha'<\omega^\mu$. Indeed, let $\alpha=\omega^{\eta_1}\cdot n_1+\omega^{\eta_2}\cdot n_2+...+\omega^{\eta_k}\cdot n_k$ and $\alpha=\omega^{\eta'_1}\cdot m_1+\omega^{\eta'_2}\cdot m_2+...+\omega^{\eta'_l}\cdot m_l$ be canonical form of ordinals $\alpha$ and $\alpha'$, respectively, where ordinals $\eta_1>\eta_2>...>\eta_k>0$, $0<n_1,n_2,...,n_k<\omega$ and $\eta'_1>\eta'_2>...>\eta'_k>0$, $0<m_1,m_2,...,m_l<\omega$, respectively. Since $\alpha<\omega^\mu$ and $\alpha'<\omega^\mu$ then $\eta_1<\mu$ and $\eta'_1<\mu$, moreover, $\omega^{\eta_1}\cdot n_1+\omega^{\eta'_1}\cdot m_1<\omega^\mu$ and hence $\alpha+\alpha'<\omega^\mu$.

Then, by Proposition 15, the isometric mapping $i:[0,+\infty_\lambda)\rightarrow{\bf No}$, given by formula $i(x)=\lambda-x$, maps $[0,+\infty_\lambda)$ into $(+\infty_\lambda,\lambda]\subset{\bf No}$.

{\bf Proof of \grqq only if\grqq.} If the isometric mapping $i:[0,+\infty_\lambda)\rightarrow{\bf No}$, given by formula $i(x)=\lambda-x$, maps $[0,+\infty_\lambda)$ into $(+\infty_\lambda,\lambda]\subset{\bf No}$, then for arbitrary ordinal $\alpha'\in[0,+\infty_\lambda)$ we have an inequality $\lambda-\alpha'>\alpha$ for all ordinals $0\leq\alpha<\lambda$ and what means that $\alpha+\alpha'<\lambda$ for each pair of ordinals  $0\leq\alpha,\alpha'<\lambda$.

One can easily see that $\lambda=\omega^\mu$ for some ordinal $\mu\geq 1$. Otherwise, the limit odinal $\lambda$ should be of a canonical form $\lambda=\omega^{\zeta_1}\cdot n_1+\omega^{\zeta_2}\cdot n_2+...+\omega^{\zeta_k}\cdot n_k$, $\zeta_1>\zeta_2>..>.\zeta_k>0$ and $0<n_1,n_2,...,n_k<\omega$. Clearly, $\alpha'=\omega^{\zeta_2}\cdot n_2+...+\omega^{\zeta_k}\cdot n_k<\lambda$ and $\alpha=\omega^{\zeta_1}\cdot n_1<\lambda$ but $\alpha+\alpha'=\lambda$. Contradiction. Thus $\alpha'=0$ and $\lambda=\omega^{\zeta_1}\cdot n_1$, where $\zeta_1=\mu$ and $n_1=1$. If it is not so, i.e., $n_1>1$, then put $\alpha=\omega^{\zeta_1}<\omega^\mu$ and $\alpha'=\omega^{\zeta_1}\cdot(n_1-1)<\lambda$ and thus $\alpha+\alpha'=\lambda$. Contradiction. Therefore, $\lambda=\omega^\mu$. $\Box$

{\bf Proposition 17.} {\it Let   $(\alpha)_{0\leq\alpha<\lambda}$ be a  $\lambda$-sequence with limit ordinal $\lambda$. Then the isometric mapping $i:[0,+\infty_\lambda)\rightarrow[a,+\infty_\lambda)\subset{\bf No}$, given by formula $i(x)=x+a$ for any fixed number $0\leq a\in[0,+\infty_\lambda)$ maps $[0,+\infty_\lambda)$ onto $[a,+\infty_\lambda)$ if and only if   $\lambda=\omega^\mu$ for some ordinal number $\mu\geq 1$.}

{\bf Proof of \grqq if\grqq.} If the isometric mapping $i:[0,+\infty_\lambda)\rightarrow[a,+\infty_\lambda)\subset{\bf No}$, given by formula $i(x)=x+a$ for any fixed number $0\leq x\in[0,+\infty_\lambda)$ maps $[0,+\infty_\lambda)$ onto $[a,+\infty_\lambda)$, then for any ordinals 
$0\leq\alpha\leq\alpha'<\lambda$ we obtain $i(\alpha)=\alpha+\alpha'<\lambda$ when $a=\alpha'$ and $x=\alpha$. Then, by Proposition 16 (Proof of \grqq only if\grqq), we obtain $\lambda=\omega^\mu$ for some ordinal $\mu\geq 1$.

{\bf Proof of \grqq only if\grqq.} If $\lambda=\omega^\mu$ for some ordinal $\mu\geq 1$, then, by Proposition 16 (Proof of \grqq if\grqq), $\alpha+\alpha'<\lambda$ for any ordinals $0\leq\alpha<\alpha'<\lambda$. Let $a\in[0,+\infty_\lambda)$ be a fixed number and $x$ be an arbitrary number in $[0,+\infty_{\omega^\mu})$. If we consider the smallest ordinal $\alpha$ such that $a\leq\alpha$ and the smallest ordinal $\alpha'$ such that $x\leq\alpha'$, then $i(x)=x+a\leq \alpha'+\alpha<\omega^\mu$. Thus, $i$ is an injection. For any number $x'\in[a,+\infty_{\omega^\mu})$ we have $0\leq x'-a<+\infty_{\omega^\mu}$. Assuming $x=x'-a$ we obtain $i(x)=x'$. Thus, $i$ is a surjection. $\Box$

Proposition 17 is a more general statement than the well-known ordinal number theory statement that we used above when studying self-similar skands in Proposition 4, which is the well-known ordinal number theory statement is the following:

{\it In order for every remainder of $\alpha\not=0$ to be equal to $\alpha$, it is necessary and sufficient that $\alpha=\sigma+\rho$.} 

And we remind the reader that
{\it an ordinal number $\rho$ is called the remainder of an ordinal $\alpha$ if $\rho\not=0$ and there exists an ordinal number $\sigma$ such that $\alpha=\sigma+\rho$.}

In this case, the ordinal numbers introduced by Cantor are considered as a type of well-ordered sets, in particular, an ordinal $\alpha$ is a type of the well-ordered set $\{0,1,2,...,\alpha',...\}$, where $0\leq\alpha'<\alpha$, $\rho$ is a type of the well-ordered set $\{0,1,2,...,\alpha',...\}$, where $0\leq\alpha'<\rho$ and $\sigma$ is a type of the well-ordered set $\{0,1,2,...,\alpha',...\}$, where $0\leq\alpha'<\sigma$.

Moreover, $\sigma+\rho$ is a usual sum of ordinal numbers and differs from sum of Conway numbers, i.e., $\sigma+\rho$ above means Cantor's arithmetic and in general $\sigma+\rho\not=\rho+\sigma$ while in Conway's arithmetic always $\sigma+\rho=\rho+\sigma$; $\sigma+\rho\not<\rho$ when $\sigma<\rho$ in Cantor's arithmetic and $\sigma+\rho=\max\{\sigma+\rho,\rho+\sigma\}$ in Conway's arithmetic.

The approach to ordinal numbers may seem \grqq similar\grqq\, to that of Cantor and Conway, but Conway's approach is quite different in form Cantor's approach, although they are similar in content, and significantly different in the form and thus in arithmetic.

 Indeed, if Conway's canonical form $\alpha=\{0,1,2,...,\alpha',...\,|\,\}$, where $0\leq\alpha'<\alpha$, of an ordinal number $\alpha$ is similar to Cantor's definition of $\alpha$ as a type of the well-ordered set $\{0,1,2,...,\alpha',...\}$, where  $0\leq\alpha'<\alpha$, then Conway's other forms of the same number  $\alpha=\{\alpha_0,\alpha_0+1,\alpha_0+2,...,\alpha_0+\alpha',...\,|\,\}$ for every fixed $\alpha_0<\alpha$ and all $\alpha'$ such that $\alpha_0\leq\alpha_0+\alpha'<\alpha$, do not match Cantor's definition in the same way since,  in general,  types of the well-ordered sets $\{0,1,2,...,\alpha',...\}$ and $\{\alpha_0,\alpha_0+1,\alpha_0+2,...,\alpha_0+\alpha',...\}$ are different, in general, e.g., $\{0,1,2,...,n-1\}$ is according to  Cantor an ordinal $n$ and $\{0,1,2,...,n-1\,|\,\}$ is according to Conway an ordinal $n$, but $\{2,3,4,...,n-1\}$ is according to Conway  an ordinal $n-2$ although $\{2,3,4,...,n-1\,|\,\}$ is still $n$. On the other hand, $\omega=\{0,1,2,...,\alpha,...\,|\,\}=\{\alpha_0,\alpha_0+1,\alpha_0+2,...,\alpha_0+\alpha,...\,|\,\}$, for every fixed $0<\alpha_0<\omega$ and both sets $\{0,1,2,...,\alpha,...\}$ and $\{\alpha_0,\alpha_0+1,\alpha_0+2,...,\alpha_0+\alpha,...\}$ are of the same type. This is the meaning of the above statement, which we will prove now in terms of Conway's number theory.

Note that despite the difference in the form of the definition of ordinal numbers and their arithmetic by Cantor and Conway, the ordinal numbers themselves are the same and maintain the same order, i.e., say with conventional designations $\alpha_{can}=\{0,1,2,...,\alpha',...\}=\{0,1,2,...,\alpha',...\,|\,\}=\alpha_{con}$ and $\alpha_{can}<\alpha_{can}\Leftrightarrow\alpha_{con}<\alpha_{con}$.

{\bf Proposition 18.} {\it Let   $(\alpha)_{0\leq\alpha<\lambda}$ be a  $\lambda$-sequence with limit ordinal $\lambda$. Then for any ordinal number $0<\alpha_0<\lambda$ sets $\{0,1,2,...,\alpha,..\}_{0\leq\alpha<\lambda}$ and $\{\alpha_0,\alpha_0+1,\alpha_0+2,...,\alpha_0+\alpha,...\}_{\alpha_0\leq\alpha<\lambda}$ are of the same type as well-ordered sets, i.e., there is an isomorphism $f:\{0,1,2,...,\alpha,..\}_{0\leq\alpha<\lambda}\rightarrow\{\alpha_0,\alpha_0+1,\alpha_0+2,...,\alpha_0+\alpha,...\}_{0\leq\alpha<\lambda}$ if and only if $\lambda=\omega^\mu$ for some ordinal $0<\mu<\Omega$.}

{\bf Proof of \grqq only if\grqq.} Let $f:\{0,1,2,...,\alpha,..\}_{0\leq\alpha<\lambda}\rightarrow\{\alpha_0,\alpha_0+1,\alpha_0+2,...,\alpha_0+\alpha,...\}_{0\leq\alpha<\lambda}$ be an isomorphism  of well-ordred sets for each ordinal $0<\alpha_0<\lambda$. 

If $\lambda=\omega^{\eta_1}\cdot n_1+\omega^{\eta_2}\cdot n_2+...+\omega^{\eta_k}\cdot n_k$ were a canonical form of ordinals $\lambda$, where ordinals $\eta_1>\eta_2>...>\eta_k>0$, $0<n_1,n_2,...,n_k<\omega$, then $\alpha=\omega^{\eta_2}\cdot n_2+...+\omega^{\eta_k}\cdot n_k<\omega^{\eta_1}\cdot n_1=\alpha'<\lambda$.

Thus $\alpha$ must be equal to $0$ in order not to contradict the assumption of Proposition 18. If now $\lambda=\omega^{\eta_1}\cdot n_1$ and $n_1\not=1$, then putting $\alpha=\omega^{\eta_1}<\omega^{\eta_1}\cdot (n_1-1)=\alpha'<\lambda$ we obtain $\omega^{\eta_1}\cdot (n_1-1)=0$ and thus $\lambda=\omega^\eta+1$. $\Box$

{\bf Proof of \grqq  if\grqq.} Let $\lambda$ be equal to $\omega^\mu$ for some ordinal $0<\mu<\Omega$. We will show that there exists an isomorphism  of well-ordered sets $f:\{0,1,2,...,\alpha,..\}_{0\leq\alpha<\lambda}\rightarrow\{\alpha_0,\alpha_0+1,\alpha_0+2,...,\alpha_0+\alpha,...\}_{0\leq\alpha<\lambda}$  for each ordinal $0<\alpha_0<\lambda$.

Let $\alpha_0=\omega^{\eta_1}\cdot n_1+\omega^{\eta_2}\cdot n_2+...+\omega^{\eta_k}\cdot n_k$ be a canonical form of $\alpha_0$, where ordinals $\eta_1>\eta_2>...>\eta_k\geq 0$, $0<n_1,n_2,...,n_k<\omega$. Since $0<\alpha_0<\omega^\mu$ for some ordinal $0<\mu<\Omega$ then $\eta_1<\mu$ and $\eta_1+1\leq\mu$. Indeed, as we saw above for each ordinals $0\leq\alpha\leq\alpha'<\omega^\mu$ there is inequality $\alpha+\alpha'<\omega^\mu$ we conclude that $m\cdot\omega^{\eta_1}\cdot n_1<\omega^\mu$ for all natural numbers $0< m<\omega$ because for $m=1$ we have $\omega^{\eta_1}\cdot n_1<\omega^\mu$ 
and thus $\omega^{\eta_1}\cdot n_1+(m-1)\cdot\omega^{\eta_1}\cdot n_1<\omega^\mu$ for all $1<m<\omega$, by putting $\alpha=\cdot\omega^{\eta_1}\cdot n_1<\omega^\mu$ and $\alpha'=(m-1)\cdot\omega^{\eta_1}\cdot n_1<\omega^\mu$.

Clearly, $\{\omega^{\eta_1}\cdot n_1,2\cdot\omega^{\eta_1}\cdot n_1,...,m\cdot\omega^{\eta_1}\cdot n_1,...\,|\,\}\leq\omega^\mu$, $0<m<\omega$, 
and thus $\omega^{\eta_1}\cdot n_1+(m-1)\cdot\omega^{\eta_1}\cdot n_1<\omega^\mu$ for all $1<m<\omega$, by putting $\alpha=\cdot\omega^{\eta_1}\cdot n_1<\omega^\mu$ and $\alpha'=(m-1)\cdot\omega^{\eta_1}\cdot n_1<\omega^\mu$,
because $\{ordinals\,\,\, \alpha<\omega^\mu\,|\,\}=\omega^\mu$ and $\{m\cdot\omega^{\eta_1}\cdot n_1\}_{0\leq m<\omega}$ is a part of all ordinals $\{\alpha\}_{0\leq\alpha<\omega^\mu}$.

So, if $\eta_1+1<\mu$, then we put $f(\alpha)=\alpha+\alpha_0$ when $0\leq<\omega^{\eta_1+1}$ and $f(\alpha)=\alpha$ when $\omega^{\eta_1+1}\leq\alpha<\omega^\mu$. And if $\eta_1+1=\mu$, then we put $f(\alpha)=\alpha+\alpha_0$ when $0\leq\alpha<\omega^\mu$. The mapping $f(\alpha)=\alpha+\alpha_0$ is correct because the sequence $(m\cdot\omega^{\eta_1}\cdot n_1)_{0\leq m<\omega}$ is cofinal to the sequence 
$(\alpha)_{0\leq\alpha<\omega^\mu}$ and the sequence $(\alpha+\alpha_0)_{0\leq\alpha<\omega^\mu}$ is cofinal to $(m\cdot\omega^{\eta_1}\cdot n_1)_{0\leq m<\omega}$. $\Box$

{\bf Proposition 19.} {\it For any ordinal $0<\delta<\Omega$, there exists an ordinal $\omega\leq\lambda<\Omega$ $(\lambda=\omega^\mu, 0<\mu<\Omega)$ such that the \grqq jump\grqq\, $(+\infty_\lambda,\lambda]$ of a sequence $(\alpha)_{0\leq\alpha<\lambda}$ is significantly larger than $\delta$, i.e., there is an isometric embedding $i:[0,+\infty_\lambda)\rightarrow(+\infty_\lambda,\lambda]$ and thus $i([0,\delta])\subset(+\infty_\lambda,\lambda]$; moreover, for any number $b=\sum\limits_{0\leq\beta<\gamma}\omega^{y_\beta}r_\beta$ $($normal form of $b)$ in which $\gamma$ denotes some ordinal, the numbers $r_\beta$ $(\beta<\gamma)$ are non-zero reals, and the numbers $y_\beta$ form a strictly decreasing sequence of numbers such that   for all $0\leq\beta<\gamma$ numbers $\mu-1<y_\beta<\mu$, when $\mu$ is not a limit ordinal, and  $\kappa<y_\beta<\mu$ for all ordinals $\kappa<\mu$, when $\mu$ is a limit ordinal, the interval $(+\infty_\lambda,+\infty_{\frac{\lambda}{2}})$ is a discrete sum of intervals $(+\infty_{c_{-}},+\infty_{c_{+}})$, where $+\infty_{c_{+}}$ is an infinity symbol with an index (a label) $c_{+}$defined by the sequence $(b+\alpha)_{0\leq\alpha<\lambda}$, and $+\infty_{c_{-}}$ is an infinity symbol with an index (a label) $c_{-}$, defined by the sequence $(b-\alpha)_{0\leq<\alpha<\lambda}$, over all such numbers $b\in B$ $(B$ is a proper class of all such numbers $b)$, i.e., $(+\infty_\lambda,+\infty_{\frac{\lambda}{2}})=\bigcup\limits_{b\in B}(+\infty_{c_{-}},+\infty_{c_{+}})$.}

{\bf Proof.} Let  $\delta$ be an ordinal $0<\delta<\Omega$. Consider the smallest ordinal $\mu$ such that $\delta<\omega^\mu=\lambda$. Then, by Proposition 16, there is an isometric embedding $i:[0,\lambda)\rightarrow (+\infty_\lambda,\lambda]$, given by formula $i(x)=\lambda-\alpha$ for each $0\leq\alpha<\lambda,$ and clearly, $i([0,\delta])\subset[0,\lambda)$.

Let $b\in(+\infty_\lambda,+\infty_{\frac{\lambda}{2}})$ be any number, whose  normal form is  $b=\sum\limits_{0\leq\beta<\gamma}\omega^{y_\beta}r_\beta$  in which $\gamma$ denotes some ordinal, the numbers $r_\beta$ $(\beta<\gamma)$ are non-zero reals, and the numbers $y_\beta$ form a strictly decreasing sequence of numbers such that   for all $0\leq\beta<\gamma$ numbers $\mu-1<y_\beta<\mu$, when $\mu$ is not a limit ordinal, and  $\kappa<y_\beta<\mu$ for all ordinals $\kappa<\mu$, when $\mu$ is a limit ordinal. Then $b-\alpha>\alpha'$ for all ordinals $0\leq\alpha'<\lambda$ and thus $b-\alpha\in(+\infty_\lambda,+\infty_{\frac{\lambda}{2}})$. Indeed, if this were not the case and the inequality were to hold $b-\alpha_0<\alpha_1$ for some ordinals $0<\alpha_0\leq\alpha_1<\lambda$, then there would be   the inequality would be $b-\alpha_0<\alpha_1$ or what is the same $b-\alpha_0\in[0,+\infty_\lambda)$. Since $i:[0,b]\rightarrow[0,b]$, given by formula $i(x)=b-x$, is an isometric mapping it is clear that $b-\alpha\in[0,b]$ for all ordinal $0\leq\alpha<\lambda$ because $\alpha<b$ for all ordinal $0\leq\alpha<\lambda$ or what is the same $b\in(+\infty_\lambda,+\infty_{\frac{\lambda}{2}})$ we conclude that $[b-\alpha_0-\alpha_1,b-\alpha_0]\subset[0,\alpha_1]$, but it is impossible because the length of $[\lambda-\alpha_0-\alpha_1,\lambda-\alpha_0]$ strictly less than the length of $[0,\alpha_1]$ since $(\lambda-\alpha_0)-(\lambda-\alpha_0-\alpha_1)=\alpha_1$ and thus $\alpha_1<\alpha_1$ what impossible because for any ordinal $\alpha$ there is no such inequality as $\alpha<\alpha$.

 The same arguments lead to the fact that for any number $b\in(+\infty_\lambda,+\infty_{\frac{\lambda}{2}})$ numbers $b+\alpha\in(+\infty_\lambda,+\infty_{\frac{\lambda}{2}})$ for all $0\leq\alpha<\lambda$.

Clearly, for different $b,b'\in B$ we obtain $(+\infty_{c_{-}},+\infty_{c_{+}})\cap(+\infty_{c'_{-}},+\infty_{c'_{+}})=\emptyset$ $(+\infty_\lambda,+\infty_{\frac{\lambda}{2}})=\bigcup\limits_{b\in B}(+\infty_{c_{-}},+\infty_{c_{+}})$.

Note that if $b<b'$ and $b''=b'-b=\sum\limits_{0\leq\beta<\gamma}\omega^{y_\beta}r_\beta$ $($normal form of $b'')$ in which $\gamma$ denotes some ordinal, the numbers $r_\beta$ $(\beta<\gamma)$ are non-zero reals, and the numbers $y_\beta$ form a strictly decreasing sequence of numbers such that   for all $0\leq\beta<\gamma$ numbers $y_\beta<\kappa$, for some $\kappa<\mu$, then $(+\infty_{c_{-}},+\infty_{c_{+}})=(+\infty_{c'_{-}},+\infty_{c'_{+}})$. What is the same if $b'\in(+\infty_{c_{-}},+\infty_{c_{+}})$, then $i:(+\infty_{c_{-}},+\infty_{c_{+}})\rightarrow(+\infty_{c'_{-}},+\infty_{c'_{+}})$, given by $i(x)=x+b'-b$ when $b<b'$ and  $i(x)=x+b'-b$ when $b'<b$. In other words, if $(+\infty_{c_{-}},+\infty_{c_{+}})\cap(+\infty_{c'_{-}},+\infty_{c'_{+}})\not=\emptyset$, then $(+\infty_{c_{-}},+\infty_{c_{+}})=(+\infty_{c'_{-}},+\infty_{c'_{+}})$. $\Box$

{\bf Remark 17.}  Proposition 19 is true not only in the framework of the $NBG$ theory, which we implied when we talked about proper classes of interval numbers $(+\infty_{c_{-}},+\infty_{c_{+}})$, but also in the theory of super-classes. In this case, every number, say $b$ , in the super-class $(+\infty_\lambda,+\infty_{\frac{\lambda}{2}})$ has $(+\infty_{c_{-}},+\infty_{c_{+}})$ as its proper super-subclass and so it is impossible for a human to imagine the magnitude of this \grqq jump\grqq\, $(+\infty_\lambda,\lambda]$ of a sequence $(\alpha)_{0\leq\alpha<\lambda=\omega^\mu}$, i.e., there is a \grqq perpetuum mobile\grqq\, increasing and decreasing sequence $(b_\alpha)_{0\leq\alpha}$, $b_\alpha\in(+\infty_\lambda,+\infty_{\frac{\lambda}{2}})$ such these intervals of numbers $(+\infty_{b_{\alpha-}},+\infty_{b_{\alpha+}})\subset(+\infty_\lambda,+\infty_{\frac{\lambda}{2}})$ are linear ordered and disjoint. However, for the ordinal $\lambda=\omega^{\mu+1}$, the \grqq jump\grqq\, $(+\infty_\lambda,\lambda]$ of a sequence $(\alpha)_{0\leq\alpha<\lambda=\omega^{\mu+1}}$ is in the same sense incredibly larger than the previous \grqq jump\grqq.

We note that a certain law of \grqq jumps\grqq\, emerges. The first \grqq jump\grqq\, $(+\infty_\omega,\omega]$ is repeated $\omega$-times, i.e., $(+\infty_{2\cdot\omega},2\cdot\omega]$, ..., $(+\infty_{n\cdot\omega},n\cdot\omega]$, ...\,, $1<n<\omega$, and then a greater \grqq jump\grqq\, $(+\infty_{\omega^2},\omega^2]$ occurs. Next, the second \grqq jump\grqq\, is repeated $\omega$-times along with the $\omega$-times of the first \grqq jump\grqq\, , i.e., $(+\infty_{2\cdot\omega^2},2\cdot\omega^2]$, ..., $(+\infty_{n\cdot\omega^2},n\cdot\omega^2]$, ...\,, $1<n<\omega$, and a third \grqq jump\grqq\, is made $(+\infty_{\omega^3},\omega^3]$ , greater than all the previous ones combined, and so on. In a word, there is a kind of {\it periodically increasing pulsating flow} of \grqq jumps\grqq, shortly called here \grqq pulsars\grqq.

There is a certain cyclicity or repetition of \grqq jumps\grqq\, between the limit ordinals of the form $\omega^\mu$, $0<\mu<\Omega$.
 Additionally, all \grqq jumps\grqq\, of the same magnitude can be listed using ordinal numbers, as they are well-ordered.

{\bf Proposition 20.} {\it All \grqq jumps\grqq\, $(+\infty_\lambda,\lambda]$, for all limit ordinals $\omega\leq...<\lambda<\omega^\mu$,  for a fixed ordinal $0<\mu<\Omega$, are well-ordered and moreover, all equal \grqq jumps\grqq\, between them can be  calculated by   types of the well-ordered sets.}

{\bf Proof.} It is clear that all \grqq jumps\grqq\, $(+\infty_\lambda,\lambda]$, for all limit ordinals $\omega\leq...<\lambda<\omega^\mu$,  for a fixed ordinal $0<\mu<\Omega$, are well-ordered by the ends $\lambda$ of intervals $(+\infty_\lambda,\lambda]$.

The second part of Proposition 20 can be proved by the transfinite induction. We start with the beginning of induction $\lambda=\omega$. There the only one \grqq jumps\grqq\, $(+\infty_\omega]$. The second \grqq jumps\grqq\, $(+\infty_{2\cdot\omega},2\cdot]$. The $n$th \grqq jumps\grqq\, is $(+\infty_{n\cdot\omega},n\cdot\omega]$. The $\omega$th \grqq jumps\grqq\, is $(+\infty_{\omega^2},\omega^2]$. Moreover, we correspond to each \grqq jumps\grqq\, is $(+\infty_{\alpha\cdot\omega},\alpha\cdot\omega]$ an ordinal $\alpha\cdot\omega$, $0<\alpha<\omega$, and for  $(+\infty_{\omega^2},\omega^2]$ we correspond an ordinal  $\omega^2.$ In total, we have $\omega$ \grqq jumps\grqq\,  \grqq length\grqq\, $(+\infty_\omega,\omega]$ and one \grqq jump\grqq\, length\grqq\, $(+\infty_{\omega^2},\omega^2]$ or, what is the same $\omega^\nu$, $2>\nu\geq 0$. Similar for $\lambda=\omega^n$, i.e., $n>\nu\geq 0$ or what is the same, we have $\omega$ \grqq jumps\grqq\,  \grqq length\grqq\, $(+\infty_\omega,\omega]$ and one \grqq jump\grqq\, length\grqq\, $(+\infty_{\omega^2},\omega^2]$. $\Box$

\bigskip

A few brief remarks about the \grqq compressions\grqq\, dual to the notion of  \grqq jumps\grqq. First about the simplest \grqq compression\grqq\, $[\frac{1}{\omega},+\infty_{\frac{1}{\omega}})$ of the $\omega$-sequence $(\frac{1}{\alpha})_{0<\alpha<\omega}$. Naturally, that interval  $[\frac{1}{\omega},+\infty_{\frac{1}{\omega}})\subset{\bf No}$ of Conway numbers is much less than the interval $(+\infty_{\frac{1}{\omega}},1]$, but the interval $[0,\frac{1}{\omega}]$ is incomparably less than the interval $[\frac{1}{\omega},+\infty_{\frac{1}{\omega}})$, in the language of mathematics it is  anti-Archimedean, i.e., for every number $x\in[\frac{1}{\omega},+\infty_{\frac{1}{\omega}})$ there is a number $y\in[\frac{1}{\omega},+\infty_{\frac{1}{\omega}})$ such that $x<y$ and for all natural numbers $n$ there are inequalities $n\cdot x<y$. This can be expressed more precisely here by means of astrophysical terminology conventional for mathematics, as \grqq black hole\grqq, \grqq event horizon\grqq\, and \grqq singularity\grqq\, as  the \grqq compression\grqq\, $[\frac{1}{\omega},+\infty_{\frac{1}{\omega}})$, the symbol of infinity $+\infty_{\frac{1}{\omega}}$ and the number $\frac{1}{\omega}$, respectively. Since the minimum number $\frac{1}{\omega}$ cannot be exceeded by any multiple addition or multiplication of natural numbers $n$ from the interval $[\frac{1}{\omega},+\infty_{\frac{1}{\omega}})$, as well as all other numbers $x\in [\frac{1}{\omega},+\infty_{\frac{1}{\omega}})$ in this interval, despite the fact that they exceed this number  $\frac{1}{\omega}$.

\grqq Compressions\grqq\, are interesting in their own right from the point of view of algebraic structures, namely, subgroups, subrings, subalgebras, submodules and vector subspaces of Conway's Field ${\bf No}$. One can easily see that the interval $(-\infty_{-\frac{1}{\omega}},+\infty_{\frac{1}{\omega}})\subset{\bf No}$ is  an addition group, actually a ring as well as a module over the ring of integers ${\mathbb Z}$ and the field of all finite Conway numbers ${\bf F}\stackrel{def}{=}(-\infty_{-\omega},+\infty_\omega)\cap{\bf No}$  and even a vector space  over the field of real numbers ${\bf R}$ as well as over ${\bf F}$.

Moreover, there is a filtration of the ring ${\bf No}$,  turning it into a filtered ring
 ${\bf No}\supset...\supset {\bf A}_{\zeta'}\supset   {\bf A}_\zeta\supset...\supset{\bf A}_{\frac{1}{\zeta}}.\supset   {\bf A}_{\frac{1}{\zeta'}}\supset...$, where ${\bf A}_{\zeta'}\stackrel{def}{=}(-\infty_{-\zeta'},+\infty_{\zeta'})\cap{\bf No}$, ${\bf A}_\zeta\stackrel{def}{=}(-\infty_{-\zeta},+\infty_\zeta)\cap{\bf No}$ and ${\bf A}_{\frac{1}{\zeta}}\stackrel{def}{=}(-\infty_{-\frac{1}{\zeta}},+\infty_{\frac{1}{\zeta}})\cap{\bf No}$ and ${\bf A}_{\frac{1}{\zeta'}}\stackrel{def}{=}(-\infty_{-\frac{1}{\zeta'}},+\infty_{\frac{1}{\zeta'}})\cap{\bf No}$, $\zeta=\omega^{\omega^\mu}$ and $\zeta'=\omega^{\omega^{\mu+1}}$,  where  $0\leq\mu<\Omega$,  and additive groups ${\bf A}_{\frac{1}{\zeta}}$ is not $\zeta$-Archimedean over  rings ${\bf A}_{\zeta}$, i.e., for every numbers $x\in{\bf A}_{\frac{1}{\zeta}}$ there is a number $y\in{\bf A}_{\frac{1}{\zeta}}$ such that $x<y$ and for all natural numbers $z\in{\bf A}_\zeta$ there are inequalities $z\cdot x<y$.

There are also many multiplicative subgroups of the groups ${\bf B}_\Omega=(+\infty_{\frac{1}{\Omega}},+\infty_\Omega)\cap{\bf No}^+$ and ${\bf C}_\Omega \stackrel{def}{=}{\bf No}\setminus\{0\}=[(-\infty_{-\Omega},-\infty_{-\frac{1}{\Omega}})\cup(+\infty_{\frac{1}{\Omega}}+\infty_{\Omega})]\cap{\bf No}$, respectively.

 Clearly, additive groups ${\bf A}_\zeta$ define, for each  $\zeta=\omega^{\omega^\mu}$,  where  $0\leq\mu<\Omega$, multiplicative subgroups  ${\bf B}_\zeta \stackrel{def}{=}{\bf A}^+_\zeta\setminus{\bf A}_{\frac{1}{\zeta}}=(+\infty_{\frac{1}{\zeta}},+\infty_\zeta)\cap{\bf No}$ and ${\bf C}_\zeta \stackrel{def}{=}{\bf A}_\zeta\setminus{\bf A}_{\frac{1}{\zeta}}=[(-\infty_{-\zeta},-\infty_{-\frac{1}{\zeta}})\cup(+\infty_{\frac{1}{\zeta}}+\infty_{\zeta})]\cap{\bf No}$,  ${\bf B}_{\frac{1}{\zeta}}\stackrel{def}{=}\{1+x\}_{x\in{\bf A}_{\frac{1}{\zeta}}}$ and ${\bf C}_{\frac{1}{\zeta}}\stackrel{def}{=}\{-1+x\}_{x\in{\bf A}_{\frac{1}{\zeta}}}\cup\{1+x\}_{x\in{\bf A}_{\frac{1}{\zeta}}}$ of the multiplicative groups  
${\bf B}_\zeta$ and 
${\bf C}_\zeta$, respectively.

Thus we obtain a double-filtered group ${\bf C}_\Omega$
 
 ${\bf C}_\Omega\supset...\supset {\bf C}_{\zeta'}\supset   \,\,\,{\bf C}_\zeta\supset...\supset   {\bf C}_\omega\supset{\bf C}_{\frac{1}{\omega}}\supset...\supset {\bf C}_{\frac{1}{\zeta'}}\supset   {\bf C}_{\frac{1}{\zeta}}\supset\,\,\,\,...$

$\,\,\bigcup\,\,\,\,\,\,\,\,\,...\,\,\,\,\,\,\,\,\bigcup\,\,\,\,\,\,\,\,\,\,\,\,\,\,\bigcup\,\,\,\,\,\,\,\,\,...\,\,\,\,\,\,\,\,\,\,\bigcup\,\,\,\,\,\,\,\,\,\,\bigcup\,\,\,\,\,\,\,\,\,\,...\,\,\,\,\,\,\,\bigcup\,\,\,\,\,\,\,\,\,\,\,\,\,\,\bigcup\,\,\,\,\,\,\,\,\,\,\,\,\,\,\,...$

${\bf B}_\Omega\supset\,...\,\supset {\bf B}_{\zeta'}\supset  \,\, {\bf B}_\zeta\,\,\supset...\supset   {\bf B}_\omega\supset{\bf B}_{\frac{1}{\omega}}\supset...\supset {\bf B}_{\frac{1}{\zeta'}}\supset   {\bf B}_{\frac{1}{\zeta}}\supset\,\,\,\,...$

where $\zeta=\omega^{\omega^\mu}$, $\zeta'=\omega^{\omega^{\mu+1}}$,    $0\leq\mu<\Omega$.

It is known that additive group of the real numbers ${\bf R}_\omega={\bf R}$ is isomorphic to the multiplicative group ${\bf R}^+_\omega={\bf R}^+$; moreover, the field of all {\it finite} Conway numbers ${\bf A}_\omega$  is isomorphic to ${\bf B}_\omega$, i.e., $\varphi(x)=e^x$ is such isomorphism (see \cite{l22}, p. 43), in particular, $\varphi |_{{\bf A}_{\frac{1}{\zeta}}}:{\bf A}_{\frac{1}{\zeta}}\rightarrow{\bf B}_{\frac{1}{\zeta}}$ is an isomorphism of groups when $\zeta=\omega^{\omega^\mu}$, $0\leq\mu<\Omega$, (see below), but whether there is an isomorphism $\varphi:{\bf A}_\zeta\rightarrow{\bf B}_\zeta$, $\zeta>\omega$, of groups ${\bf A}_\zeta$ and ${\bf B}_\zeta$, in particular, ${\bf No}$ and ${\bf No}^+$,  is an {\it open problem}. 

The answer is most likely to be {\it no}.
 Although there is an isomorphism of ${\bf No}$ to the discrete subgroup ${\bf E}\stackrel{def}{=}\{\omega^x\}_{x\in{\bf No}}$ of the multiplicative group ${\bf No}^+$ and even to a larger multiplicative subgroup ${\bf D}^+_\Omega\subset{\bf No}^+$ (see next Proposition 21).
\bigskip
\bigskip
 \bigskip
{\bf Proposition 21.} {\it There is an additive  double-filtered group ${\bf L}_\Omega$

${\bf L}_\Omega\supset...\,\supset {\bf L}_{\zeta'}\supset   \,\,\,{\bf L}_\zeta\supset...\supset   {\bf L}_\omega\supset{\bf L}_{\frac{1}{\omega}}\supset...\supset {\bf L}_{\frac{1}{\zeta'}}\,\,\supset   {\bf L}_{\frac{1}{\zeta}}\supset\,\,\,\,...$

$\,\,\bigcup\,\,\,\,\,\,\,\,\,...\,\,\,\,\,\,\,\,\,\bigcup\,\,\,\,\,\,\,\,\,\,\,\,\,\,\bigcup\,\,\,\,\,\,\,\,...\,\,\,\,\,\,\bigcup\,\,\,\,\,\,\,\,\,\,\,\,\bigcup\,\,\,\,\,\,\,\,\,...\,\,\,\,\,\,\,\,\,\,\bigcup\,\,\,\,\,\,\,\,\,\,\,\,\bigcup\,\,\,\,\,\,\,\,\,\,\,\,\,\,\,...$

${\bf A}_\Omega\supset...\supset \,\,{\bf A}_{\zeta'}\supset   {\bf A}_\zeta\supset...\supset   {\bf A}_\omega\supset{\bf A}_{\frac{1}{\omega}}\supset...\supset {\bf A}_{\frac{1}{\zeta'}}\supset   {\bf A}_{\frac{1}{\zeta}}\supset\,\,\,\,...\,\,\,\,\,\,\,\,\,\,\,\,\,\,\,\,\,\,\,\,\,\,\,\,\,\,$
where $\zeta=\omega^{\omega^\mu}$, $\zeta'=\omega^{\omega^{\mu+1}}$,    $0\leq\mu\leq\Omega$, and there exists an isomorphism $\varphi_\Omega: {\bf L}_\Omega\rightarrow {\bf D}_\Omega$ of groups   
such that $\varphi_{\Omega}({\bf L}_\zeta)\subset{\bf C}_\zeta$ and  $\varphi_\Omega({\bf A}_\zeta)\subset{\bf B}_\zeta$ for all $\zeta=\omega^{\omega^\mu}$, $0<\mu\leq\Omega$,   $\varphi_\Omega({\bf L}_\omega)={\bf C}_\omega$, $\varphi_\Omega({\bf L}_{\frac{1}{\omega}})={\bf C}_{\frac{1}{\omega}}$ and  $\varphi_\Omega({\bf A}_{\frac{1}{\zeta}})={\bf B}_{\frac{1}{\zeta}}$ for all $\zeta=\omega^{\omega^\mu}$, $0<\mu\leq\Omega$.}

{\bf Proof.} Let's define  a mapping $\varphi:{\bf No}\rightarrow{\bf No}^+$ by the following forlula 
\begin{equation}
\label{f203}
\varphi(x)=e^{x'+x''}\stackrel{def}{=}e^{\omega\cdot{\frac{x'}{\omega}}}\cdot e^{x''}=\omega^{{\frac{x'}{\omega}}}\cdot e^{x''}
\end{equation}
 by setting $e^\omega=\omega$ and consider for each $x\in{\bf No}$ its parts   $x=x'+x''$, where  $x'=\sum\limits_{0\leq\beta<\gamma'}\omega^{y'_\beta}r'_\beta$ (normal form of $x'$) in which $\gamma'$ denotes some ordinal, the numbers $r'_\beta$ ($\beta<\gamma'$) are non-zero reals, and the numbers $y'_\beta$ form a strictly decreasing sequence of numbers such that $y'_\beta>0$, for all $0\leq\beta<\gamma'$ and a normal form of $x''=\sum\limits_{0\leq\beta<\gamma''}\omega^{y''_\beta}r''_\beta$ such that $y''_\beta\leq 0$, for all $0\leq\beta<\gamma''$ (actually, $x''\in{\bf A}_\omega$) and $\omega^{\frac{x'}{\omega}}$ is given by formula $(\ref{f314})$, $e^{x''}=1+x''+\frac{x''^2}{2!}+...+\frac{x''^n}{n!}...$, where
\begin{equation}
\label{f2037}
e^{x''}-1=x''+\frac{x''^2}{2!}+\frac{x''^n}{n!}...=\{x'',x''+\frac{x''^2}{2!},x''+\frac{x''^2}{2!}{2!}+...\frac{x''^n}{n!},...\,|\,1,\frac{1}{2},...,\frac{1}{n},...\}
\end{equation}
when $x''>0$ and
\begin{equation}
\label{f2038}
e^{x''}-1=x''+\frac{x''^2}{2!}+\frac{x''^n}{n!}...=\{-1,=\frac{1}{2},...,-\frac{1}{n},..\,|\,x'',x''+\frac{x''^2}{2!},x''+\frac{x''^2}{2!}{2!}+...\frac{x''^n}{n!},...\}
\end{equation}
when $x''<0$ (see details in \cite{l1332}, p. 25).

(A number $x'$ is called an {\it initial} part of $x$, $x''$ is called a {\it finite} part of $x$ and 
$y=\varphi(x)=e^x$ is called an exponential function, where $e\in{\bf R}$ is the Euler number.)

The expotential function $y=e^x$ does indeed behave like $x$th power of $e$  taking into account  that the Conway function has the following properties: $\omega^0=1$, $\omega^{-x}=\frac{1}{\omega^x}$, $\omega^{x+y}=\omega^x\cdot\omega^y$ (see \cite{l22}, Theorem 20, p. 32).

Notice  that for any two numbers $x_1,x_2\in{\bf No}$ such that $e^{x_1}=e^{x_2}$, then $\frac{e^{x_1}}{e^{x_2}}=1$ and also $e^{x_1-x_2}=1$ thus $x_1-x_2=0$ and $x_1=x_2$. Thus the exponential function $y=a^x$ is injective.

Denote by ${\bf D}^+_\Omega=\varphi({\bf No})$,  ${\bf D}^+_\zeta=\varphi({\bf A}_\zeta)$ and by ${\bf D^+}_{\frac{1}{\zeta}}=\varphi({\bf A}_{\frac{1}{\zeta}})$, 
where $\zeta=\omega^{\omega^\mu}$, $\zeta'=\omega^{\omega^{\mu+1}}$,    $0\leq\mu\leq\Omega$. Now we define ${\bf D}_\Omega={\bf D}_\Omega\cup{\bf D}^-_\Omega$, where ${\bf D}^-_\Omega\stackrel{def}{=}\{-x\}_{x\in{\bf D}^+_\Omega}.$

Thus we obtain a filtered multiplicative group ${\bf D}_\Omega$:

${\bf D}_\Omega\supset...\,\supset {\bf D}_{\zeta'}\supset   \,\,\,{\bf D}_\zeta\supset...\supset   {\bf D}_\omega\supset{\bf D}_{\frac{1}{\omega}}\supset...\supset {\bf D}_{\frac{1}{\zeta'}}\,\,\supset   {\bf D}_{\frac{1}{\zeta}}\supset\,\,\,\,...$

Now we can see that if a number $y$ has a canonical form $(\ref{f0201})$ $y=\sum\limits_{0\leq\beta<\alpha}\omega^{z_\beta}\cdot r_\beta$ such that $r_0>0$,   a canonical form of the number $z_0=\sum\limits_{0\leq\gamma<\alpha_0}\omega^{t_\gamma}\cdot s_\gamma$ has the following property: $t_\gamma>-1$,  for all $0\leq\gamma<\alpha_0$, then there is a number $x\in{\bf No}$ such that $e^x=y$ and hence $y\in{\bf D}^+_\Omega$. Thus $x=\ln y$ is a logarithm of $y$.

 Really, if $y=\sum\limits_{0\leq\beta<\alpha}\omega^{z_\beta}\cdot r_\beta$ has a canonical form as above, i.e.,  
$t_\beta> -1$, then obviously
\begin{equation}
\label{f4748}
y=\omega^{z_0}\cdot r_0\cdot(1+\delta)
\end{equation}
where 
\begin{equation}
\label{f4749}
\delta=\sum\limits_{0<\beta<\alpha}\omega^{z_\beta-z_0}\cdot \frac{r_\beta}{r_0}
\end{equation}
and $\delta$ is an infinitesimal number, i.e., Conway numbers $\delta$ such that $|\delta|<\frac{1}{\alpha}$ for all natural numbers $0<\alpha<\omega$.

 Now we can define $x$ as a value of the logarithmic function $\ln y$  by the following formulas:
\begin{equation}
\label{f133}
x=\ln y\stackrel{def}{=}\omega\cdot{z_0}+\ln r_0+\ln (1+\delta),
\end{equation}
 which is well-defined because $\ln r_0$ is a real-valued logarithm of real number $r_0>0$ and $\ln(1+\delta)$ is given by the following formulas:

\begin{equation}
\label{f4203}
\ln(1+\delta)=\delta-\frac{\delta^2}{2}+...+(-1)^{n-1}\frac{\delta^n}{n}+...=\{\delta,\delta-\frac{\delta^2}{2},...,\delta-\frac{\delta^2}{2}+...+(-1)^{n-1}\frac{\delta^n}{n},...\,|\,1,\frac{1}{2},...,\frac{1}{n},...\}
\end{equation}
when $\delta>0$ and
\begin{equation}
\label{f4203}
\ln(1+\delta)=\delta-\frac{\delta^2}{2}+...+(-1)^{n-1}\frac{\delta^n}{n}+...=\{-1,-\frac{1}{2},...,-\frac{1}{n},...\,|\,\delta,\delta-\frac{\delta^2}{2},...,\delta-\frac{\delta^2}{2}+...+(-1)^{n-1}\frac{\delta^n}{n},...\}
\end{equation}
when $\delta<0$.

It is really an exponent whose power $e^x=y$. Indeed, 
\begin{equation}
\label{f4750}
e^x=e^{\omega\cdot{z_0}+\ln r_0+\ln (1+\delta)}=e^{\omega\cdot z_0}\cdot e^{\ln r_0}\cdot e^{\ln (1+\delta)}=
\omega^{z_0}\cdot r_0\cdot(1+\delta)=y,
\end{equation}
because $e^{\omega\cdot z_0}=\omega^{z_0}$ since $\omega\cdot z_0=\sum\limits_{0\leq\gamma<\alpha_0}\omega^{t_\gamma+1}\cdot s_\gamma$ and $t_\gamma+1>0$ for all $0\leq\gamma<\alpha_0$, and thus the formula $e^{\omega\cdot z_0}=\omega^{y_0}$ is correct because the number $\omega\cdot z_0=z'+z''$, where $z''=0$.

Moreover, every number $y\in{\bf D}^+_\Omega$ has such canonical form. Indeed, let
$y=\sum\limits_{0\leq\beta<\alpha}\omega^{z_\beta}\cdot r_\beta$ be a canonical form of $y$. Since $y\in{\bf D}^+_\Omega$ there is a number $x\in{\bf No}$ such that $e^x=y$. Consider a canonical form $x=\sum\limits_{0\leq\alpha<\gamma}\omega^{y_\alpha}\cdot r_\alpha$ in which $\gamma$ denotes some ordinal, the numbers $r_\alpha$ $(0\leq\alpha<\gamma)$ are non-zero reals, and the numbers $y_\alpha$ form a descending sequence of numbers.  

Let $x=x'+x''$  be the decomposition of a number  $x$ into two terms: an initial infinite number  $x'=\sum\limits_{0\leq\alpha<\gamma'}\omega^{y_\alpha}\cdot r_\alpha$ 
in which $\gamma'\leq\gamma$, $y_\alpha>0$ for all $0\leq\alpha<\gamma'$,  $y_{\gamma'}\leq 0$ , and
a finite number $x''=\sum\limits_{\gamma'\leq\alpha<\gamma'}\omega^{y_\alpha}\cdot r_\alpha$. Then $y=\sum\limits_{0\leq\beta<\alpha}\omega^{z_\beta}\cdot r_\beta=e^x=e^{x'}\cdot e^{x''}= e^{\omega\cdot\frac{x'}{\omega}}\cdot e^{x''}=\omega^{\frac{x'}{\omega}}\cdot e^{x''}=\omega^{\frac{x'}{\omega}}\cdot(1+\delta) =\omega^{\frac{x'}{\omega}}+\omega^{\frac{x'}{\omega}}\cdot\delta$, where $\delta$ is a infinitesimal. Thus $z_0=(\sum\limits_{0\leq\alpha<\gamma'}\omega^{y_\alpha}\cdot r_\alpha)\cdot\omega^{-1}=\sum\limits_{0\leq\alpha<\gamma'}\omega^{y_\alpha-1}\cdot r_\alpha=\sum\limits_{0\leq\alpha<\gamma'}\omega^{t_\alpha}\cdot r_\alpha$ and $t_\alpha>-1$ because $y_\alpha>0$.

Now let's define  an additive group ${\bf L}_\Omega$ in the following way: 
$${\bf L}_\Omega=\ln[{\bf D}_\Omega]\bigoplus{\bf i}\pi\cdot{\mathbb Z}_2,$$ i.e., a direct sum of additive groups $\ln[{\bf D}_\Omega]$ and ${\bf i}\pi\cdot{\mathbb Z}_2$, where $\ln[{\bf D}_\Omega]={\bf A}_\Omega$ and ${\mathbb Z}_2=\{0,1\}$ with $0+1=1+0=1$, $1+1=0$, $\pi\cdot 1=\pi$, is an additive group  with the following adding: $x+y=(\ln x'+{\bf i}\pi\cdot x'')+(\ln y'+{\bf i}\pi\cdot y'')=(\ln x'+\ln y')+{\bf i}\pi\cdot( x''+ y'')$, where $x',y'\in{\bf D}_\Omega$ and $x'',y''\in{\mathbb Z}_2$, where ${\bf i}^2=-1$ and $\pi\in{\bf R}$. 
Actually, ${\bf L}_\Omega=({\bf No}\bigoplus{\bf i}\pi\cdot{\mathbb Z})/{\bf i}(\pi\cdot{\mathbb Z})$, where ${\bf No}\bigoplus{\bf i}\pi\cdot{\mathbb Z}\subset{\bf C}[{\bf i}]$ and ${\bf C}[{\bf i}]$ is a proper class of all Conway complex numbers $z=x+{\bf i}y$, $x,y\in{\bf A}_\Omega$,  ${\bf i}^2=-1$ and $\pi\cdot{\mathbb Z}$ is an additive subgroup of real numbers ${\bf R}$.

Then $\varphi_\Omega:\ln[{\bf D}_\Omega]\bigoplus{\bf i}\pi\cdot{\mathbb Z}_2\rightarrow{\bf D}_\Omega$ is given by the following formula $\varphi_\Omega(x)=e^x=e^{\ln x'}\cdot e^{{\bf i}x''}=\varphi(\ln x')\cdot(\cos x''+{\bf i}\sin x'')$.

 Indeed, for any $x\in\ln[{\bf D}_\Omega]\bigoplus{\bf i}\pi\cdot{\mathbb Z}_2$ there is a number $x'\in{\bf D}_\Omega$ and $x''\in{\mathbb Z}_2$ such that $x=\ln x'+{\bf i}\pi\cdot x''$, where $x''\in{\mathbb Z_2}$. 
Then $e^x=e^{\ln x'+{\bf i}\pi\cdot x''}=e^{\ln x'}\cdot e^{{\bf i}\pi\cdot x''}= x'\cdot(\cos \pi\cdot x''+{\bf i}\sin \pi\cdot x'')=x'\cdot \cos x''$. 
Clearly, if $x''=0$, then $\varphi(x)\in{\bf D}_\Omega$;  if $x''=\pi$, then $\varphi(x)\in-{\bf D}_\Omega\stackrel{def}{=}\{-x\}_{x\in{\bf D}_\Omega}$.

Moreover, $\varphi_\Omega$ is an isomorphism of groups. Indeed, if $x,y\in\ln[{\bf D}_\Omega]\bigoplus{\bf i}\pi\cdot{\mathbb Z}_2$, then $x+y=(\ln x'+{\bf i}\pi\cdot x'') +(\ln y' +{\bf i}\pi\cdot y'')=\ln x' +\ln y' +{\bf i}\pi\cdot(x''+y'')$ and $\varphi_\Omega(x+y)=e^{\ln x' +\ln y' +{\bf i}\pi\cdot(x''+y'')}=e^{\ln x'}\cdot e^{\ln y'}\cdot e^{{\bf i}\pi\cdot x''}\cdot e^{{\bf i}\pi\cdot y''}=e^{\ln x'}\cdot e^{{\bf i}\pi\cdot x''}\cdot e^{\ln y'}\cdot e^{{\bf i}\pi\cdot y''}=e^{\ln x'+{\bf i}\pi\cdot x''}\cdot e^{\ln y'+{\bf i}\pi\cdot y''}=\varphi_\Omega(x)\cdot\varphi_\Omega(y).$ 

One can see that
$\varphi_{\Omega}({\bf L}_\zeta)\subset{\bf C}_\zeta$ and  $\varphi_\Omega({\bf A}_\zeta)\subset{\bf B}_\zeta$ for all $\zeta=\omega^{\omega^\mu}$, $0<\mu\leq\Omega$,   $\varphi_\Omega({\bf L}_\omega)={\bf C}_\omega$ and  $\varphi_\Omega({\bf A}_{\frac{1}{\zeta}})={\bf B}_{\frac{1}{\zeta}}$ for all $\zeta=\omega^{\omega^\mu}$, $0\leq\mu\leq\Omega$ because each number $y\in{\bf C}_\omega$ has a canonical form $y=\sum\limits_{0\leq\beta<\gamma}\omega^{z_\beta}r_\beta$ such that $z_0= \omega^0=1>0$ and hence there exists $\ln y\in{\bf A}_\omega$.

And we obtain the following filtrations:

${\bf C}_\Omega\supset...\supset \,\,{\bf C}_{\zeta'}\supset   {\bf C}_\zeta\,\,\,\supset...\supset   {\bf C}_\omega\supset{\bf C}_{\frac{1}{\omega}}\supset...\,\,\supset {\bf C}_{\frac{1}{\zeta'}}\supset   {\bf C}_{\frac{1}{\zeta}}\supset\,\,\,...$

$\,\,\bigcup\,\,\,\,\,\,\,\,\,...\,\,\,\,\,\,\,\,\,\,\,\bigcup\,\,\,\,\,\,\,\,\,\,\,\,\bigcup\,\,\,\,\,\,\,\,\,\,\,...\,\,\,\,\,\,\,\,\Vert\,\,\,\,\,\,\,\,\,\,\,\,\,\,\Vert\,\,\,\,\,\,\,\,\,\,\,\,\,...\,\,\,\,\,\,\,\,\,\,\,\,\Vert\,\,\,\,\,\,\,\,\,\,\,\,\,\,\,\,\Vert\,\,\,\,\,\,\,\,\,\,\,\,\,\,\,...$

${\bf D}_\Omega\supset...\,\,\supset \,\,{\bf D}_{\zeta'}\supset   {\bf D}_\zeta\supset...\supset   {\bf D}_\omega\supset{\bf D}_{\frac{1}{\omega}}\supset...\,\,\,\supset {\bf D}_{\frac{1}{\zeta'}}\supset   {\bf D}_{\frac{1}{\zeta}}\supset\,\,\,\,...\,\,\,\,\,\,\,\,\,\,\,\,\,\,\,\,\,\,\,\,\,\,$
and

${\bf B}_\Omega\supset...\,\,\supset \,\,{\bf B}_{\zeta'}\,\,\,\supset   {\bf B}_\zeta\supset...\supset   {\bf B}_\omega\supset{\bf B}_{\frac{1}{\omega}}\,\,\,\supset...\,\,\,\,\,\,\supset {\bf B}_{\frac{1}{\zeta'}}\supset   {\bf B}_{\frac{1}{\zeta}}\supset\,\,\,\,...$

$\,\,\bigcup\,\,\,\,\,\,\,\,\,...\,\,\,\,\,\,\,\,\,\,\,\bigcup\,\,\,\,\,\,\,\,\,\,\,\,\bigcup\,\,\,\,\,\,\,\,\,\,\,...\,\,\,\,\,\,\,\,\,\,\Vert\,\,\,\,\,\,\,\,\,\,\,\,\,\,\Vert\,\,\,\,\,\,\,\,\,\,\,\,\,\,\,\,...\,\,\,\,\,\,\,\,\,\,\,\,\,\,\Vert\,\,\,\,\,\,\,\,\,\,\,\,\,\,\,\,\Vert\,\,\,\,\,\,\,\,\,\,\,\,\,\,\,...$

${\bf D}^+_\Omega\supset...\,\,\supset \,\,{\bf D}^+_{\zeta'}\supset   {\bf D}^+_\zeta\supset...\supset   {\bf D}^+_\omega\supset{\bf D^+}_{\frac{1}{\omega}}\supset\,...\,\,\,\supset {\bf D}^+_{\frac{1}{\zeta'}}\supset   {\bf D}^+_{\frac{1}{\zeta}}\supset\,\,\,\,...\,\,\,\,\,\,\,\,\,\,\,\,\,\,\,$
where $\zeta=\omega^{\omega^\mu}$, $\zeta'=\omega^{\omega^{\mu+1}}$,    $0\leq\mu<\Omega$.
 $\Box$

{\bf Remark 18.} Note that the proof of Proposition 21 is not so simple and short. The fact is that defining the functions $e^x$  for a finite Conway number $x$ that is not a real number and $\ln x$ when $x$ is not a positive real number (note also that always $x$ is not a infinitesimal number) is not so simple and it will take several pages to identify them (see \cite{l1332}, pp. 25-31).

{\bf Remark 19.} Of course, there is a more general mapping $\bar\varphi:{\bf C}[{\bf i}]\rightarrow{\bf C}[{\bf i}]\setminus\{0\}$ of groups  given by the following formula: $\bar\varphi(x+{\bf i}y)=e^z=e^x\cdot e^{{\bf i}y}=e^x\cdot(\cos y +{\bf i}\sin y)$, where $z=x+{\bf i}y\in{\bf C}[{\bf i}]$, $x,y\in{\bf No}$ and $\cos y=\cos y''$ and $\sin y=\sin y''$, $y=y'+y''$ (see paragraph: ${\bf 2. 2}$ above), but we will omit this extension, which is the same as the previous one, except that this homomorphism $\bar\varphi$ is not a monomorphism because it is periodic along the imaginary axis, i.e., $e^{z+{\bf i}2\pi}=e^z$.

{\bf Proposition 22.} {\it For each number $y\in{\bf D}^+_\Omega$ all numbers $\bar y$ which are commensurate $y$ are also in ${\bf D}^+_\Omega$.}

{\bf Proof.} If  $y\in{\bf D}^+_\Omega$, then as we saw above  a canonical form  $y=\sum\limits_{0\leq\beta<\alpha}\omega^{z_\beta}\cdot r_\beta$ of $y$ has the property: a canonical form of the number $z_0=\sum\limits_{0\leq\gamma<\alpha_0}\omega^{t_\gamma}\cdot s_\gamma$ is such that $t_\gamma>-1$,  for all $0\leq\gamma<\alpha_0$.

If numbers $y$ and $\bar y$ are commensurate, then for some natural number $n$ we have inequalities $y<n\cdot\bar y$ and $\bar y<n\cdot y$ which imply $\omega^{z_0}\cdot r_0<n\cdot \omega^{\bar z_0}\cdot \bar r_0$ and $\omega^{\bar z_0}\cdot \bar r_0<n\cdot\omega^{ z_0}\cdot r_0$ and we obtain $z_0=\bar z_0$. Note that $\omega ^{z_0}$ is a unique simplest number in this commensurate class and called a {\it leader}.

{\bf Corollary.} ${\bf D}^+_\Omega=\bigcup\limits_{\omega^{z_0}\in{\bf D}^+_\Omega}[\omega^{z_0}]$, where $[\omega^{z_0}]$ is a commensurate class of number $\omega^{z_0}$, $z_0=\sum\limits_{0\leq\gamma<\alpha_0}\omega^{t_\gamma}\cdot s_\gamma$ is such that $t_\gamma>-1$,  for all $0\leq\gamma<\alpha_0$ and the class ${\bf No}^+\setminus{\bf D}^+_\Omega$ is the class of all numbers $y$ such that there is no magnitude like $\ln y$.

This Corollary shows that class ${\bf D}^+_\Omega$ is quite extensive,  in comparison with class ${\bf A}^+_\omega$ and the same cannot be said about Class ${\bf No}^+\setminus{\bf D}^+_\Omega$. Let's give some illustrative examples.

{\bf 1.} $(+\infty_{\omega^{-2}},+\infty_{\omega^{-\frac{1}{2}}})\cup(+\infty_{\omega^{\frac{1}{2}}},+\infty_{\omega^2})\subset{\bf D}^+_\Omega$ and  $(+\infty_{\frac{1}{c_1}},+\infty_{\frac{1}{\omega}})\cup(+\infty_\omega,+\infty_{c_1})\subset{\bf No}^+\setminus{\bf D}^+_\Omega$, where  $c_1=\{\alpha\cdot\omega^{\omega^{-1}}\,|\,\omega^{\alpha^{-1}}\}_{0<\alpha<\omega}$.

Note that there are similar intervals, e.g., $(+\infty_{\omega^{\frac{1}{4}}},+\infty_{\omega^{\frac{3}{4}}})$ and $(+\infty_{\omega^{\frac{1}{4}+\frac{1}{\omega}}},+\infty_{\omega^{\frac{3}{4}+\frac{1}{\omega}}})$ between intervals $(+\infty_\omega,+\infty_{c_1})$  and $(+\infty_{\omega^{-2}},+\infty_{\omega^2})$  such that $(+\infty_{\omega^{\frac{1}{4}}},+\infty_{\omega^{\frac{3}{4}}})\subset{\bf D}^+_\Omega$ and $(+\infty_{\omega^{\frac{1}{4}+\frac{1}{\omega}}},+\infty_{\omega^{\frac{3}{4}+\frac{1}{\omega}}})\subset{\bf No}^+\setminus{\bf D}^+_\Omega$, respectively.

{\bf 2.} $(+\infty_{\omega^{-n}},+\infty_{\omega^{-n+1-\frac{1}{2}}})\cup(+\infty_{\omega^{n-1+\frac{1}{2}}},+\infty_{\omega^n})\subset{\bf D}^+_\Omega$ and  $(+\infty_{\frac{1}{c_n}},+\infty_{\frac{1}{\omega^{n-1}}})\cup(+\infty_{\omega^{n-1}},+\infty_{c_n})\subset{\bf No}^+\setminus{\bf D}^+_\Omega$, where  $c_n=\{\alpha\cdot\omega^{n-1+\omega^{-1}}\,|\,\omega^{n-1+\alpha^{-1}}\}_{0<\alpha<\omega}$

{\bf 3}. $(+\infty_{\frac{1}{\omega^{\omega+1}}},+\infty_{\omega^{-\frac{1}{\omega^{\omega+1}}}})\cup
(+\infty_{\omega^{\frac{1}{\omega^{\omega+1}}}}, +\infty_{\omega^{\omega+1}})\subset{\bf D}^+_\Omega$. And more general, for all ordinals $0<\beta<\Omega$   intervals $(+\infty_{\frac{1}{{c_\beta}}}, 
+\infty_{\frac{1}{{c_{-\beta}}}})\cup(+\infty_{c_{-\beta}},
+\infty_{c_\beta})\subset{\bf D}^+_\Omega$, where  $c_\beta=\{\alpha\cdot\omega^{\omega^{\beta^{-1}}}\,|\,\omega^{\omega^{\alpha^{-1}}}\}_{0<\alpha<\omega}$ and 
$c_{-\beta}=\{\frac{\omega^{\omega^{\beta^{-1}}}}{\alpha}\,|\,\omega^{\omega^{\alpha^{-1}}}\}_{0<\alpha<\omega}$ and

   intervals $(+\infty_{\frac{1}{{c_\beta}}}, 
+\infty_{\frac{1}{{c_{-\beta}}}})\cup(+\infty_{c_{-\beta}},
+\infty_{c_\beta})\subset{\bf No}^+\setminus{\bf D}^+_\Omega$, where  $c_\beta=\{\alpha\cdot\omega^{\omega^{\beta^{-1}}}\cdot\omega^{\omega^{-1}}\,|\,\omega^{\omega^{\alpha^{-1}}}\cdot\omega^{\omega^{-1}}\}_{0<\alpha<\omega}$ and 
$c_{-\beta}=\{\frac{\omega^{\omega^{\beta^{-1}}}\cdot\omega^{\omega^{-1}}}{\alpha}\,|\,\omega^{\omega^{\alpha^{-1}}}\}_{0<\alpha<\omega}$

 And so on and so forth.

\bigskip
\bigskip
In conclusion, I would like to say two things. The form of Conway's numbers is the most successful form I know, i.e., $x=\{L\,|\,R\}$, \cite{l22}, p. 4,  allowing us to see the unique feature of ordinal numbers $\lambda=\{L\,|\,\}$, \cite{l22}, p. 26, i.e., when $R=\emptyset$, with their pulsating periodic \grqq jumps\grqq, as well as the role of real numbers $x=\{x-(1/n)\,|\,x-(1/n)\}_{n>0}$, \cite{l22}, p. 24,  in defining  Conway's exponential function, \cite{l22}, p. 31, and the normal form of number $x=\sum\limits_{\beta<\alpha}\omega^{y_\beta}\cdot r_\beta$, \cite{l22}, p. 33. 

And of course Conway's principle of induction: \grqq Objects may be created from earlier objects in any reasonably constructive fashion\grqq, \cite{l22}, p. 66, allows us to see  general laws not only for Cantor's \grqq consistent multiplicity\grqq,  allowed by formal systems, but also for cases of undefinable multiplicity, such as eternal sequences $(\alpha)_{\alpha>0}$ that Cantor considered as \grqq inconsistent multiplicity\grqq.

 \end{document}